\documentclass[aos,preprint]{imsart}

\RequirePackage[OT1]{fontenc}
\RequirePackage{amsthm,amsmath,amsfonts,mathrsfs,amssymb,bm,dsfont,enumitem,graphicx}
\RequirePackage[numbers]{natbib}
\RequirePackage[colorlinks,citecolor=blue,urlcolor=blue]{hyperref} 
\RequirePackage{subfig}

\pdfoutput=1

\arxiv{arXiv:1804.04916}

\startlocaldefs
\numberwithin{equation}{section}
\theoremstyle{plain}
\newtheorem{thm}{Theorem}[section]
\newtheorem{lem}{Lemma}[section]
\newtheorem{assumption}{Assumption}
\endlocaldefs

\theoremstyle{definition}
\newtheorem{remark_tmp}{Remark}[section]
\newenvironment{remark}
{ \begin{remark_tmp} 	}
	{ 
		\medskip\hfill{$\blacklozenge$}
	\end{remark_tmp} 
}

\theoremstyle{definition}
\newtheorem{definition_tmp}{Definition}[section]
\newenvironment{definition}
{ \begin{definition_tmp} 	}
	{ 
		\medskip\hfill{$\blacksquare$}
	\end{definition_tmp} 
}

\allowdisplaybreaks[1]

\DeclareMathOperator*{\argmin}{arg\,min}

\DeclareMathOperator*{\esssup}{ess\,sup}
\DeclareMathOperator*{\tr}{trace}
\DeclareMathOperator*{\sign}{sign}
\DeclareMathOperator*{\supp}{supp}
\DeclareMathOperator*{\diam}{diam}
\DeclareMathOperator*{\vol}{vol}
\DeclareMathOperator*{\diag}{diag}

\DeclareMathOperator*{\clo}{clo}

\renewcommand{\P}{\mathbb{P}}
\newcommand{\E}{\mathbb{E}}
\newcommand{\V}{\mathbb{V}}
\newcommand{\Cov}{\mathbb{C}\mathrm{ov}}
\newcommand{\I}{\mathds{1}}
\newcommand{\G}{\mathbb{G}}

\newcommand{\bH}{\mathbf{H}}
\newcommand{\bI}{\mathbf{I}}
\newcommand{\bM}{\mathbf{M}}
\newcommand{\bN}{\mathbf{N}}
\newcommand{\bP}{\mathbf{P}}
\newcommand{\bQ}{\mathbf{Q}}
\newcommand{\bR}{\mathbf{R}}
\newcommand{\bX}{\mathbf{X}}
\newcommand{\bZ}{\mathbf{Z}}
\newcommand{\bb}{\mathbf{b}}
\newcommand{\be}{\mathbf{e}}
\newcommand{\boldm}{\mathbf{m}}
\newcommand{\bp}{\mathbf{p}}
\newcommand{\bq}{\mathbf{q}}
\newcommand{\bt}{\mathbf{t}}
\newcommand{\bv}{\mathbf{v}}
\newcommand{\bx}{\mathbf{x}}
\newcommand{\by}{\mathbf{y}}

\newcommand{\bbeta}{\boldsymbol{\beta}}
\newcommand{\bkappa}{\boldsymbol{\kappa}}

\newcommand{\bgamma}{\boldsymbol{\gamma}}
\newcommand{\bphi}{\boldsymbol{\phi}}
\newcommand{\btau}{\boldsymbol{\tau}}

\newcommand{\bSigma}{\boldsymbol{\Sigma}}
\newcommand{\bUpsilon}{\boldsymbol{\Upsilon}}
\newcommand{\bvsig}{\boldsymbol{\varsigma}}

\newcommand{\bgTmp}{\widehat{\bm{\gamma}}}
\newcommand{\bg}[1]{\bm{\bgTmp}_{\bq,#1}(\bx)}

\newcommand{\bPiTmp}{\bm{\Pi}}
\newcommand{\bPi}[1]{\bm{\bPiTmp}_{#1}}

\newcommand{\bOTmp}{\Omega}
\newcommand{\bO}[1]{\bOTmp_{#1}(\bx)}
\newcommand{\bOhat}[1]{\widehat{\bOTmp}_{#1}(\bx)}
\newcommand{\bOstar}[1]{\widehat{\bOTmp}^*_{#1}(\bx)}  

\newcommand{\bSigTmp}{\bm{\Sigma}}
\newcommand{\bSig}[1]{\bm{\bSigTmp}_{#1}}
\newcommand{\bSighat}[1]{\bm{\widehat{\bSigTmp}}_{#1}}
\newcommand{\bSigbar}[1]{\bm{\bar{\bSigTmp}}_{#1}}

\newcommand{\e}{\varepsilon} 

\newcommand{\bpt}{\bm{\tilde{\bp}}}

\setenumerate[1]{label=\bf(\alph*)}

\begin{document}

\begin{frontmatter}
\title{Large Sample Properties of\\ Partitioning-Based Series Estimators}
\runtitle{Partitioning-Based Series Estimators}

\begin{aug}
\author{\fnms{Matias D.} \snm{Cattaneo}\thanksref{t1}\ead[label=e1]{cattaneo@princeton.edu}},
\author{\fnms{Max H.} \snm{Farrell}\ead[label=e2]{max.farrell@chicagobooth.edu}}
\and
\author{\fnms{Yingjie} \snm{Feng}\ead[label=e3]{yjfeng@umich.edu}}

\thankstext{t1}{Financial support from the National Science Foundation (SES 1459931) is gratefully acknowledged.}
\runauthor{Cattaneo, Farrell, and Feng}

\affiliation{Princeton University, University of Chicago, and Princeton University}

\address{Matias D. Cattaneo\\
	Department of Operations Research\\ and Financial Engineering\\
	Princeton University\\
	Princeton, NJ 08544\\
	\printead{e1}
}

\address{Max H. Farrell\\
	Booth School of Business\\
	University of Chicago\\
	Chicago, IL 60637\\
	\printead{e2}
}

\address{Yingjie Feng\\
	Department of Politics\\
	Princeton University\\
	Princeton, NJ 08544\\
	\printead{e3}
}
\end{aug}

\begin{abstract}
We present large sample results for partitioning-based least squares nonparametric regression, a popular method for approximating conditional expectation functions in statistics, econometrics, and machine learning. First, we obtain a general characterization of their leading asymptotic bias. Second, we establish integrated mean squared error approximations for the point estimator and propose feasible tuning parameter selection. Third, we develop pointwise inference methods based on undersmoothing and robust bias correction. Fourth, employing different coupling approaches, we develop uniform distributional approximations for the undersmoothed and robust bias-corrected $t$-statistic processes and construct valid confidence bands. In the univariate case, our uniform distributional approximations require seemingly minimal rate restrictions and improve on approximation rates known in the literature. Finally, we apply our general results to three partitioning-based estimators: splines, wavelets, and piecewise polynomials. The supplemental appendix includes several other general and example-specific technical and methodological results. A companion \textsf{R} package is provided.
\end{abstract}



\end{frontmatter}

\section{Introduction}\label{section: Introduction}

We study the standard nonparametric regression setup, where $\{(y_i,\bx_i'), i = 1, \ldots n\}$ is a random sample from the model 
\begin{equation}\label{equation: regression model}
y_i=\mu(\bx_i)+\e_i ,    	 	\qquad  	  \E[\e_i|\bx_i]=0  ,    	 	\qquad  	  \E[\e_i^2|\bx_i]=\sigma^2(\bx_i)  ,
\end{equation}
for a scalar response  $y_i$ and a $d$-vector of continuously distributed covariates $\bx_i = (x_{1,i}, \ldots, x_{d,i})'$ with compact support $\mathcal{X}$. The object of interest is the unknown regression function $\mu(\cdot)$ and its derivatives. We focus on \emph{partitioning-based}, or locally-supported, series (linear sieve) least squares regression estimators, which are characterized by two features. First, the support $\mathcal{X}$ is partitioned into non-overlapping cells, which are then used to form a set of basis functions. Second, the final fit is determined by a least squares regression using these bases. The key distinguishing characteristic is that each basis function is nonzero only on a small, contiguous set of cells of the partition. Popular examples include splines, compactly supported wavelets, and piecewise polynomials. For this class of estimators, we develop novel bias approximations, integrated mean squared error (IMSE) expansions useful for tuning parameter selection, and pointwise and uniform estimation and inference results, with and without bias correction techniques.

A partitioning-based estimator is made precise by the partition of $\mathcal{X}$ and basis expansion used. Let $\Delta=\{\delta_{l} \subset \mathcal{X}: 1 \leq l \leq \bar{\kappa}\}$ be a collection of $\bar{\kappa}$ open and disjoint sets, the closure of whose union is $\mathcal{X}$ (or, more generally, covers $\mathcal{X}$). We restrict $\delta_{l}$ to be polyhedral, which allows for tensor products of (marginally-formed) intervals as well as other popular partitioning shapes. Based on this partition, the dictionary of $K$ basis functions, each of order $m$ (e.g., $m=4$ for cubic splines) is denoted by $\bx_i \mapsto \bp(\bx_i) := \bp(\bx_i;\Delta,m)=(p_1(\bx_i;\Delta,m),\cdots,p_K(\bx_i;\Delta,m))'$. For $\bx \in \mathcal{X}$ and $\bq = (q_1, \ldots, q_d)'\in\mathbb{Z}_+^d$, the partial derivative $\partial^\bq\mu(\bx)$ is estimated by least squares regression
\begin{equation}\label{equation: partitioning estimator}
\widehat{\partial^\bq\mu}(\bx)=\partial^\bq \bp(\bx)'\widehat{\bbeta},\qquad\qquad   		  
\widehat{\bbeta}\in\underset{\mathbf{b}\in\mathbb{R}^K}{\argmin}\,\sum_{i=1}^{n}(y_i - \bp(\bx_i)'\mathbf{b})^2,
\end{equation}
where $\partial^\bq\mu(\bx) = \partial^{q_1 + \cdots + q_d}  \mu(\bx) / \partial^{q_1}x_1\cdots\partial^{q_d}x_d$ (for boundary points defined from the interior of $\mathcal{X}$ as usual) and $\mu(\bx):=\partial^\mathbf{0}\mu(\bx)$. 

The approximation power of this class of estimators comes from two user-specified parameters: the granularity of the partition $\Delta$ and the order $m\in\mathbb{Z}_+$ of the basis. The choice $m$ is often fixed in practice, and hence we regard $\Delta$ as the tuning parameter. Under our assumptions, $\bar{\kappa} \to \infty$ as the sample size $n\to\infty$, and the volume of each $\delta_l$ shrinks proportionally to $h^d$, where $h=\max\{\diam(\delta):\delta\in\Delta\}$ serves as a universal measure of the granularity. Thus, as $\bar{\kappa} \to \infty$, $h^d$ vanishes at the same rate, and with each basis being supported only on a finite number of cells, $K$ diverges proportionally as well. Complete, detailed examples of bases and partitioning schemes are discussed in the online supplement for brevity.

Our first contribution, in Section \ref{sec:bias}, is a general characterization of the bias of partitioning-based estimators, which we then use for both tuning parameter selection and robust bias correction. In the supplement, we specialize our generic bias approximation to splines, wavelets, and piecewise polynomials over different partitioning schemes, leading to novel results.

Our second contribution, in Section \ref{sec:imse}, is a general integrated mean squared error (IMSE) expansion for partitioning-based estimators. These results lead to IMSE-optimal partitioning choices, and hence deliver IMSE-optimal point estimators of the regression function and its derivatives. We show that the IMSE-optimal choice of partition granularity obeys $h_{\mathtt{IMSE}} \asymp n^{-1/(2m+d)}$, which translates to the familiar $K_{\mathtt{IMSE}} \asymp n^{-d/(2m+d)}$, and give a precise characterization of the leading constant. For simple cases on tensor-product partitions, some results exist for splines \cite{Agarwal-Studden_1980_AoS,Zhou-Shen-Wolfe_1998_AoS,Zhou-Wolfe_2000_SS} and piecewise polynomials \cite{Cattaneo-Farrell_2013_JoE}. In addition to generalizing these results substantially, our characterization for compactly supported wavelets appears to be new.

The IMSE-optimal partitioning scheme, and consistent implementations thereof, cannot be used directly to form valid pointwise or uniform (in $\bx \in \mathcal{X}$) inference procedures. Undersmoothing (employing a finer partition than would be IMSE-optimal) is theoretically valid for inference, but difficult to implement in a principled way. Inspired by results proving that undersmoothing is never optimal relative to bias correction for kernel-based nonparametrics \cite{Calonico-Cattaneo-Farrell_2018_JASA}, we develop three robust bias-corrected inference procedures using our new bias characterizations of partitioning-based estimators. These methods are more involved than their kernel-based counterparts, but are still based on least squares regression using partitioning-based estimation. Specifically, we show that the conventional partitioning-based estimator $\widehat{\partial^\bq\mu}(\bx)$ and the three bias-corrected estimators we propose have a common structure, which we exploit to obtain general pointwise and uniform distributional approximations under weak (sometimes minimal) conditions. These robust bias correction results for partitioning-based estimators, both pointwise and uniform in $\bx$, are practically useful because they allow for mean squared error minimizing tuning parameter choices, thus offering a data-driven method combining optimal point estimation and valid inference on the same partitioning scheme. 

Section \ref{sec:pointwise} establishes pointwise in $\bx\in\mathcal{X}$ distributional approximations for both conventional and robust bias-corrected $t$-statistics based on partitioning-based estimators. These pointwise distributional results are made uniform in Section \ref{sec:coupling}, where we establish a strong approximation for the whole $t$-statistic processes, indexed by the point $\bx \in \mathcal{X}$, covering both conventional and robust bias-corrected inference. To illustrate, Section \ref{sec:bands} constructs valid confidence bands for (derivatives of) the regression function using our uniform distributional approximations. When compared to the current literature, we obtain a strong approximation to the \textit{entire} $t$-statistic process under either weaker or seemingly minimal conditions on the tuning parameter $h$ (i.e., on $K$ or $\bar{\kappa}$), depending on the case under consideration.

Section \ref{sec:numerical} summarizes Monte Carlo results, Section \ref{sec:proofs} gives selected proofs, and Section \ref{sec:conclusion} concludes. A supplemental appendix (SA hereafter) gives complete proofs, several new technical and methodological results, further Monte Carlo evidence, and applies our general results to splines, wavelets, and piecewise polynomials. A companion {\sf R} package \cite{Cattaneo-Farrell-Feng_2019_lspartition} is provided.

\subsection{Related Literature}

This paper contributes primarily to two literatures, nonparametric regression and strong approximation. There is a vast literature on nonparametric regression, summarized in many textbook treatments \cite[e.g.,][and references therein]{Fan-Gijbels_1996_Book,Gyorfi-etal_2002_book,Ruppert-Wand-Carroll_2009_book}. Of particular relevance are treatments of series (linear sieve) methods, which offer some results concerning partitioning-based estimators in particular, many times limited to splines, wavelets, or piecewise polynomials considered separately \cite{Newey_1997_JoE,Huang_1998_AoS,Zhou-Shen-Wolfe_1998_AoS,Huang_2003_AoS,Cattaneo-Farrell_2013_JoE,Belloni-Chernozhukov-Chetverikov-Kato_2015_JoE,Chen-Christensen_2015_JOE,Chen-Christensen_2018_QE,Belloni-Chernozhukov-Chetverikov-FernandezVal_2019_JoE}. Piecewise polynomial fits on partitions have a long and ongoing tradition in statistics, dating at least to the regressogram of Tukey \cite{Tukey1961_Berkeley}, continuing through \cite{Stone_1982_AoS} (named local polynomial regression therein) and \cite{Cattaneo-Farrell_2013_JoE}, and up to modern, data-driven partitioning techniques such as regression trees \cite{Breiman-etal1984_book}, trend filtering \cite{Tibshirani2014_AOS}, and related methods \cite{Zhang-Singer_2010_Book}. Partitioning-based methods have also featured as inputs or preprocessing in treatment effects \cite{Cattaneo-Farrell_2011_BookCh,Calonico-Cattaneo-Titiunik_2015_JASA}, empirical finance \cite{Cattaneo-etal2019_sorting}, ``binscatter'' analysis \cite{Cattaneo-Crump-Farrell-Feng_2019_Binscatter}, and other settings. The bias corrections we develop for series estimation and uniform inference follow recent work on kernel-based nonparametric inference \cite{Calonico-Cattaneo-Titiunik_2014_ECMA,Calonico-Cattaneo-Farrell_2018_JASA,Calonico-Cattaneo-Farrell_2019_wp}. Our coupling and strong approximation results relate to early work discussed in \cite[][Chapter 22]{Eggermont-LaRiccia_2009_Book} and the more recent work in \cite{Chernozhukov-Lee-Rosen_2013_ECMA,Chernozhukov-Chetverikov-Kato_2014a_AoS,Chernozhukov-Chetverikov-Kato_2014b_AoS,Chernozhukov-Chetverikov-Kato_2015_PTRF,Chernozhukov-Chetverikov-Kato_2016_SPA} and \cite{Zhai_2018_TPRF}, as well as with the results for series estimators in \cite{Belloni-Chernozhukov-Chetverikov-Kato_2015_JoE} and \cite{Belloni-Chernozhukov-Chetverikov-FernandezVal_2019_JoE}. See also \cite{Zaitsev_2013_RMSurveys} for a review on strong approximation methods.

\subsection{Notation}
For a $d$-tuple $\bq=(q_1,\ldots,q_d)\in\mathbb{Z}^d_+$, define $[\bq]=\sum_{j=1}^{d}q_j$, $\bx^\bq=x_1^{q_1}x_2^{q_2}\cdots x_d^{q_d}$ and $\partial^\bq \mu(\bx)=\partial^{[\bq]} \mu(\bx)/\partial x_1^{q_1}\ldots\partial x_d^{q_d}$. Unless explicitly stated otherwise, whenever $\bx$ is a boundary point of some closed set, the partial derivative is understood as the limit with $\bx$ ranging within it. Let $\mathbf{0}=(0,\cdots,0)'$ be the length-$d$ zero vector. We set $\mu(\bx):=\partial^\mathbf{0}\mu(\bx)$ and $\widehat{\mu}_j(\bx):=\widehat{\partial^\mathbf{0}\mu}_j(\bx)$ for $j=0,1,2,3$ and collect the covariates as $\bX = [\bx_1, \ldots, \bx_n]'$. The tensor product or Kronecker product operator is $\otimes$. The smallest integer greater than or equal to $u$ is $\lceil u \rceil$. For two random variables $X$ and $Y$, $X=_d Y$ denotes that they have the same probability law.

We use several norms. For a vector $\mathbf{v}=(v_1,\ldots,v_M)\in\mathbb{R}^M$, we write $\|\mathbf{v}\|=(\sum_{i=1}^Mv_i^2)^{1/2}$ and $\dim(\mathbf{v})=M$. For a matrix $\mathbf{A}\in\mathbb{R}^{M\times N}$, $\|\mathbf{A}\|=\max_i\sigma_i(\mathbf{A})$ and $\|\mathbf{A}\|_\infty=\max_{1\leq i\leq M}\sum_{j=1}^N|a_{ij}|$ for operator norms induced by $L_2$ and $L_\infty$ norms, where $\sigma_i(\mathbf{A})$ is the $i$-th singular value of $\mathbf{A}$, and $\lambda_{\min}(\mathbf{A})$ is the minimum eigenvalue of $\mathbf{A}$. 

We use empirical process notation: $\E_n[g(\bx_i)]=\frac{1}{n}\sum_{i=1}^ng(\bx_i)$ and $\G_n[g(\bx_i)]=\frac{1}{\sqrt{n}}\sum_{i=1}^n(g(\bx_i)-\E[g(\bx_i)])$. For sequences, $a_n\lesssim b_n$ denotes $\limsup_n|a_n/b_n|$ is finite, $a_n=O_\P(b_n)$ denotes $\limsup_{\epsilon\to\infty}\limsup_n\P[|a_n/b_n|\geq\epsilon]=0$, $a_n=o(b_n)$ denotes $a_n/b_n\to 0$, $a_n=o_\P(b_n)$ denotes $a_n/b_n\to_\P 0$, where $\to_\P$ is convergence in probability, and $a_n\asymp b_n$ denotes $a_n\lesssim b_n$ and $b_n\lesssim a_n$. Limits are taken as $n\to\infty$ ($h\to0$, $K\to\infty$, as appropriate), unless otherwise stated.

Finally, throughout the paper, $r_n >0$ denotes a non-vanishing sequence and $\nu>0$ denotes a fixed constant used to characterize moment bounds.

\section{Setup}
\label{sec:setup}

We first make precise our setup and assumptions. Our first assumption restricts the data generating process.

\begin{assumption}[Data Generating Process]\label{assumption: DGP} \leavevmode
	\begin{enumerate}
		\item $\{(y_i, \bx_i'):1\leq i\leq n\}$ are i.i.d.\ satisfying (\ref{equation: regression model}), where $\bx_i$ has compact connected support $\mathcal{X}\subset\mathbb{R}^d$ and an absolutely continuous distribution function. The density of $\bx_i$, $f(\cdot)$, and the conditional variance of $y_i$ given $\bx_i$, $\sigma^2(\cdot)$, are bounded away from zero and continuous.
		
		\item $\mu(\cdot)$ is $S$-times continuously differentiable, for $S \geq [\bq]$, and all $\partial^{\bvsig}\mu(\cdot)$, $[\bvsig]=S$, are H\"older continuous with exponent $\varrho>0$.
	\end{enumerate}
\end{assumption}

The next two assumptions specify a set of high-level conditions on the partition and basis: we require that the partition is ``quasi-uniform'' and the basis is ``locally'' supported.

\begin{assumption}[Quasi-Uniform Partition] \label{assumption: quasi-uniform partition}
	The ratio of the sizes of inscribed and circumscribed balls of each $\delta\in\Delta$ is bounded away from zero uniformly in $\delta\in\Delta$, and
	\[\frac{\max\{\diam(\delta): \delta\in\Delta\}}{\min\{\diam(\delta):\delta\in\Delta\}}\lesssim 1,\]
	where $\diam(\delta)$ denotes the diameter of $\delta$. Further, for $h=\max\{\diam(\delta):\delta\in\Delta\}$, assume $h=o(1)$.
\end{assumption}

This condition implies that the size of each $\delta\in\Delta$ can be well characterized by the diameter of $\delta$, where we use $h$ as a universal measure of mesh sizes of elements in $\Delta$. In the univariate case, it reduces to a bounded mesh ratio. A special case of a quasi-uniform partition is one formed via a tensor product of univariate marginal partitions on each dimension of $\bx\in\mathcal{X}$, with appropriately chosen knot positions. The SA (\S SA-3) gives details and discusses this special example of partitioning scheme. If $\Delta$ covers only strict subset of $\mathcal{X}$, then our results hold on that subset. 

We focus on nonrandom partitions. Data-dependent partitioning can be accommodated by sample splitting: estimating the partition configuration in one subsample and performing inference in the other. In this way, quite general partitions can be used with our results, including data-driven methods such as regression trees and other modern machine learning techniques. In fact, these modern methods would typically generate non-tensor-product partitioning schemes. We defer studying general data-dependent partitioning for future work, but note that a few specific results are available \cite{Breiman-etal1984_book,Nobel_1996_AoS,Calonico-Cattaneo-Titiunik_2015_JASA}.

The second assumption on the partitioning-based estimators employs generalized notions of \textit{stable local basis} \citep{Davydov_2001_JAT} and \textit{active basis} \citep{Huang_2003_AoS}. We say a function $p(\cdot)$ on $\mathcal{X}$ is \textit{active} on $\delta\in\Delta$ if it is not identically zero on $\delta$.

\begin{assumption}[Local Basis]\label{assumption: Local Basis} \leavevmode
	\begin{enumerate}
		\item For each basis function $p_k$, $k=1,\ldots, K$, the union of elements of $\Delta$ on which $p_k$ is active is a connected set, denoted by $\mathcal{H}_k$. For all $k=1,\ldots, K$, both the number of elements of $\mathcal{H}_k$ and the number of basis functions which are active on $\mathcal{H}_k$ are bounded by a constant.
		
		\item For any $\mathbf{a}=(a_1,\cdots, a_K)'\in\mathbb{R}^{K}$,
		\[\mathbf{a}'\int_{\mathcal{H}_k}\bp(\bx;\Delta, m)\bp(\bx; \Delta, m)'\,d\bx\, \mathbf{a}\gtrsim a_k^2h^d,
		\qquad k=1,\ldots,K.
		\]
		
		\item Let $[\bq] < m$. For an integer $\varsigma \in [[\bq], m)$, for all $\bvsig, [\bvsig]\leq \varsigma$, 
		\[h^{-[\bvsig]} \lesssim \inf_{\delta\in\Delta}\inf_{\bx\in\clo({\delta})}\|\partial^{\bvsig}\bp(\bx;\Delta,m)\|
		\leq     \sup_{\delta\in\Delta}\sup_{\bx\in\clo({\delta})}\|\partial^{\bvsig}\bp(\bx;\Delta,m)\|
		\lesssim h^{-[\bvsig]}
		\]
		where $\clo({\delta})$ is the closure of $\delta$, and for $[\bvsig]=\varsigma+1$,
		\[\sup_{\delta\in\Delta}\sup_{\bx\in\clo({\delta})}\|\partial^{\bvsig}\bp(\bx;\Delta,m)\|\lesssim h^{-\varsigma-1}.\]
	\end{enumerate}
\end{assumption}

Assumption \ref{assumption: Local Basis} imposes conditions ensuring the stability of the $L_2$ projection operator onto the approximating space. Condition \ref{assumption: Local Basis}(a) requires that each basis function in $\bp(\bx;\Delta, m)$ be supported by a region consisting of a finite number of cells in $\Delta$. Therefore, as $\bar{\kappa}\to\infty$ (and $h\to 0$), each element of $\Delta$ shrinks and all the basis functions are ``locally supported'' relative to the whole support of the data. Another common assumption in least squares regression is that the regressors are not too co-linear: the minimum eigenvalue of $\E[\bp(\bx_i)\bp(\bx_i)']$ is usually assumed to be bounded away from zero. Since the local support condition in Assumption \ref{assumption: Local Basis}(a) implies a banded structure for this matrix, it suffices to require that the basis functions are not too co-linear locally, as stated in Assumption \ref{assumption: Local Basis}(b). These two assumptions are very similar to Conditions A.2 and Conditions A.3 in the Appendix of \cite{Huang_2003_AoS}, and therefore they could also be used to establish theoretical results analogous to those discussed in that appendix (such results are not needed herein because our proofs are different). Finally, Assumption \ref{assumption: Local Basis}(c) controls the magnitude of the local basis in a uniform sense.

Assumptions \ref{assumption: quasi-uniform partition} and \ref{assumption: Local Basis} implicitly relate the number of approximating series terms, the number of knots used and the maximum mesh size: $K\asymp \bar{\kappa}\asymp h^{-d}$. By restricting the growth rate of these tuning parameters, the least squares partitioning-based estimator satisfying the above conditions is well-defined in large samples. We next state a high-level requirement that gives explicit expression of the leading approximation error. For each $\bx\in\mathcal{X}$, let $\delta_{\bx}$ be the element of $\Delta$ whose closure contains $\bx$ and $h_\bx$ be the diameter of this $\delta_{\bx}$.

\begin{assumption}[Approximation Error]\label{assumption: Misspecification Error}
	Let $S \geq m$. For all $\bvsig$ satisfying $[\bvsig]\leq \varsigma$ in Assumption \ref{assumption: Local Basis}, there exists $s^*\in\mathcal{S}_{\Delta,m}$, the linear span of $\bp(\bx;\Delta,m)$, and \[\mathscr{B}_{m,\bvsig}(\bx)=-\sum_{\bm{u}\in\Lambda_m}\partial^{\bm{u}}\mu(\bx)h_\bx^{m-[\bvsig]}B_{\bm{u},\bvsig}(\bx)\]
	such that
	\begin{equation}\label{eq: Approximation Power with Bias}
	\underset{\bx\in\mathcal{X}}{\sup}\;|\partial^{\bvsig} \mu(\bx)-\partial^{\bvsig} s^*(\bx)+\mathscr{B}_{m,\bvsig}(\bx)| \lesssim h^{m+\varrho-[\bvsig]}
	\end{equation}
	and
	\begin{equation}
	\sup_{\delta\in\Delta}\sup_{\bx_1,\bx_2\in\clo({\delta})}\frac{|B_{\bm{u},\bvsig}(\bx_1)-B_{\bm{u},\bvsig }(\bx_2)|}{\|\bx_1-\bx_2\|}\lesssim h^{-1}
	\end{equation}
	where $B_{\bm{u},\bvsig}(\cdot)$ is a known function that is bounded uniformly over $n$, and 
	$\Lambda_m$ is a multi-index set, which depends on the basis, with $[\bm{u}]=m$ for $\bm{u}\in\Lambda_m$.
\end{assumption}

The usual rate-only assumption, $\sup_{\bx\in\mathcal{X}}|\partial^\bq \mu(\bx)-\partial^\bq s^*(\bx)|\lesssim h^{m-[\bq]}$, which is implied by Assumptions \ref{assumption: Misspecification Error}, will not suffice for our bias correction and IMSE expansion results: \eqref{eq: Approximation Power with Bias} is needed. The terms $B_{\bm{u},\bvsig}(\bx)$ in $\mathscr{B}_{m,\bvsig}(\bx)$ are known functions of $\bx$ which depend on the particular partitioning scheme and bases used. The only unknowns are the higher-order derivatives of $\mu(\cdot)$. In the SA (\S SA-6) we verify this (and other assumptions) for splines, wavelets, and piecewise polynomials, including explicit formulas for the leading error in \eqref{eq: Approximation Power with Bias} and precise characterizations of $\Lambda_m$. We assume enough smoothness exists to characterize these terms: see \cite{Calonico-Cattaneo-Farrell_2018_JASA} for a discussion when smoothness constrains inference.

The function $\mathscr{B}_{m,\bvsig}$ is understood as the approximation error in $L_\infty$ norm, and is not in general the misspecification (or smoothing) bias of a series estimator. In least squares series regression settings, the leading smoothing bias is described by two terms in general: $\mathscr{B}_{m,\bvsig}$ and the accompanying error from the linear projection of $\mathscr{B}_{m,\mathbf{0}}$ onto $\mathcal{S}_{\Delta,m}$. We formalize this result in Lemma \ref{lem:bias} below. The second bias term is often ignored because in several cases the leading approximation error $\mathscr{B}_{m,\mathbf{0}}$ is \textit{approximately orthogonal} to $\bp$ with respect to the Lebesgue measure, that is, 
\begin{equation}\label{eq: approx orthogonal}
	\max_{1\leq k\leq K} \int_{\mathcal{H}_k} p_k(\bx;\Delta,m)\mathscr{B}_{m, \mathbf{0}}(\bx)\,d\bx=o(h^{m+d}),
\end{equation}
under Assumptions \ref{assumption: DGP}--\ref{assumption: Misspecification Error}. In some simple cases, \eqref{eq: approx orthogonal} is automatically satisfied if one constructs the leading error based on a basis representing the orthogonal complement of $\mathcal{S}_{\Delta, m}$. When \eqref{eq: approx orthogonal} holds, the leading term in $L_\infty$ approximation error coincides with the leading misspecification (or smoothing) bias of a partitioning-based series estimator. When a stronger quasi-uniformity condition holds (i.e., neighboring cells are of the same size asymptotically), a sufficient condition for \eqref{eq: approx orthogonal} is simply the orthogonality between $B_{\bm{u},\mathbf{0}}$ and $\bp$ in $L_2$ with respect to the Lebesgue measure, for all $\bm{u}\in\Lambda_m$.

For general partitioning-based estimators this orthogonality need not hold. For example, \eqref{eq: approx orthogonal} is hard to verify when the partitioning employed is sufficiently uneven, as is usually the case when employing machine learning methods. All our main results hold when this orthogonality fails, and importantly, our bias correction methods and IMSE expansion explicitly account for the $L_2$ projection of $\mathscr{B}_{m,\mathbf{0}}$ onto the approximating space spanned by $\bp$.


\section{Characterization and Correction of Bias}\label{sec:bias}

We now precisely characterize the bias of $\widehat{\partial^\bq\mu}(\bx)$ under Assumptions \ref{assumption: DGP}--\ref{assumption: Misspecification Error}, but not assuming \eqref{eq: approx orthogonal}. Then, using this result, we develop valid IMSE expansions and three robust bias-corrected inference procedures. This section focuses on bias correction, and Section \ref{sec:pointwise} presents the associated robust Studentization adjustments for inference, following the ideas in \cite{Calonico-Cattaneo-Farrell_2018_JASA} for kernel-based nonparametrics.

Given our assumptions, the estimator $\widehat{\partial^\bq\mu}(\bx)$ of \eqref{equation: partitioning estimator} can be written as
\begin{equation}
\label{eqn:mu hat US}
\widehat{\partial^\bq\mu}_0(\bx) := \bg{0}'\E_n[\bPi{0}(\bx_i) y_i],
\end{equation}
where
\[\widehat{\bgamma}_{\bq,0}(\bx)' := \partial^\bq \bp(\bx)'  \E_n[ \bp(\bx_i) \bp(\bx_i)']^{-1}
  \quad \text{ and } \quad
  \bPi{0}(\bx_i) := \bp(\bx_i).\]
The subscript ``0'' differentiates this estimator from the bias-corrected versions below. We give our first result, proven in the SA, \S SA-10.2.

\begin{lem}[Conditional Bias]\label{lem:bias}
	Let Assumptions \ref{assumption: DGP}, \ref{assumption: quasi-uniform partition}, \ref{assumption: Local Basis}, and \ref{assumption: Misspecification Error} hold. If $\frac{\log n}{nh^d}=o(1)$, then
	\begin{align}
	&\E[\widehat{\partial^\bq\mu}_0(\bx)|\bX ] - \partial^\bq \mu(\bx) \nonumber\\
	& \qquad = \bg{0}'\E_n[\bPi{0}(\bx_i)\mu(\bx_i)]-\partial^\bq \mu(\bx) \label{eqn:pre-bias}\\
	& \qquad = \mathscr{B}_{m,\bq}(\bx) - \bg{0}'\E_n[\bPi{0}(\bx_i)\mathscr{B}_{m,\mathbf{0}}(\bx_i)] + O_\P(h^{m+\varrho-[\bq]}).\label{eqn:asymp-bias}
	\end{align}
\end{lem}

The proof of this lemma generalizes an idea in \cite[Theorem 2.2]{Zhou-Shen-Wolfe_1998_AoS} to handle partitioning-based series estimators beyond the specific example of $B$-Splines on tensor-product partitions. The first component $\mathscr{B}_{m,\bq}(\bx)$ is the leading term in the asymptotic error expansion and depends on the function space generated by the series employed. The second component comes from the least squares regression, and it can be interpreted as the projection of the leading approximation error onto the space spanned by the basis employed. Because the approximating basis $\bp(\bx)$ is locally supported (Assumption \ref{assumption: Local Basis}), the orthogonality condition in \eqref{eq: approx orthogonal}, when it holds, suffices to guarantee that the projection of leading error is of smaller order (such as for $B$-splines on a tensor-product partition). In general the bias will be $O(h^{m-[\bq]})$ and further, in finite samples both terms may be important even if \eqref{eq: approx orthogonal} holds.

We consider three bias correction methods to remove the leading bias terms of Lemma \ref{lem:bias}. All three methods rely, in one way or another, on a higher-order basis: for some $\tilde{m} > m$, let $\bpt(\bx) := \bpt(\bx; \tilde{\Delta}, \tilde{m})$ be a basis of order $\tilde{m}$ defined on a partition $\tilde{\Delta}$ which has maximum mesh $\tilde{h}$. Objects accented with a tilde always pertain to this secondary basis and partition for bias correction. In practice, a simple choice is $\tilde{m}=m+1$ and $\tilde{\Delta}=\Delta$.

The first approach is to use a higher-order basis in place of the original basis. This is thus named \emph{higher-order-basis bias correction} and numbered as approach $j=1$. In complete parallel to \eqref{eqn:mu hat US} define
\begin{equation}\label{eqn:naive bc}
\widehat{\partial^\bq\mu}_1(\bx) := \bg{1}'\E_n[\bPi{1}(\bx_i) y_i ],	
\end{equation}
where
\[\quad \bg{1}' := \partial^\bq \bpt(\bx)'  \E_n[ \bpt(\bx_i) \bpt(\bx_i)']^{-1}
  \quad \text{ and } \quad
  \bPi{1}(\bx_i) := \bpt(\bx_i).\]
This approach can be viewed as a bias correction of the original point estimator because, trivially, $\widehat{\partial^\bq\mu}_1(\bx) = \widehat{\partial^\bq\mu}_0(\bx) - (\widehat{\partial^\bq\mu}_0(\bx) - \widehat{\partial^\bq\mu}_1(\bx))$. Valid inference based on $\widehat{\partial^\bq\mu}_1(\bx)$ can be viewed as ``undersmoothing'' applied to the higher-order point estimator, but is distinct from undersmoothing $\widehat{\partial^\bq\mu}_0(\bx)$ (i.e., using a finer partition $\Delta$ and keeping the order fixed). \cite{Huang_2003_AoS} used this idea to remove the asymptotic bias of splines estimators.

Our second approach makes use of the generic expression of the least squares bias in \eqref{eqn:pre-bias}. The unknown objects in this expression are $\mu$ and $\partial^\bq \mu$, both of which can be estimated using the higher-order estimator \eqref{eqn:naive bc}. By plugging these into \eqref{eqn:pre-bias} and subtracting the result from $\widehat{\partial^\bq\mu}_0(\bx)$, we obtain the \emph{least squares bias correction}, numbered as approach 2:
\begin{align}\label{eqn:generic bc}
\widehat{\partial^\bq\mu}_2(\bx)
&:= \widehat{\partial^\bq\mu}_0(\bx) - \left(\bg{0}'\E_n[\bPi{0}(\bx_i)\widehat{\mu}_1(\bx_i)] - \widehat{\partial^\bq\mu}_1(\bx)\right)\\
&:= \bg{2}'\E_n[\bPi{2}(\bx_i) y_i]\nonumber
\end{align}
where
\begin{equation*}
\begin{split}
&\quad
\widehat{\bgamma}_{\bq,2}(\bx)' := \Big(\bg{0}',-\bg{0}'\E_n[\bp(\bx_i)\bpt(\bx_i)'] \E_n[\bpt(\bx_i)\bpt(\bx_i)']^{-1} + \bg{1}'\Big)\\
&\text{ and } \qquad
\bPi{2}(\bx_i) := \left( \bp(\bx_i)', \bpt(\bx_i)'\right)',
\end{split}
\end{equation*}
which is exactly of the same form as $\widehat{\partial^\bq\mu}_0(\bx)$ and $\widehat{\partial^\bq\mu}_1(\bx)$ (cf., \eqref{eqn:mu hat US} and \eqref{eqn:naive bc}), except for the change in $\widehat{\bgamma}_{\bq,j}(\bx)$ and $\bPi{j}(\bx_i)$.

Finally, approach number 3 targets the leading terms in Equation \eqref{eqn:asymp-bias}. We dub this approach \emph{plug-in bias correction}, as it specifically estimates the leading bias terms, in fixed-$n$ form, of $\widehat{\partial^\bq\mu}_0(\bx)$ according to Assumption \ref{assumption: Misspecification Error}. To be precise, we employ the explicit plug-in bias estimator
\[\widehat{\mathscr{B}}_{m,\bq}(\bx)= -\sum_{\bm{u}\in\Lambda_{m}} \Big(\widehat{\partial^{\bm{u}}\mu}_1(\bx)\Big) h_\bx^{m-[\bq]} B_{\bm{u},\bq}(\bx),
\]
with $[\bq]< m$ and $\Lambda_m$ as in Assumption \ref{assumption: Misspecification Error}, leading to
\begin{align}\label{eqn:asymptotic bc}
\widehat{\partial^\bq\mu}_3(\bx)
&:= \widehat{\partial^\bq\mu}_0(\bx) - \left(\widehat{\mathscr{B}}_{m,\bq}(\bx) - \bg{0}'\E_n[\bPi{0}(\bx_i)\widehat{\mathscr{B}}_{m,\mathbf{0}}(\bx_i)]\right)\\
&:= \bg{3}'\E_n[\bPi{3}(\bx_i) y_i]\nonumber
\end{align}
where
\begin{equation*}
\begin{split}
&\widehat{\bgamma}_{\bq,3}(\bx)' := \Big(\bg{0}', \sum_{\bm{u}\in\Lambda_{m}} \Big\{\widehat{\bgamma}_{\bm{u},1}(\bx)' h_\bx^{m-[\bq]}B_{\bm{u},\bq}(\bx)\\
&\hspace{1.7in} - \widehat{\bgamma}_{\bq,0}(\bx)'\E_n[\bp(\bx_i)h_{\bx_i}^{m}B_{\bm{u},\mathbf{0}}(\bx_i)\widehat{\bgamma}_{\bm{u},1}(\bx_i)']\Big\}\Big),\\
&\text{ and } \qquad \bPi{3}(\bx_i) := \left( \bp(\bx_i)', \bpt(\bx_i)'\right)'.
\end{split}
\end{equation*}
When the orthogonality condition \eqref{eq: approx orthogonal} holds, the second correction term in $\widehat{\partial^\bq\mu}_3(\bx)$ is asymptotically negligible relative to the first. However, in finite samples both terms can be important, so we consider the general case.

Our results employing bias correction will require the following conditions on the higher-order basis used for bias estimation.
 
\begin{assumption}[Bias Correction]\label{assumption: BC bases}
	The partition $\tilde{\Delta}$ satisfies Assumption \ref{assumption: quasi-uniform partition}, with maximum mesh $\tilde{h}$, and the basis $\bpt(\bx;\tilde{\Delta}, \tilde{m})$, $\tilde{m} > m$, satisfies Assumptions \ref{assumption: Local Basis} and \ref{assumption: Misspecification Error} with $\tilde{\varsigma} = \tilde{\varsigma}(\tilde{m}) \geq m$ in place of $\varsigma$. Let $\rho := h/\tilde{h}$, which obeys $\rho \to \rho_0\in(0,\infty)$. In addition, for $j=3$, either (i) $\bpt(\bx;\tilde{\Delta}, \tilde{m})$ spans a space containing the span of $\bp(\bx;\Delta, m)$, and for all $\bm{u}\in\Lambda_{m}$, $\partial^{\bm{u}}\bp(\bx;\Delta, m)=\mathbf{0}$; or (ii) both $\bp(\bx;\Delta, m)$ and $\bpt(\bx;\tilde{\Delta}, \widetilde{m})$ reproduce polynomials of degree 
	$[\bq]$.
\end{assumption}

In addition to removing the leading bias, Assumption \ref{assumption: BC bases} requires that the asymptotic variance of bias-corrected estimators is properly bounded from below in a uniform sense, which is critical for inference. Additional conditions are needed for plug-in bias correction ($j=3$) due to the more complicated covariance between $\widehat{\partial^\bq\mu}_0$ and the estimated leading bias. Orthogonality properties due to the projection structure of the least squares bias correction ($j=2$) removes these ``covariance'' components in the variance of $\widehat{\partial^\bq\mu}_2$. The natural choice of $\tilde{\Delta} = \Delta$ and $\tilde{m} = m+1$ will obey this restriction under intuitive conditions. In the SA, Assumption \ref{assumption: BC bases} is verified for splines, wavelets, and piecewise polynomials (\S SA-6), and we also compare theoretically the alternative bias correction strategies (\S SA-7.2).

\section{IMSE and Convergence Rates}\label{sec:imse}

We establish two main results related to the point estimator $\widehat{\partial^\bq\mu}_0(\bx)$: a valid IMSE expansion for the estimator, which gives as a by-product an estimate of its $L_2$ convergence rate, and its uniform convergence rate.

\subsection{IMSE-Optimal Point Estimation}

We first give a general IMSE approximation, which then is specialized for tensor-product partitions. These expansions are used to obtain optimal choices of partition size from a point estimation perspective.

Our first result holds for any partition $\Delta$ satisfying Assumption \ref{assumption: quasi-uniform partition}.

\begin{thm}[IMSE]\label{thm:IMSE}
	Let Assumptions \ref{assumption: DGP}, \ref{assumption: quasi-uniform partition}, \ref{assumption: Local Basis}, and \ref{assumption: Misspecification Error} hold. If $\frac{\log n}{nh^d}=o(1)$, then for a weighting function $w(\bx)$ that is continuous and bounded away from zero on $\mathcal{X}$,
	\begin{align*}
		&\int_\mathcal{X} \E[(\widehat{\partial^\bq\mu}_0(\bx)-\partial^\bq \mu(\bx))^2|\bX]w(\bx)\,d\bx\\
		&\qquad = \frac{1}{n}\Big(\mathscr{V}_{\Delta,\bq} + o_\P(h^{-d-2[\bq]})\Big) + \Big(\mathscr{B}_{\Delta,\bq} + o_\P(h^{2m-2[\bq]})\Big)
	\end{align*}
	where 
	\begin{align*}
		\mathscr{V}_{\Delta,\bq}
		& = \tr\Big( \bSigma_{0}  \int_\mathcal{X}\bgamma_{\bq, 0}(\bx) \bgamma_{\bq, 0}(\bx)' w(\bx)d\bx\Big)
		\asymp h^{-d-2[\bq]},\\
		\mathscr{B}_{\Delta,\bq}
		& = \int_\mathcal{X} \Big(\mathscr{B}_{m,\bq}(\bx)
		-  \bgamma_{\bq, 0}(\bx)' \E[\bp(\bx_i)\mathscr{B}_{m,\mathbf{0}}(\bx_i)]\Big)^2 w(\bx)d\bx
		\lesssim h^{2m-2[\bq]},
	\end{align*}
	$\bSigma_{0}:= \E[ \bPi{0}(\bx_i)\bPi{0}(\bx_i)' \sigma^2(\bx_i)]$, and $\bgamma_{\bq,0}(\bx)' := \partial^\bq \bp(\bx)'  \E[ \bp(\bx_i) \bp(\bx_i)']^{-1}$.
\end{thm}	

This theorem, proven in the SA, \S SA-10.5, shows that the leading term in the integrated (and pointwise) variance of $\widehat{\partial^\bq\mu}_0(\bx)$ is of order $n^{-1}h^{-d-2[\bq]}$. For the bias term, on the other hand, the theorem  only establishes an upper bound: to bound the bias component from below, stronger conditions on the regression function would be needed. This rate bound is sharp in general.

The quantities $\mathscr{V}_{\Delta,\bq}$ and $\mathscr{B}_{\Delta,\bq}$ are nonrandom sequences depending on the partitioning scheme $\Delta$ in a complicated way, and need not converge as $h\to0$. Nevertheless, when the integrated squared bias does not vanish ($\mathscr{B}_{\Delta,\bq}\neq0$), Theorem \ref{thm:IMSE} implies that the IMSE-optimal mesh size $h_{\mathtt{IMSE}}$ is proportional to $n^{-1/(2m+d)}$, or equivalently, the IMSE-optimal number of series terms $K_{\mathtt{IMSE}}\asymp n^{d/(2m+d)}$. Furthermore, because the IMSE expansion is obtained for a given partition scheme, the result in Theorem \ref{thm:IMSE} can be used to evaluate different partitioning schemes altogether, and to select the ``optimal'' one in an IMSE sense. 

Theorem \ref{thm:IMSE} generalizes prior work substantially. Existing results cover only special cases, such as piecewise polynomials \cite{Cattaneo-Farrell_2013_JoE} or splines \cite{Agarwal-Studden_1980_AoS,Zhou-Shen-Wolfe_1998_AoS,Zhou-Wolfe_2000_SS} on tensor-product partitions only, and often restricting to $d = 1$ or $[\bq] = 0$.

We now consider the special case of a tensor-product partition where the ``tuning parameter'' $\Delta$ reduces to the vector of partitioning knots $\bkappa=(\kappa_1,\dots,\kappa_d)'$, where $\kappa_\ell$ is the number of subintervals used for the $\ell$-th covariate. We assume that $\Delta$ and $\bp(\cdot)$ obey the following regularity conditions, so that the limiting constants in the IMSE approximation can be characterized.

\begin{assumption}[Regularity for Asymptotic IMSE]\label{asmpt:rect}
	Suppose that $\mathcal{X}=\otimes_{\ell=1}^d\mathcal{X}_\ell\subset\mathbb{R}^d$, which is normalized to $[0,1]^d$ without loss of generality, and $\Delta$ is a tensor-product partition. For $\bx \in [0,1]^d$, denote $\delta_\bx=\{t_{\ell ,l_\bx }\leq x_\ell \leq t_{\ell ,l_\bx +1},  1\leq \ell \leq d \}$, where $l_\bx < \kappa_\ell$ (see SA, \S SA-3 for details). Let $\mathbf{b}_\bx = (b_{\bx,1}, \ldots, b_{\bx,d})$ collect the interval lengths $b_{\bx,\ell} = |t_{\ell ,l_\bx +1} - t_{\ell ,l_\bx }|$. In addition:
	
	\begin{enumerate}
		\item For $\ell = 1, \ldots, d$, $\sup_{\bx\in[0,1]^d}|b_{\bx,\ell}-\kappa_\ell^{-1}g_\ell(\bx)^{-1}|=o(\kappa_\ell^{-1})$, where $g_\ell(\cdot)$ is bounded away from zero continuous.
		
		\item For all $\delta\in\Delta$ and $\bm{u}_1,\bm{u}_2\in\Lambda_{m}$, there exist constants  $\eta_{\bm{u}_1,\bm{u}_2,\bq}$ such that
		\[
		\int_{\delta}\frac{h_\bx^{2m-2[\bq]}}{\mathbf{b}_\bx^{\bm{u}_1+\bm{u}_2-2\bq}}B_{\bm{u}_1,\bq}(\bx)B_{\bm{u}_2,\bq}(\bx)\,d\bx
		= \eta_{\bm{u}_1,\bm{u}_2,\bq}\vol(\delta)
		\]
		where $\vol(\delta)$ denotes the volume of $\delta$.
		
		\item There exists a set of points $\{\boldsymbol{\tau}_k\}_{k=1}^K$ such that $\btau_k\in\supp(p_k(\cdot))$ for each $k=1,\ldots, K$, and $\{\btau_k\}_{k=1}^K$ can be assigned into $J+\breve{J} < \infty$ groups such that $\{\boldsymbol{\tau}_{s,k_s}\}_{k_s=1}^{K_s}$, $s=1,\ldots,J+\breve{J}$, $\sum_{s=1}^{J+\breve{J}}K_s=K$, and the following conditions hold: (i) For all $1\leq s\leq J$, $\{\delta_{\boldsymbol{\tau}_{s,k_s}}\}_{k_s=1}^{K_s}$ are pairwise disjoint and $\vol\big([0,1]^d\setminus\bigcup_{k_s=1}^{K_s}\delta_{\boldsymbol{\tau}_{s,k_s}}\big)=o(1)$; and (ii) for all $J+1\leq s\leq J+\breve{J}$, $\vol\big(\bigcup_{k_s=1}^{K_s}\delta_{\boldsymbol{\tau}_{s,k_s}}\big)=o(1)$. 
	\end{enumerate}   
\end{assumption}

Part (a) slightly strengthens the quasi-uniform condition imposed in Assumption \ref{assumption: quasi-uniform partition}, but allows for quite general transformations of the knot location. Part (b) ensures that the ``local'' integral of the product between any two $B_{\bm{u},\bq}(\cdot)$ for $\bm{u}\in\Lambda_m$, which depend on the basis but not $\mu(\bx)$, is proportional to the volume of the cell. The scaling factor is due to the use of the lengths of intervals on each axis (denoted by $\bb_\bx$) to characterize the approximation error for a tensor-product partition, instead of the more general diameter used in Section \ref{sec:setup}. Finally, part (c) describes how the supports of the basis functions cover the whole support of data. Specifically, it requires that the approximating basis $\bp$ can be divided into $J+\breve{J}$ groups. The supports of functions in each of the first $J$ groups constitute ``almost'' complete covers of $\mathcal{X}$. In contrast, the supports of functions in other groups are negligible in terms of volume. In such a case, we refer to $J$ as the number of complete covers generated by the supports of basis functions. For tensor product $B$-splines (with simple knots) and wavelets, each subrectangle in $\Delta$ can be associated with one basis function in $\bp$ and the supports of the remaining functions are asymptotically negligible in terms of volume. Thus, $J=1$ in these two examples. For piecewise polynomials of total order $m$, within each subrectangle the unknown function is approximated by a multivariate polynomial of degree $m-1$, and thus $J=\binom{d+m-1}{m-1}$. This condition is used to ensure that the summation over the number of basis functions converges to a well-defined integral as $K\asymp h^{-d}\to\infty$.

We then have the following result for $\widehat{\mu}_0(\bx)$, proven in the SA, \S SA-10.7.

\begin{thm}[Asymptotic IMSE]\label{thm:Limiting IMSE}
	Suppose that the conditions in Theorem \ref{thm:IMSE} and Assumption \ref{asmpt:rect} hold. Then, for $[\bq]=0$,
	\[\mathscr{V}_{\bkappa,\mathbf{0}} = \Big(\prod_{\ell=1}^d \kappa_\ell\Big) \mathscr{V}_{\mathbf{0}} \; + o(h^{-d}), \qquad
	\mathscr{V}_{\mathbf{0}} =  J \int_{[0,1]^d} \frac{\sigma^2(\bx)}{f(\bx)} \Big(\prod_{\ell=1}^{d} g_\ell(\bx)\Big)w(\bx)\,d\bx ,
	\]
	and, provided that \eqref{eq: approx orthogonal} holds, 
	\[\mathscr{B}_{\bkappa,\mathbf{0}}
	= \sum_{\bm{u}_1,\bm{u}_2\in\Lambda_m}
	\boldsymbol{\kappa}^{-(\bm{u}_1+\bm{u}_2)} \mathscr{B}_{\bm{u}_1,\bm{u}_2,\mathbf{0}} + o(h^{2m}),
	\]
	\[\mathscr{B}_{\bm{u}_1,\bm{u}_2,\mathbf{0}}
	= \eta_{\bm{u}_1,\bm{u}_2,\mathbf{0}} \int_{[0,1]^d}\frac{\partial^{\bm{u}_1}\mu(\bx)\partial^{\bm{u}_2}\mu(\bx)}{\mathbf{g}(\bx)^{\bm{u}_1+\bm{u}_2}}w(\bx)d\bx.
	\]
\end{thm}

The bias approximation requires the approximate orthogonality condition \eqref{eq: approx orthogonal} which is satisfied by $B$-splines, wavelets, and piecewise polynomials. It appears to be an open question whether $\mathscr{V}_{\bkappa,\bq}$ and $\mathscr{B}_{\bkappa,\bq}$ converge to a well-defined limit when general basis functions are considered. \cite{Cattaneo-Farrell_2013_JoE} showed convergence to well defined limits for piecewise polynomials, but their result is not easy to extend to cover other bases functions without imposing $\bq=\mathbf{0}$ and the approximate orthogonality condition \eqref{eq: approx orthogonal}. The SA (\S SA-3) contains more details and other results.

Theorem \ref{thm:Limiting IMSE} justifies the IMSE-optimal choice of number of knots:
\[\bkappa_{\mathtt{IMSE},\mathbf{0}} 
= \argmin_{\bkappa\in\mathbb{Z}^d_{++}} \left\{\frac{1}{n} \Big(\prod_{\ell=1}^d \kappa_\ell\Big) \mathscr{V}_{\mathbf{0}}
+ \sum_{\bm{u}_1,\bm{u}_2\in\Lambda_m}
\boldsymbol{\kappa}^{-(\bm{u}_1+\bm{u}_2)} \mathscr{B}_{\bm{u}_1,\bm{u}_2,\mathbf{0}} \right\},
\]
and, in particular, when the same number of knots is used in all margins,
\[\kappa_{\mathtt{IMSE},\mathbf{0}}
= \bigg\lceil \left(\frac{2m\sum_{\bm{u}_1,\bm{u}_2\in\Lambda_m} \mathscr{B}_{\bm{u}_1,\bm{u}_2,\mathbf{0}}}
{d\mathscr{V}_{\mathbf{0}}}\right)^{\frac{1}{2m+d}} 
n^{\frac{1}{2m+d}} \bigg\rceil
\]

Data-driven versions of this IMSE-optimal choice, and extensions to derivative estimation, are discussed in the SA (\S SA-8) and fully implemented in our companion general-purpose {\sf R} package \texttt{lspartition} \cite{Cattaneo-Farrell-Feng_2019_lspartition}. While beyond the scope of this paper, it would be of interest to study the theoretical properties of cross-validation methods as an alternative way of constructing IMSE-optimal tuning parameter selectors for partitioning-based estimators.

\subsection{Convergence Rates}

Theorem \ref{thm:IMSE} immediately delivers the $L_2$ convergence rate for the point estimator $\widehat{\partial^\bq\mu}_0(\bx)$. For completeness, we also establish its uniform convergence rate. Recall that $\nu>0$.

\begin{thm}[Convergence Rates]\label{thm:Convergence Rates}
	Let Assumptions \ref{assumption: DGP}, \ref{assumption: quasi-uniform partition} and \ref{assumption: Local Basis} hold. Assume also that $\sup_{\bx\in\mathcal{X}}|\partial^\bq \mu(\bx)-\partial^\bq s^*(\bx)|\lesssim h^{m-[\bq]}$ with $s^*$ defined in Assumption \ref{assumption: Misspecification Error}. Then, if $\frac{\log n}{nh^d}=o(1)$,
	\[\int_\mathcal{X} \Big(\widehat{\partial^\bq\mu}_0(\bx)-\partial^\bq \mu(\bx)\Big)^2w(\bx)\,d\bx
	  \lesssim_\P \frac{1}{nh^{d+2[\bq]}} + h^{2(m-[\bq])}
	\]
	If, in addition,
	\begin{description}
		\item[(i)] $\sup_{\bx\in \mathcal{X}}\E[|\varepsilon_i|^3\exp(|\varepsilon_i|)]<\infty$ and $\frac{(\log n)^3}{n h^{d}}\lesssim 1$, or
		\item[(ii)] $\sup_{\bx\in \mathcal{X}}\E[|\varepsilon_i|^{2+\nu}]<\infty$ and $\frac{n^{\frac{2}{2+\nu}}(\log n)^{\frac{2\nu}{4+2\nu}}}{nh^d}\lesssim 1$,
	\end{description}
	then
	\[\sup_{\bx\in\mathcal{X}}\,\Big|\widehat{\partial^\bq\mu}_0(\bx)-\partial^\bq \mu(\bx)\Big|^2
	  \lesssim_\P \frac{\log n}{nh^{d+2[\bq]}} + h^{2(m-[\bq])}.
	\]
\end{thm}

This theorem, proven in the SA, \S SA-10.11, shows that the partitioning-based estimators can attain the optimal $L_2$ and uniform convergence rate \cite{Stone_1982_AoS} by proper choice of partitioning scheme, under our high-level assumptions. (The full force of Assumption \ref{assumption: Misspecification Error} is not needed for this result.) \cite{Cattaneo-Farrell_2013_JoE} were the first to show existence of a series estimator (in particular, piecewise polynomials) attaining the optimal uniform convergence rate, a result that was later generalized to other series estimators in \cite{Belloni-Chernozhukov-Chetverikov-Kato_2015_JoE,Chen-Christensen_2015_JOE}.

\section{Pointwise Inference}\label{sec:pointwise}

We give pointwise inference based on classical undersmoothing and all three bias correction methods. All four point estimators take the form $\widehat{\partial^\bq\mu}_j(\bx) = \bg{j}'\E_n[\bPi{j}(\bx_i) y_i]$, where $j=0$ corresponds to the conventional partitioning estimator, and $j=1,2,3$ refer to the three distinct bias correction strategies. Infeasible inference would be based on the standardized $t$-statistics
\[T_j(\bx) = \frac{ \widehat{\partial^\bq\mu}_j(\bx) - \partial^\bq\mu(\bx)}{ \sqrt{\Omega_j(\bx) /n}}, \qquad
  \Omega_j(\bx) = \bgamma_{\bq,j}(\bx)' \bSigma_{j} \bgamma_{\bq,j}(\bx),
\]
where, for each $j=0,1,2,3$, $\bgamma_{\bq,j}(\bx)$ are defined as $\widehat{\bgamma}_{\bq,j}$ in \eqref{eqn:mu hat US}, \eqref{eqn:naive bc}, \eqref{eqn:generic bc}, and \eqref{eqn:asymptotic bc}, respectively, but with sample averages and other estimators replaced by their population counterparts, and $\bSigma_{j} := \E[ \bPi{j}(\bx_i)\bPi{j}(\bx_i)' \sigma^2(\bx_i)]$. These $t$-statistics are infeasible, but they nonetheless capture the additional variability introduced by the bias correction approach when $j=1,2,3$, the key idea behind robust bias-corrected inference \cite{Calonico-Cattaneo-Titiunik_2014_ECMA,Calonico-Cattaneo-Farrell_2018_JASA}. We also discuss below Studentization, that is, replacing $\Omega_j(\bx)$ with a consistent estimator.

\subsection{Distributional Approximation} Our first result establishes the limiting distribution of the standardized $t$-statistics $T_j(\bx)$.

\begin{thm}[Asymptotic Normality]\label{thm:pointwise inference}
	Let Assumptions \ref{assumption: DGP}, \ref{assumption: quasi-uniform partition}, \ref{assumption: Local Basis}, and \ref{assumption: Misspecification Error} hold. Assume $\sup_{\bx\in\mathcal{X}} \E[\e_i^2 \I\{|\e_i|>M\}|\bx_i=\bx]\to 0$ as $M\to\infty$, and $\frac{\log n}{nh^d}=o(1)$. Furthermore, for $j=0$, assume $nh^{2m+d}=o(1)$; and for $j=1,2,3$, assume Assumption \ref{assumption: BC bases} holds and $nh^{2m+d}\lesssim 1$.
	
	Then, for each $j=0,1,2,3$ and $\bx\in\mathcal{X}$, $\sup_{u\in\mathbb{R}} |\P[T_j(\bx)\leq u] - \Phi(u) | = o(1)$, where $\Phi(u)$ denotes the cumulative distribution function of $\mathsf{N}(0,1)$.
\end{thm}

This theorem, proven in the SA, \S SA-10.9, gives a valid Gaussian approximation for the $t$-statistics $T_j(\bx)$, pointwise in $\bx\in\mathcal{X}$. The regularity conditions imposed are extremely mild, and in perfect quantitative agreement with those used in \cite{Belloni-Chernozhukov-Chetverikov-Kato_2015_JoE} for $j=0$ (undersmoothing). For $j=1,2,3$ (robust bias correction), the result is new and the restrictions are in qualitative agreement with those obtained for kernel-based nonparametrics.

\subsection{Implementation}

To make the results in Theorem \ref{thm:pointwise inference} feasible, we replace $\Omega_j(\bx)$ with a consistent estimator. Specifically, we consider the four feasible $t$-statistics, $j=0,1,2,3$,
\begin{equation}
\begin{split}\label{eqn:t stat}
& \widehat{T}_j(\bx) = \frac{ \widehat{\partial^\bq\mu}_j(\bx) - \partial^\bq\mu(\bx)}{ \sqrt{ \bOhat{j} /n}}  ,   	\qquad \qquad  	  \bOhat{j} =  \bg{j}' \bSighat{j} \bg{j}   ,   		\\
& \widehat{\bSigma}_j = \E_n[ \bPi{j}(\bx_i)\bPi{j}(\bx_i)' \widehat{\e}_{i,j}^2]  ,    	\qquad \qquad  	  \widehat{\e}_{i,j} = y_i - \widehat{\mu}_j(\bx_i),
\end{split}
\end{equation}
Once the basis functions and partitioning schemes are chosen, the statistic $\widehat{T}_j(\bx)$ is readily implementable. The following theorem gives sufficient conditions for valid pointwise inference.

\begin{thm}[Variance Consistency]\label{thm:variance consistency}
	Let Assumptions \ref{assumption: DGP}, \ref{assumption: quasi-uniform partition}, \ref{assumption: Local Basis}, and \ref{assumption: Misspecification Error} hold. If $j=1,2,3$, also let Assumption \ref{assumption: BC bases} hold. In addition, assume one of the following holds:
	\begin{description}
		\item[(i)] $\sup_{\bx\in \mathcal{X}}\E[|\varepsilon_i|^{2+\nu}]<\infty$ and $\frac{n^{\frac{2}{2+\nu}}(\log n)^{\frac{2\nu}{4+2\nu}}}{nh^d}=o(1)$, or
		\item[(ii)] $\sup_{\bx\in \mathcal{X}}\E[|\varepsilon_i|^3\exp(|\varepsilon_i|)]<\infty$ and $\frac{(\log n)^3}{n h^{d}}=o(1)$.
	\end{description}
	 Then, for each $j=0,1,2,3$, $|\bOhat{j} - \bO{j} | = o_\P(h^{-d-2[\bq]})$.
\end{thm}

This result, proven in the SA, \S SA-10.12, together with Theorem \ref{thm:pointwise inference}, delivers feasible inference. Valid $100(1-\alpha)\%$, $\alpha\in(0,1)$, confidence intervals for $\partial^\bq\mu(\bx)$ are formed in the usual way:
\[ \Big[ \; \widehat{\partial^\bq\mu}_j(\bx) \pm \Phi^{-1}(1-\alpha/2) \cdot \sqrt{\bOhat{j}/n} \; \Big], \qquad j=0,1,2,3.\]
For $j=1,2,3$, the IMSE-optimal partitioning scheme choice derived in Section \ref{sec:imse} can be used directly, while for $j=0$ the partitioning has to be undersmoothed (i.e., made finer than the IMSE-optimal choice) in order to obtain valid confidence intervals. Our results generalize, under weaker conditions, prior work on univariate regression splines \cite{Zhou-Shen-Wolfe_1998_AoS,Zhou-Wolfe_2000_SS,Huang_2003_AoS}.

\section{Uniform Inference}\label{sec:coupling}

We next give a valid distributional approximation for the \emph{whole} process $\{\widehat{T}_j(\bx): \bx \in \mathcal{X}\}$, for each $j=0,1,2,3$. We establish this approximation using two distinct coupling strategies. We then propose a simulation-based feasible implementation of the result. We close by applying our results to construct valid confidence bands for $\partial^\bq \mu(\cdot)$.

\subsection{Strong Approximations}

The stochastic processes $\{\widehat{T}_j(\bx): \bx \in \mathcal{X}\}$ are not asymptotically tight, and therefore do not converge weakly in $\mathcal{L}^\infty(\mathcal{X})$, where $\mathcal{L}^\infty(\mathcal{X})$ denotes the set of all (uniformly) bounded real functions on $\mathcal{X}$ equipped with uniform norm. Nevertheless, their finite sample distribution can be approximated by carefully constructed Gaussian processes (in a possibly enlarged probability space).

We first employ the following lemma to simplify the problem. Recall that $r_n$ is some non-vanishing positive sequence and $\nu>0$.

\begin{lem}[Hats Off]\label{lem:Strong Approximation Hats Off}
	Let Assumptions \ref{assumption: DGP}, \ref{assumption: quasi-uniform partition}, \ref{assumption: Local Basis}, and \ref{assumption: Misspecification Error} hold. Assume one of the following holds:
	\begin{description}
		\item[(i)] $\sup_{\bx\in \mathcal{X}}\E[|\varepsilon_i|^{2+\nu}|\bx_i=\bx]<\infty$ and $\frac{n^{\frac{2}{2+\nu}}(\log n)^{\frac{2+2\nu}{2+\nu}}}{nh^d} = o(r_n^{-2})$; or
		\item[(ii)] $\sup_{\bx\in \mathcal{X}}\E[|\varepsilon_i|^3\exp(|\varepsilon_i|)|\bx_i=\bx]<\infty$ and $\frac{(\log n)^4}{n h^{d}} = o(r_n^{-2})$.
	\end{description}
	Furthermore, if $j=0$, assume $nh^{d+2m}=o(r_n^{-2})$; and, if $j=1,2,3$, assume Assumption \ref{assumption: BC bases} holds and $nh^{d+2m+2\varrho} = o(r_n^{-2})$. Then
	\begin{equation*}\label{eq:T-hatsoff}
	\sup_{\bx\in\mathcal{X}} \Big| \widehat{T}_j(\bx)  - t_j(\bx) \Big| = o_\P(r_n^{-1}), \quad
	t_j(\bx) = \frac{\bgamma_{\bq,j}(\bx)'}{ \sqrt{ \Omega_{j}(\bx)}}
	\G_n[\bPi{j}(\bx_i)\e_i].
	\end{equation*}
\end{lem}

Lemma \ref{lem:Strong Approximation Hats Off} requires that the estimation and sampling uncertainty of $\widehat{\bgamma}_{\bq,j}$ and $\bOhat{j}$, as well as the smoothing bias of $\widehat{\partial^\bq\mu}_j(\bx)$, be negligible uniformly over $\bx\in\mathcal{X}$. The proof is given in \S\ref{proof:lem:Strong Approximation Hats Off} and relies on a new technical lemma stated in \S\ref{proof:techlemmas}. This technical approximation step allows us to focus on developing a distributional approximation for the infeasible stochastic processes $\{t_j(\bx): \bx \in \mathcal{X}\}$, $j=0,1,2,3$. We make precise our uniform distributional approximation in the following definition.

\begin{definition}[Strong Approximation]\label{def: strong approx}
	For each $j=0,1,2,3$, the law of the stochastic process $\{t_j(\bx), \bx \in \mathcal{X}\}$ is approximated by that of a Gaussian process $\{Z_j(\bx), \bx \in \mathcal{X}\}$ in $\mathcal{L}^\infty(\mathcal{X})$ if the following condition holds: in a sufficiently rich probability space, there exists a copy $t'_j(\cdot)$ of $t_j(\cdot)$ and a standard Normal random vector $\bN_{K_j}\thicksim \mathsf{N}(\mathbf{0},\bI_{K_j})$ with $K_j=\dim(\bPi{j}(\bx))$ such that
	\[\sup_{\bx\in\mathcal{X}} \Big| t'_j(\bx) - Z_j(\bx) \Big| = o_\P(r_n^{-1}), \qquad
	  Z_j(\bx) = \frac{\bgamma_{\bq,j}(\bx)' \bSigma_{j}^{1/2}}{ \sqrt{ \Omega_{j}(\bx)}}  \bN_{K_j}.
	\]
	This approximation is denoted by $t_j(\cdot) =_{d} Z_j(\cdot) + o_{\P}(r_n^{-1})$ in $\mathcal{L}^\infty(\mathcal{X})$.
\end{definition}

This definition gives the precise meaning of uniform distributional approximation of $t_j(\cdot)$ by a Gaussian process $Z_j(\cdot)$, and also provides the explicit characterization of such Gaussian process. We establish this strong approximation in two distinct ways. For $d=1$, we develop a novel two-step coupling approach based on the classical Koml\'os-Major-Tusn\'ady (KMT) construction \cite{Komlos-Major-Tusnady_1975_Zeitsch,Komlos-Major-Tusnady_1976_Zeitsch}. For $d>1$, we apply an improved version of the classical Yurinskii construction \cite{Yurinskii_1978_TPA}.

\subsubsection{Unidimensional Regressor}

Let $d=1$. The following theorem gives a valid distributional approximation for $\{\widehat{T}_j(x): x\in\mathcal{X}\}$ using the Gaussian process $\{Z_j(x): x\in\mathcal{X}\}$, for $j=0,1,2,3$, in the sense of Definition \ref{def: strong approx}.

\begin{thm}[Strong Approximation: KMT]\label{thm:strong approximation KMT}
	Let the assumptions and conditions of Lemma \ref{lem:Strong Approximation Hats Off} hold with $d=1$. If $j=2,3$, also assume 
	$\frac{(\log n)^{3/2}}{\sqrt{nh}} = o(r_n^{-2})$. Then, for each $j=0,1,2,3$, $t_j(\cdot) =_d Z_j(\cdot) + o_\P(r_n^{-1})$ in $\mathcal{L}^\infty(\mathcal{X})$, where $Z_j(\cdot)$ is given in Definition \ref{def: strong approx}.
\end{thm}

The proof of this result, in \S\ref{proof:thm:strong approximation KMT}, employs a two-step coupling approach:

\begin{description}
	\item[Step 1.] On a sufficiently rich probability space, there exists a copy $t'_j(\cdot)$ of $t_j(\cdot)$, and an i.i.d.\ sequence $\{\zeta_i: 1\leq i \leq n\}$ of standard Normal random variables, such that
	\begin{equation*}\label{eq:Z-condcoupling}
	\sup_{x\in\mathcal{X}} \Big| t'_j(x)-z_j(x) \Big|= o_\P(r_n^{-1}), \quad
	z_j(x) = \frac{\bgamma_{\bq,j}(x)'}{ \sqrt{ \Omega_{j}(x)}} \G_n[\bPi{j}(x_i)\sigma(x_i)\zeta_i].
	\end{equation*}

	\item[Step 2.] On a sufficiently rich probability space, there exists a copy $z'_j(\cdot)$ of $z_j(\cdot)$, and the standard Normal random vector $\bN_{K_j}$ from Definition \ref{def: strong approx} such that $z'_j(\cdot) =_d \bar{Z}_j(\cdot)$ conditional on $\bX$, where
	\[\bar{Z}_j(x) = \frac{\bgamma_{\bq,j}(x)'\bSigbar{j}^{1/2}}{ \sqrt{ \Omega_{j}(x)}} \bN_{K_j}, \quad
	\bSigbar{j} := \E_n[ \bPi{j}(x_i) \bPi{j}(x_i)' \sigma^2(x_i)],
	\]
	and
	\begin{equation*}\label{eq:Z-coupling}
	\sup_{x\in\mathcal{X}} \Big| \bar{Z}_j(x)-Z_j(x) \Big|= o_\P(r_n^{-1}).
	\end{equation*}		
\end{description}

These two steps summarize our strategy for constructing the unconditionally Gaussian process $\{Z_j(x), x \in \mathcal{X}\}$ approximating the distribution of the whole $t$-statistic processes $\{t_j(x):x\in\mathcal{X}\}$: we first couple $t_j(\cdot)$ to the process $z_j(\cdot)$, which is Gaussian only conditionally on $\mathbf{X}$ but not unconditionally (Step 1), and we then show that the unconditionally Gaussian process  $Z_j(\cdot)$ approximates the distribution of $z_j(\cdot)$ (Step 2). 

To complete the first coupling step, we employ a version of the classical KMT inequalities that applies to independent but non-identically distributed random variables \cite{Sakhanenko_1985_Advances, Sakhanenko_1991_SAM}. We do this because the processes $\{t_j(x):x\in\mathcal{X}\}$ are characterized by a sum of independent but not identically distributed random variables conditional on $\bX$. This part of our proof is inspired by, but is distinct from, the one given in \cite[Chapter 22]{Eggermont-LaRiccia_2009_Book}, where a conditional strong approximation for smoothing splines is established. Our proof relies instead on the new general coupling Lemma \ref{lem:general cond coupling} in \S\ref{proof:thm:strong approximation KMT}.

The intermediate coupling result in Step 1 has the obvious drawback that the process $\{z_j(x):x\in\mathcal{X}\}$ is Gaussian only conditionally on $\mathbf{X}$ but not unconditionally. Step 2 addresses this shortcoming by establishing an unconditional coupling, that is, approximating the distribution of the stochastic process $z_j(\cdot)$ by that of the (unconditional) Gaussian process $Z_j(\cdot)$. As shown in \S \ref{proof:thm:strong approximation KMT}, verifying the second coupling step boils down to controlling the supremum of a Gaussian random vector of increasing dimension, and in particular the crux is to prove precise (rate) control on $\big\| \bSigbar{j}^{1/2} - \bSig{j}^{1/2} \big\|$, $ j=0,1,2,3$. Both $\bSigbar{j}$ and $\bSig{j}$ are symmetric and positive \emph{semi}-definite. Further, for $j=0,1$, $\lambda_{\min}(\bSigma_{j}) \gtrsim h^{d}$ for generic partitioning-based estimators under our assumptions, and therefore we use the bound 
\begin{equation}
\label{eq: sqrt-matrix-eigenval}
\|\mathbf{A}_1^{1/2}-\mathbf{A}_2^{1/2}\|\leq\lambda_{\min}(\mathbf{A}_2)^{-1/2}\|\mathbf{A}_1-\mathbf{A}_2\|,
\end{equation}
which holds for symmetric positive semi-definite $\mathbf{A}_1$ and symmetric positive definite $\mathbf{A}_2$ \cite[Theorem X.3.8]{Bhatia_2013_Book}. Using this bound we obtain unconditional coupling from conditional coupling without additional rate restrictions.

However, for $j=2,3$ the bound \eqref{eq: sqrt-matrix-eigenval} cannot be used in general because $\bp$ and $\bpt$ are typically not linearly independent, and hence $\bSig{j}$ will be singular. To circumvent this problem, we employ the weaker bound \cite[Theorem X.1.1]{Bhatia_2013_Book}: if $\mathbf{A}_1$ and $\mathbf{A}_2$ are symmetric positive semi-definite matrices, then
\begin{equation}\label{eq: sqrt-matrix-generic}
\|\mathbf{A}_1^{1/2}-\mathbf{A}_2^{1/2}\|\leq\|\mathbf{A}_1-\mathbf{A}_2\|^{1/2}.
\end{equation}
This bound can be used for any partitioning-based estimator, with or without bias correction, at the cost of slowing the approximation error rate $r_n$ when constructing the unconditional coupling, and hence leading to the stronger side rate condition as shown in the Theorem \ref{thm:strong approximation KMT} below. When $r_n=1$, there is no rate penalty, while the penalty is only in terms of $\log n$ terms when $r_n=\sqrt{\log n}$ (as in Theorem \ref{thm:Confidence Bands} further below). Furthermore, for certain partitioning-based series estimators it is still possible to use \eqref{eq: sqrt-matrix-eigenval} even when $j=2,3$, as the following remark discusses.

\begin{remark}[Square-root Convergence and Improved Rates]\label{rmk: sqrt-rates}
	The additional restriction imposed in Theorem \ref{thm:strong approximation KMT} for $j=2,3$, that $(\log n)^{3/2}/\sqrt{nh} = o(r_n^{-2})$, can be dropped in some special cases. For some bases it is possible to find a transformation matrix $\bUpsilon$, with $\|\bUpsilon\|_\infty\lesssim 1$, and a basis $\bm{\check{\bp}}$, which obeys Assumption \ref{assumption: Local Basis}, such that $( \bp(\cdot)', \bpt(\cdot)')' =  \bUpsilon\bm{\check{\bp}}(\cdot)$. In other words, the two bases $\bp$ and $\bpt$ can be expressed in terms of another basis $\bm{\check{\bp}}$ without linear dependence. Then, a positive lower bound holds for $\lambda_{\min}(\bSigma_{j}), j=2,3$, implying that the bound \eqref{eq: sqrt-matrix-eigenval} can be used instead of \eqref{eq: sqrt-matrix-generic}. For example, for piecewise polynomials and $B$-splines with equal knot placements for $\bp$ and $\bpt$, a natural choice of $\bm{\check{\bp}}$ is simply a higher-order  polynomial basis on the same partition. Since each function in $\bp$ and $\bm{\check{\bp}}$ is a polynomial on each $\delta\in\Delta$ and nonzero on a fixed number of cells, the ``local representation'' condition $\|\bUpsilon\|_\infty\lesssim 1$ automatically holds. See the SA (\S SA-6) for more details.
\end{remark}

An alternative unconditional strong approximation for general series estimators was obtained by \cite{Belloni-Chernozhukov-Chetverikov-Kato_2015_JoE} for the case of undersmoothing inference ($j=0$). Their proof employs the classical Yurinskii's coupling inequality that controls the convergence rate of partial sums in terms of Euclidean norm, leading to the rate restriction $r_n^6K^5/n\to0$, up to $\log n$ terms, which does not depend on $\nu>0$. In contrast, Theorem \ref{thm:strong approximation KMT} employs a (conditional) KMT-type coupling and then a second (unconditional) coupling approximation, and make use of the banded structure of the Gram matrix formed by local bases, to obtain weaker restrictions. Under bounded polynomial moments, we require only $r_n^6K^3/n^{3\nu/(2+\nu)}\to0$, up to $\log n$ terms. For example, when $\nu=2$ and $r_n=\sqrt{\log n}$ this translates to $K^2/n\to0$, up to $\log n$ terms, which is weaker than previous results in the literature. Under the  sub-exponential conditional moment restriction, the rate condition can be relaxed all the way to $K/n\to0$, up to $\log n$ terms, which appears to be a minimal condition. This is for the entire $t$-statistic process. In addition, Theorem \ref{thm:strong approximation KMT} gives novel strong approximation results for robust bias-corrected $t$-statistic processes.

\begin{remark}[Strong Approximation: KMT for Haar Basis]\label{rmk:KMT Haar}
	Our two-step coupling approach builds on the new coupling Lemma \ref{lem:general cond coupling}, which appears to be hard to extend to $d>1$, except for the important special case the undersmoothed ($j=0$) $t$-statistic process $\{\widehat{T}_0(x): x\in\mathcal{X}\}$ constructed using Haar basis, which is a spline, wavelet, and piecewise polynomial with $m=1$. In the SA, \S SA-5.1, we establish $t_0(\cdot) =_d Z_0(\cdot) + o_\P(r_n^{-1})$ in $\mathcal{L}^\infty(\mathcal{X})$ for any $d\geq1$ under the same conditions of Lemma \ref{lem:Strong Approximation Hats Off}.	
\end{remark}

\subsubsection{Multidimensional Regressors}

Let $d\geq1$. The method of proof employed to establish Theorem \ref{thm:strong approximation KMT} does not extend easily to multivariate regressors ($d>1$) in general. Therefore, we present an alternative strong approximation result based on an improved version of the classical Yurinskii's coupling inequality, recently developed by \cite{Belloni-Chernozhukov-Chetverikov-FernandezVal_2019_JoE}.

\begin{thm}[Strong Approximation: Yurinskii]\label{thm:strong approximation Yurinskii}
	Let the assumptions and conditions of Lemma \ref{lem:Strong Approximation Hats Off} hold. Furthermore, assume $\nu\geq1$ and $\frac{(\log n)^4}{nh^{3d}} = o(r_n^{-6})$. Then, for each $j=0,1,2,3$, $t_j(\cdot) =_d Z_j(\cdot) + o_\P(r_n^{-1})$ in $\mathcal{L}^\infty(\mathcal{X})$, where $Z_j(\cdot)$ is given in Definition \ref{def: strong approx}.
\end{thm}

This strong approximation result, proven in \S\ref{proof:thm:strong approximation Yurinskii}, does not have optimal (i.e. minimal) restrictions, but nonetheless improves on previous results by exploiting the specific structure of the partitioning-based estimators, while also allowing for any $d\geq1$. Specifically, the result sets $\nu=1$ and requires $r_n^6K^3/n\to0$, up to $\log n$ terms. While not optimal when $\nu>3$ (see Remark \ref{rmk:KMT Haar} for a counterexample), it improves on the condition $r_n^6K^5/n\to0$, up to $\log n$ terms, mentioned previously. In addition, Theorem \ref{thm:strong approximation Yurinskii} gives novel strong approximation results for robust bias-corrected $t$-statistic processes.

\subsection{Implementation}

We present a simple plug-in approach that gives a (feasible) approximation to the infeasible standardized Gaussian processes $\{Z_j(\bx):\bx\in\mathcal{X}\}$, in order to conduct inference using the results in Theorem \ref{thm:strong approximation KMT} or Theorem \ref{thm:strong approximation Yurinskii}. In the SA (\S SA-5.2), we also give another plug-in approach and one based on the wild bootstrap. The following definition gives a precise description of how the approximation works.

\begin{definition}[Simulation-Based Strong Approximation]\label{definition: strong approx simul}
	Let $\P^*[\cdot]=\P[\cdot|\by, \bX]$ denote the probability operator conditional on the data. For each $j=0,1,2,3$, the law of the Gaussian process $\{Z_j(\bx): \bx \in \mathcal{X}\}$ is approximated by a (feasible) Gaussian process $\{\widehat{Z}_j(\bx):\bx\in\mathcal{X}\}$, with known distribution conditional on the data $(\by, \bX)$, in $\mathcal{L}^\infty(\mathcal{X})$, if the following condition holds: on a sufficiently rich probability space there exists a copy $\widehat{Z}'_j(\cdot)$ of $\widehat{Z}_j(\cdot)$ such that $\widehat{Z}'_j(\cdot)=_d Z_j(\cdot)$ conditional on the data, and
	\[\P^*\left[\sup_{\bx\in\mathcal{X}} |\widehat{Z}'_j(\bx) - Z_j(\bx)| \geq \eta r_n^{-1} \right] = o_\P(1), \qquad \forall \eta>0,\]
	where, for a $\bN_{K_j}\thicksim \mathsf{N}(\mathbf{0},\bI_{K_j})$ with $K_j=\dim(\bPi{j}(\bx))$,
	\begin{equation*}\label{eq:Z-coupling-hat-unconditional}
	\widehat{Z}_j(\bx) = \frac{\bg{j}' \bSighat{j}^{1/2}}{ \sqrt{ \bOhat{j}} } \bN_{K_j},
	\qquad \bx\in\mathcal{X}, \qquad j=0,1,2,3.
	\end{equation*}
	This approximation is denoted by $\widehat{Z}_j(\cdot) =_{d^*} Z_j(\cdot) + o_{\P^*}(r_n^{-1})$ in $\mathcal{L}^\infty(\mathcal{X})$.
\end{definition}

From a practical perspective, Definition \ref{definition: strong approx simul} implies that sampling from $\widehat{Z}_j(\cdot)$, conditional on the data, is possible and provides a valid distributional approximation of $Z_j(\cdot)$, for each $j=0,1,2,3$. The feasible process $\widehat{Z}_j(\cdot)$ relies on a direct plug-in approach, where all the unknown quantities are replaced by consistent estimators already used in forming $\widehat{T}_j(\bx)$. Resampling is from a multivariate standard Gaussian of dimension $K_j$, not $n$. 

\begin{thm}[Plug-in Approximation]\label{thm:Plug-in Aprox Unconditional}
	Let the assumptions and conditions of Lemma \ref{lem:Strong Approximation Hats Off} hold. Furthermore, for $j=2,3$:
	\begin{description}
		\item[(i)] when $\sup_{\bx\in \mathcal{X}}\E[|\varepsilon_i|^{2+\nu}|\bx_i=\bx]<\infty$, assume $\frac{n^{\frac{1}{2+\nu}}(\log n)^{\frac{4+3\nu}{4+2\nu}}}{\sqrt{nh^d}} = o(r_n^{-2})$; or
		\item[(ii)] when $\sup_{\bx\in \mathcal{X}}\E[|\varepsilon_i|^3\exp(|\varepsilon_i|)|\bx_i=\bx]<\infty$, assume $\frac{(\log n)^{5/2}}{\sqrt{n h^{d}}} = o(r_n^{-2})$.
	\end{description}
	Then, for each $j=0,1,2,3$, $\widehat{Z}_j(\cdot) =_{d^*} Z_j(\cdot) + o_{\P^*}(r_n^{-1})$ in $\mathcal{L}^\infty(\mathcal{X})$, where $\widehat{Z}_j(\cdot)$ is given in Definition \ref{definition: strong approx simul}.
\end{thm}

This result, proven in \S\ref{proof:thm:Plug-in Aprox Unconditional}, strengthens the rate condition for $j=2,3$ compared to Theorems \ref{thm:strong approximation KMT} ($d=1$) and \ref{thm:strong approximation Yurinskii} ($d\geq1$) only by logarithmic factors when $r_n=\sqrt{\log n}$. Moreover, if the structure discussed in Remark \ref{rmk: sqrt-rates} holds, then this additional condition can be dropped.

\subsection{Application: Confidence Bands}
\label{sec:bands}

A natural application of Theorems \ref{thm:strong approximation KMT}, \ref{thm:strong approximation Yurinskii} and \ref{thm:Plug-in Aprox Unconditional} is to construct confidence bands for the regression function or its derivatives. Specifically, for $j=0,1,2,3$ and $\alpha\in(0,1)$, we seek a quantile $q_j(\alpha)$ such that
\[\P\left[\sup_{\bx\in\mathcal{X}} |\widehat{T}_j(\bx)| \leq q_j(\alpha) \right] =1-\alpha + o(1),
\]
which then can be used to construct uniform $100(1-\alpha)$-percent confidence bands for $\partial^\bq \mu(\bx)$ of the form
\[\Big[\widehat{\partial^\bq \mu}_j(\bx)
\pm q_j(\alpha)\sqrt{\bOhat{j}/n}
\;:\; \bx\in\mathcal{X} 
\Big].
\]

The following theorem, proven in \S\ref{proof:thm:Confidence Bands}, establishes a valid distributional approximation for the suprema of the $t$-statistic processes $\{\widehat{T}_j(\bx):\bx\in\mathcal{X}\}$ using \cite[Lemma 2.4]{Chernozhukov-Chetverikov-Kato_2014b_AoS} to convert our strong approximation results into convergence of distribution functions in terms of Kolmogorov distance.

\begin{thm}[Confidence Bands]\label{thm:Confidence Bands}
	Let the conditions of Theorem \ref{thm:strong approximation KMT} or Theorem \ref{thm:strong approximation Yurinskii} hold with $r_n=\sqrt{\log n}$. 
	If the corresponding conditions of Theorem \ref{thm:Plug-in Aprox Unconditional} hold for each $j=0,1,2,3$, then 
	\[\sup_{u\in\mathbb{R}} \left|\P\left[ \sup_{\bx\in \mathcal{X}} |\widehat{T}_j(\bx)| \leq u \right]
	- \P^*\left[ \sup_{\bx\in \mathcal{X}} |\widehat{Z}_j(\bx)| \leq u \right]
	\right| = o_\P(1).
	\]
\end{thm}

\cite{Chernozhukov-Chetverikov-Kato_2014a_AoS,Chernozhukov-Chetverikov-Kato_2014b_AoS} recently showed that if one is only interested in the supremum of an empirical process rather than the \textit{whole} process, then the sufficient conditions for distributional approximation could be weakened compared to earlier literature. Their result applied Stein's method for Normal approximation to show that suprema of general empirical processes can be approximated by a sequence of suprema of Gaussian processes, under the usual undersmoothing conditions (i.e., $j=0$). They illustrate their general results by considering $t$-statistic processes for both kernel-based and series-based nonparametric regression: \cite[Remark 3.5]{Chernozhukov-Chetverikov-Kato_2014b_AoS} establishes a result analogous to Theorem \ref{thm:Confidence Bands} under the side rate condition $K/n^{1-2/(2+\nu)}=o(1)$, up to $\log n$ terms (with $q=2+\nu$ in their notation). In comparison, our result for $j=0$ and $d=1$ in Theorem \ref{thm:Confidence Bands}, under the same moment conditions, requires exactly the same side condition, up to $\log n$ terms. Theorems \ref{thm:strong approximation KMT} and \ref{thm:Confidence Bands} show that the \emph{whole} $t$-statistic process for partitioning-based series estimators, and not just the suprema thereof, can be approximated under the same weak conditions when $d=1$. The same result holds for sub-exponential moments, where the rate condition becomes minimal: $K/n=o(1)$, up to $\log n$ factors. In addition, Theorem \ref{thm:Confidence Bands} gives new inference results for bias-corrected estimators ($j=1,2,3$).

For the case of special univariate regression splines, \cite{Zhou-Shen-Wolfe_1998_AoS} constructs conservative confidence bands under the assumption of normal errors and the rate restriction $K^2/n=o(1)$. In comparison, the confidence band constructed using Theorem \ref{thm:Confidence Bands} has asymptotically exact coverage rate, and requires substantially weaker tuning parameter rate restrictions.

\section{Simulations}
\label{sec:numerical}

We conducted a Monte Carlo investigation of the finite sample performance of our methods. Only a summary is given here, while the SA (\S SA-9) contains complete results and details.

We considered three univariate ($d=1$), two bivariate ($d=2$), and two trivariate ($d=3$) data generating processes. We shall summarize one univariate design here for brevity. We set $\mu(x)=\sin(\pi x-\pi/2) / (1 + 2(2x-1)^2(\sign(2x-1)+1))$, with $\sign(\cdot)$ denoting the sign function. We generate samples $\{(y_i, x_i): i=1,\ldots, n\}$ from $y_i=\mu(x_i) + \e_i$, where $x_i\thicksim \mathsf{U}[0,1]$ and $\e_i\thicksim \mathsf{N}(0,1)$, independent of each other. We consider $5,000$ simulated datasets with $n=1,000$ each. Results based on splines and wavelets are presented. We use linear splines or Daubechies (father) wavelets of order $2$ ($m=2$) to form the point estimator $\widehat{\mu}_0(x)$, and quadratic splines or Daubechies wavelets of order $3$ ($\tilde{m}=3$) for bias correction, on the same evenly spaced partitioning for point estimation and bias correction ($\Delta=\tilde{\Delta}$).

Table \ref{table: simuls} reports (simulated) root mean squared error for point estimators (column ``RMSE''), coverage rate and average interval length of pointwise $95\%$ nominal confidence intervals at $x=0.5$ (columns ``CR'' and ``AL''), and analogous uniform inference results (columns under ``Uniform''). For $B$-splines, $\kappa$ is set using either an infeasible IMSE-optimal choice ($\kappa_\mathtt{IMSE}$), a rule-of-thumb estimate ($\hat{\kappa}_\mathtt{ROT}$), or a direct plug-in estimate ($\hat{\kappa}_\mathtt{DPI}$). For wavelets, the tuning parameter is instead the resolution level (resp., $s_{\mathtt{IMSE}}$, $\hat{s}_{\mathtt{ROT}}$, or $\hat{s}_{\mathtt{DPI}}$), which is the logarithm of the number of subintervals (to base $2$). Finally, Table \ref{table: simuls} reports all four (estimation and) inference methods discussed in this paper ($j=0,1,2,3$), except for wavelets with plug-in bias correction ($j=3$) because of the lack of smoothness of low-order wavelet bases.

Numerical findings are consistent with our theoretical results. Robust bias correction seems to perform quite well, always delivering close-to-correct coverage, both pointwise and uniformly. In addition, our rule-of-thumb (ROT) and direct plug-in (DPI) knot selection procedures for tensor-product partitions exhibited good numerical performance. 

\section{Main Technical Lemmas and Proofs}\label{sec:proofs}

\subsection{Technical Lemma}\label{proof:techlemmas}

Let $\widehat{\bQ}_m=\E_n[\bp(\bx_i)\bp(\bx_i)']$, $\widehat{\bQ}_{\tilde{m}}=\E_n[\bpt(\bx_i)\bpt(\bx_i)']$, $\bQ_m=\E[\bp(\bx_i)\bp(\bx_i)']$, and  $\bQ_{\tilde{m}}=\E[\bpt(\bx_i)\bpt(\bx_i)']$.

\begin{lem} \label{lem:technical lemma}
    Let Assumptions \ref{assumption: DGP}, \ref{assumption: quasi-uniform partition}, \ref{assumption: Local Basis}, and \ref{assumption: BC bases} hold. If $\frac{\log n}{nh^d}=o(1)$, then:
    \emph{(i)} $\|\widehat{\bQ}_m-\bQ_m\|\lesssim_\P h^d\sqrt{\frac{\log n}{nh^d}}$, $\|\widehat{\bQ}_m-\bQ_m\|_\infty\lesssim_\P h^d\sqrt{\frac{\log n}{nh^d}}$;
    \emph{(ii)} $\|\widehat{\bQ}_m\|\lesssim_\P h^d$, $\|\widehat{\bQ}_m^{-1}\|_\infty\lesssim_\P h^{-d}$; 
    \emph{(iii)} for each $j=0,1,2,3$, $\sup_{\bx\in \mathcal{X}}\|\bgamma_{\bq, j}(\bx)'\|_\infty
    \lesssim h^{-d-[\bq]}$,  $\sup_{\bx\in\mathcal{X}}\|\bg{j}'-\bgamma_{\bq, j}(\bx)'\|_\infty
    \lesssim h^{-d-[\bq]}\sqrt{\frac{\log n}{nh^d}}$, $\inf_{\bx\in\mathcal{X}}\|\bgamma_{\bq,j}(\bx)'\|\gtrsim h^{-d-[\bq]}$; and
    \emph{(iv)} for $j=0,1,2,3$, $\sup_{\bx\in\mathcal{X}}\Omega_j(\bx)\lesssim h^{-d-2[\bq]}$ and $\inf_{\bx\in\mathcal{X}} \Omega_j(\bx)
    \gtrsim h^{-d-2[\bq]}$.\smallskip
    
    \noindent Proof: SA, Section SA-10.\qed
\end{lem} 
These results for $\widehat{\bQ}_m$ and $\bQ_m$ also hold for $\widehat{\bQ}_{\tilde{m}}$ and $\bQ_{\tilde{m}}$ under Assumption \ref{assumption: BC bases}. See the SA (\S SA-2) for details and other related results.

\subsection{Proof of Lemma \ref{lem:Strong Approximation Hats Off}}\label{proof:lem:Strong Approximation Hats Off}

First, suppose condition \textbf{(i)} holds. Theorem SA-4.2 of the SA shows
$\sup_{\bx\in\mathcal{X}}\big|\widehat{\Omega}_0(\bx)-\Omega_0(\bx)\big| \lesssim_\P n^{-\frac{1}{2}}h^{-\frac{3d}{2}-2[\bq]}\big[(\log n)^{\frac{1}{2}}+n^{\frac{1}{2+\nu}}(\log n)^{\frac{\nu}{4+2\nu}}+\sqrt{n}h^{\frac{d}{2}+m}\big]$ and, for $j=1,2,3$, $\sup_{\bx\in\mathcal{X}}\big|\widehat{\Omega}_j(\bx)-\Omega_j(\bx)\big| \lesssim_\P n^{-\frac{1}{2}}h^{-\frac{3d}{2}-2[\bq]}
\big[(\log n)^{\frac{1}{2}}+n^{\frac{1}{2+\nu}}(\log n)^{\frac{\nu}{4+2\nu}}+\sqrt{n}h^{\frac{d}{2}+m+\varrho}\big]$.
Then, for $j=0,1,2,3$,
$
\sup_{\bx\in\mathcal{X}}
        |(\widehat{\partial^\bq \mu}_j(\bx)-\partial^\bq \mu(\bx))/\sqrt{\Omega_j(\bx)/n}
            - (\widehat{\partial^\bq \mu}_j(\bx)-\partial^\bq \mu(\bx))/\sqrt{\widehat{\Omega}_j(\bx)/n}|
\lesssim_\P \sqrt{n}h^{3d/2+3[\bq]}
             \sup_{\bx\in\mathcal{X}}\big|\widehat{\partial^\bq \mu}_j(\bx)-\partial^\bq \mu(\bx)\big|
             \sup_{\bx\in\mathcal{X}}\big|\widehat{\Omega}_j(\bx)-\Omega_j(\bx)\big|
= o_\P(r_n^{-1}),
$
where the result follows from Lemma \ref{lem:technical lemma}, Theorem SA-4.1, the uniform convergence rate of $\widehat{\Omega}_j(\bx)$, and the rate conditions imposed.

The result under the conditions in \textbf{(ii)} follows similarly. \qed

\subsection{Proof of Theorem \ref{thm:strong approximation KMT}}\label{proof:thm:strong approximation KMT}
We first prove the following general lemma. Let $TV_{\mathcal{X}}(g(\cdot))$ denote the total variation of $g(\cdot)$ on $\mathcal{X}\subseteq\mathbb{R}$.

\begin{lem}[Kernel-Based KMT Coupling] \label{lem:general cond coupling}
	Suppose $\{(x_i,\e_i):1\leq i \leq n \}$ are i.i.d., with $x_i\in\mathcal{X}\subseteq\mathbb{R}$ and $\sigma_i^2:=\sigma^2(x_i)=\E[\e_i^2|x_i]$. Let $\{A(x):=\G_n[\mathscr{K}(x, x_i)\e_i], x\in\mathcal{X}\}$ be a stochastic process with $\mathscr{K}(\cdot, \cdot):\mathbb{R} \times \mathbb{R} \mapsto \mathbb{R}$ an $n$-varying kernel function possibly depending on $\bX$. Assume one of the following holds:
	\begin{description}
		\item[(i)] $\sup_{x\in\mathcal{X}}\E[|\e_i|^{2+\nu}|x_i=x]<\infty$, for some $\nu>0$, and 
		\[
		\begin{split}
		&\sup_{x\in\mathcal{X}}\max_{1\leq i\leq n}|\mathscr{K}(x, x_i)| = o_\P(r_n^{-1}n^{-\frac{1}{2+\nu}+\frac{1}{2}}),\\
		&\sup_{x\in \mathcal{X}}TV_{\mathcal{X}}(\mathscr{K}(x,\cdot))
		=o(r_n^{-1}n^{-\frac{1}{2+\nu}+\frac{1}{2}});
		\quad or
		\end{split}
		\]
		
		\item[(ii)] $\sup_{x\in \mathcal{X}} \E[|\e_i|^3\exp(|\e_i|)|x_i=x] < \infty$ and
		\[
		\begin{split}
		&\sup_{x\in\mathcal{X}}\max_{1\leq i\leq n}|\mathscr{K}(x, x_i)|	= o_\P(r_n^{-1}(\log n)^{-1}\sqrt{n}),\\
		&\sup_{x\in \mathcal{X}} TV_{\mathcal{X}}(\mathscr{K}(x,\cdot))
		=o(r_n^{-1}(\log n)^{-1}\sqrt{n}).
		\end{split}
		\]
	\end{description}
	Then, on a sufficiently rich probability space, there exists a copy $A'(\cdot)$ of $A(\cdot)$, and an i.i.d. sequence $\{\zeta_i: 1\leq i \leq n\}$ of standard Normal random variables such that
	$A(x)=_d \G_n\big[\mathscr{K}(x, x_i)\sigma_i\zeta_i\big] + o_\P(r_n^{-1})$ in $\mathcal{L}^\infty(\mathcal{X})$.
	
\end{lem}

\textit{Proof}. Suppose condition \textbf{(i)} holds. Let $\{x_{(i)}:1\leq i\leq n\}$ be the order statistics of $\{x_i:1\leq i\leq n\}$ such that $x_{(1)}\leq x_{(2)}\leq\cdots\leq x_{(n)}$, which also induces the concomitants $\{\e_{[i]}:1\leq i\leq n\}$ and $\{\sigma^2_{[i]}=\sigma^2(x_{(i)}):1\leq i\leq n\}$.
Conditional on $\mathbf{X}$, $\{\e_{[i]}:1\leq i\leq n\}$ is an independent mean zero sequence with $\V[\e_{[i]}|\bX]=\sigma^2_{[i]}$. By \cite[Corollary 5]{Sakhanenko_1991_SAM}, there exists a sequence of i.i.d standard normal random variables $\{\zeta_{[i]}:1\leq i\leq n\}$ such that $\max_{1\leq l\leq n}|S_{l,n}|\lesssim_\P n^{\frac{1}{2+\nu}}$, where $S_{l,n}:=\sum_{i=1}^l (\e_{[i]} - \sigma_{[i]}\zeta_{[i]})$. Then, using summation by parts,
$\sup_{x\in\mathcal{X}}|\sum_{i=1}^{n}\mathscr{K}(x, x_{(i)})(\e_{[i]}-\sigma_{[i]}\zeta_{[i]})|
= \sup_{x\in\mathcal{X}}|\mathscr{K}(x, x_{(n)})S_{n,n}
   - \sum_{i=1}^{n-1}S_{i,n}\left(\mathscr{K}(x,x_{(i+1)}) - \mathscr{K}(x,x_{(i)})\right)|
\leq (\sup_{x\in\mathcal{X}}\\
\max_{1\leq i\leq n}|\mathscr{K}(x,x_i)|
      + \sup_{x\in\mathcal{X}}\sum_{i=1}^{n-1}|\mathscr{K}(x,x_{(i+1)})-\mathscr{K}(x,x_{(i)})|)\max_{1\leq l \leq n}|S_{l,n}|
$.
Since $\sum_{i=1}^{n-1}\left|\mathscr{K}(x,x_{(i+1)})-\mathscr{K}(x,x_{(i)})\right|\leq TV_{\mathcal{X}}(\mathscr{K}(x,\cdot))$, we have $A(x)=_d\G_n[\mathscr{K}(x,x_i)\sigma_{i}\zeta_i]+ o_\P(r_n^{-1})$.

When \textbf{(ii)} holds, the proof is the same except that under the stronger moment restriction, $\max_{1\leq l\leq n}|S_{l,n}|\lesssim_\P \log n$ by \cite[Theorem 1]{Sakhanenko_1985_Advances}.
\qed\medskip

To prove Theorem \ref{thm:strong approximation KMT}, for each $j=0,1,2,3$, let $\mathscr{K}(x, u) = \bgamma_{\bq,j}(x)'\bPi{j}(u)/\\
\sqrt{\Omega_j(x)}$ and observe that $\sup_{x\in \mathcal{X}}\sup_{u\in\mathcal{X}} |\mathscr{K}(x, u)|
\lesssim h^{-d/2}$. By Lemma \ref{lem:technical lemma}, the uniform bound on the total variation of $\mathscr{K}(x,u)$ can be verified easily. Alternatively, simply note that $\big|\sum_{i=1}^{n-1}S_{i,n} \left(\mathscr{K}(x,x_{(i+1)}) -\mathscr{K}(x,x_{(i)})\right)\big|
\leq\big\|\bgamma_{\bq,j}(x)'/\sqrt{\Omega_j(x)}\big\|_\infty\big\|\sum_{i=1}^{n-1}S_{i,n}(\bPi{j}(x_{(i+1)})-\bPi{j}(x_{(i)}))\big\|_\infty$. By Assumption \ref{assumption: Local Basis} and Lemma \ref{lem:technical lemma},
$\sup_{x\in\mathcal{X}}\|\bgamma_{\bq,j}(x)'/\sqrt{\Omega_j(x)}\|_\infty \lesssim h^{-d/2}$. Denote the $l$th element of $\bPi{j}(\cdot)$ by $\pi_{j,l}(\cdot)$. Then, $\max_{1\leq l \leq K_j} \big|\sum_{i=1}^{n-1}\big(\pi_{j,l}(x_{(i+1)})-\pi_{j,l}(x_{(i)})\big)S_{l,n}\big|\\
\leq \max_{1\leq l \leq K_j} \sum_{i=1}^{n-1}\big|\pi_{j,l}(x_{(i+1)})-\pi_{j,l}(x_{(i)})\big| \max_{1\leq \ell \leq n}\big|S_{\ell,n}\big|$. By Assumption \ref{assumption: Local Basis} and \ref{assumption: BC bases}, $\max_{1\leq l \leq K_j}
\sum_{i=1}^{n-1}|\pi_{j,l}(x_{(i+1)})-\pi_{j,l}(x_{(i)})|\lesssim 1$. Thus, using Lemma \ref{lem:general cond coupling}, under the corresponding moment conditions and rate restrictions, there exists an independent standard normal sequence $\{\zeta_i:1\leq i\leq n\}$ such that $\G_n[\mathscr{K}(x,x_i)\e_i]=_d z_j(x)+o_\P(r_n^{-1})$.

Next, note that
$
z_j(\bx) =_{d|\bX}
            \bgamma_{\bq,j}(\bx)'\bSig{j}^{1/2}\bN_{K_j}/\sqrt{\Omega_j(\bx)}
          + \bgamma_{\bq,j}(\bx)'\big(\bSigbar{j}^{1/2}-\bSig{j}^{1/2}\big)\bN_{K_j}/\sqrt{\Omega_j(\bx)}
$
where $\bN_{K_j}$ is a $K_j$-dimensional standard normal vector (independent of $\bX$) and ``$=_{d|\bX}$'' denotes that two processes have the same conditional distribution given $\bX$. For the second term, by the Gaussian maximal inequality \cite[Lemma 13]{Chernozhukov-Lee-Rosen_2013_ECMA}, $\E\big[ \big\|\big(\bSigbar{j}^{1/2}-\bSig{j}^{1/2}\big)\bN_{K_j}\big\|_\infty\,\big|\bX\big]
\lesssim \sqrt{\log n} \,\big\|\bSigbar{j}^{1/2}-\bSig{j}^{1/2}\big\|$. By the same argument as that in the proof of Lemma SA-2.1 in \S SA-10.1, $\big\|\bSigbar{j}-\bSig{j}\big\| \lesssim_\P h^d (\log n/(nh^d))^{1/2}$.
Then by \cite[Theorem X.1.1]{Bhatia_2013_Book}, $\big\|\bSigbar{j}^{1/2}-\bSig{j}^{1/2}\big\| \lesssim_\P h^{d/2}(\log n/(nh^d))^{1/4}.$ For $j=0,1$, a sharper bound is available: by \cite[Theorem X.3.8]{Bhatia_2013_Book} and Lemma \ref{lem:technical lemma},
$\|\bSigbar{j}^{1/2}-\bSig{j}^{1/2}\|
  \leq \lambda_{\min}(\bSig{j})^{-1/2} \|\bSigbar{j}-\bSigma_{j}\|
  \lesssim_\P h^{d/2}\sqrt{\frac{\log n}{nh^d}}.
$
Thus,
$
\E\big[\sup_{\bx\in\mathcal{X}}\big|\bgamma_{\bq,j}(\bx)'
\big(\bSigbar{j}^{1/2}-\bSig{j}^{1/2}\big)\bN_{K_j}/\sqrt{\Omega_j(\bx)}\big|\;\big|\bX \big]
\lesssim_\P h^{-\frac{d}{2}}\sqrt{\log n}\,\big\|\bSigbar{j}^{1/2}-\bSig{j}^{1/2}\big\|=o_\P(r_n^{-1}).
$
The results now follow from Markov inequality and Dominated Convergence.
\qed

\subsection{Proof of Theorem \ref{thm:strong approximation Yurinskii}}\label{proof:thm:strong approximation Yurinskii}
It suffices to verify the conditions in Lemma 39 of \cite{Belloni-Chernozhukov-Chetverikov-FernandezVal_2019_JoE}.
For $j=0,1,2,3$, let $\boldsymbol{\xi}_i=\frac{1}{\sqrt{n}}\bPi{j}(\bx_i)\e_i$. $\{\boldsymbol{\xi}_i:1\leq i \leq n\}$ is an i.i.d. sequence of $K_j$-dimensional random vectors, and
$\sum_{i=1}^{n}\E[\|\boldsymbol{\xi}_i\|^2\|\boldsymbol{\xi}_i\|_\infty]
 =\E\big[\big\|\bPi{j}(\bx_i)\e_i\big\|^2\big\|\bPi{j}(\bx_i)\e_i\big\|_\infty\big]/\sqrt{n}
 \lesssim \E\big[\bPi{j}(\bx_i)'\bPi{j}(\bx_i)|\e_i|^3 \big]/\sqrt{n}
 \lesssim n^{-1/2}
$
by Assumption \ref{assumption: Local Basis}, the moment condition imposed, and Lemma \ref{lem:technical lemma}.
On the other hand, let $\{\mathbf{g}_i:1\leq i\leq n\}$ be a sequence of independent Gaussian vectors with mean zero and variance $\bSig{j}/n$. By properties of Gaussian random variables and Lemma \ref{lem:technical lemma},
$(\E[\|\mathbf{g}_i\|_\infty^2])^{1/2}\lesssim \sqrt{\frac{\log n}{n}}$, $\sum_{i=1}^{n}(\E[\|\mathbf{g}_i\|^4])^{1/2}\lesssim
\tr\big(\sum_{i=1}^{n}\E[\boldsymbol{\xi}_i\boldsymbol{\xi}_i']\big)\lesssim 1$, and thus
$L_n:=\sum_{i=1}^{n}\E[\|\boldsymbol{\xi}_i\|^2 \|\boldsymbol{\xi}_i\|_\infty]+\sum_{i=1}^{n}\E[\|\mathbf{g}_i\|^2\\
\|\mathbf{g}_i\|_\infty]\lesssim \sqrt{\frac{\log n}{n}}$.
Then, there exists a $K_j$-dimensional normal vector $\mathbf{N}_{K_j}$ with variance $\bSig{j}$ such that for any $t>0$,
$
\P\big(\|\sum_{i=1}^{n}\boldsymbol{\xi}_i-\mathbf{N}_{K_j}\|_\infty>3r_n^{-1}h^{\frac{d}{2}}t\big) \leq \min_{\tau\geq 0}
\big(2\P(\|\bZ\|_\infty>\tau)+h^{-\frac{3d}{2}}t^{-3}r_n^3L_n\tau^2\big)
\lesssim t^{-3}r_n^3(\log n)^{\frac{3}{2}}/\sqrt{nh^{3d}}
$
where $\bZ$ is a $K_j$-dimensional standard Gaussian vector, and the second inequality follows by setting $\tau=C\sqrt{\log n}$ for a sufficiently large $C>0$. Using $\sup_{x\in\mathcal{X}}\|\bgamma_{\bq,j}(x)'/\sqrt{\Omega_j(x)}\|_\infty \lesssim h^{-d/2}$ again, the result follows.
\qed

\subsection{Proof of Theorem \ref{thm:Plug-in Aprox Unconditional}}\label{proof:thm:Plug-in Aprox Unconditional}
For each $j=0,1,2,3$,\newline
$
\widehat{Z}_j(\bx)-Z_j(\bx)=
\Big(\frac{\bg{j}}{\widehat{\Omega}^{1/2}_j(\bx)}-
 \frac{\bgamma_{\bq,j}(\bx)'}{\Omega^{1/2}_j(\bx)}\Big)\bSighat{j}^{1/2}\bN_{K_j}
 + \frac{\bgamma_{\bq,j}(\bx)'}{\sqrt{\Omega_j(\bx)}}[\bSighat{j}^{1/2}-\bSig{j}^{1/2}]\bN_{K_j}
$.\newline
Each term on the RHS is a mean-zero Gaussian process conditional on the data. The desired results then follow by Lemma \ref{lem:technical lemma}, Theorem SA-4.2 in \S SA-4, and applying the Gaussian maximal inequality to each term as in Section \ref{proof:thm:strong approximation KMT}.
\qed

\subsection{Proof of Theorem \ref{thm:Confidence Bands}}\label{proof:thm:Confidence Bands}
By Theorem \ref{thm:strong approximation KMT} or Theorem \ref{thm:strong approximation Yurinskii}, there exists a sequence of constants $\eta_n$ such that $\eta_n=o(1)$ and
$\P(|\sup_{\bx\in\mathcal{X}}|\widehat{T}_j(\bx)|-\sup_{\bx\in\mathcal{X}}|Z_j(\bx)||>\eta_n/r_n)=o(1)$.
Therefore, for any $u\in\mathbb{R}$,
$
\P\big[\sup_{\bx\in\mathcal{X}}|\widehat{T}_j(\bx)|\leq u\big]
\leq \P\big[\{\sup_{\bx\in\mathcal{X}}|\widehat{T}_j(\bx)|\leq u\} 
            \cap \{|\sup_{\bx\in\mathcal{X}}|\widehat{T}_j(\bx)|-\sup_{\bx\in\mathcal{X}}|Z_j(\bx)||\leq\eta_n/r_n\}\big]
     + \P\big[\{|\sup_{\bx\in\mathcal{X}}|\widehat{T}_j(\bx)|-\sup_{\bx\in\mathcal{X}}|Z_j(\bx)||>\eta_n/r_n\}\big]
\leq\P\left[\sup_{\bx\in\mathcal{X}}|Z_j(\bx)|\leq u+\eta_n/r_n\right]+o(1)
\leq\P\left[\sup_{\bx\in\mathcal{X}}|Z_j(\bx)|\leq u\right] + Cr_n^{-1}\eta_n\E[\sup_{\bx\in\mathcal{X}}|Z_j(\bx)|]+ o(1)
$
for some constant $C>0$, where the last inequality holds by the anti-concentration inequality due to \cite{Chernozhukov-Chetverikov-Kato_2014b_AoS}. By the Gaussian maximal inequality, $\E[\sup_{\bx\in\mathcal{X}}|Z_j(\bx)|]\lesssim \sqrt{\log n}$. Since $r_n=\sqrt{\log n}$, the second term on the right of the last inequality is $o(1)$. The reverse of the inequality follows similarly, and thus
$\sup_{u\in\mathbb{R}}\big|\P\big[\sup_{\bx\in\mathcal{X}}|\widehat{T}_j(\bx)|\leq u\big]-\P\big[\sup_{\bx\in\mathcal{X}}|Z_j(\bx)|
\leq u\big]\big|=o(1)$. In addition, by Theorem \ref{thm:Plug-in Aprox Unconditional}, $\widehat{Z}_j(\cdot)$ is approximated by the same Gaussian process conditional on the data. The result then follows by the same argument.
\qed

\section{Conclusion}\label{sec:conclusion}

We presented new asymptotic results for partitioning-based least squares regression estimators. The first main contribution gave a general IMSE expansion for the point estimators. The second set of contributions were pointwise and uniform distributional approximations, with and without robust bias correction, for $t$-statistic processes indexed by $\bx\in\mathcal{X}$, with improvements in rate restrictions and convergence rates. For the case of $d=1$, our uniform approximation results rely on a new coupling approach, which delivered seemingly minimal rate restrictions. Furthermore, we apply our general results to three popular special cases: $B$-splines, compactly supported wavelets, and piecewise polynomials. Finally, we provide a general purpose {\sf R} package \texttt{lspartition} \cite{Cattaneo-Farrell-Feng_2019_lspartition}.

\section*{Acknowledgements}
We thank Victor Chernozhukov, Denis Chetverikov, Michael Jansson, Xinwei Ma, and Whitney Newey for useful discussions. We also thank the co-Editors, Edward George and Ming Yuan, an associate editor, and the reviewers for thoughtful comments that significantly improved this paper. See \ref{suppA} for supplementary materials.

\begin{supplement}
\sname{Supplement A}\label{suppA}
\stitle{Additional Technical Results, Omitted Proofs, Implementation Details, and Further Simulation Results}
\slink[url]{http://arxiv.org/pdf/1804.04916}
\sdescription{The SA gives omitted proofs and additional technical results that may be of independent interest, including pointwise and uniform stochastic linearization useful in semiparametric settings (\S SA-4; see, in particular, Remark SA-4.1), theoretical comparisons between bias correction approaches, and a discussion of the relationship between $B$-Splines and polynomials. Details on implementation, specific examples, and further simulation evidence are also reported.}
\end{supplement}

\bibliography{Cattaneo-Farrell-Feng_2019_AoS--Bibliography}
\bibliographystyle{imsart-number}

\begin{table}
	\renewcommand{\arraystretch}{1.1}
	\caption{Simulation Evidence}\label{table: simuls}
	\subfloat[B-Splines ($m=2$, $\tilde{m}=3$, $\Delta=\tilde{\Delta}$, Evenly Spaced Partition)]
	         {\resizebox{.8\textwidth}{!}{
\begin{tabular}{lrcrcrrcrr}
\hline\hline
\multicolumn{1}{l}{\bfseries }&\multicolumn{1}{c}{\bfseries }&\multicolumn{1}{c}{\bfseries }&\multicolumn{1}{c}{\bfseries }&\multicolumn{1}{c}{\bfseries }&\multicolumn{2}{c}{\bfseries Pointwise}&\multicolumn{1}{c}{\bfseries }&\multicolumn{2}{c}{\bfseries Uniform}\tabularnewline
\cline{6-7} \cline{9-10}
\multicolumn{1}{l}{}&\multicolumn{1}{c}{$\kappa$}&\multicolumn{1}{c}{}&\multicolumn{1}{c}{RMSE}&\multicolumn{1}{c}{}&\multicolumn{1}{c}{CR}&\multicolumn{1}{c}{AL}&\multicolumn{1}{c}{}&\multicolumn{1}{c}{UCR}&\multicolumn{1}{c}{AW}\tabularnewline
\hline
{\bfseries $j=0$}&&&&&&&&&\tabularnewline
~~$\kappa_{\mathtt{IMSE}}$&$3.0$&&$0.046$&&$91.6$&$0.328$&&$79.9$&$0.384$\tabularnewline
~~$\hat{\kappa}_{\mathtt{ROT}}$&$4.1$&&$0.002$&&$94.9$&$0.254$&&$90.1$&$0.433$\tabularnewline
~~$\hat{\kappa}_{\mathtt{DPI}}$&$4.7$&&$0.008$&&$93.8$&$0.311$&&$91.6$&$0.460$\tabularnewline
\hline
{\bfseries $j=1$}&&&&&&&&&\tabularnewline
~~$\kappa_{\mathtt{IMSE}}$&$3.0$&&$0.003$&&$94.8$&$0.226$&&$93.8$&$0.426$\tabularnewline
~~$\hat{\kappa}_{\mathtt{ROT}}$&$4.1$&&$0.008$&&$94.8$&$0.297$&&$93.5$&$0.473$\tabularnewline
~~$\hat{\kappa}_{\mathtt{DPI}}$&$4.7$&&$0.007$&&$94.9$&$0.294$&&$93.2$&$0.497$\tabularnewline
\hline
{\bfseries $j=2$}&&&&&&&&&\tabularnewline
~~$\kappa_{\mathtt{IMSE}}$&$3.0$&&$0.004$&&$94.7$&$0.268$&&$94.1$&$0.443$\tabularnewline
~~$\hat{\kappa}_{\mathtt{ROT}}$&$4.1$&&$0.008$&&$94.8$&$0.312$&&$93.4$&$0.497$\tabularnewline
~~$\hat{\kappa}_{\mathtt{DPI}}$&$4.7$&&$0.004$&&$94.8$&$0.330$&&$93.6$&$0.526$\tabularnewline
\hline
{\bfseries $j=3$}&&&&&&&&&\tabularnewline
~~$\kappa_{\mathtt{IMSE}}$&$3.0$&&$0.016$&&$90.0$&$0.320$&&$89.0$&$0.413$\tabularnewline
~~$\hat{\kappa}_{\mathtt{ROT}}$&$4.1$&&$0.007$&&$94.2$&$0.275$&&$93.0$&$0.463$\tabularnewline
~~$\hat{\kappa}_{\mathtt{DPI}}$&$4.7$&&$0.005$&&$94.2$&$0.322$&&$93.2$&$0.490$\tabularnewline
\hline
\end{tabular}
}}\\
	\subfloat[Wavelets ($m=2$, $\tilde{m}=3$, $\Delta=\tilde{\Delta}$, Evenly Spaced Partition)]
	         {\resizebox{.8\textwidth}{!}{
\begin{tabular}{lrcrcrrcrr}
\hline\hline
\multicolumn{1}{l}{\bfseries }&\multicolumn{1}{c}{\bfseries }&\multicolumn{1}{c}{\bfseries }&\multicolumn{1}{c}{\bfseries }&\multicolumn{1}{c}{\bfseries }&\multicolumn{2}{c}{\bfseries Pointwise}&\multicolumn{1}{c}{\bfseries }&\multicolumn{2}{c}{\bfseries Uniform}\tabularnewline
\cline{6-7} \cline{9-10}
\multicolumn{1}{l}{}&\multicolumn{1}{c}{$s$}&\multicolumn{1}{c}{}&\multicolumn{1}{c}{RMSE}&\multicolumn{1}{c}{}&\multicolumn{1}{c}{CR}&\multicolumn{1}{c}{AL}&\multicolumn{1}{c}{}&\multicolumn{1}{c}{UCR}&\multicolumn{1}{c}{AW}\tabularnewline
\hline
{\bfseries $j=0$}&&&&&&&&&\tabularnewline
~~$s_{\mathtt{IMSE}}$&$3.0$&&$0.002$&&$94.2$&$0.476$&&$91.1$&$0.509$\tabularnewline
~~$\hat{s}_{\mathtt{ROT}}$&$2.0$&&$0.002$&&$94.2$&$0.476$&&$91.1$&$0.509$\tabularnewline
~~$\hat{s}_{\mathtt{DPI}}$&$2.8$&&$0.002$&&$94.2$&$0.476$&&$91.1$&$0.509$\tabularnewline
\hline
{\bfseries $j=1$}&&&&&&&&&\tabularnewline
~~$s_{\mathtt{IMSE}}$&$3.0$&&$0.036$&&$93.6$&$0.449$&&$89.9$&$0.504$\tabularnewline
~~$\hat{s}_{\mathtt{ROT}}$&$2.0$&&$0.036$&&$93.6$&$0.449$&&$89.9$&$0.504$\tabularnewline
~~$\hat{s}_{\mathtt{DPI}}$&$2.8$&&$0.036$&&$93.6$&$0.449$&&$89.9$&$0.504$\tabularnewline
\hline
{\bfseries $j=2$}&&&&&&&&&\tabularnewline
~~$s_{\mathtt{IMSE}}$&$3.0$&&$0.009$&&$94.2$&$0.523$&&$91.4$&$0.576$\tabularnewline
~~$\hat{s}_{\mathtt{ROT}}$&$2.0$&&$0.009$&&$94.2$&$0.523$&&$91.4$&$0.576$\tabularnewline
~~$\hat{s}_{\mathtt{DPI}}$&$2.8$&&$0.009$&&$94.2$&$0.523$&&$91.4$&$0.576$\tabularnewline
\hline
\end{tabular}
}}
	
	\flushleft\textbf{Notes}:\newline
	(i) Pointwise = pointwise inference at $x=0.5$, Uniform = uniform inference.\newline
	(ii) RMSE = root MSE of point estimator, CR = coverage rate of $95\%$ nominal confidence intervals, AL = average interval length of $95\%$ nominal confidence intervals.\newline
	(iii) UCR = uniform coverage rate of $95\%$ nominal confidence band, AW = average width of $95\%$ nominal confidence band.\newline
	(iv) $\kappa_\mathtt{IMSE}$ and $s_\mathtt{IMSE}$ = infeasible IMSE-optimal number of partitions, $\hat{\kappa}_\mathtt{ROT}$ and $\hat{s}_\mathtt{ROT}$ = feasible rule-of-thumb (ROT) implementation of $\kappa_\mathtt{IMSE}$, $\hat{\kappa}_\mathtt{DPI}$ and $\hat{s}_\mathtt{DPI}$ = feasible direct plug-in (DPI) implementation of $\kappa_\mathtt{IMSE}$. See \S SA-8 and \S SA-9 in supplemental appendix for more details.
	
\end{table}

\end{document}


\begin{frontmatter}
\title{Large Sample Properties of \\
	Partitioning-Based Series Estimators\medskip\\
       Supplemental Appendix
\bigskip }
\runtitle{Supplemental: Partitioning-Based Series Estimators}
\begin{aug}
	\author{\fnms{Matias D.} \snm{Cattaneo}\thanksref{t1}\ead[label=e1]{cattaneo@princeton.edu}},
	\author{\fnms{Max H.} \snm{Farrell}\ead[label=e2]{max.farrell@chicagobooth.edu}}
	\and
	\author{\fnms{Yingjie} \snm{Feng}\ead[label=e3]{yjfeng@umich.edu}}
	
	\thankstext{t1}{Financial support from the National Science Foundation (SES 1459931) is gratefully acknowledged.}
	\runauthor{Cattaneo, Farrell, and Feng}
	
	\affiliation{Princeton University, University of Chicago, and Princeton University}
	
	\address{Matias D. Cattaneo\\
		Department of Operations Research\\ and Financial Engineering\\
		Princeton University\\
		Princeton, NJ 08544\\
		\printead{e1}
	}
	
	\address{Max H. Farrell\\
		Booth School of Business\\
		University of Chicago\\
		Chicago, IL 60637\\
		\printead{e2}
	}
	
	\address{Yingjie Feng\\
		Department of Politics\\
		Princeton University\\
		Princeton, NJ 08544\\
		\printead{e3}
	}
	
\end{aug}

\bigskip
\begin{abstract}
This supplement gives omitted theoretical proofs of the results discussed in the main paper, additional technical results and methodological discussions, which may be of independent interest, and further simulation evidence. Section \ref{section: setup} repeats the setup and assumptions. Section \ref{section: lemmas} states a few important lemmas. Then, Section \ref{section: IMSE rec} provides an integrated mean squared error expansion for the case of general tensor-product partitions, building on the discussion in the main text. Section \ref{section: convergence rates} proves pointwise and uniform stochastic linearizations. Section \ref{section: uniform} states all uniform inference results, including two feasible inference methods omitted from the main text. Section \ref{section: examples} discusses leading examples of partitioning-based estimators, where the main assumptions used in the paper are verified for each case, including splines, wavelets, and piecewise polynomials. Some interesting conceptual and technical asides are given in Section \ref{section: extension}. Finally, implementation and other numerical issues are discussed in Section \ref{section: implementation}, and complete results from a simulation study are presented in Section \ref{section: numerical}. See also the companion {\sf R} package {\tt lspartition} detailed in \cite{Cattaneo-Farrell-Feng_2019_lspartition} and available at
\begin{center}
	\url{http://sites.google.com/site/nppackages/lspartition/}.
\end{center}
\end{abstract}

\begin{keyword}[class=MSC]
	\kwd[Primary ]{62H10}\kwd{62M99}\kwd{57R12}
	\kwd[; secondary ]{62M99}
\end{keyword}

\begin{keyword}
	\kwd{nonparametric regression}\kwd{series methods}\kwd{sieve methods}\kwd{robust bias correction}\kwd{uniform inference}\kwd{strong approximation}\kwd{tuning parameter selection}
\end{keyword}

\end{frontmatter}


\clearpage
\tableofcontents 
\clearpage


\section{Setup and Assumptions} \label{section: setup}

For a $d$-tuple $\bq=(q_1,\ldots,q_d)\in\mathbb{Z}^d_+$, define $[\bq]=\sum_{\ell=1}^{d}q_\ell$, $\bx^\bq=x_1^{q_1}x_2^{q_2}\cdots x_d^{q_d}$, $\bq!=q_1!\cdots q_d!$ and $\partial^\bq \mu(\bx)=\partial^{[\bq]} \mu(\bx)/\partial x_1^{q_1}\ldots\partial x_d^{q_d}$. Unless explicitly stated otherwise, whenever $\bx$ is a boundary point of some closed set, the partial derivative is understood as the limit with $\bx$ ranging within it. Let $\mathbf{0}=(0,\cdots,0)'$ and $\mathbf{1}=(1,\cdots, 1)'$ be the length-$d$ vectors of zeros and ones respectively.
The tensor product or Kronecker product operator is $\otimes$, and the entrywise division operator (Hadamard division) is $\oslash$. The smallest integer greater than or equal to $u$ is $\lceil u \rceil$.

We use several norms. For a column vector $\mathbf{v}=(v_1,\ldots,v_M)'\in\mathbb{R}^M$, we let $\|\mathbf{v}\|=(\sum_{i=1}^M v_i^2)^{1/2}$, $\|\mathbf{v}\|_\infty=\max_{1\leq i\leq M}|v_i|$ and $\dim(\mathbf{v})=M$. For a matrix $\mathbf{A}\in\mathbb{R}^{M\times N}$,
$\|\mathbf{A}\|_1=\max_{1\leq j \leq N}\sum_{i=1}^{M}|a_{ij}|$, 
$\|\mathbf{A}\|=\max_i\sigma_i(\mathbf{A})$ and $\|\mathbf{A}\|_\infty=\max_{1\leq i\leq M}\sum_{j=1}^N|a_{ij}|$ for operator norms induced by $L_1$, $L_2$, and $L_\infty$ norms respectively, where $\sigma_i(\mathbf{A})$ is the $i$-th singular value of $\mathbf{A}$. $\lambda_{\max}(\mathbf{A})$ and $\lambda_{\min}(\mathbf{A})$ denote the maximum and minimum eigenvalues of $\mathbf{A}$. For a real-valued function $g(\bx)$, $\|g\|_{L_p(\mathcal{X})}=(\int_\mathcal{X}|g(x)|^p\,d\bx)^{1/p}$ and $\|g\|_{L_\infty(\mathcal{X})}=\esssup_{\bx\in\mathcal{X}}|g(\bx)|$, and let
$\mathcal{C}^s(\mathcal{X})$ denote the space of $s$-times continuously differentiable functions on $\mathcal{X}$.

We also use standard empirical process notation:  $\E_n[g]=\E_n[g(\bx_i)]=\frac{1}{n}\sum_{i=1}^ng(\bx_i)$, $\E[g]=\E[g(\bx_i)]$, and $\G_n[g]=\G_n[g(\bx_i)]=\frac{1}{\sqrt{n}}\sum_{i=1}^n(g(\bx_i)-\E[g(\bx_i)])$. For sequences of numbers or random variables, we use $a_n \lesssim b_n$ to denote $\limsup_n |a_n/b_n|$ is finite, and $a_n \lesssim_\P b_n$ or $a_n=O_\P(b_n)$ to denote $\limsup_{\epsilon\rightarrow\infty}\limsup_n\P[|a_n/b_n|\geq\epsilon]=0$. $a_n=o(b_n)$ implies $a_n/b_n\rightarrow 0$, and $a_n=o_\P(b_n)$ implies that $a_n/b_n\rightarrow_\P 0$ where $\rightarrow_\P$ denotes convergence in probability. $a_n\asymp b_n$ implies that $a_n\lesssim b_n$ and $b_n\lesssim a_n$. $\rightsquigarrow$ denotes convergence in distribution, and for two random variables $X$ and $Y$, $X=_d Y$ implies that they have the same probability distribution.

We set $\mu(\bx):=\partial^\mathbf{0}\mu(\bx)$ and likewise $\widehat{\mu}_j(\bx):=\widehat{\partial^\mathbf{0}\mu}_j(\bx)$ for $j=0,1,2,3$. Let $\bX=[\bx_1, \ldots, \bx_n]'$. Finally, let $C, C_1, C_2, \ldots$ denote universal constants which may be different in different uses.

\subsection{Setup}

As a complete reference, this section repeats the setup, assumptions, and notation used in the main paper. 

We study the standard nonparametric regression setup,  where $\{(y_i,\bx_i'), i = 1, \ldots n\}$ is a random sample from the model 
\begin{equation}\label{equation: regression model}
y_i=\mu(\bx_i)+\e_i ,    	 	\qquad  	  \E[\e_i|\bx_i]=0  ,    	 	\qquad  	  \E[\e_i^2|\bx_i]=\sigma^2(\bx_i)  ,
\end{equation}
for a scalar response $y_i$ and a $d$-vector of continuously distributed covariates $\bx_i = (x_{1,i}, \ldots, x_{d,i})'$ with compact support $\mathcal{X}$. The object of interest is the unknown regression function $\mu(\cdot)$ and its derivatives.

We focus on \emph{partitioning-based}, or locally-supported, series (linear sieve) least squares regression estimators. Let $\Delta=\{\delta_l: 1\leq l \leq \bar{\kappa}\}$ be a collection of open and disjoint subsets of $\mathcal{X}$ such that the closure of their union is $\mathcal{X}$.
Every $\delta_l\in\Delta$ is restricted to be polyhedral.
An important special case is tensor-product partitions, in which case the support of the regressors is of tensor product form and each dimension of $\mathcal{X}$ is partitioned marginally into intervals.

Based on a particular partition, the dictionary of $K$ basis functions, each of order $m$ (e.g., $m=4$ for cubic splines) is denoted by
\[\bx_i \mapsto \bp(\bx_i) := \bp(\bx_i;\Delta,m)=(p_1(\bx_i;\Delta,m),\cdots,p_K(\bx_i;\Delta,m))'.\]
For a point $\bx \in \mathcal{X}$ and $\bq = (q_1, \ldots, q_d)'\in\mathbb{Z}_+^d$, the partial derivative $\partial^\bq\mu(\bx)$ is estimated by least squares regression
\[\widehat{\partial^\bq\mu}(\bx)=\partial^\bq \bp(\bx)'\widehat{\bbeta},\qquad\qquad   		  
\widehat{\bbeta}\in\underset{\mathbf{b}\in\mathbb{R}^K}{\argmin}\,\sum_{i=1}^{n}(y_i - \bp(\bx_i)'\mathbf{b})^2,\]
Our assumptions will guarantee that the sample matrix $\E_n[\bp(\bx_i) \bp(\bx_i)']$ is nonsingular with probability approaching one in large samples, and thus we write the estimator as
\begin{equation}
\begin{split}
\label{eqn:mu hat US}
& \widehat{\partial^\bq\mu}_0(\bx) := \bg{0}'\E_n\left[\bPi{0}(\bx_i) y_i\right]  	,   \quad \text{ where }  \\
 &  \bg{0}' := \partial^\bq \bp(\bx)'  \E_n[ \bp(\bx_i) \bp(\bx_i)']^{-1}
\quad \text{ and } \quad
\bPi{0}(\bx_i) := \bp(\bx_i).	
\end{split}
\end{equation}
The order $m$ of the basis is usually fixed in practice, and thus the tuning parameter for this class of nonparametric estimators is $\Delta$. As $n\rightarrow\infty$, $\bar{\kappa}\rightarrow\infty$, and the volume of each $\delta_l$ shrinks proportionally to $h^d$, where $h=\max\{\diam(\delta): \delta\in\Delta\}$.

\subsection{Assumptions}

We list our main assumptions for completeness. Detailed discussion can be found in the main paper. 
The first assumption concerns the data generating process.

\begin{assumption}[Data Generating Process]\label{assumption: DGP} \leavevmode
	\begin{enumerate}
		\item $\{(y_i, \bx_i'):1\leq i\leq n\}$ are i.i.d.\ satisfying (\ref{equation: regression model}), where $\bx_i$ has compact connected support $\mathcal{X}\subset\mathbb{R}^d$ and an absolutely continuous distribution function. The density of $\bx_i$, $f(\cdot)$, and the conditional variance of $y_i$ given $\bx_i$, $\sigma^2(\cdot)$, are bounded away from zero and continuous.
		
		\item $\mu(\cdot)$ is $S$-times continuously differentiable, for $S \geq [\bq]$, and all $\partial^{\bvsig}\mu(\cdot)$, $[\bvsig]=S$, are H\"older continuous with exponent $\varrho>0$.
	\end{enumerate}
\end{assumption}

\begin{assumption}[Quasi-Uniform Partition] \label{assumption: quasi-uniform partition}
	The ratio of the sizes of inscribed and circumscribed balls of each $\delta\in\Delta$ is bounded away from zero uniformly in $\delta\in\Delta$, and
	\[\frac{\max\{\diam(\delta): \delta\in\Delta\}}{\min\{\diam(\delta):\delta\in\Delta\}}\lesssim 1,\]
	where $\diam(\delta)$ denotes the diameter of $\delta$. Further, for $h=\max\{\diam(\delta):\delta\in\Delta\}$, assume $h=o(1)$.
\end{assumption}

The next assumption requires that the basis is \textit{locally} supported. We employ the notion of \emph{active basis}: a function $p(\cdot)$ on $\mathcal{X}$ is \emph{active} on $\delta\in\Delta$ if it is not identically zero on $\delta$. 
\begin{assumption}[Local Basis]\label{assumption: Local Basis} \leavevmode
	\begin{enumerate}
		\item For each basis function $p_k$, $k=1,\ldots, K$, the union of elements of $\Delta$ on which $p_k$ is active is a connected set, denoted by $\mathcal{H}_k$. For all $k=1,\ldots, K$, both the number of elements of $\mathcal{H}_k$ and the number of basis functions which are active on $\mathcal{H}_k$ are bounded by a constant.
		
		\item For any $\mathbf{a}=(a_1,\cdots, a_K)'\in\mathbb{R}^{K}$,
		\[\mathbf{a}'\int_{\mathcal{H}_k}\bp(\bx;\Delta, m)\bp(\bx; \Delta, m)'\,d\bx\, \mathbf{a}\gtrsim a_k^2h^d,
		\qquad k=1,\ldots,K.
		\]
		
		\item Let $[\bq] < m$. For an integer $\varsigma \in [[\bq], m)$, for all $\bvsig, [\bvsig]\leq \varsigma$, 
		\[h^{-[\bvsig]} \lesssim \inf_{\delta\in\Delta}\inf_{\bx\in\clo({\delta})}\|\partial^{\bvsig}\bp(\bx;\Delta,m)\|
		\leq     \sup_{\delta\in\Delta}\sup_{\bx\in\clo({\delta})}\|\partial^{\bvsig}\bp(\bx;\Delta,m)\|
		\lesssim h^{-[\bvsig]}
		\]
		where $\clo({\delta})$ is the closure of $\delta$, and for $[\bvsig]=\varsigma+1$,
		\[\sup_{\delta\in\Delta}\sup_{\bx\in\clo({\delta})}\|\partial^{\bvsig}\bp(\bx;\Delta,m)\|\lesssim h^{-\varsigma-1}.\]
	\end{enumerate}
\end{assumption}
We remind readers that Assumption \ref{assumption: quasi-uniform partition} and \ref{assumption: Local Basis} implicitly relate the number of approximating series terms, the number of cells in $\Delta$, and the maximum mesh size: $K\asymp \bar{\kappa}\asymp h^{-d}$.

The next assumption gives an explicit high-level expression of the leading approximation error which is needed for bias correction and integrated mean squared error (IMSE) expansion. To simplify notation, for each $\bx\in\mathcal{X}$ define $\delta_{\bx}$ as the element of $\Delta$ whose closure contains $\bx$ and $h_\bx$ for the diameter of this $\delta_{\bx}$.
\begin{assumption}[Approximation Error]\label{assumption: Misspecification Error}
	Let $S \geq m$. For all $\bvsig$ satisfying $[\bvsig]\leq \varsigma$ in Assumption \ref{assumption: Local Basis}, there exists $s^*\in\mathcal{S}_{\Delta,m}$, the linear span of $\bp(\bx;\Delta,m)$, and \[\mathscr{B}_{m,\bvsig}(\bx)=-\sum_{\bm{u}\in\Lambda_m}\partial^{\bm{u}}\mu(\bx)h_\bx^{m-[\bvsig]}B_{\bm{u},\bvsig}(\bx)\]
	such that
	\begin{equation}\label{eq: Approximation Power with Bias}
	\underset{\bx\in\mathcal{X}}{\sup}\;|\partial^{\bvsig} \mu(\bx)-\partial^{\bvsig} s^*(\bx)+\mathscr{B}_{m,\bvsig}(\bx)| \lesssim h^{m+\varrho-[\bvsig]}
	\end{equation}
	and
	\begin{equation}
	\sup_{\delta\in\Delta}\sup_{\bx_1,\bx_2\in\clo({\delta})}\frac{|B_{\bm{u},\bvsig}(\bx_1)-B_{\bm{u},\bvsig }(\bx_2)|}{\|\bx_1-\bx_2\|}\lesssim h^{-1}
	\end{equation}
	where $B_{\bm{u},\bvsig}(\cdot)$ is a known function that is bounded uniformly over $n$, and 
	$\Lambda_m$ is a multi-index set, which depends on the basis, with $[\bm{u}]=m$ for $\bm{u}\in\Lambda_m$.
\end{assumption}

Our last assumption concerns the basis used for bias correction. Specifically, for some $\tilde{m} > m$, let $\bpt(\bx) := \bpt(\bx; \tilde{\Delta}, \tilde{m})$ be a basis of order $\tilde{m}$  defined on a partition $\tilde{\Delta}$ with maximum mesh $\tilde{h}$. Objects accented with a tilde always pertain to this secondary basis and partition for bias correction. See the main text for complete discussion.

\begin{assumption}[Bias Correction]\label{assumption: BC bases}
	The partition $\tilde{\Delta}$ satisfies Assumption \ref{assumption: quasi-uniform partition}, with maximum mesh $\tilde{h}$, and the basis $\bpt(\bx;\tilde{\Delta}, \tilde{m})$,  $\tilde{m} > m$, satisfies Assumptions \ref{assumption: Local Basis} and \ref{assumption: Misspecification Error} with $\tilde{\varsigma} = \tilde{\varsigma}(\tilde{m}) \geq m$ in place of $\varsigma$. Let $\rho := h/\tilde{h}$, which obeys $\rho \rightarrow \rho_0\in(0,\infty)$. In addition, for $j=3$, either (i) $\bpt(\bx;\tilde{\Delta}, \tilde{m})$ spans a space containing the span of $\bp(\bx;\Delta, m)$, and for all $\bm{u}\in\Lambda_{m}$, $\partial^{\bm{u}}\bp(\bx;\Delta, m)=\mathbf{0}$; or (ii) both $\bp(\bx;\Delta, m)$ and $\bpt(\bx;\tilde{\Delta}, \widetilde{m})$ reproduce polynomials of degree 
	$[\bq]$.
\end{assumption}


\section{Technical Lemmas}\label{section: lemmas}
We present a series of technical lemmas that will be used for pointwise and uniform analysis. To begin with, we introduce additional notation to simplify our expressions. The following are frequently used outer-product matrices:
\begin{align*}
&\bQ_m =\E[\bp(\bx_i)\bp(\bx_i)'],\qquad \widehat{\bQ}_m = \E_n[\bp(\bx_i)\bp(\bx_i)'],\\
&\bQ_{\tilde{m}} =\E[\bpt(\bx_i)\bpt(\bx_i)'],\qquad \widehat{\bQ}_{\tilde{m}}=\E_n[\bpt(\bx_i)\bpt(\bx_i)'],\\
&\bQ_{m,\tilde{m}}=\E[\bp(\bx_i)\bpt(\bx_i)'],\quad\,
\widehat{\bQ}_{m,\tilde{m}} = \E_n[\bp(\bx_i)\bpt(\bx_i)'].
\end{align*}

In the main paper, we define four partitioning-based nonparametric estimators, numbered as $j=0,1,2,3$, which can be written in exactly the same form:
\[ \widehat{\partial^\bq\mu}_j(\bx) := \bg{j}'\E_n[\bPi{j}(\bx_i) y_i ].
\]
We summarize the definitions of these quantities here. For $j=0$ (classical estimator without bias correction),
\begin{equation} \label{eq: classical est}
\bg{0}' := \partial^\bq\bp(\bx)'\widehat{\bQ}_m^{-1} \quad
\text{and}\quad
\bPi{0}(\bx_i) := \bp(\bx_i) .
\end{equation}

For $j=1$ (high-order-basis bias correction),
\begin{equation} \label{eq: higher-order bc}
\bg{1}':=\partial^\bq\bpt(\bx)'\widehat{\bQ}_{\tilde{m}}^{-1}, \quad
\text{and} \quad
\bPi{1}(\bx_i) := \bpt(\bx_i).
\end{equation}

For $j=2$ (least squares bias correction),
\begin{equation}\label{eq: ls bc}
\begin{split}
&\widehat{\bgamma}_{\bq,2}(\bx)' := \Big(\bg{0}',-\bg{0}'\widehat{\bQ}_{m, \tilde{m}} \widehat{\bQ}_{\tilde{m}}^{-1} + \bg{1}'\Big) \\
&\text{ and } \quad
\bPi{2}(\bx_i) := \left( \bp(\bx_i)', \bpt(\bx_i)'\right)'.
\end{split}
\end{equation}

For $j=3$ (plug-in bias correction),
\begin{equation}\label{eq: plug in bc}
\begin{split}
&\widehat{\bgamma}_{\bq,3}(\bx)'
= \Big(\bg{0}',
\sum_{\bm{u}\in\Lambda_{m}} \Big\{
\widehat{\bgamma}_{\bm{u},1}(\bx)' h_\bx^{m-[\bq]}B_{\bm{u},\bq}(\bx) \\
&\hspace{1.3in}- \widehat{\bgamma}_{\bq,0}'\E_n[\bp(\bx_i)h_{\bx_i}^{m}B_{\bm{u},\mathbf{0}}(\bx_i)\widehat{\bgamma}_{\bm{u},1}(\bx_i)']
\Big\}
\Big),\\
&\text{ and } \qquad
\bPi{3}(\bx_i) := \left( \bp(\bx_i)', \bpt(\bx_i)'\right)'.
\end{split}
\end{equation}

For $j=0, 1, 2, 3$, $\bgamma_{\bq, j}(\bx)$ is defined as $\bg{j}$ in \eqref{eq: classical est}, \eqref{eq: higher-order bc}, \eqref{eq: ls bc} and \eqref{eq: plug in bc} but with sample averages replaced by their population counterparts. Finally, define
\begin{equation} \label{eq: Sig}
\begin{split}
&\bSig{j}=\E[\bPi{j}(\bx_i)\bPi{j}(\bx_i)'\sigma^2(\bx_i)], \qquad
\bSigbar{j}= \E_n[\bPi{j}(\bx_i)\bPi{j}(\bx_i)'\sigma^2(\bx_i)] \\
&\bSighat{j}=\E_n[\bPi{j}(\bx_i)\bPi{j}(\bx_i)'\widehat{\e}_{i,j}^2], \qquad
\widehat{\e}_{i,j}=y_i-\widehat{\mu}_j(\bx_i),  \\
& \Omega_j(\bx)=\bgamma_{\bq, j}(\bx)'\bSig{j}\bgamma_{\bq, j}(\bx), \quad \text{and} \quad 
\widehat{\Omega}_j(\bx)=\bg{j}'\bSighat{j}\bg{j}.
\end{split}
\end{equation}

The next lemma establishes the rate of convergence and boundedness for $\widehat{\bQ}_m$ and other relevant matrices, which will be crucial for both pointwise and uniform analysis. Since the orders of bases used to construct bias corrected estimators are fixed and mesh size ratio $\rho\rightarrow \rho_0\in(0,\infty)$ by Assumption \ref{assumption: BC bases}, the same conclusions also apply to $\widehat{\bQ}_{\tilde{m}}$, $\bQ_{\tilde{m}}$ and other related quantities, though we do not make it explicit in the statement.

\begin{lem}\label{lem: Q matrix}
	Under Assumptions \ref{assumption: DGP}--\ref{assumption: Local Basis},
	\begin{equation*}
	h^d\lesssim\lambda_{\min}(\bQ_m)\leq\lambda_{\max}(\bQ_m)\lesssim h^d,\quad \|\bQ_m^{-1}\|_\infty\lesssim h^{-d}.
	\end{equation*}
	Moreover, when $\frac{\log n}{nh^d}=o(1)$, we have
	\begin{align*}
	&\|\widehat{\bQ}_m-\bQ_m\|_\infty\lesssim_\P h^d\sqrt{\frac{\log n}{nh^d}}=o_\P(h^d),\quad
	\|\widehat{\bQ}_m-\bQ_m\|\lesssim_\P h^{d}\sqrt{\frac{\log n}{nh^d}}=o_\P(h^d),\\
	&\|\widehat{\bQ}_m\|\lesssim_\P h^d,\quad
	\|\widehat{\bQ}_m^{-1}\|\lesssim_\P h^{-d}, \quad
	\|\widehat{\bQ}_m^{-1}\|_\infty\lesssim_\P h^{-d},\\
	&\|\widehat{\bQ}_m^{-1}-\bQ_m^{-1}\|_\infty \lesssim_\P h^{-d}\sqrt{\frac{\log n}{nh^d}}=o_\P(h^{-d}),\\
	&\|\widehat{\bQ}_m^{-1}-\bQ_m^{-1}\|\lesssim_\P h^{-d}\sqrt{\frac{\log n}{nh^d}}=o_\P(h^{-d}),\quad \text{and}\\
	&\|\bSigbar{0}-\bSig{0}\| \lesssim_\P h^d\sqrt{\frac{\log n}{nh^d}}=o_\P(h^d).
	\end{align*}
\end{lem}

The next lemma shows that the asymptotic error expansion specified in Assumption \ref{assumption: Misspecification Error} translates into the bias expansion for the classical point estimator $\widehat{\partial^\bq\mu}_0(\cdot)$, which is (partly) reported in Lemma 3.1 of the main paper. Importantly, the following lemma also gives a \textit{uniform} bound on the conditional bias which is crucial for least squares bias correction ($j=2$) and uniform inference.
 
\begin{lem}[Conditional Bias] \label{lem: bias}
	Let Assumptions \ref{assumption: DGP}--\ref{assumption: Misspecification Error} hold. If $\frac{\log n}{nh^d}=o(1)$, then,
	\[\E[\widehat{\partial^\bq\mu}_0(\bx)|\bX]-\partial^\bq\mu(\bx)
	=\mathscr{B}_{m,\bq}(\bx) - 
	\bg{0}'\E_n[\bPi{0}(\bx_i)\mathscr{B}_{m,\mathbf{0}}(\bx_i)]+O_\P(h^{m+\varrho-[\bq]}),
	\]
	and
	\[\sup_{\bx\in\mathcal{X}}\Big|\E[\widehat{\partial^\bq\mu}_0(\bx)|\bX]-\partial^\bq\mu(\bx)\Big| \lesssim_\P h^{m-[\bq]}.\]
	Moreover, if the following condition is also satisfied,
	\begin{equation}\label{eq: approx orth}
	\max_{1\leq k\leq K}\int_{\mathcal{H}_k} p_k(\bx;\Delta, m)\mathscr{B}_{m,\mathbf{0}}(\bx)\,d\bx=o(h^{m+d}),
	\end{equation}
	then
	$\|\bg{0}'\E_n[\bPi{0}(\bx_i)\mathscr{B}_{m,\mathbf{0}}(\bx_i)]\|_{L_\infty(\mathcal{X})}=o_\P(h^{m-[\bq]})$.
\end{lem}
Equation \eqref{eq: approx orth} implies that the leading approximation error is approximately orthogonal to $\bp(\cdot)$, in which case the $L_\infty$ approximation error coincides with the leading smoothing bias of $\widehat{\partial^\bq\mu}_0(\cdot)$.

The next lemma shows that $\bg{j}$ in \eqref{eq: classical est}, \eqref{eq: higher-order bc}, \eqref{eq: ls bc}, and \eqref{eq: plug in bc} converge to their population counterpart $\bgamma_{\bq, j}(\bx)$ in a proper sense, and all the $\bgamma_{\bq, j}(\bx)$ are sufficiently controlled in terms of both $L_\infty$- and $L_2$-operator norms. 

\begin{lem}\label{lem: basis}
	Let Assumptions \ref{assumption: DGP}, \ref{assumption: quasi-uniform partition}, \ref{assumption: Local Basis}, and \ref{assumption: BC bases} hold. If $\frac{\log n}{nh^d}=o(1)$, then,
	\begin{equation*}
	\begin{split}
	&\sup_{\bx\in\mathcal{X}}\|\bgamma_{\bq,j}(\bx)'\|_\infty\lesssim h^{-d-[\bq]},\quad \sup_{\bx\in\mathcal{X}}\|\bg{j}'-\bgamma_{\bq,j}(\bx)'\|_\infty\lesssim_\P h^{-d-[\bq]}\sqrt{\frac{\log n}{nh^d}},\\
	&\sup_{\bx\in\mathcal{X}}\|\bgamma_{\bq,j}(\bx)'\|\lesssim h^{-d-[\bq]},\quad\;\;\,
	\inf_{\bx\in\mathcal{X}}\|\bgamma_{\bq,j}(\bx)'\|\gtrsim h^{-d-[\bq]}, \qquad \text{and}\\
	&\sup_{\bx\in\mathcal{X}}\|\bg{j}'-\bgamma_{\bq,j}(\bx)'\|\lesssim_\P h^{-d-[\bq]}\sqrt{\frac{\log n}{nh^d}}.
	\end{split}
	\end{equation*}	
\end{lem}

The last lemma in this section proves that the asymptotic variance of classical and bias-corrected estimators is properly bounded. Importantly, we need not only an upper bound, but also a lower bound since the variance term appears in the denominator of $t$-statistics.

\begin{lem} \label{lem: Avar}
	Let Assumptions \ref{assumption: DGP}, \ref{assumption: quasi-uniform partition}, \ref{assumption: Local Basis}, and \ref{assumption: BC bases} hold. If $\frac{\log n}{nh^d}=o(1)$,
	\[\sup_{\bx\in\mathcal{X}}
	\Omega_j(\bx) \lesssim h^{-d-2[\bq]} \quad\text{and}\quad
	\inf_{\bx\in\mathcal{X}}\Omega_j(\bx)\gtrsim h^{-d-2[\bq]}.
	\]
\end{lem}
This lemma also implies that bias correction does not change the asymptotic order of the variance.

\section{Integrated Mean Squared Error} \label{section: IMSE rec}
We now state an integrated mean squared error (IMSE) result that specializes Theorem 4.1 in the main paper to the case of general tensor-product partitions but is less restricted than Theorem 4.2. In particular, we allow for any $\bq$ and do not rely on Assumption 6(c). This result gives the key ingredient for characterizing the limit $\mathscr{V}_{\Delta,\bq}$ and $\mathscr{B}_{\Delta,\bq}$ in some generality.

To be specific, throughout this section, we assume the support of the regressors is of tensor-product form, that each dimension of $\mathcal{X}$ is partitioned marginally into intervals, and $\Delta$ is the tensor product of these intervals. Let $\mathcal{X}_\ell = [\underline{x}_\ell, \overline{x}_\ell]$ be the support of covariate $\ell = 1, 2, \ldots, d$, and partition this into $\kappa_\ell$ disjoint subintervals defined by $\{\underline{x}_\ell=t_{\ell,0}<t_{\ell,1}<\cdots<t_{\ell,\kappa_\ell-1}<t_{\ell,\kappa_\ell}=\overline{x}_\ell\}$. If this partition of $\mathcal{X}_\ell$ is $\Delta_\ell$, then a complete partition of $\mathcal{X}$ can be formed by tensor products of the one-dimensional partitions: $\Delta=\otimes_{\ell=1}^d\Delta_\ell$, with $\bkappa=(\kappa_1,\kappa_2,\dots,\kappa_d)'$ collecting the number of subintervals in each dimension of $\bx_i$, and $\bar{\kappa} = \kappa_1\kappa_2\cdots\kappa_d$. A generic cell of this partition is the rectangle
\begin{equation}
	\label{eqn:rect}
	\delta_{l_1\ldots l_d}=\{\bx:t_{\ell ,l_\ell }< x_\ell < t_{\ell ,l_\ell +1},\quad 0\leq l_\ell \leq\kappa_\ell -1 \quad 1\leq \ell \leq d \}.
\end{equation}
Assumption \ref{assumption: quasi-uniform partition} can be verified by choosing the knot positions appropriately, often dividing $\mathcal{X}_\ell$ uniformly or by quantiles.

As before, $\delta_{\bx}$ is the subrectangle in $\Delta$ containing $\bx$, and we write $\bb_{\bx}$ for the vector collecting the interval lengths of $\delta_{\bx}$ (as defined in Assumption 6). $\bt_{\bx}^L$ 
denotes the start point of $\delta_{\bx}$. For $\ell=1, \ldots, d$, for a generic cell $\delta_{l_1\ldots l_d}$ as in \eqref{eqn:rect}, we write $b_{\ell, l_\ell}=t_{\ell, l_\ell+1} -
t_{\ell, l_\ell}$, $l_\ell=0, \ldots, \kappa_\ell-1$, $b_\ell=\max_{0\leq l_\ell\leq\kappa_\ell-1} b_{\ell,l_\ell}$ and $b=\max_{1\leq \ell \leq d} b_\ell$ ($b\asymp h$ by Assumption \ref{assumption: quasi-uniform partition}).

\begin{thm}[IMSE for Tensor-Product Partitions]\label{thm: IMSE-R}
	Suppose that the conditions in Theorem 4.1 and Assumptions 6(a)--(b) hold. Then, for any arbitrary  sequence of points $\{\btau_k\}_{k=1}^K$ such that $\btau_k\in\supp(p_k(\cdot))$ for each $k=1,\ldots, K$,
	\[\mathscr{V}_{\Delta,\bq} = \mathscr{V}_{\bkappa,\bq} + o(h^{-d-2[\bq]})
	\qquad\text{and}\qquad
	\mathscr{B}_{\Delta,\bq} = \mathscr{B}_{\bkappa,\bq} + o(h^{2m-2[\bq]}), \quad \text{with}
	\]
	\begin{align*}
		\mathscr{V}_{\bkappa,\bq}
		&= \bkappa^{\mathbf{1}+2\bq} \sum_{k=1}^K \Big[\frac{\sigma^2(\btau_k)w(\btau_k)}{f(\btau_k)} \prod_{\ell=1}^d g_\ell(\btau_k)\Big] 
		\be_k'\bH_\mathbf{0}^{-1}\bH_\bq\be_k \vol(\delta_{\btau_k})
		\asymp h^{-d-2[\bq]},\\
		\mathscr{B}_{\bkappa,\bq}
		&= \sum_{\bm{u}_1,\bm{u}_2 \in\Lambda_m} \bkappa^{-(\bm{u}_1+\bm{u}_2-2\bq)}
		\Big\{\mathscr{B}_{\bm{u}_1,\bm{u}_2,\bq}
		+ \bv_{\bm{u}_1,\mathbf{0}}'\bH_\mathbf{0}^{-1}\bH_\bq\bH_\mathbf{0}^{-1}\bv_{\bm{u}_2,\mathbf{0}}
		- 2\mathbf{v}_{\bm{u}_1,\bq}'\bH_\mathbf{0}^{-1}\bv_{\bm{u}_2,\mathbf{0}}\Big\}\\
		&\lesssim h^{2m-2[\bq]},
	\end{align*}
	where the $K$-dimensional vector $\be_k$ is the $k$-th unit vector (i.e., $\be_k$ has a $1$ in the $k$-th position and $0$ elsewhere), $\mathbf{g}(\bx)=(g_1(\bx),\dots,g_d(\bx))'$ is discussed below, 
	\[\bH_\bq = \bkappa^{-2\bq}\int_{[0,1]^d}\partial^\bq\bp(\bx)\partial^\bq\bp(\bx)'d\bx,\]
	\[\mathscr{B}_{\bm{u}_1,\bm{u}_2,\bq}
	= \eta_{\bm{u}_1,\bm{u}_2,\bq} \int_{[0,1]^d}\frac{\partial^{\bm{u}_1}\mu(\bx)\partial^{\bm{u}_2}\mu(\bx)}{\mathbf{g}(\bx)^{\bm{u}_1+\bm{u}_2-2\bq}}w(\bx)d\bx,
	\]
	and the $K$-dimensional vector $\bv_{\bm{u},\bq}$ has $k$-th typical element
	\[\frac{\sqrt{w(\btau_k)}\partial^{\bm{u}}\mu(\btau_k)}{\bkappa^{\bq}\mathbf{g}(\btau_k)^{\bm{u}-\bq}}
	\int_{\mathcal{X}}\frac{h_\bx^{m-[\bq]}}{\bb_{\bx}^{\bm{u}-\bq}}\partial^\bq p_k(\bx)B_{\bm{u},\bq}(\bx)d\bx.
	\]
\end{thm}

This theorem takes advantage of the assumed tensor-product structure of the partition $\Delta$ to express the leading bias and variance as proportional to the number of subintervals used for each regressor. The accompanying constants are expressed as sums over the local contributions of the basis function $\bp(\cdot)$ used, and are easy to be shown bounded in general (they may still not converge to any well-defined limit at this level of generality). 

Assumption 6(a) used in Theorem \ref{thm: IMSE-R} (and in Theorem 4.2) slightly strengthens the quasi-uniform condition imposed in Assumption \ref{assumption: quasi-uniform partition}. It can be verified by choosing knot positions appropriately. Specifically, let $\Delta$ be a tensor-product partition (see Equation \eqref{eqn:rect}) with (marginal) knots chosen as
\[t_{\ell,l} = G_\ell^{-1}\left(\frac{l}{\kappa_\ell}\right), \qquad l=0,1,\dots,\kappa_\ell,\qquad \ell=1,2,\dots,d,\]
where $G_\ell(\cdot)$ is a univariate continuously differentiable distribution function and $G_\ell^{-1}(v)=\inf\{x\in\mathbb{R}:G_\ell(x)\geq v\}$. In this case the function $g_\ell(\cdot)$ in Assumption 6(a) is simply the density of $G_\ell(\cdot)$. Two examples commonly used in practice are:
(i) evenly-spaced partitions, denoted by $\Delta_\mathtt{ES}$, where $G_\ell(x)=x$ and $g_\ell(x)=1$, and (ii) quantile-spaced partitions, denoted by $\Delta_\mathtt{QS}$, where $G_\ell(x)=\widehat{F}_\ell(x)$ with $\widehat{F}_\ell(x)$ the empirical distribution function for the $\ell$-th covariate, $\ell=1,2,\dots,d$. For the case of quantile-spaced partitioning, if $\widehat{F}_\ell(x)$ convergences to $F_\ell(x)=\P[x_\ell\leq x]$ in a suitable sense, $g_\ell(x)=dF_\ell(x)/dx$, i.e., the marginal density of $x_\ell$. See \cite{Agarwal-Studden_1980_AoS} and \cite{Barrow-Smith_1978_QAM} for a slightly more high-level condition in the context of univariate $B$-splines.

Assumption 6(b) in Theorem \ref{thm: IMSE-R} (and in Theorem 4.2) involves a scaling factor $\frac{h_\bx^{2m-2[\bq]}}{\bb_\bx^{\bm{u}_1+\bm{u}_2-2[\bq]}}$. In Section \ref{section: setup} we do not add specific restrictions on the shape of cells in $\Delta$ and thus the diameter of a cell is used to conveniently express the order of approximation error (denoted by $h$), but for tensor-product partitions the approximation error is usually characterized by the lengths of intervals on each axis (denoted by $\bb_\bx$), as illustrated in Section \ref{section: examples}. Therefore, the scaling factor makes the above theorem immediately apply to the three leading examples in Section \ref{section: examples} based on tensor-product partitioning schemes (splines, wavelets, and piecewise polynomials).

Using Theorem \ref{thm: IMSE-R}, the IMSE-optimal tuning parameter selector for partitioning-based series estimators becomes
\[\bkappa_{\mathtt{IMSE}, \bq} = \argmin_{\bkappa\in\mathbb{Z}^d_{++}} \left\{ \frac{1}{n} \mathscr{V}_{\bkappa,\bq} + \mathscr{B}_{\bkappa,\bq} \right\},
\]
which is still given in implicit form because the leading constants are not known to converge at this level of generality. It follows that $\bkappa_{\mathtt{IMSE},\bq}=(\kappa_{\mathtt{IMSE},\bq,1},\dots,\kappa_{\mathtt{IMSE},\bq,d})'$ with $\kappa_{\mathtt{IMSE},\bq,\ell} \asymp n^{\frac{1}{2m+d}}$, $\ell=1,2,\cdots,d$.

\section{Rates of Convergence} \label{section: convergence rates}

In this section, we discuss rates of convergence for the classical estimator and three bias-corrected estimators. Some results, such as the $L_2$ and uniform convergence rates for $\widehat{\partial^{\bq}\mu}_0(\bx)$, are also reported in Theorem 4.3 of the main paper.

To begin with, we first report pointwise linearization which is an important intermediate step towards pointwise asymptotic normality.

\begin{lem}[Pointwise Linearization]\label{lem: pointwise linearization}
	Let Assumptions \ref{assumption: DGP}--\ref{assumption: BC bases} hold. If $\frac{\log n}{nh^d}=o(1)$, then
	for each $\bx\in\mathcal{X}$ and $j=0,1,2,3$,
	\begin{align*}
	&\widehat{\partial^\bq \mu}_j(\bx)-\partial^\bq \mu(\bx)=\bgamma_{\bq,j}(\bx)'\E_n[\bPi{j}(\bx_i)\e_i]+R_{1n,\bq}(\bx)+R_{2n,\bq}(\bx),
	\quad \text{where}\\
	&R_{1n,\bq}(\bx):=\left(\bg{j}-\bgamma_{\bq,j}(\bx)'\right)
	\E_n[\bPi{j}(\bx_i)\e_i] \lesssim_\P 
	\frac{\sqrt{\log n}}{nh^{d+[\bq]}},\\
	&R_{2n,\bq}(\bx):=\E[\widehat{\partial^\bq \mu}_j(\bx)|\bX]-\partial^\bq \mu(\bx)\lesssim_\P h^{m-[\bq]}.
	\end{align*}
	Furthermore, for $j=1,2,3$, $R_{2n,\bq}(\bx)\lesssim_\P h^{m+\varrho-[\bq]}$.
\end{lem}
We can use Lemma \ref{lem: pointwise linearization} to establish the limiting distribution of the standardized $t$-statistics as given in Theorem 5.1 of the main paper. The proof is available in Section \ref{proof: Theorem 5.1}. 
\begin{remark}
	The term denoted by $R_{2n,\bq}(\bx)$ captures the conditional bias of the point estimator. We establish a sharp bound on it by applying Bernstein's maximal inequality to control the \emph{largest} element of $\E_n[\bPi{0}(\bx_i)\mathscr{B}_{m, \mathbf{0}}(\bx_i)]$ and employing the bound on the \emph{uniform norm} of $\widehat{\bQ}_m^{-1}$ derived in Lemma \ref{lem: Q matrix}. This result improves on previous results in the literature and, in particular, confirms a conjecture posed by \cite[][Comment 4.2(ii), p. 352]{Belloni-Chernozhukov-Chetverikov-Kato_2015_JoE}, for partitioning-based series estimation in general. To be precise, in their setup, $\ell_k c_k$ can be understood as an uniform bound on the $L_2$ approximation error (orthogonal to the approximating basis). Thus, using our results, we obtain 
	$$\bp(\bx)'(\widehat{\bQ}^{-1}-\bQ^{-1})\G_n[\bp(\bx_i)\ell_kc_k]\lesssim_\P
	\sqrt{\frac{\log n}{nh^d}}\cdot\sqrt{\frac{\log n}{h^d}}\ell_kc_k,$$
	which coincides with Equation (4.15) of \cite{Belloni-Chernozhukov-Chetverikov-Kato_2015_JoE}, up to a normalization, thereby improving on the approximation established in their Equation (4.12). See our proof of Lemmas \ref{lem: bias} and \ref{lem: pointwise linearization} for more details.
\end{remark}

For uniform convergence, we need to upgrade the above lemma to a uniform version at the cost of some stronger conditions:
\begin{lem}[Uniform Linearization]\label{lemma: uniform linearization}
	Let Assumptions \ref{assumption: DGP}--\ref{assumption: BC bases} hold. 
	For each $j=0,1,2,3$, define $R_{1n,\bq}(\bx)$ and $R_{2n,\bq}(\bx)$ as in Lemma \ref{lem: pointwise linearization}.
	Assume one of the following holds:
	\begin{enumerate}[label=(\roman*)]
		\item $\E[|\e_i|^{2+\nu}]<\infty$ for some $\nu>0$, and $\frac{n^{\frac{2}{2+\nu}}(\log n)^{\frac{2\nu}{4+2\nu}}}
		{nh^d} \lesssim 1$; or
		\item $\E[|\e_i|^3\exp(|\e_i|)]<\infty$, and $\frac{(\log n)^3}{nh^d}\lesssim 1$.
	\end{enumerate} 	
	Then,
	$\sup_{\bx\in\mathcal{X}}\,|R_{1n,\bq}(\bx)|\lesssim_\P
	\frac{\log n}{nh^{d+[\bq]}} =: \bar{R}_{1n,\bq}$.
	For $j=0$, $\sup_{\bx\in\mathcal{X}}\,|R_{2n,\bq}(\bx)|\lesssim_\P h^{m-[\bq]} =: \bar{R}_{2n,\bq}$, and for $j=1,2,3$,
	$\sup_{\bx\in\mathcal{X}}\,|R_{2n,\bq}(\bx)|\lesssim_\P h^{m+\varrho-[\bq]} =: \bar{R}_{2n,\bq}$. 
\end{lem}

Given Lemma \ref{lem: pointwise linearization} and \ref{lemma: uniform linearization}, we are ready to derive the desired rates of convergence.
\begin{thm}[Convergence Rates] \label{thm: uniform convergence}
	Let Assumption \ref{assumption: DGP}, \ref{assumption: quasi-uniform partition}, and \ref{assumption: Local Basis} hold. Also assume that 
	$\sup_{\bx\in \mathcal{X}}|\partial^{\bq}\mu(\bx)-\partial^{\bq}s^*(\bx)|\lesssim h^{m-[\bq]}$ with $s^*$ defined in Assumption \ref{assumption: Misspecification Error}. Then, if $\frac{\log n}{nh^d}=o(1)$,
	\[ \int_{\mathcal{X}}\Big(\widehat{\partial^{\bq}\mu}_0(\bx)
	-\partial^{\bq}\mu(\bx)\Big)^2w(\bx)d\bx\lesssim_\P 
	\frac{1}{nh^{d+2[\bq]}}+h^{2(m-[\bq])}.
	\]  
	In addition, under the conditions of Lemma \ref{lemma: uniform linearization}, for each $j=0,1,2,3$,
	\begin{align*}
	\sup_{\bx\in\mathcal{X}}\,\Big|\widehat{\partial^\bq \mu}_j(\bx)-\partial^\bq \mu(\bx)\Big|&\lesssim_\P h^{-(d/2+[\bq])}\sqrt{\frac{\log n}{n}}
	+\bar{R}_{1n,\bq}+\bar{R}_{2n,\bq}=: R_{\bq, j}^{\mathtt{uc}}
	\end{align*}
	where $\bar{R}_{1n,\bq}$ and $\bar{R}_{2n,\bq}$ are defined in Lemma \ref{lemma: uniform linearization}.
\end{thm}

Relying on the uniform convergence results, we can further show the consistency of variance estimates $\widehat{\Omega}_j(\bx)$ defined in Section \ref{section: lemmas}. The results in the following theorem are partly reported in Theorem 5.2 in the main paper.

\begin{thm}[Variance Estimate] \label{thm: variance}
	Let Assumptions \ref{assumption: DGP}--\ref{assumption: BC bases} hold. 
	\begin{enumerate}[label=(\roman*)]
		\item Assume that $\E[|\e_i|^{2+\nu}]<\infty$ for some $\nu>0$, and $\frac{n^{\frac{2}{2+\nu}}(\log n)^{\frac{2\nu}{4+2\nu}}}
		{nh^d} = o(1)$. Then for each $j=0,1,2,3$,
		\begin{align*}
		&\|\bSighat{j}-\bSigma_j\|\lesssim_\P h^d\left(R^{\mathtt{uc}}_{\mathbf{0},j}+
		\frac{n^{\frac{1}{2+\nu}}(\log n)^{\frac{\nu}{4+2\nu}}}
		{\sqrt{nh^d}}\right)=o_\P(h^d),\;\text{ and}\\
		&\sup_{\bx\in\mathbf{X}}\Big|\widehat{\Omega}_j(\bx)
		-\Omega_j(\bx)\Big|\lesssim_\P h^{-d-2[\bq]}\left(R^{\mathtt{uc}}_{\mathbf{0},j}
		+\frac{n^{\frac{1}{2+\nu}}(\log n)^{\frac{\nu}{4+2\nu}}}
		{\sqrt{nh^d}}\right)=o_\P(h^{-d-2[\bq]})
		\end{align*}
		where $R_{\mathbf{0},j}^{\mathtt{uc}}$ is the uniform convergence rate given in Theorem \ref{thm: uniform convergence} with $\bq=\mathbf{0}$.
		
		\item Assume that $\E[|\e_i|^3\exp(|\e_i|)]<\infty$, and $\frac{(\log n)^3}{nh^d} =o(1) $. Then for each $j=0,1,2,3$,
		\begin{align*}
		&\|\bSighat{j}-\bSigma_j\|\lesssim_\P h^d\left(R^{\mathtt{uc}}_{\mathbf{0},j}+
		\frac{(\log n)^{3/2}}{\sqrt{nh^d}}\right)=o_\P(h^d),\;\text{ and}\\
		&\sup_{\bx\in\mathcal{X}}\Big|\widehat{\Omega}_j(\bx)
		-\Omega_j(\bx)\Big|\lesssim_\P h^{-d-2[\bq]}\left(R^{\mathtt{uc}}_{\mathbf{0},j}+
		\frac{(\log n)^{3/2}}{\sqrt{nh^d}}\right)=o_\P(h^{-d-2[\bq]}).
		\end{align*}
	\end{enumerate} 
\end{thm}


\section{Uniform Inference}\label{section: uniform}
\subsection{Strong Approximation}

Now we move on to uniform inference. The following $t$-statistic processes are of interest:
\[ \widehat{T}_j(\bx)=\frac{\widehat{\partial^\bq\mu}_j(\bx)-
\mu(\bx)}{\sqrt{\widehat{\Omega}_j(\bx)/n}}, \quad \bx\in\mathcal{X},\quad j=0,1,2,3.
\]

Let $r_n$ be a positive non-vanishing sequence, which will be used to denote the approximation error rate in the following analysis. As the first step, we employ Lemma \ref{lem: basis} and Theorem \ref{thm: variance} to show that the sampling and estimation uncertainty of $\bg{j}$ and $\widehat{\Omega}_j(\bx)$ are negligible uniformly over $\bx\in\mathcal{X}$. Specifically, define
\[t_j(\bx) = \frac{\bgamma_{\bq,j}(\bx)'}{\sqrt{\Omega_j(\bx)}}\G_n[
\bPi{j}(\bx_i)\e_i], \quad \bx\in\mathcal{X}, \quad j=0,1,2,3.
\]
The following lemma, which appears as Lemma 6.1 in the main paper, shows that $\widehat{T}_j(\cdot)$ can be approximated by $t_j(\cdot)$ uniformly. 

\begin{lem}[Hats Off]\label{lem: hats off}
	Let Assumptions \ref{assumption: DGP}--\ref{assumption: BC bases} hold. Assume one of the following holds:
	\begin{enumerate}[label=(\roman*)]
		\item $\E[|\e_i|^{2+\nu}]<\infty$ for some $\nu>0$, and $\frac{n^{\frac{2}{2+\nu}}(\log n)^{\frac{2+2\nu}{2+\nu}}}
		{nh^d}=o(r_n^{-2})$; or
		\item $\E[|\e_i|^3\exp(|\e_i|)]<\infty$, and $\frac{(\log n)^4}{nh^d}=o(r_n^{-2}) $.
	\end{enumerate} 
	Furthermore, for $j=0$, assume $nh^{d+2m}=o(r_n^{-2})$; and for $j=1,2,3$, assume $nh^{d+2m+2\varrho}=o(r_n^{-2})$. Then for each $j=0,1,2,3$,
	$\sup_{\bx\in\mathcal{X}}|\widehat{T}_j(\bx) - t_j(\bx)|=o_\P(r_n^{-1})$. 
\end{lem}
The proof can be found in Section 8.2 of the main paper.
Now our task reduces to construct valid distributional approximation to $t_j(\cdot)$ for each $j=0,1,2,3$, in a proper sense. Specifically, we want to show that on a sufficiently rich enough probability space there exists a copy $t'_j(\cdot)$ of $t_j(\cdot)$ and a Gaussian process $Z_j(\cdot)$ such that
$\sup_{\bx\in\mathcal{X}}|t'_j(\bx)-Z_j(\bx)|=o_\P(r_n^{-1})$.
When such a construction is possible, we can use $Z_j(\cdot)$ to approximate the distribution of $t_j(\cdot)$, as well as $\widehat{T}_j(\cdot)$ in view of Lemma \ref{lem: hats off}. To save our notation, we denote this strong approximation by $\widehat{T}_j(\cdot)=_d Z_j(\cdot)+o_\P(r_n^{-1})$ in $\mathcal{L}^\infty(\mathcal{X})$ where $\mathcal{L}^\infty(\mathcal{X})$ refers to the set of all uniformly bounded real functions on $\mathcal{X}$ equipped with uniform norm.

We will show in the following that there exists an unconditional Gaussian approximating process
\[Z_j(\bx)=\frac{\bgamma_{\bq,j}(\bx)'\bSig{j}^{1/2}}
{\sqrt{\Omega_j(\bx)}}\bN_{K_j}, \quad \bx\in\mathcal{X},
\quad \text{for}\quad j=0,1,2,3,
\]
where $\bN_{K_j}$ is a $K_j$-dimensional standard Normal vector with $K_j=\dim(\bPi{j}(\cdot))$.
As discussed in the main paper, we may have two strategies to construct strong approximations. The next theorem employs KMT coupling techniques.

\begin{thm}[Strong Approximation: KMT Coupling]
	\label{thm: KMT coupling}
	Let the conditions of Lemma \ref{lem: hats off} hold, and for $j=2,3$, also assume that $\frac{(\log n)^{3/2}}{\sqrt{nh^d}}=o_\P(r_n^{-2})$. In addition, assume one of the following holds:
	\begin{enumerate}[label=(\roman*)]
		\item $\sup_{x\in\mathcal{X}}\E[|\e_i|^{2+\bar{\nu}}|\bx_i=\bx]<\infty$ for some $\bar{\nu}>0$, and  
		$\frac{n^{\frac{2}{2+\bar{\nu}}}}{nh^d}=o(r_n^{-2})$;
		\item 
		$\sup_{\bx\in\mathcal{X}}
		\E[|\e_i|^3\exp(|\e_i|)|\bx_i=\bx] < \infty$.
	\end{enumerate}
Then for $d=1$ and $j=0,1,2,3$,
$\widehat{T}_j(\cdot)=_d Z_j(\cdot) +o_\P(r_n^{-1})$ in  $\mathcal{L}^\infty(\mathcal{X})$.
When $d>1$ and $\bp(\cdot)$ is a Haar basis, the same result still holds for $j=0$.
\end{thm} 
  
Unfortunately, extending this to general multi-dimensional cases needs some non-trivial results which are not available in current strong approximation literature, to the best of our knowledge. Thus, when $d\geq 2$, we employ an improved version of Yurinskii's inequality due to \cite{Belloni-Chernozhukov-Chetverikov-FernandezVal_2019_JoE}, which allows us to replace the Euclidean norm in the original version of Yurinskii's inequality by a sup-norm. Since basis functions considered in this paper are locally supported, a metric of distance between random vectors in terms of sup-norm leads to a weaker rate restriction than that previously used in literature, though still stronger than what we obtained in Theorem \ref{thm: KMT coupling}. 
The result is reported in Theorem 6.2 of the main paper. We restate it in the following for completeness. The proof can be found in Section 8.4 of the main paper.
\begin{thm}[Strong Approximation: Yurinskii's coupling]
	\label{thm: Yurinskii coupling}
	Let the conditions of Lemma \ref{lem: hats off} hold. In addition, assume that $\sup_{\bx\in\mathcal{X}}\E[|\e_i|^{3}|\bx_i=\bx]<\infty$ and $\frac{(\log n)^4}{nh^{3d}}=o(r_n^{-6})$. Then for $j=0,1,2,3$,
	$\widehat{T}_j(\cdot)=_d Z_j(\cdot) +o_\P(r_n^{-1})$ in $\mathcal{L}^\infty(\mathcal{X})$.
\end{thm}

\subsection{Implementation}
$\{Z_j(\bx):\bx\in\mathcal{X}\}$ is still an infeasible approximation of $\{\widehat{T}_j(\bx):\bx\in\mathcal{X}\}$. Our next objective is to construct practicably feasible distributional approximations. To be precise, we construct Gaussian processes $\{\widehat{Z}_j(\bx):\bx\in\mathcal{X}\}$, with distributions known conditional on the data $(\by, \bX)$, and there exists a copy $\widehat{Z}'_j(\cdot)$ of $\widehat{Z}_j(\cdot)$ in a sufficiently rich probability space such that (i) $\widehat{Z}'_j(\cdot)=_d \widehat Z_j(\cdot)$ conditional on the data $(\by, \bX)$ and (ii) for some positive non-vanishing sequence $r_n$ and for all $\eta>0$,
\[\P^*\left[\sup_{\bx\in\mathcal{X}} |\widehat{Z}'_j(\bx) - Z_j(\bx)| \geq \eta r_n^{-1} \right] = o_\P(1),\]
where $\P^*[\cdot]=\P[\cdot|\by, \bX]$ denotes the probability operator conditional on the data. When such a feasible process exists, we write $\widehat{Z}_j(\cdot) =_{d^*} Z_j(\cdot) + o_{\P^*}(r_n^{-1})$ in $\mathcal{L}^\infty(\mathcal{X})$. From a practical perspective, sampling from $\widehat{Z}_j(\cdot)$, conditional on the data, is possible and provides a valid distributional approximation.

Our first construction is a direct plug-in approach using the conclusion of Theorem \ref{thm: KMT coupling} and \ref{thm: Yurinskii coupling}. All unknown objects are replaced by consistent estimators already used in the feasible $t$-statistics:
\begin{equation*}
\widehat{Z}_j(\bx) = \frac{\bg{j}' \bSighat{j}^{1/2}}{ \sqrt{ \bOhat{j}} } \bN_{K_j},
\qquad \bx\in\mathcal{X}, \qquad j=0,1,2,3
\end{equation*}
The validity of this method is shown in Theorem 6.3 of the main paper. We restate it in the following whose proof is available in Section 8.5 of the main paper. 

\begin{thm}[Plug-in Approximation] \label{thm: plug-in approx}
	Let the conditions in Lemma \ref{lem: hats off} hold. Furthermore, for $j=2,3$, 
	\begin{enumerate}[label=(\roman*)]
		\item when $\E[|\e_i|^{2+\nu}]<\infty$ for some $\nu>0$, assume that $\frac{n^{\frac{1}{2+\nu}}(\log n)^{\frac{4+3\nu}{4+2\nu}}}
	{\sqrt{nh^d}}=o(r_n^{-2})$;
	\item when $\E[|\e_i|^3\exp(|\e_i|)]<\infty$, assume that $\frac{(\log n)^{5/2}}{\sqrt{nh^d}}=o(r_n^{-2})$.
	\end{enumerate}
	Then for each $j=0,1,2,3$, $\widehat{Z}_j(\cdot)=_{d^*} Z_j(\cdot) + o_{\P^*}(r_n^{-1})$
	in $\mathcal{L}^\infty(\mathcal{X})$.
\end{thm}

Our second construction, which is not reported in the main paper, is inspired by the intermediate approximation in Lemma 8.2 of the main paper.
Again, $\bgamma_{\bq,j}(\bx)$ and $\Omega_j(\bx)$ can be simply replaced by their sample analogues, and $\zeta_i$'s can be generated by sampling from an $n$-dimensional standard Gaussian vector independent of the data. The key unknown quantity is the conditional variance function $\sigma^2(\bx)=\E[\e_i^2|\bx_i=\bx]$. If the errors are homoskedastic, one can simply let $\sigma^2(\cdot)=1$. If not, we need an estimator $\widehat{\sigma}(\cdot)$ of the conditional variance satisfying a mild condition on uniform convergence, and then we can construct
\[\widehat{z}_j(\bx)=\frac{\bg{j}'}{\sqrt{\bOhat{j}}}
\G_n[\bPi{j}(\bx_i)\zeta_i\widehat{\sigma}(\bx_i)],\quad \bx\in\mathcal{X}, \quad j=0,1,2,3.
\]
In practice one may simply construct a nonparametric estimator of $\sigma^2(\cdot)$ by using, for example, regression splines or smoothing splines. See \cite[Chapter 22.4]{Eggermont-LaRiccia_2009_Book} for more details. We do not elaborate on this issue, but the next theorem shows the validity of this approach.

\begin{thm} \label{thm: plug-in approx cond}
	Let the conditions in Lemma \ref{lem: hats off} hold, and $\hat{\sigma}^2(\cdot)$ is an estimator of $\sigma^2(\cdot)$ such that 
	$\max_{1\leq i\leq n}|\hat{\sigma}^2(\bx_i)-\sigma^2(\bx_i)|=o_\P(r_n^{-1}(\log n)^{-1/2})$. For $j=2,3$, further assume
	$\frac{(\log n)^{3/2}}{\sqrt{nh^d}}=o(r_n^{-2})$. Then for each $j=0,1,2,3$, $\widehat{z}_j(\cdot) =_{d^*} Z_j(\cdot)+o_{\P^*}(r_n^{-1})$ in $\mathcal{L}^\infty(\mathcal{X})$.
\end{thm}

Our third approach to approximating the infeasible $Z_j(\cdot)$ employs an easy-to-implement wild bootstrap procedure. Specifically, for i.i.d. bounded random variables $\{\omega_i:1\leq i\leq n\}$ with $\E[\omega_i]=0$ and $\E[\omega_i^2]=1$ independent of the data, we construct bootstrapped $t$-statistics:
\begin{equation*}
\widehat{z}^*_j(\bx) = \frac{\bg{j}'} {\sqrt{\bOstar{j}}}
\G_n[\bPi{j}(\bx_i) \omega_i \widehat{\e}_{i,j}],
\qquad \bx\in\mathcal{X}, \qquad j=0,1,2,3,
\end{equation*}
where $\widehat{\e}_{i,j}$ is defined in \eqref{eq: Sig}, and the bootstrap Studentization $\bOstar{j}$ is constructed using $\bSighat{j}^* = \E_n[\bPi{j}(\bx_i) \bPi{j}(\bx_i)' \omega_i^2 \widehat{\e}_{i,j}^2]$.

In comparison with previous plug-in approximations, the bootstrap method requires an additional Gaussian approximation step as in either Theorem \ref{thm: KMT coupling} or Theorem \ref{thm: Yurinskii coupling} (conditional on the data). These results are stated in the following two theorems. $\E^*[\cdot]$ denotes the expectation conditional on the data.

\begin{thm}[Wild Bootstrap: KMT] \label{thm: wild bootstrap KMT}
	Let the conditions of Lemma \ref{lem: hats off} hold. In addition,  
	\begin{enumerate}[label=(\roman*)]
		\item when $\E[|\e_i|^{2+\nu}]<\infty$ for some $\nu>0$, assume that for $j=2,3$, $\frac{n^{\frac{1}{2+\nu}}(\log n)^{\frac{4+3\nu}{4+2\nu}}}
		{\sqrt{nh^d}}=o(r_n^{-2})$;
		\item when $\E[|\e_i|^3\exp(|\e_i|)]<\infty$, assume that there exists some constant $C>0$ such that for $1\leq i\leq n$,
		$$C\E^*[|\e_i\omega_i|^3\exp(|\e_i\omega_i|)]\leq \E^*[|\e_i\omega_i|^2]\quad \text{a.s.},$$ and for $j=2,3$, $\frac{(\log n)^{5/2}}{\sqrt{nh^d}}=o(r_n^{-2})$.
	\end{enumerate}
	Then for each $j=0,1,2,3$, $\widehat{z}^*_j(\cdot)=_{d^*} Z_j(\cdot)+o_{\P^*}(r_n^{-1})$ in $\mathcal{L}^\infty(\mathcal{X})$.
\end{thm}
The additional condition in (ii) of the above theorem is generally difficult to verify in practice unless some restrictive conditions, e.g. bounded support of $\e_i$, are imposed. 

\begin{thm}[Wild Bootstrap: Yurinskii] \label{thm: wild bootstrap Yurinskii}
	Let the conditions in Theorem \ref{thm: Yurinskii coupling} hold. In addition, for $j=2,3$, 
    \begin{enumerate}[label=(\roman*)]
    	\item when $\E[|\e_i|^{2+\nu}]<\infty$ for some $\nu>0$, assume that $\frac{n^{\frac{1}{2+\nu}}(\log n)^{\frac{4+3\nu}{4+2\nu}}}
    	{\sqrt{nh^d}}=o(r_n^{-2})$;
    	\item when $\E[|\e_i|^3\exp(|\e_i|)]<\infty$, assume that $\frac{(\log n)^{5/2}}{\sqrt{nh^d}}=o(r_n^{-2})$.
    \end{enumerate}
    Then for each $j=0,1,2,3$, $\widehat{z}_j^*(\cdot)=_{d^*}Z_j(\cdot)+o_{\P^*}(r_n^{-1})$ in $\mathcal{L}^\infty(\mathcal{X})$.
\end{thm}

\subsection{Application: Confidence Bands}
An immediate application of Theorems \ref{thm: KMT coupling}--\ref{thm: wild bootstrap Yurinskii} is to construct confidence bands for the regression function or its derivatives. Specifically, for $j=0,1,2,3$, and $\alpha\in(0,1)$, we seek a quantile $q_j(\alpha)$ such that
\[\P\left[\sup_{\bx\in\mathcal{X}} |\widehat{T}_j(\bx)| \leq q_j(\alpha) \right] =1-\alpha + o(1),
\]
which then can be used to construct uniform $100(1-\alpha)$-percent confidence bands for $\partial^\bq \mu(\bx)$ of the form
\[\left[\widehat{\partial^\bq \mu}_j(\bx)
\pm q_j(\alpha)\sqrt{\bOhat{j}/n}
\;:\; \bx\in\mathcal{X} 
\right].
\]
The following theorem establishes a valid distributional approximation for the suprema of the $t$-statistic processes $\{\widehat{T}_j(\bx):\bx\in\mathcal{X}\}$ using \cite[Lemma 2.4]{Chernozhukov-Chetverikov-Kato_2014b_AoS} to convert our strong approximation results into convergence of distribution functions in terms of Kolmogorov distance.

\begin{thm}[Confidence Band] \label{thm: confidence band}
	Let the conditions of Theorem \ref{thm: KMT coupling} or Theorem \ref{thm: Yurinskii coupling} hold with $r_n=\sqrt{\log n}$. If the corresponding conditions of Theorem \ref{thm: plug-in approx} for each $j=0,1,2,3$ hold, then 
	\[ 
	\sup_{u\in\mathbb{R}}\left|
	\P\left[\sup_{\bx\in\mathcal{X}}|\widehat{T}_j(\bx)|\leq u\right]
	-\P^*\left[\sup_{\bx\in\mathcal{X}}|\widehat{Z}_j(\bx)|\leq u\right]\right|=o_\P(1).
	\]
	If the corresponding conditions in Theorem \ref{thm: plug-in approx cond} hold for each $j=0,1,2,3$, then
	\[ 
	\sup_{u\in\mathbb{R}}\left|
	\P\left[\sup_{\bx\in\mathcal{X}}|\widehat{T}_j(\bx)|\leq u\right]
	-\P^*\left[\sup_{\bx\in\mathcal{X}}|\widehat{z}_j(\bx)|\leq u\right]\right|=o_\P(1).
	\]
	If the corresponding conditions in Theorem \ref{thm: wild bootstrap KMT} or \ref{thm: wild bootstrap Yurinskii} hold for each $j=0,1,2,3$, then
	\[ 
	\sup_{u\in\mathbb{R}}\left|
	\P\left[\sup_{\bx\in\mathcal{X}}|\widehat{T}_j(\bx)|\leq u\right]
	-\P^*\left[\sup_{\bx\in\mathcal{X}}|\widehat{z}_j^*(\bx)|\leq u\right]\right|=o_\P(1).
	\]
\end{thm}


\section{Examples} \label{section: examples}
We now illustrate how Assumptions \ref{assumption: Local Basis}, \ref{assumption: Misspecification Error}, and \ref{assumption: BC bases} are verified for popular basis choices, including $B$-splines, wavelets, and piecewise polynomials, and show when \eqref{eq: approx orth} holds. Note well that even in the cases where \eqref{eq: approx orth} holds, we will not make use of this property neither herein nor in our software implementation \cite{Cattaneo-Farrell-Feng_2019_lspartition}, since both bias terms of Lemma \ref{lem: bias} may be important in finite samples.

The first three subsections illustrate each of these assuming rectangular support, $\mathcal{X}=\otimes_{\ell=1}^d\mathcal{X}_\ell\subset\mathbb{R}^d$, and correspondingly use tensor-product partitions. We note that these three in fact overlap in the special cases of $m=1$ on a tensor-product partition (see Equation \ref{eqn:rect}), leading to so-called zero-degree splines, Haar wavelets, and regressograms (piecewise constant). The final subsection considers a general partition. Recall that whenever $\Delta$ covers only strict subset of $\mathcal{X}$, all our results hold on that subset; the first three subsections may also be interpreted in this light. 

\subsection{B-Splines on Tensor-Product Partitions}
A univariate spline is a piecewise polynomial satisfying certain smoothness constraints. For some integer $m_\ell\geq 1$, let $\mathcal{S}_{\Delta_\ell,m_\ell}$ be the set of splines of order $m_\ell$ with univariate partition $\Delta_\ell$, i.e.,
\[\small
\mathcal{S}_{\Delta_\ell,m_\ell}=\Big\{s\in \mathcal{C}^{m_\ell-2}(\mathcal{X}_\ell)\colon s(x)\text{ is a
	polynomial of degree}\; (m_\ell-1)\;
 \text{on each }
[t_{\ell,l_\ell},t_{\ell,l_\ell+1}]\Big\}.\]
If $m_\ell=1$, it is interpreted as the set of piecewise constant functions that are discontinuous on $\Delta_\ell$.
$\mathcal{S}_{\Delta_\ell, m_\ell}$ is a vector space and can be spanned by many equivalent representing bases. $B$-splines as a local basis are well studied in literature and enjoy many nice properties \cite{
	Micula-Micula_1999_Handbook,Schumaker_2007_Book}. 

Define an extended knot sequence $\Delta_{\ell,e}$ such that
\[t_{\ell,-m_\ell+1}=t_{\ell,-m_\ell+2}=\cdots=t_{\ell,0}<t_{\ell,1}<\cdots<t_{\ell,\kappa_{\ell}-1}<t_{\ell,\kappa_{\ell}}=\cdots=t_{\ell,\kappa_{\ell}+m_\ell-1},\]
Then, the $m_\ell$-th order $B$-spline with knot sequence $\Delta_{\ell,e}$ is
\[p_{l_\ell,m_\ell}(x_\ell)=(t_{\ell,l_\ell}-t_{\ell,l_\ell-m_\ell})[t_{\ell,l_\ell-m_\ell}, \cdots,
t_{\ell,l_\ell}]_t(t-x_\ell)_+^{m_\ell-1},\; l_\ell=1,\ldots,K_\ell=\kappa_{\ell}+m_\ell-1,\]
where $(a)_+=\max\{a,0\}$ and $[t_1,t_2,\cdots,t_l]_t\,\mu(t,x)$ denotes the divided difference of $\mu(t,x)$ with respect to $t$, given a sequence of knots $t_1\leq t_2\leq\cdots\leq t_l$. 

When there is no chance of confusion, we shall write $p_{l_\ell}(x_\ell)$ instead of $p_{l_\ell, m_\ell}(x_\ell)$. Accordingly, the space of tensor-product polynomial splines of order $\boldm=(m_1,\ldots,m_d)$ with partition $\Delta$ is spanned by tensor products of univariate $B$-splines
\[\mathcal{S}_{\Delta,\boldm}=\otimes_{\ell=1}^d
\mathcal{S}_{\Delta_\ell,m_\ell}=\text{span}\{p_{l_1}(x_1)p_{l_2}(x_2)\cdots
p_{l_d}(x_d)\}_{l_1=1,\ldots,l_d=1}^{K_1,\ldots,K_d}.\]
We have a total of $K=\prod_{\ell=1}^{d}K_\ell$ basis functions. The order of univariate basis could vary across dimensions, but for simplicity we assume that $m_1=\cdots=m_d=m$ and hence write $\mathcal{S}_{\Delta, m}=\mathcal{S}_{\Delta, \boldm}$. Also, let $p_{l_1\ldots l_d}(\bx)=p_{l_1}(x_1)\cdots p_{l_d}(x_d)$ to further simplify notation.

Arrange this set of basis functions by first increasing $l_d$ from $1$ to $K_d$ with other $l_\ell$'s fixed at $1$ and then increasing $l_\ell$ sequentially. As a result, we construct a one-to-one correspondence $\varphi$ mapping from $\{(l_1,\ldots,l_d): 1\leq l_\ell\leq K_\ell, \ell=1,\ldots,d\}$ to $\{1,\ldots, K\}$. Then, we can write $p_k(\bx)=p_{\varphi^{-1}(k)}(\bx)$, $k=1,\ldots, K$, which is consistent with our notation in Section \ref{section: setup}.

The following lemma shows that Assumptions \ref{assumption: Local Basis} and \ref{assumption: Misspecification Error} hold for $B$-splines.

\begin{lem}[$B$-Splines Estimators]\label{lemma: B-Splines}
	Let $\bp(\bx)$ be a tensor-product $B$-Spline basis of order $m$, and suppose Assumptions \ref{assumption: DGP} and \ref{assumption: quasi-uniform partition} hold with $m\leq S$. Then:
	\begin{enumerate}
		\item $\bp(\bx)$ satisfies Assumption \ref{assumption: Local Basis}.
		
		\item If, in addition,
		\begin{equation}\label{eq: quasi-uniform partition for Splines}
		\max_{0\leq l_\ell\leq\kappa_\ell-2}|b_{\ell,l_\ell+1}-b_{\ell,l_\ell}|=O(b^{1+\varrho}), \qquad \ell=1,\ldots, d,
		\end{equation}
		then Assumption \ref{assumption: Misspecification Error} holds with $\varsigma = m-1$ and
		\[\mathscr{B}_{m,\bvsig}(\bx)
		=-\sum_{\bm{u}\in\Lambda_m} \frac{\partial^{\bm{u}} \mu(\bx)
			h_\bx^{m-[\bvsig]}}{(\bm{u}-\bvsig)!} \frac{\mathbf{b}_{\bx}^{\bm{u}-\bvsig}}{h_{\bx}^{m-[\bvsig]}} 
		B^\mathtt{S}_{\bm{u}-\bvsig}\Big((\bx-\bt_\bx^L)\oslash\mathbf{b}_\bx\Big)
		\]
		where $\Lambda_m=\{\bm{u}\in\mathbb{Z}_+^d:[\bm{u}]=m,\text{  and  } u_\ell=m \text{ for some }1\leq \ell\leq d\}$ and $B^\mathtt{S}_{\bm{u}}(\bx)$ is the product of univariate Bernoulli polynomials; that is, $B^\mathtt{S}_{\bm{u}}(\bx):=\prod_{\ell=1}^d B_{u_\ell}(x_\ell)$ with $B_{u_\ell}(\cdot)$ being the $u_\ell$-th Bernoulli polynomial and $B^\mathtt{S}_{\bm{u}}(\cdot)=0$ if $\bm{u}$ contains negative elements. Furthermore, Equation \eqref{eq: approx orth} holds.
		
		\item Let $\bpt(\bx)$ be a tensor-product $B$-Spline basis of order $\tilde{m}>m$ on the same partition $\Delta$, and assume $\tilde{m}\leq S$. Then Assumption \ref{assumption: BC bases} is satisfied.
		
	\end{enumerate} 
\end{lem}

Equation \eqref{eq: quasi-uniform partition for Splines} gives a precise definition of the strong quasi-uniform condition on the partition scheme. Assumption \ref{assumption: quasi-uniform partition} requires that the volume of all cells vanish at the same rate but allows for any constant proportionality between neighboring cells. Presently, cells are further restricted to be of the same volume asymptotically, and further, a specific rate is required that is related to the smoothness of $\mu(\cdot)$. Note that, for example, equally spaced knots satisfy this conditional trivially. For other schemes, additional work may be needed. Under \eqref{eq: quasi-uniform partition for Splines}, \cite{Barrow-Smith_1978_QAM} obtained an expression of the leading asymptotic error of univariate splines, which was later used by \cite{Zhou-Shen-Wolfe_1998_AoS} and \cite{Zhou-Wolfe_2000_SS}, among others. Lemma \ref{lemma: B-Splines} extends previous results to the multi-dimensional case, in addition to showing that the high-level conditions in Assumption \ref{assumption: Local Basis} and \ref{assumption: Misspecification Error} hold for $B$-Splines.

When \eqref{eq: quasi-uniform partition for Splines} fails, it may be possible to obtain results, but additional cumbersome notation would be needed and the results may be less useful. However, it is straightforward to verify that without \eqref{eq: quasi-uniform partition for Splines}, there still exists $s^*\in\mathcal{S}_{\Delta,m}$ such that for all $[\bvsig]\leq \varsigma$,
\begin{equation}\label{eq: Approximation Power}
\sup_{\bx\in\mathcal{X}}|\partial^{\bvsig} \mu
-\partial^{\bvsig} s^*|\lesssim h^{m-[\bvsig]},
\end{equation}
where recall that $\mathcal{S}_{\Delta,m}$ denotes the linear span of $\bp(\bx;\Delta,m)$. However, this result does not allow for bias correction or IMSE expansion. 


Finally, $\Lambda_m$ contains only the multi-indices corresponding to $m$-th order partial derivatives of $\mu(\bx)$.
This is due to the fact that, as a variant of polynomial approximation, the total order of tensor-product $B$-splines is not fixed at $m$, i.e.\ some higher-order components are retained that approximate terms with $m$-th order cross partial derivatives. This feature distinguishes tensor-product splines from multivariate splines which do control the total order of bases.

\subsection{Wavelets on Tensor-Product Partitions}

Our results apply to compactly supported wavelets, such as Cohen-Daubechies-Vial wavelets  \cite{Cohen_1993_ACHA}. For more background details see \cite{Meyer_1995_Book,Hardle-Kerkyacharian-Picard-Tsybakov_2012_Book,Chui_2016_Book}, and references therein.

To describe these estimators, we first introduce the general definition of wavelets and then show that a large class of orthogonal wavelets satisfy our high-level assumptions. For some integer $m\geq 1$, we call $\phi$ a (univariate) scaling function or ``father wavelet'' of degree $m-1$ if (i) $\int_\mathbb{R}\phi(x_\ell)\,dx_\ell=1$, (ii) $\phi$ and all its derivatives up to order $m-1$ decrease rapidly as $|x_\ell|\rightarrow \infty$, and (iii) $\{\phi(x_\ell-l):l\in\mathbb{Z}\}$ forms a Riesz basis for a closed subspace of $L_2(\mathbb{R})$. A real-valued function $\psi$ is called a (univariate) ``mother wavelet'' of degree $m-1$ if (i) $\int_{\mathbb{R}}x_\ell^v\psi(x_\ell)dx_\ell=0$ for $0\leq v\leq m-1$, (ii) $\psi$ and all its derivatives decrease rapidly as $|x_\ell|\rightarrow\infty$, and (iii) $\{2^{s/2}\psi(2^sx_\ell-l):s,l\in\mathbb{Z}\}$ forms a Riesz basis of $L_2(\mathbb{R})$.

We restrict our attention to orthogonal wavelets with compact support. The father wavelet $\phi$ and the mother wavelet $\psi$ are both compactly supported, and for any integer $s_0\geq 0$, any function in $L_2(\mathbb{R})$ admits the following $(m-1)$-regular wavelet multiresolutional expansion:
\[\mu(x_\ell)=\sum_{l=-\infty}^{\infty}a_{s_0l}\phi_{s_0l}(x_\ell)+\sum_{s=s_0}^{\infty}\sum_{l=-\infty}^{\infty}b_{sl}\psi_{sl}(x_\ell),\;x_\ell\in\mathbb{R},\quad \text{where}\]
\begin{align*}
a_{sl} &=\int_\mathbb{R}\mu(x_\ell)\phi_{sl}(x_\ell)dx_\ell, \qquad \phi_{sl}(x_\ell)=2^{s/2}\phi(2^sx_\ell-l),\\
b_{sl} &=\int_\mathbb{R}\mu(x_\ell)\psi_{sl}(x_\ell)dx_\ell, \qquad \psi_{sl}(x_\ell)=2^{s/2}\psi(2^sx_\ell-l),
\end{align*}
and $\{\phi_{s_0l},l\in\mathbb{Z};\psi_{sl},s\geq s_0,l\in\mathbb{Z}\}$ is an orthonormal basis of $L_2(\mathbb{R})$. Accordingly, to construct an orthonormal basis of $L_2([0,1])$, one can pick those basis functions supported in the interior of $[0,1]$ and add some boundary correction functions. For details of construction of such a basis see, for example, \cite{Cohen_1993_ACHA}. With a slight abuse of notation, in what follows we write $\{\phi_{s_0l},l\in\mathcal{L}_{s_0};\psi_{sl},s\geq s_0,l\in\mathcal{L}_s\}$ for an orthonormal basis of $L_2([0,1])$ rather than $L_2(\mathbb{R})$ where $\mathcal{L}_{s_0}$ and $\mathcal{L}_s$ are some proper index sets depending on $s_0$ and $s$ respectively. Again, we use tensor product to form a multidimensional wavelet basis, and then it fits into the context of our analysis. 

The compactness of the father wavelet is needed to verify Assumption \ref{assumption: Local Basis}. To see this, consider a standard multiresolutional analysis setting. A large space, say $L_2([0,1])$, is decomposed into a sequence of nested subspaces 
$\{0\}\subset\cdots\subset V_{-1}\subset V_0\subset V_1\cdots\subset L_2([0,1])$. Generally, $\{\phi_{sl},l\in\mathcal{L}_s\}$ constitutes a basis for $V_s$, and $\{\psi_{sl},l\in\mathcal{L}_s\}$ forms a basis for the orthogonal complement $W_s$ of $V_s$ in $V_{s+1}$. Since the support of $\phi$ is compact and $\{\phi_{sl}\}$ is generated simply by dilation and translation of $\phi$, $[0,1]$ can be viewed as implicitly partitioned. Specifically, by increasing the resolution level from $s$ to $s+1$, the length of the support of $\phi_{sl}$ is halved. Hence, it is equivalent to placing additional partitioning knots at the midpoint of each subinterval. In addition, as the sample size $n$ grows, we allow the resolution level $s$ to increase (thus written as $s_n$ when needed), but each basis function $\phi_{s_nl}$ is supported by only a finite number of subintervals since the generating scaling function is compactly supported.

The above discussion associates the generating process of wavelets at different levels with a partitioning scheme. It is easy to see that Assumption \ref{assumption: quasi-uniform partition} is automatically satisfied in this case, and given a resolution level $s_n$, the mesh width $b=2^{-s_n}$. On the other hand, given a starting level $s_0$, $\{\phi_{s_0l},l\in\mathcal{L}_{s_0};\psi_{sl},s_0\leq s\leq s_n,l\in\mathcal{L}_s\}$ does not satisfy Assumption \ref{assumption: Local Basis}, since some functions in such a basis have increasing supports as the resolution level $s_n\to\infty$. Nevertheless, they span the same space as $\{\phi_{s_n,l},l\in\mathcal{L}_{s_n}\}$ and hence there exists a linear transformation that connects the two equivalent bases. Least squares estimators are invariant to such a transformation. Therefore, we will focus on the following tensor-product (father) wavelet basis
\begin{equation}\label{eq10}
\bp(\bx):=\otimes_{\ell=1}^d 2^{-s_n/2} \bphi_{s_n}(x_\ell)
\end{equation}
where $\bphi_{s_n}$ is a vector containing all functions in $\{\phi_{s_n,l},l\in\mathcal{L}_{s_n}\}$. By multiplying the basis by $2^{-s_n/2}$, we drop the normalizing constants that originally appear in the construction of orthonormal basis. As the next lemma shows, this large class of orthogonal wavelet bases satisfy our assumptions.

\begin{lem}[Wavelets Estimators]\label{lemma-w}
	Let $\phi$ and $\psi$ be a scaling function and a wavelet function of degree $m-1$ with $q+1$ continuous derivatives, $\bp(\bx)$ be the tensor-product orthogonal (father) wavelet basis of degree $m-1$ generated by $\phi$, and suppose Assumption \ref{assumption: DGP} holds with $m\leq S$.
	\begin{enumerate}
		\item $\bp(\bx)$ satisfies Assumption \ref{assumption: Local Basis}.
		
		\item Assumption \ref{assumption: Misspecification Error} holds with $\varsigma=q$ and
		\[\mathscr{B}_{m,\bvsig}(\bx)
		= -\sum_{\bm{u}\in\Lambda_m} \frac{\partial^{\bm{u}} \mu(\bx)
			h^{m-[\bvsig]}}{\bm{u}!}\frac{b^{m-[\bvsig]}}{h^{m-[\bvsig]}}
		B^\mathtt{W}_{\bm{u},\bvsig}(\bx/b)
		\]
		where $\Lambda_m=\{\bm{u}\in\mathbb{Z}_+^d:[\bm{u}]=m,\text{  and  } u_\ell=m \text{ for some } 1\leq \ell\leq d\}$ and $B^\mathtt{W}_{\bm{u},\bvsig}(\bx)=\sum_{s\geq0} \partial^{\bvsig} \xi_{\bm{u},s}(\bx)$ converges uniformly, with $\xi_{\bm{u},s}(\cdot)$ being a linear combination of products of univariate father wavelet $\phi$ and the mother wavelet $\psi$; its exact form is notationally cumbersome and is given in Equation \eqref{eq: wavelet error}. Furthermore, Equation \eqref{eq: approx orth} holds.
		
		\item Let $\tilde{\phi}$ be a scaling function of degree $\tilde{m}-1$ with $m+1$ continuous derivatives for some $\tilde{m}>m$, $\bpt(\bx)$ be the tensor-product orthogonal wavelet basis generated by $\tilde{\phi}$ that has the same resolution level as $\bp(\bx)$, and assume $\tilde{m}\leq S$. Then Assumption \ref{assumption: BC bases} is satisfied.
		
	\end{enumerate}
\end{lem}

In addition to verifying that our high-level assumptions hold for wavelets, this result gives a novel asymptotic error expansion for multidimensional compactly supported wavelets. Our derivation employs the ideas in \cite{Sweldens-Piessens_1994_NM} and exploits the tensor-product structure of the wavelet basis. The end result is similar to tensor-product splines, in the sense that the total order of the approximating basis is not fixed at $m$ and thus $\Lambda_m$ is the same as that for $B$-splines.

\subsection{Generalized Regressograms on Tensor-Product Partitions}
\label{section: Partitioning Estimators}

To construct the generalized regressogram, we first define the piecewise polynomials. What will distinguish these from splines is that (i) each polynomial is supported on exactly one cell, and relatedly (ii) no continuity is assumed between cells. First, for some fixed integer $m\in\mathbb{Z}_+$, let $\mathbf{r}(x_\ell)=(1,x_\ell,\ldots,x_\ell^{m-1})'$ denote a vector of powers up to degree $m-1$. To extend it to a multidimensional basis, we take the tensor product of $\mathbf{r}(x_\ell)$, denoted by a column vector $\mathbf{R}(\bx)$. The total order of such a basis is not fixed, and its behavior is more similar to tensor-product $B$-splines. Following \cite[][and references therein]{Cattaneo-Farrell_2013_JoE}, we exclude all terms with degree greater than $m-1$ in $\mathbf{R}(\bx)$. Hence the remaining elements in $\mathbf{R}(\bx)$ are given by $\bx^{\boldsymbol{\alpha}}=x_1^{\alpha_1}\cdots x_d^{\alpha_d}$ for a unique $d$-tuple $\boldsymbol{\alpha}$ such that $[\boldsymbol{\alpha}]\leq m-1$. Then we ``localize'' this basis by restricting it to a particular subrectangle $\delta_{l_1\cdots l_d}$. Specifically, we write $\bp_{l_1\ldots l_d}(\bx)=\I_{\delta_{l_1\ldots l_d}}(\bx)\mathbf{R}(\bx)$, where $\I_{\delta_{l_1\ldots l_d}}(\bx)$ is equal to 1 if $\bx\in \delta_{l_1\ldots l_d}$ and $0$ otherwise. Finally, we rotate the basis by centering each basis function at the start point of the corresponding cell and scale it by interval lengths.
We can arrange the basis functions according to a particular ordering $\varphi$: first order the subrectangles the same way as that for $B$-splines, and then within each subrectangle $\delta_{l_1\ldots l_d}$ the basis functions in $\mathbf{R}(\bx)$ are ordered ascendingly in $\boldsymbol{\alpha}$ and $\ell=1,\ldots, d$. Using the same notation as in Section \ref{section: setup}, we still write the basis as $\bp(\bx)=(p_1(\bx),\cdots, p_K(\bx))'$.

The following lemma shows that Assumptions \ref{assumption: Local Basis} and \ref{assumption: Misspecification Error} hold for generalized regressograms.

\begin{lem}[Generalized Regressograms]\label{lemma: Partitioning Estimators}
	Let $\bp(\bx)$ be the rotated piecewise polynomial basis of degree $m-1$ based on Legendre polynomials, and suppose that Assumptions \ref{assumption: DGP} and \ref{assumption: quasi-uniform partition} hold with $m\leq S$. Then,
	\begin{enumerate}
		\item $\bp(\bx)$ satisfies Assumption \ref{assumption: Local Basis}.
		
		\item Assumption \ref{assumption: Misspecification Error} holds with $\varsigma = m-1$ and
		\[\mathscr{B}_{m,\bvsig}(\bx)
		= -\sum_{\bm{u}\in\Lambda_m}  \frac{\partial^{\bm{u}} \mu(\bx)
			h_\bx^{m-[\bvsig]}}{(\bm{u}-\bvsig)!} \frac{\mathbf{b}_\bx^{\bm{u}-\bvsig}}{h_{\bx}^{m-[\bvsig]}}
		B^\mathtt{P}_{\bm{u}-\bvsig}\Big((\bx-\bt^L_\bx) \oslash \mathbf{b}_\bx\Big),
		\]
		where $\Lambda_m = \{\bm{u}:[\bm{u}]=m\}$ and
		\[B^\mathtt{P}_{\bm{u}}(\bx) := \prod_{\ell=1}^{d}\binom{2 u_\ell}{u_\ell}^{-1}P_{u_\ell}(x_\ell),\]
		with $P_{u_\ell}(\cdot)$ being the $u_\ell$-th shifted Legendre polynomial orthogonal on $[0,1]$, and $B^\mathtt{P}_{\bm{u}}(\cdot)=0$ if $\bm{u}$ contains negative elements. Furthermore, Equation \eqref{eq: approx orth} holds.
		
		\item Let $\bpt(\bx)$ be a piecewise polynomial basis of degree $\tilde{m}-1$ on the same partition $\Delta$ for some $\tilde{m}>m$, and assume $\tilde{m} \leq S$. Then Assumption \ref{assumption: BC bases} is satisfied.
		
	\end{enumerate}
\end{lem}

The leading asymptotic error obtained in Lemma \ref{lemma: Partitioning Estimators} differs from the one in \cite{Cattaneo-Farrell_2013_JoE} because it is expressed in terms of orthogonal polynomials. Here we employ Legendre polynomials $\bar{P}_m(x)$ orthogonal on $[-1,1]$ with respect to the Lebesgue measure, and then shift them to $P_m(x)=\bar{P}_m(2x-1)$. Thus the shifted Legendre polynomials are orthogonal on $[0,1]$. Further, we allow for more general partitioning schemes.

\subsection{Generalized Regressograms on General Partitions}

Moving away from tensor-product partitions will impede verification of Assumption \ref{assumption: Misspecification Error} in general. A more typical result, such as Equation \eqref{eq: Approximation Power} can be more easily proven in generality. In practice, however, many data-driven partitioning selection procedures, such as those induced by regression trees, do not necessarily lead to a tensor-product partition. There is also a large literature in approximation theory discussing bases constructed on triangulations or other general partition schemes \citep[e.g.,][]{Micula-Micula_1999_Handbook}. To demonstrate the power of our theory, we show presently that for the generalized regressogram we can obtain concrete results on general partitions. 

The basis is as described above, but now allowing for cells of general shape. Let $\mathcal{X}$ be convex and $\bt_\delta$ be an arbitrary point in each $\delta \in \Delta$. Hence, $\bt_{\delta_\bx}$ is simply a point in the cell containing $\bx$. Construct the rotated piecewise polynomial basis as in Section \ref{section: Partitioning Estimators}, but centering the basis at $\bt_\delta$ and rescaling it by the diameter of $\delta$. 

\begin{lem}\label{lemma-pp}
	Let $\bp(\bx)$ be the rotated piecewise polynomial basis of degree $m-1$ on a general partition $\Delta$, and suppose that Assumptions \ref{assumption: DGP} and \ref{assumption: quasi-uniform partition} hold with $m\leq S$. Then,
	\begin{enumerate}
		\item $\bp(\bx)$ satisfies Assumption \ref{assumption: Local Basis}.
		
		\item Assumption \ref{assumption: Misspecification Error} holds with $\varsigma = m-1$ and
		\[\mathscr{B}_{m,\bvsig}(\bx)
		= -\sum_{\bm{u}\in\Lambda_m} \frac{\partial^{\bm{u}} \mu(\bx)
			h_\bx^{m-[\bvsig]}}{(\bm{u}-\bvsig)!}
		\Big(\frac{\bx-\bt_{\delta_\bx}}{h_\bx}\Big)^{\bm{u}-\bvsig}
		\I(\bm{u}\geq\bvsig),
		\]
		where $\Lambda_m = \{\bm{u}:[\bm{u}]=m\}$ and $\bm{u}\geq\bvsig\Leftrightarrow u_1\geq \varsigma_1, \cdots, u_d\geq \varsigma_d$.
		
		\item Let $\bpt(\bx)$ be a piecewise polynomial of degree $\tilde{m}-1$ on the same partition $\Delta$ for some $\tilde{m}>m$, and assume $\tilde{m}\leq S$. Then Assumption \ref{assumption: BC bases} is satisfied.
	\end{enumerate}
\end{lem} 

It is a challenging task to construct orthogonal polynomial bases on non-rectangular domains, which makes Equation \eqref{eq: approx orth} hard to be satisfied when employing partitioning estimators on general partitions. Thus, our more general characterization (and correction) of the bias is quite useful in this case.

\section{Discussion and Extensions} \label{section: extension}
\subsection{Connecting Splines and Piecewise Polynomials}
There is an important connection between splines and piecewise polynomial bases. Essentially, the former can be viewed as a piecewise polynomial basis with certain continuity restrictions. Therefore, the estimators based on the two bases are linked by utilizing the well-known results about regressions with linear constraints. See \cite{Buse_1977_JASA} for an illustration of cubic splines as a special case of restricted least squares. 

Formally, let us start with the (rotated) piecewise polynomials $\bp$ discussed in Section \ref{section: examples}. For expositional simplicity, we only discuss tensor-product partitions here. Clearly, $\bp$ spans a vector space $\mathcal{P}$ containing all ``piecewise'' polynomials with degree no greater than $m-1$ on $\mathcal{X}$:
\[\mathcal{P}:=\Big\{s(\cdot):\;s(\cdot)=\sum_{k=1}^{K}a_k p_k(\cdot), \; a_k\in\mathbb{R}\Big\}.\]
Functions in this space are continuous within each subrectangle, but might have ``jumps'' along boundaries of cells. Then by imposing certain continuity restrictions on functions in $\mathcal{P}$, we can construct a subspace $\mathcal{S}\subset\mathcal{P}$
\[\mathcal{S}:=\Big\{s(\cdot):\; s(\cdot) = \sum_{k=1}^{K}a_k p_k(\cdot), \; a_k\in\mathbb{R}, \text{ and}\; \partial^{\bvsig}s(\cdot)\in \mathcal{C}(\mathcal{X}),\; \forall [\bvsig] \leq \iota\Big\}.\]
where $\iota\leq m-2$ is a positive integer that controls the smoothness of the basis (and thus the smoothness of the estimated function). Since the derivatives of polynomials are linear in coefficients, the continuity constraints are linear.

Now consider the following restricted least squares:
\begin{equation*}
\widehat{\bbeta}_{\mathtt{r}}=\underset{\mathbf{b}\in\mathbb{R}^K}{\arg\min} \,\|\by-\bP\mathbf{b}\|^2\quad\text{ s.t. }\; \bR\mathbf{b}=0
\end{equation*}   
where $\bP=(\bp(\bx_1),\ldots, \bp(\bx_n))'$ and $\bR$ is a $\vartheta_\iota\times K$ restriction matrix. $\vartheta_\iota$ denotes the number of restrictions depending on the required smoothness $\iota$. If there are no redundant constraints, $\bR$ has full row rank. It is well known that given the restriction matrix $\bR$, the least squares estimator can be written as
\[\widehat{\bbeta}_{\mathtt{r}}=[\bI-(\bP'\bP)^{-1}\bR'(\bR(\bP'\bP)^{-1}\bR')^{-1}\bR](\bP'\bP)^{-1}\bP' \by.\]
Since the unrestricted least squares estimator $\widehat{\bbeta}_{\mathtt{ur}}=(\bP'\bP)^{-1}\bP'\by$, the above equation implies that
$\widehat{\bbeta}_{\mathtt{r}}=[\bI-(\bP'\bP)^{-1}\bR'(\bR(\bP'\bP)^{-1}\bR')^{-1}\bR]\widehat{\bbeta}_{\mathtt{ur}}
=:(\bI-\mathbf{U})\widehat{\bbeta}_{\mathtt{ur}}$.
Therefore,
\begin{equation}\label{eq: restricted ls}
\widehat{\mu}_{\mathtt{r}}(\bx)=\bp(\bx)'(\bI-\mathbf{U})
\widehat{\bbeta}_{\mathtt{ur}}
=\widehat{\mu}_{\mathtt{ur}}(\bx)-\bp(\bx)'\mathbf{U}
\widehat{\bbeta}_{\mathtt{ur}}
\end{equation}
where $\widehat{\mu}_{\mathtt{ur}}(\bx):=\bp(\bx)'\widehat{\bbeta}_{\mathtt{ur}}$ is the unrestricted estimator.
With this relation we can derive expressions of bias and variance for the restricted estimator. Clearly, the conditional variance $\V[\widehat{\mu}_{\mathtt{r}}(\bx)|\bX]=\bp(\bx)'(\bI-\mathbf{U})
\V[\widehat{\mu}_{\mathtt{ur}}(\bx)|\bX](\bI-\mathbf{U})'\bp(\bx)$. On the other hand, as shown in Lemma \ref{lemma: Partitioning Estimators}, there exists $s^*(\bx):=\bp(\bx)'\bbeta^*$ such that $\|\mu-s^*+\mathscr{B}_{m,\mathbf{0}}\|_{L_\infty(\mathcal{X})}
=o(h^m)$. Hence
\begin{align*}
&\E[\widehat{\mu}_{\mathtt{r}}(\bx)|\bX]-\mu(\bx)\\
&=\bp(\bx)'(\bI-\mathbf{U})(\bP'\bP)^{-1}\bP'[\bP\bbeta^*-
\pmb{\mathscr{B}}_{m,\mathbf{0}}+o(h^m)]-\mu(\bx)\\
&=\mathscr{B}_{m,\mathbf{0}}(\bx)-\bp(\bx)'(\bI-\mathbf{U})(\bP'\bP)^{-1}
\bP'[\pmb{\mathscr{B}}_{m,\mathbf{0}}+o(h^m)]-\bp(\bx)'\mathbf{U}\bbeta^*+o(h^m)
\end{align*}
where $\pmb{\mathscr{B}}_{m, \mathbf{0}}=(\mathscr{B}_{m, \mathbf{0}}(\bx_1), \ldots, \mathscr{B}_{m, \mathbf{0}}(\bx_n))'$.
$\bp(\bx)'\mathbf{U}\bbeta^*$ can be viewed as a measure of to what extent the continuity constraints are satisfied. When $s^*(\bx)$ is continuous up to order $\iota$, $\bR\boldsymbol{\beta}^*$ is exactly $0$ and thus this term vanishes.

To repeat the previous analysis for such a restricted estimator, the main challenge is to analyze the asymptotic properties of the outer product of the restriction matrix. Once we know its limiting eigenvalue distributions (e.g., bounds on extreme eigenvalues),  $\mathbf{U}$ can be properly bounded, and then the conclusions in previous sections may be established with similar proofs. Unfortunately, for general multidimensional cases, it is difficult and tedious to specify a non-redundant set of continuity constraints and analyze the eigenvalue distributions of $\bR\bR'$, thus damping the usefulness of such a method, whereas when $d=1$, the restriction matrix is quite straightforward and well bounded. 

Formally, let $\bR$ be a restriction matrix corresponding to $\bar{\iota}=\iota+1$ continuity constraints at each partitioning knot, implying that $s\in\mathcal{S}$ has $\iota$ continuous derivatives. We write the $\ell$th restriction at the $k$th knot as $\mathbf{r}_{k\ell}$, corresponding to the requirement that the $(\ell-1)$th derivatives of functions in $\mathcal{S}$ are continuous. Explicitly, the entire restriction matrix admits the following structure:
\begin{equation}\label{eq: R matrix}
\bR=
\begin{array}{c@{\!\!\!}l}
\left[ \begin{array}[c]{c}
\mathbf{r}_{11}'\\
\vdots\\
\mathbf{r}_{1\bar{\iota}}'\\
\vdots\\
\mathbf{r}_{\kappa 1}'\\
\vdots\\
\mathbf{r}_{\kappa\bar{\iota}}'\\
\end{array}  \right]
&
\begin{array}[c]{@{}l@{\,}l}
\left. \begin{array}{c} \vphantom{\mathbf{r}_{11}'}  \\ 
\vphantom{\vdots}\\ \vphantom{\mathbf{r}_{1\bar{\iota}}'} \end{array} \right\} & \text{Restrictions at the 1st knot} \\
\vphantom{\vdots}&\\
\left. \begin{array}{c} \vphantom{\mathbf{r}_{\kappa 1}'}  \\ \vphantom{\vdots}
\\ \vphantom{\mathbf{r}_{\kappa\bar{\iota}}'} \end{array} \right\} & \text{Restrictions at the $\kappa$th knot}
\end{array}
\end{array}
\end{equation}
As $\kappa$ increases, the dimension of $\bR$ also grows. Moreover, $\iota$ cannot exceed $m-2$ since when $m$ continuity constraints are imposed, $\mathcal{S}$ degenerates to the space of global polynomials of degree no greater than $m-1$. 

The asymptotic behavior of the restricted estimator is closely related to the extreme eigenvalues of $\bR\bR'$. For a general restriction matrix $\bR$ with fixed dimensions, if it contains non-redundant constraints only, the minimum eigenvalue of $\bR\bR'$ is nonzero. When the number of constraints $\vartheta_\iota\rightarrow\infty$, however, the limit of the minimum eigenvalue does not have to be nonzero, and its limiting behavior depends on the specific structure of constraints. The next lemma shows that for the particular restricted estimators considered here, the eigenvalues of $\bR\bR'$ are indeed bounded and bounded away from zero uniformly over the number of knots for $\iota\leq m-2$.
 
\begin{lem}[Restriction Matrix]\label{lem: restriction matrix}
	Let $\mathbf{R}$ be the restriction matrix described as \eqref{eq: R matrix} with $\iota\leq m-2$. Then
	$$1\lesssim\lambda_{\min}(\mathbf{RR}')\leq\lambda_{\max}(\mathbf{RR}')\lesssim 1.$$
\end{lem}

The proof of this lemma employs the specific structure of the restriction matrix.
Generally, the outer product of $\bR$ takes the following form:
	\begin{equation}\label{eq: RR'}
	\bR\bR'=
	\begin{bmatrix}
	\mathbf{A} & \mathbf{B} &\mathbf{0} &\cdots & \mathbf{0}\\
	\mathbf{B}' &\mathbf{A} & \mathbf{B}&\cdots & \mathbf{0}\\
	\vdots& &\ddots&& \vdots\\
	\mathbf{0}&\cdots&\mathbf{B}'&\mathbf{A}&\mathbf{B}\\ 
	\mathbf{0}&\cdots&&\mathbf{B}'&\mathbf{A}
	\end{bmatrix}
	\end{equation}
	where
	\begin{equation*}
	\mathbf{A}=\begin{bmatrix}
	\mathbf{r}_{11}'\mathbf{r}_{11}&\mathbf{r}_{11}'\mathbf{r}_{12}&\cdots&\mathbf{r}_{11}'\mathbf{r}_{1\bar{\iota}}\\
	\mathbf{r}_{12}'\mathbf{r}_{11}&\mathbf{r}_{12}'\mathbf{r}_{12}&\cdots&\mathbf{r}_{12}'\mathbf{r}_{1\bar{\iota}}\\
	\vdots&\cdots&&\vdots\\
	\mathbf{r}_{1\bar{\iota}}'\mathbf{r}_{11}&\mathbf{r}_{1\bar{\iota}}'\mathbf{r}_{12}&\cdots&\mathbf{r}_{1\bar{\iota}}'\mathbf{r}_{1\bar{\iota}}
	\end{bmatrix},\quad
	\mathbf{B}=\begin{bmatrix}
	\mathbf{r}_{11}'\mathbf{r}_{21}&\mathbf{r}_{11}'\mathbf{r}_{22}&\cdots&\mathbf{r}_{11}'\mathbf{r}_{2\bar{\iota}}\\
	\mathbf{r}_{12}'\mathbf{r}_{21}&\mathbf{r}_{12}'\mathbf{r}_{22}&\cdots&\mathbf{r}_{12}'\mathbf{r}_{2\bar{\iota}}\\
	\vdots&\cdots&&\vdots\\
	\mathbf{r}_{1\bar{\iota}}'\mathbf{r}_{21}&\mathbf{r}_{1\bar{\iota}}'\mathbf{r}_{22}&\cdots&\mathbf{r}_{1\bar{\iota}}'\mathbf{r}_{2\bar{\iota}}
	\end{bmatrix}.
	\end{equation*}
	
	Importantly, the form described in \eqref{eq: RR'} is usually referred to as a (tridiagonal) block Toeplitz matrix, meaning that it is a tridiagonal block matrix containing blocks repeated down the diagonals. It is well known that its asymptotic eigenvalue distribution is characterized by the Fourier transform
	$$\pmb{\mathscr{T}}_{\bar{\iota}}(\omega)=\mathbf{A}+(\mathbf{B}+\mathbf{B}')\cos\omega,\quad \omega\in[0, 2\pi].$$
	As $\kappa\rightarrow\infty$, $\lambda_{\min}(\bR\bR')$ converges to the minimum attained by the minimum eigenvalue of $\pmb{\mathscr{T}}_{\bar{\iota}}(\omega)$ as a function of $\omega$ on $[0,2\pi]$. Similarly, the limit of $\lambda_{\max}(\bR\bR')$ is the maximum attained by the maximum eigenvalue of $\pmb{\mathscr{T}}_{\bar{\iota}}(\omega)$.

\subsection{Comparison of Bias Correction Approaches}
We make some comparison of the three bias correction approaches considered in this paper.

First, the higher-order correction ($j=1$) and least squares correction ($j=2$) are closely related. 
Let $\bP=(\bp(\bx_1), \ldots, \bp(\bx_n))'$ and $\tilde{\bP}=(\bpt(\bx_1), \ldots, \bpt(\bx_n))'$. Clearly,
\begin{align*}
\widehat{\mu}_2(\bx)&=\widehat{\mu}_1(\bx)+ \bp(\bx)'(\bP'\bP)^{-1}\bP'(\mathbf{I}-\tilde{\bP}(\tilde{\bP}'
\tilde{\bP})^{-1}\tilde{\bP})\by\\
&=\widehat{\mu}_1(\bx)+\bp(\bx)'(\bP'\bP)^{-1}\bP'\bM_{\bpt}\by
\end{align*}
where $\bM_{\bpt}:=\mathbf{I}-\tilde{\bP}(\tilde{\bP}'\tilde{\bP})^{-1}
\tilde{\bP}'$. Importantly, when $\bp$ and $\bpt$ generate nested models, i.e., there exists a transformation matrix $\bUpsilon$ such that $\bp(\cdot)=\bUpsilon\bpt(\cdot)$, it is easy to see that $\bM_{\bpt}\bP=\mathbf{0}$. Thus the higher-order and least squares bias correction approaches are equivalent. When $\bpt$ and $\bp$ are not nested bases, the two methods will typically differ in variance and bias. 

To compare their variance, we generally have
\begin{align*}
\V[\widehat{\mu}_2(\bx)|\bX]&=\V[\widehat{\mu}_1(\bx)|\bX]
+\V[\bp(\bx)'(\bP'\bP)^{-1}\bP'\bM_{\bpt}\by|\bX]\\
&\quad+2\Cov[\bpt(\bx)'(\tilde{\bP}'\tilde{\bP})^{-1}\tilde{\bP}\by,\;\bp(\bx)'(\bP'\bP)^{-1}\bP'\bM_{\bpt}\by|\bX].
\end{align*}
When $\sigma^2(\bx)=\sigma^2$, the covariance term is $0$, and thus $\widehat{\mu}_2$ has variance no less than that of $\widehat{\mu}_1$. The same conclusion is true for asymptotic variance as shown in the proof of Lemma \ref{lem: Avar}.

Regarding their bias, let $\mathfrak{B}_{\tilde{m}, \mathbf{0}}(\bx)
=\E[\widehat{\mu}_1(\bx)|\bX] - \mu(\bx)$ denote the conditional bias of $\widehat{\mu}_1(\bx)$, and $\boldsymbol{\mathfrak{B}}_{\tilde{m},\mathbf{0}}:=
(\mathfrak{B}_{\tilde{m},\mathbf{0}}(\bx_1),\ldots,\mathfrak{B}_{\tilde{m},\mathbf{0}}(\bx_n))'$. Then
$
\E[\widehat{\mu}_2(\bx)|\bX]-\mu(\bx)
=\mathfrak{B}_{\tilde{m},\mathbf{0}}(\bx)-
\bp(\bx)'(\bP'\bP)^{-1}\bP'\boldsymbol{\mathfrak{B}}_{\tilde{m},\mathbf{0}}
$. Clearly, the second term will asymptotically get close to the projection of $\mathfrak{B}_{\tilde{m},\mathbf{0}}(\bx)$ onto the space spanned by $\bp$: $$\mathscr{L}_{\bp}[\mathfrak{B}_{\tilde{m},\mathbf{0}}](\bx)
:=\bp(\bx)'(\E[\bp(\bx)\bp(\bx)'])^{-1}\E[\bp(\bx)\mathfrak{B}_{\tilde{m},\mathbf{0}}(\bx)]$$
where $\mathscr{L}_{\bp}[\cdot]$ denotes the projection operator.
Therefore, when $\bp$ and $\bpt$ are not nested bases, we typically have the bias of $\widehat{\mu}_2$ no greater than that of $\widehat{\mu}_1$ in terms of $\|\cdot\|_{F,L_2(\mathcal{X})}$, where for a real-valued function $g(\cdot)$ on $\mathcal{X}$, $\|g\|_{F,L_2(\mathcal{X})}=(\int_{\mathcal{X}}|g(\bx)|^2dF(\bx))^{1/2}$.

According to the discussion above, the higher-order and least squares bias correction approaches do not dominate each other in general, and whether one is preferred to the other depends on the data generating process and the relation between the two approximation spaces (or more precisely, the approximation power of $\bpt$ for functions in the linear span of $\bp$). Then a natural question follows: is there an optimal weighting scheme when $\bp$ and $\bpt$ are not nested? Again, assume $\e_i$'s are homoskedastic for simplicity, and we take a weighted average of $\widehat{\mu}_1$ and $\widehat{\mu}_2$:
$\widehat{\mu}_{w,\mathtt{bc}}:=w\widehat{\mu}_2+(1-w)\widehat{\mu}_1$
where $w\in[0,1]$. Then using a conclusion in the proof of Lemma \ref{lem: Avar}, the change in the integrated
asymptotic variance (weighted by the design density $f(\bx)$) is
\[w^2\sigma^2\int_{\mathcal{X}}\bp(\bx)' \bQ_m^{-1}(\bQ_m-\bQ_{m,\tilde{m}}\bQ_{\tilde{m}}^{-1}\bQ_{m, \tilde{m}}')\bQ_m^{-1}\bp(\bx)f(\bx)d\bx=: w^2\bar{\mathscr{V}}.
\]
On the other hand, by the property of projection operator, the change in the integrated squared bias (weighted by the design density $f(\bx)$) is
\[(w^2-2w)\int_{\mathcal{X}}\left(\mathscr{L}_{\bp}[\mathfrak{B}_{\tilde{m},\mathbf{0}}](\bx)\right)^2f(\bx)d\bx=:(w^2-2w)\bar{\mathscr{B}}.
\]
It is easy to see the optimal weight is $w^*=\bar{\mathscr{B}}/(\bar{\mathscr{B}}+\bar{\mathscr{V}})$. Clearly, when variance is very small (e.g., $\sigma^2$ is small), $w^*$ is close to $1$ and $\widehat{\mu}_2$ is preferred, whereas when bias is small, $w^*$ is close to $0$ and one may want to use $\widehat{\mu}_1$.

Next, the comparison of plug-in bias correction with the other two is more complicated since $\widehat{\mu}_3$ generally cannot be viewed as a regression estimator with additional covariates, and the covariance between $\widehat{\mu}_0$ and the estimated bias does not vanish. For piecewise polynomials, however, all three bias correction approaches are simply equivalent under certain conditions.

To see this, suppose $\bp$ and $\bpt$ are constructed on the same partitioning scheme $\Delta$, but the order of basis increases from $m$ to $m+1$. $\widehat{\mu}_0$ and $\widehat{\mu}_1$ are linked by \eqref{eq: restricted ls} since $\widehat{\mu}_0$ can be viewed as a restricted regression estimator compared with $\widehat{\mu}_1$. Specifically, one can construct a polynomial series of order $m+1$ (with degree no greater than $m$) within each cell $\delta\in\Delta$, and then implement a local regression restricting the coefficients of the polynomial terms of degree $m$ to be $0$. Then the restriction matrix in this case can be written as
$\bR=\left[\mathbf{0} \quad \bI_\vartheta \right]$
where $\vartheta$ denotes the number of polynomial terms of degree $m$ and $\bR$ is a $\vartheta\times\tilde{K}$ matrix. Plug it in \eqref{eq: restricted ls}, and then use the formula for matrix inverse in block form to obtain $(\tilde{\bP}'\tilde{\bP})^{-1}$. It is easy to see that the second term on the RHS of \eqref{eq: restricted ls}, $\bp(\bx)'\mathbf{M}\hat{\bbeta}_{\mathtt{ur}}$, is exactly the same as the leading bias derived in \cite{Cattaneo-Farrell_2013_JoE} (see their proof of Theorem 3) with the $m$th derivative estimated by piecewise polynomials of order $m+1$. As explained in the proof of Lemma \ref{lemma: Partitioning Estimators}, the leading approximation error can be alternatively expressed in terms of Legendre polynomials. Asymptotically, the two expressions are equivalent since the ``locally'' orthogonalized polynomials of degree $m$ will converge to the $m$th Legendre polynomials given by Lemma \ref{lemma: Partitioning Estimators}. 

For splines or wavelets, we do not have the above equivalence in general, since bases of different orders do not generate nested spaces, and the relative performance of the three approaches depends on the relation between these approximation spaces.

\section{Implementation Details}\label{section: implementation}
In this section, we briefly discuss implementation details about choosing the IMSE-optimal tuning parameters. We restrict our attention to tensor-product partitions with the same number of knots used in every dimension. Thus the tuning parameter reduces to a scalar $\kappa$ which denotes the number of subintervals used in every dimension. Also, we let the weighting function $w(\bx)$ be the density of $\bx_i$. We offer two approaches: rule-of-thumb (ROT) and direct plug-in (DPI).
\subsection{Rule-of-Thumb Choice}
The rule-of-thumb choice is based on the special case considered in Theorem 4.2  of the main paper. Specifically, assume $\bq=\mathbf{0}$ and knots are evenly spaced. The implementation steps are summarized as follows.

\begin{itemize}
	\item \textbf{Preliminary regression}.
	Implement a preliminary regression to estimate $\mu(\cdot)$ by a global polynomial of order $(m+2)$. Denote this estimate by $\hat{\mu}_{\mathtt{pre}}(\cdot)$.
	
	\item \textbf{Bias constant}:
	Use this preliminary regression to obtain an estimate of the $m$th derivatives of $\mu(\cdot)$, i.e.,
	$\widehat{\partial^{\bm{u}} \mu}(\cdot)
	=\partial^{\bm{u}}\widehat{\mu}_{\mathtt{pre}}(\cdot)$, for each $\bm{u} \in \Lambda_m$.
	Then an estimate of the bias constant is
	\[\widehat{\mathscr{B}}_{\bm{u}_1, \bm{u}_2,\mathbf{0}}=
	\eta_{\bm{u}_1, \bm{u}_2, \mathbf{0}}
	\times \frac{1}{n}\sum_{i=1}^{n}\widehat{\partial^{\bm{u}_1} \mu}(\bx_i)\widehat{\partial^{\bm{u}_2}\mu}(\bx_i).\]
	
	\item \textbf{Variance constant}.
	Implement another regression using a global polynomial of order $(m+2)$ to estimate $\E[y_i^2|\bx_i=\bx]$. Combining it with $\hat{\mu}_{\mathtt{pre}}(\cdot)$, we can obtain an estimate of the conditional variance function, denoted by $\hat{\sigma}^2(\cdot)$, since $\sigma^2(\bx)=\E[y_i^2|\bx_i=\bx]-(\E[y_i|\bx_i=\bx])^2$. Then an estimate of the variance constant is
	\[
	\widehat{\mathscr{V}}_{\mathbf{0}} =
	\begin{cases}
	\frac{1}{n}\sum_{i=1}^{n}\hat{\sigma}^2(\bx_i) \quad &\text{for splines and wavelets,}\\
	\binom{d+m-1}{m-1}\times\frac{1}{n}\sum_{i=1}^{n}\hat{\sigma}^2(\bx_i) \quad &\text{for piecewise polynomials.}
	\end{cases}
	\]
	
	\item \textbf{Rule-of-thumb $\hat{\kappa}_{\mathtt{ROT}}$}.
	Using the above results, a simple rule-of-thumb choice of $\kappa$ is
	\[\hat{\kappa}_{\mathtt{ROT}}=\bigg\lceil\left(
	\frac{2m\sum_{\bm{u}_1, \bm{u}_2\in\Lambda_m}\widehat{\mathscr{B}}_{\bm{u}_1, \bm{u}_2,\mathbf{0}}}{d\widehat{\mathscr{V}}_{\mathbf{0}}}\right)^{\frac{1}{2m+d}} n^{\frac{1}{2m+d}}\bigg\rceil.\]
\end{itemize}
Rather than assume a uniform design and evenly-spaced partition, one may use a trimmed-from-below Gaussian reference model to estimate the density of $\bx_i$ and adjust the variance and bias constants based on Theorem 4.2.
Clearly, this choice of $\kappa$ is derived based on simplifying assumptions, but it still has the correct rate ($\asymp n^{\frac{1}{2m+d}}$) even in other cases.

\subsection{Direct Plug-in Choice}
A direct plug-in (DPI) procedure is summarized in the following.

\begin{itemize}
	\item \textbf{Preliminary choice of $\kappa$}:
	Implement the ROT procedure to obtain $\hat{\kappa}_{\mathtt{ROT}}$.
	
	\item \textbf{Preliminary regression}.
	Given the user-specified basis (splines, wavelets, or piecewise polynomials), knot placement scheme (``uniform" or ``quantile") and rule-of-thumb choice $\hat{\kappa}_{\mathtt{ROT}}$, implement a series regression of order $(m+1)$ to obtain derivative estimates for every $\bm{u} \in \Lambda_m$. Denote this preliminary estimate by $\widehat{\partial^{\bm{u}}\mu}_{\mathtt{pre}}(\cdot)$.
	
	\item \textbf{Bias constant}.
	Construct an estimate $\widehat{\mathscr{B}}_{m,\bq}(\cdot)$ of the leading error $\mathscr{B}_{m,\bq}(\cdot)$ simply by replacing $\partial^{\bm{u}}\mu(\cdot)$ by $\widehat{\partial^{\bm{u}}\mu}_{\mathtt{pre}}(\cdot)$.
	$\widehat{\mathscr{B}}_{m, \mathbf{0}}(\cdot)$ can be obtained similarly. Then we use the pre-asymptotic version of the conditional bias to estimate the bias constant:
	\[\widehat{\mathscr{B}}_{\kappa,\bq} =
	\frac{1}{n}\sum_{i=1}^{n}\left(\widehat{\mathscr{B}}_{m, \bq}(\bx_i)-\widehat{\bgamma}_{\bq,0}(\bx_i)'\mathbb{E}_n[\bPi{0}(\bx_i)\widehat{\mathscr{B}}_{m,\mathbf{0}}(\bx_i)]\right)^2.
	\]
	
	\item \textbf{Variance constant}.
	Implement a series regression of order $m$ with $\kappa=\hat{\kappa}_{\mathtt{ROT}}$, and then use the pre-asymptotic version of the conditional variance to obtain an estimate of the variance constant. Specifically,
	\begin{equation*}
	\widehat{\mathscr{V}}_{\kappa, \bq} =
	\frac{1}{n}\sum_{i=1}^{n}\widehat{\bgamma}_{\bq,0}(\bx_i)'\bSighat{0}
	\widehat{\bgamma}_{\bq,0}(\bx_i), \quad
	\widehat{\bSigma}_0=\mathbb{E}_n[\bPi{0}(\bx_i)\bPi{0}(\bx_i)'
	\widehat{\e}^2_{i,0}],
	\end{equation*}
	where $\widehat{\e}_{i,0}$'s are regression residuals. Different weighting schemes for residuals may be used, leading to various heteroskedasticity-consistent variance estimates. 
	
	\item \textbf{Direct plug-in $\hat{\kappa}_{\mathtt{DPI}}$}.
	Combining these results, a direct plug-in choice of $\kappa$ is
	\[\hat{\kappa}_{\mathtt{DPI}}=
	\bigg\lceil\left(\frac{2(m-[\bq])\kappa_{\mathtt{ROT}}^{2(m-[\bq])}\widehat{\mathscr{B}}_{\kappa,\bq}}{(d+2[\bq])\kappa_{\mathtt{ROT}}^{-(d+2[\bq])}\widehat{\mathscr{V}}_{\kappa,\bq}}\right)^{\frac{1}{2m+d}} n^{\frac{1}{2m+d}}\bigg\rceil.
	\]
\end{itemize}


\section{Simulations}\label{section: numerical}
In this section, we present detailed simulation results. We consider the following regression functions:
\begin{equation*}
\begin{split}
\text{Model 1: }&\mu(x) = \sin\left(\frac{\pi}{2} (2x-1)\right) /
                          (1 + 2(2x-1)^2(\sign(2x-1)+1)) \\
\text{Model 2: }&\mu(x) = \sin\left(\frac{3\pi}{2}(2x-1)\right) /
                          (1 + 18(2x-1)^2(\sign(2x-1)+1))\\
\text{Model 3: }&\mu(x) = 2x - 1 + 5\phi(20x - 10)\\
\text{Model 4: }&\mu(\bx) = \sin(5x_1)\sin(10x_2)\\
\text{Model 5: }&\mu(\bx) = (1-(4x_1-2)^2)^2\sin(5x_2)/5\\
\text{Model 6: }&\mu(\bx) = (1-(4x_1-2)^2)^2(2x_2-1)(x_3-0.5)\\
\text{Model 7: }&\mu(\bx) = \tau(x_1)\tau(x_2)\tau(x_3)
\end{split}
\end{equation*}
where $\sign(x)=-1,0$, or $1$ if $x<0$, $x=0$, or $x>0$ respectively, $\phi(\cdot)$ is the standard Normal density function, and $\tau(x)=(x-0.5)+8(x-0.5)^2+6(x-0.5)^3-30(x-0.5)^4-30(x-0.5)^5$.

Given each regression function, we generate $(y_i, \bx_i)_{i=1,\ldots, n}$ by
$y_i=\mu(\bx_i) + \e_i$ where $\bx_i\sim i.i.d.\,\mathsf{U}[0,1]^d$ and 
$\e_i\sim i.i.d. \,\mathsf{N}(0,1)$. For each model, we generate $5,000$ simulated datasets with $n=1,000$.

We first present results for (tensor-product) spline regressions in Tables \ref{table:pointwise_uni1}--\ref{table:uniform7}.
For each model, we use linear splines ($m=2$) to form the classical point estimates. For robust inference, we use quadratic splines ($\tilde{m}=3$) to implement bias correction. Both evenly-spaced and quantile-spaced knot placements are considered, and for simplicity point estimators and bias correction employ the same knot placements ($\Delta=\tilde{\Delta}$).

For each model, we present three sets of simulation evidence. First, for three fixed points, we calculate the (simulated) root mean squared error (RMSE), coverage rate (CR) and average confidence interval length (AL). The nominal coverage is set to be $95\%$. Second, to evaluate the performance of our rule-of-thumb (ROT) and direct plug-in (DPI) knot selection procedures, we show some basic summary statistics (mean, median, standard deviation, etc.) of the selected number of knots ($\hat{\kappa}_{\mathtt{ROT}}$ and $\hat{\kappa}_{\mathtt{DPI}}$). Third, for uniform confidence bands, we calculate: the proportion of values covered with probability at least $95\%$ (CP), average coverage errors (ACE), and average width of confidence band (AW), and uniform coverage rate (UCR). The quantile estimation based on the plug-in and bootstrap methods uses $1000$ random draws conditional on the data for each simulated dataset.

Next, we evaluate the performance of (Daubechies) wavelet regressions. 
Several distinctive features of wavelets should be noted. First, the tuning parameter of wavelets ($s$) is referred to as ``resolution level", which relates to the number of knots $\kappa$ by $\kappa=2^s$. Thus, as $s$ increases, the number of series terms grows very fast. This issue is even more severe for tensor-product wavelets. Given our relatively small sample size, we restrict the wavelet-based simulations to Model 1-5. Second, due to the lack of smoothness of low-order wavelet basis, the plug-in bias correction approach ($j=3$) is not feasible unless very high-order wavelet basis is used to estimate the derivatives in the leading bias, which may not be of practical interest. Thus, we only report results based on higher-order and least-squares bias correction (as well as classical estimators). 

Results for wavelets are reported in Table \ref{table:pointwise_wav1}--\ref{table:uniform_wav5}, where Daubechies (father) wavelets of order $2$ ($m=2$) are used to form classical estimators, and Daubechies wavelets of order $3$ ($\tilde{m}=3$) for bias correction. In addition, as explained above, the IMSE-optimal resolution level is related to the number of knots by  $s_{\mathtt{IMSE}}=\log_2\kappa_{\mathtt{IMSE}}$. Similarly, for the estimated resolution levels, $\hat{s}_{\mathtt{ROT}}=\log_2\hat\kappa_{\mathtt{ROT}}$, and 
$\hat{s}_{\mathtt{DPI}}=\log_2\hat\kappa_{\mathtt{DPI}}$.

Our simulation results show that the bias correction methods perform generally well in both pointwise and uniform inference. The following are some practical guidance and caveats. First, as discussed in Section \ref{section: extension}, no bias correction approach among the three dominates the others in general. When high-order bias is large, their differences may be more pronounced. Second, our simulation shows that when $d>1$, our estimators still perform relatively well if a small number of knots are used. If not, as expected, their performance could be poor since the number of regressors explodes in such cases. Third, the ROT selection procedure used in the simulation study employs a global polynomial regression of degree $m+2$ to form the estimates of derivatives in the leading bias, which may be too conservative for highly nonlinear models. In such cases, a higher-order global polynomial regression may offer a better initial choice of the tuning parameter and improve the performance of DPI selection.

\begin{figure}[!h]
	\centering
	\vspace{1in}
	\caption{Confidence Bands, Model 1, Spline}
	\includegraphics[height=0.5\textheight, width=0.8\textwidth]{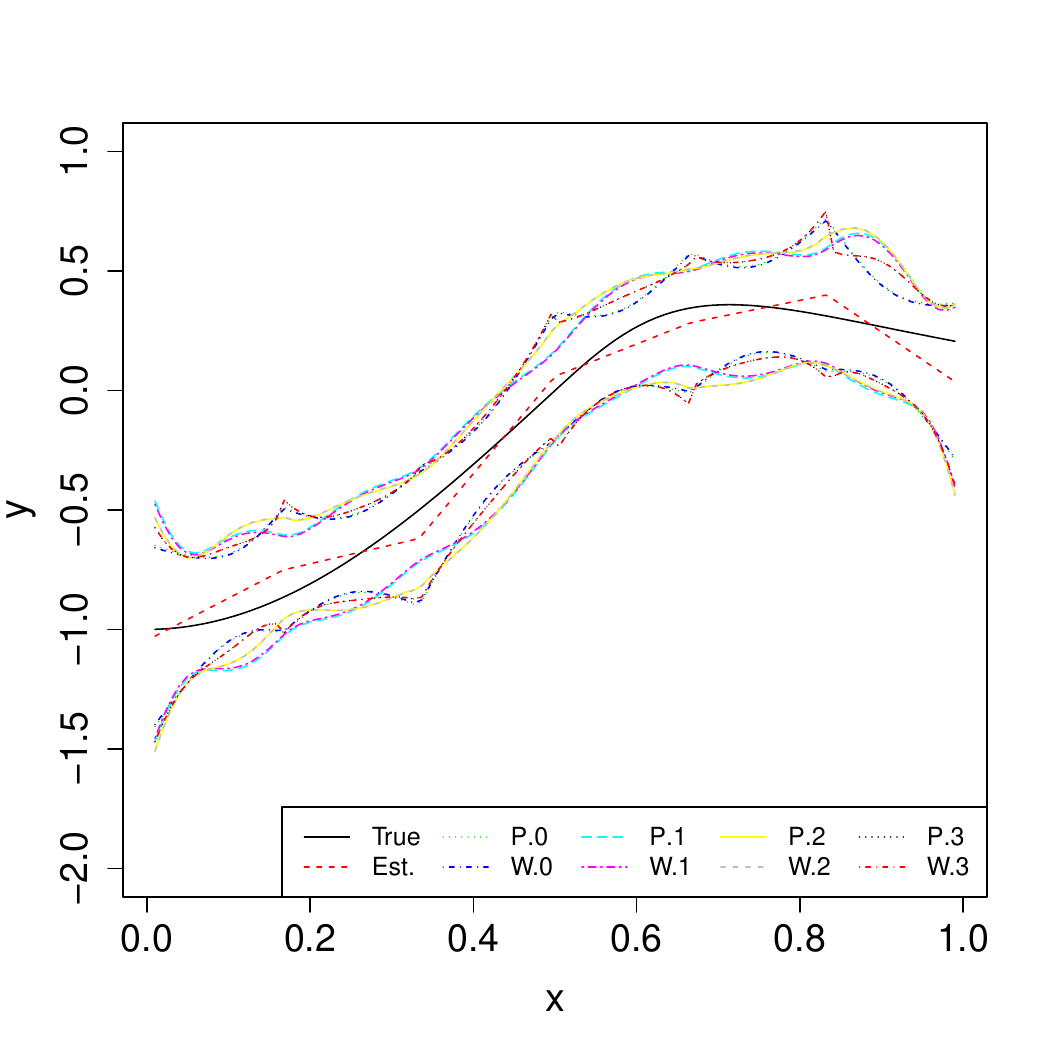}
	\label{figure:confidence band}
	\caption*{Note: ``P.$j$" denotes the confidence band based on plug-in approximation using the $j$th estimator, and ``W.$j$" based on wild bootstrap. $j=0$: the classical estimator; $j=1$: higher-order bias correction; $j=2$: least-squares bias correction; $j=3$: plug-in bias correction.}
\end{figure}

In the end, we present a figure (see Figure \ref{figure:confidence band}) as a visual illustration, which compares confidence bands based on different estimators suggested in our paper. We use Model 1 to generate a simulated dataset, and both plug-in approximation and wild bootstrap are used in constructing bands.

\newpage
\section{Proofs}\label{section: proofs-supplement}
\subsection{Proof of Lemma \ref{lem: Q matrix}}
\begin{proof}
	We first prove the boundedness of the eigenvalues of $\bQ_m$. We use the fact that
	$\lambda_{\max}(\bQ_m)=\max_{\mathbf{a}'\mathbf{a}=1}\,\mathbf{a}' \bQ_m\mathbf{a}$ and $\lambda_{\min}(\bQ_m)=\min_{\mathbf{a}'\mathbf{a}=1}\,\mathbf{a}' \bQ_m\mathbf{a}$.
	By definition of $\bQ_m$,
	\[
	\mathbf{a}' \bQ_m\mathbf{a}=\mathbf{a}'\E[\bp(\bx_i)\bp(\bx_i)']\mathbf{a}
	=\int_{\mathcal{X}} \Big(\sum_{k=1}^K a_k p_k(\bx)\Big)^2f(\bx)\,d\bx
	=: \|s(\bx)\|_{F, L_2(\mathcal{X})}^2
	\]
	where $s(\bx)=\sum_{k=1}^K a_k p_k(\bx)$.
	By Assumption \ref{assumption: DGP},
	$\|s(\bx)\|_{L_2(\mathcal{X})}^2\lesssim\|s(\bx)\|_{F, L_2(\mathcal{X})}^2 \lesssim\|s(\bx)\|_{L_2(\mathcal{X})}^2$.	
	By Assumption \ref{assumption: Local Basis}(a), the number of basis functions in $\bp(\cdot)$ which are active on a generic cell $\delta_l$, $l=1,\ldots, \bar{\kappa}$, is bounded by a constant. Denoted them by $(\bar{p}_{l,1},\cdots,\bar{p}_{l,M_l})'$, where $M_l$ may vary across $l$. It follows from Assumption \ref{assumption: Local Basis}(c) that
	$s(\bx)^2=(\sum_{k=1}^{M_l}a_k\bar{p}_{l,k}(\bx))^2
	\lesssim \sum_{k=1}^{M_l} a_k^2$ for all $\bx\in\delta_l$.
	Taking integral and summing over all $\delta_l$, we have $\|s(\bx)\|_{L_2(\mathcal{X})}^2\lesssim h^d$, and the upper bound on $\lambda_{\max}(\bQ_m)$ follows.
	On the other hand, by Assumption \ref{assumption: Local Basis}(b), $\|s(\bx)\|_{L_2(\mathcal{H}_k)}^2\gtrsim a_k^2 h^d$. Then taking sum over all $\mathcal{H}_k$,  $\|s(\bx)\|_{L_2(\mathcal{X})}^2\gtrsim h^d$, and the lower bound on $\lambda_{\min}(\bQ_m)$ follows.
	
	To derive the convergence rate of $\widehat{\bQ}_m$,
	let $\alpha_{k,l}=\frac{1}{n}\sum_{i=1}^{n}\alpha_{k,l}(i)$ be the $(k,l)$th element of $(\widehat{\bQ}_m-\bQ_m)$, where $\alpha_{k,l}(i):=p_k(\bx_i)p_l(\bx_i)-\E[p_k(\bx_i)p_l(\bx_i)]$. It follows from Assumption \ref{assumption: DGP} and \ref{assumption: Local Basis} that $\alpha_{j,l}$ is the sum of $n$ independent random variables with zero means, $|\alpha_{k,l}(i)|\lesssim 1$ uniformly over $i, k$ and $l$, and thus
	$\V[\alpha_{k,l}(i)]\lesssim h^d/n$.
	By Bernstein's inequality, for every $\vartheta>0$,
	$$\P(|\alpha_{k,l}|>\vartheta)\leq 2\exp\left(-\frac{\vartheta^2/2}{C_1h^d/n+C_2\vartheta/(3n)}\right).$$
	By Assumption \ref{assumption: Local Basis}, $(\widehat{\bQ}_m-\bQ_m)$ only has a finite number of nonzeros in any row or column, and thus for every $\vartheta>0$,
	\[\P(\max_{k,l}|\alpha_{k,l}| > \vartheta) \leq 2CK\exp\left(-\frac{\vartheta^2/2}{C_1h^d/n+C_2\vartheta/(3n)}\right).
	\]
	Then $\max_{k,l}\,|\alpha_{k,l}|\lesssim_\P h^d\sqrt{\log n/(nh^d)}$, which suffices to show that
	\[\|\widehat{\bQ}_m-\bQ_m\|_\infty\lesssim_\P h^{d}\sqrt{\log n/(nh^d)}
	\quad \text{and} \quad
	\|\widehat{\bQ}_m-\bQ_m\|_1\lesssim_\P h^d\sqrt{\log n/(nh^d)}.
	\]
	By the relation between induced operator norms,
	$\|\widehat{\bQ}_m-\bQ_m\|\lesssim_\P h^d\sqrt{\log n/(nh^d)}$.
	Hence when $\log n/(nh^d)=o(1)$, $\|\widehat{\bQ}_m-\bQ_m\|=o_\P(h^d)$.
	Notice that for any vector $\mathbf{a}\in\mathbb{R}^K$ such that $\mathbf{a}'\mathbf{a}=1$,
	$$\|\widehat{\bQ}_m-\bQ_m\|\geq|\mathbf{a}'(\widehat{\bQ}_m-\bQ_m)\mathbf{a}|=|\mathbf{a}'\widehat{\bQ}_m\mathbf{a}-\mathbf{a}' \bQ_m\mathbf{a}|.$$
	Since $h^d\lesssim\mathbf{a}' \bQ_m\mathbf{a}\lesssim h^d$, this suffices to show that $\|\widehat{\bQ}_m\|\lesssim_\P h^d$, and with probability approaching one, $\lambda_{\min}(\widehat{\bQ}_m)\gtrsim h^d$.
	
	Next, we show $\|\widehat{\bQ}_m^{-1}\|_\infty\lesssim_\P h^{-d}$. In the univariate case, it is easy to see that by Assumption \ref{assumption: Local Basis}, $\widehat{\bQ}_m$ is a banded matrix with a finite bandwidth (independent of $n$). In the multidimensional case $\widehat{\bQ}_m$ will take a ``multi-layer'' block banded matrix form.  
	Specifically, let $\mathbf{A}_{kl}$ be the generic $(k,l)$th block of $\widehat{\bQ}_m$. $\mathbf{A}_{kl}=\mathbf{0}$ if $|k-l|>L$ for some integer $L>0$, with
	each nonzero $\mathbf{A}_{kl}$ banded (``two layers'') or taking another block banded structure (``more than two layers''). To see this,	
	we first arrange the ordering of basis functions appropriately by ``rectangularizing'' the partition $\Delta$. By Assumption \ref{assumption: quasi-uniform partition}, $\Delta$ is quasi-uniform and there exists a universal measure of mesh size $h$. Construct an initial rectangular partition covering $\mathcal{X}$, which is formed as tensor products of intervals of length $h$. A generic cell in this partition is indexed by a $d$-tuple $(l_1, \ldots, l_d)$. Then for each cell in this partition, take its intersection with $\mathcal{X}$ and exclude all cells outside of $\mathcal{X}$. Thus we construct a ``trimmed'' tensor-product partition $\Delta^{\mathtt{rec}}=\{\delta^{\mathtt{rec}}_{l_1\ldots l_d}\}$.
	Clearly, each element $\delta\in\Delta$ is covered by a finite number of cells in $\Delta^{\mathtt{rec}}$. On the other hand, each $\delta^{\mathtt{rec}}\in\Delta^{\mathtt{rec}}$ is also overlapped with a finite number of cells in $\Delta$ (if not, cells in $\Delta$ overlapping with $\delta^{\mathtt{rec}}$ cannot be covered by a ball of radius $2h$). Arrange these cells by first increasing $l_d$ with other $l_\ell$'s fixed at the lowest values and then increasing $l_\ell$'s sequentially. Then arrange the basis functions in $\bp(\cdot)$ according to their supports. Specifically, start with basis functions that are active on the first cell, and then arrange those functions that are active on the second and have not been included yet. Continue this procedure until the functions active on the last cell have been included. According to this particular ordering, the Gram of this basis has the same banded structure as that of tensor-product local basis on tensor-product partitions. The nested banded structure involves at most $d$ layers. The bandwidth at each layer may be different, but is bounded by a universal constant $\bar{L}$ by Assumption \ref{assumption: Local Basis}.
	
	In the one-dimensional case, it follows from  \cite[Theorem 2.2]{Demko_1977_SIAM} that $\|\widehat{\bQ}_m^{-1}\|_\infty\lesssim_\P h^{-1}$. In the multidimensional case the original proof needs to be slightly modified. We only prove for the case when $\widehat{\bQ}_m$ is block banded with banded blocks (two layers of banded structures). The general case follows similarly.
	
	For universal constants $C_1$ and $C_2$, with probability approaching one, $\lambda_{\min}(\widehat{\bQ}_m)\geq C_1h^d$ and $\lambda_{\max}(\widehat{\bQ}_m)\leq C_2h^d$. Hence for any vector $\mathbf{a}\in\mathbb{R}^K$ such that $\mathbf{a}'\mathbf{a}=1$, there are some constants $C_3$, $C_4$ and $C_5$ such that with probability approaching one,
	\[0<C_3\leq\frac{\mathbf{a}'\widehat{\bQ}_m\mathbf{a}}{C_5h^d}\leq C_4<1.
	\]
	Hence $\boldsymbol{\Psi}:=\bI_K-\widehat{\bQ}_m/(C_5h^d)$ is a block banded matrix with banded blocks satisfying $\|\boldsymbol{\Psi}\|<1$ with probability approaching one.
	Therefore, we can write $C_5h^d\widehat{\bQ}_m^{-1}=\sum_{l=1}^{\infty}\boldsymbol{\Psi}^l$. The $(s,t)$th entry of $\widehat{\bQ}_m^{-1}$, denoted by $\alpha_{s,t}$, is an element of $\sum_{l=1}^{\infty}\boldsymbol{\Psi}^l$. Hence,
	\[|\alpha_{s,t}|\leq \sum_{l=\chi(s,t,\bar{L})}^{\infty}\|\boldsymbol{\Psi}^l\|=\frac{\|\boldsymbol{\Psi}\|^{\chi(s,t,L)}}{1-\|\boldsymbol{\Psi}\|}
	\]
	where $\chi(s, t, \bar{L})$ is a number depending on the row index $s$, column index $t$ and the upper bound on bandwidths $\bar{L}$. We further denote the block row index and column index of $\alpha_{s,t}$ as $r_s$ and $r_t$, and the row index and column index within the block containing it as $\iota_s$ and $\iota_t$. $|r_s-r_t|$ and $|\iota_s-\iota_t|$ measure how far away $\alpha_{s,t}$ is from the diagonals of the entire matrix and the block it belongs to. As in the one-dimensional case, the first few products of $\boldsymbol{\Psi}$ do not contribute to off-diagonal blocks of the inverse matrix. As $|r_s-r_t|$ increases, $\chi(s,t,\bar{L})$ also gets larger. Meanwhile, since each block of $\widehat{\bQ}_m$ is also banded, $\chi(s,t,\bar{L})$ also increase with $|\iota_s-\iota_t|$. By the same argument as in \cite{Demko_1977_SIAM}, 
	$|\alpha_{s,t}|\leq (1-\|\boldsymbol{\Psi}\|)^{-1}
	\|\boldsymbol{\Psi}\|^{|r_s-r_t|/C_6+|\iota_s-\iota_t|/C_6}$ where $C_6$ is some constant depending on $\bar{L}$. Thus $\alpha_{s,t}$ exponentially decays when $|r_s-r_t|$ or $|\iota_s-\iota_t|$ becomes large. By convergence of geometric series, $\|\widehat{\bQ}_m^{-1}\|\lesssim_\P h^{-d}$.
	$\|\bQ_m^{-1}\|_\infty\lesssim h^{-d}$ follows similarly.
	
	Finally, the above results immediately imply that
	\[
	\begin{split}
	&\|\widehat{\bQ}_m^{-1}-\bQ_m^{-1}\|_\infty\leq \|\widehat{\bQ}_m^{-1}\|_\infty\|\widehat{\bQ}_m-\bQ_m\|_\infty\|\bQ_m^{-1}\|_\infty\lesssim_\P h^{-d}\sqrt{\log n/(nh^d)}, \quad \text{and}\\
	&\|\widehat{\bQ}_m^{-1}-\bQ_m^{-1}\|\leq \|\widehat{\bQ}_m^{-1}\|\|\widehat{\bQ}_m-\bQ_m\|\|\bQ_m^{-1}\|\lesssim_\P h^{-d}\sqrt{\log n/(nh^d)}.
	\end{split}
	\]
	The bound on $\|\E_n[\bp(\bx_i)\bp(\bx_i)'\sigma^2(\bx_i)]
	-\E[\bp(\bx_i)\bp(\bx_i)'\sigma^2(\bx_i)]\|$
	follows similarly, since $\sigma^2(\bx)\lesssim 1$ uniformly over $\bx\in\mathcal{X}$.
\end{proof}

\subsection{Proof of Lemma \ref{lem: bias}}
\begin{proof}
	We first prove the uniform bound on the conditional bias. By Assumption \ref{assumption: Misspecification Error} we can find $s^*\in\mathcal{S}_{\Delta,m}$ such that
	$\sup_{\bx\in\mathcal{X}}|\partial^\bq\mu(\bx)-\partial^\bq s^*(\bx)|\lesssim h^{m-[\bq]}$. Since
	\begin{align*}
	\E[\widehat{\partial^\bq\mu}_0(\bx)|\bX] - \partial^\bq \mu(\bx)
	&=\partial^\bq\bp(\bx)'\widehat{\bQ}_m^{-1}\E_n[\bp(\bx_i) \mu(\bx_i)]-\partial^\bq \mu(\bx)\\
	&=\partial^\bq\bp(\bx)'\widehat{\bQ}_m^{-1}
	\E_n[\bp(\bx_i)(\mu(\bx_i)-s^*(\bx_i))]+
	\partial^\bq s^*(\bx)-\partial^\bq \mu(\bx),
	\end{align*}
	it suffices to show that the first term in the second line is properly bounded uniformly over $\bx\in\mathcal{X}$. It follows from Lemma \ref{lem: Q matrix} and Assumption \ref{assumption: Local Basis} that
	\begin{equation*}
	\quad \sup_{\bx\in\mathcal{X}}\,
	\Big|\partial^\bq\bp(\bx)'\widehat{\bQ}_m^{-1}
	\E_n[\bp(\bx_i)(\mu(\bx_i)-s^*(\bx_i))]\Big|
	\lesssim_\P h^{-[\bq]-d}\|\E_n[\bp(\bx_i)(\mu(\bx_i)-s^*(\bx_i))]\|_\infty.
	\end{equation*}
	By Assumption \ref{assumption: Local Basis} and \ref{assumption: Misspecification Error}, $\max_{1\leq k\leq K}|\E[p_k(\bx_i)(\mu(\bx_i)-s^*(\bx_i))]|\lesssim h^{m+d}$,
	$\V[p_k(\bx_i)(\mu(\bx_i)-s^*(\bx_i))]\lesssim h^{d+2m}$,
	and $\sup_{\bx\in\mathcal{X}}|p_k(\bx)(\mu(\bx)-s^*(\bx))|\lesssim h^m$. By Bernstein's inequality, for any $\vartheta>0$,
	\begin{align*}
	&\P\left(\max_{1\leq k\leq K}
	\Big|\E_n[p_k(\bx_i)(\mu(\bx_i)-s^*(\bx_i))]-\E[p_k(\bx_i)(\mu(\bx_i)-s^*(\bx_i))]\Big|>\vartheta\right)\\
	\leq&\,\sum_{k=1}^{K}\P\left(\Big|\E_n[p_k(\bx_i)(\mu(\bx_i)-s^*(\bx_i))]-\E[p_k(\bx_i)(\mu(\bx_i)-s^*(\bx_i))]\Big|>\vartheta\right)\\
	\leq&\, 2K\exp\left(\frac{-\vartheta^2/2}{C_1h^{m}\vartheta/(3n)+
		  C_2h^{2m+d}/n}\right).
	\end{align*}
	Therefore,
	\[
	\max_{1\leq k\leq K}\Big|\E_n[p_k(\bx_i)(\mu(\bx_i)-s^*(\bx_i))]-\E[p_k(\bx_i)(\mu(\bx_i)-s^*(\bx_i))]\Big|
	\lesssim_\P h^{m+d}\sqrt{\frac{\log n}{nh^d}},\]
	which suffices to prove the desired uniform bound since $\frac{\log n}{nh^d}=o(1)$. 
		
	Next, we prove the leading bias expansion. By Assumption \ref{assumption: Misspecification Error}, we can find $s^*\in\mathcal{S}_{\Delta, m}$ such that
	$\sup_{\bx\in\mathcal{X}}|\partial^\bq\mu(\bx)-
	\partial^\bq s^*(\bx) + \mathscr{B}_{m,\bq}(\bx)|
	\lesssim h^{m+\varrho-[\bq]}$.
	Then the conditional bias can be expanded as follows:
	\begin{align*}
	&\E[\widehat{\partial^\bq\mu}_0(\bx)|\bX]-\partial^{\bq}\mu(\bx)=\partial^\bq\bp(\bx)'\widehat{\bQ}_m^{-1}\E_n[\bp(\bx_i)\mu(\bx_i)]
	- \partial^\bq \mu(\bx) \\
	=&\,\partial^\bq\bp(\bx)'\widehat{\bQ}_m^{-1}
	\E_n[\bp(\bx_i)\mu(\bx_i)]-\partial^\bq s^*(\bx)
	+\mathscr{B}_{m,\bq}(\bx)+O(h^{m+\varrho-[\bq]}) \\
	=&\,\mathscr{B}_{m,\bq}(\bx) +\partial^\bq\bp(\bx)'\widehat{\bQ}_m^{-1}
	\E_n[\bp(\bx_i)(\mu(\bx_i)-s^*(\bx_i))]+O(h^{m+\varrho-[\bq]}).
	\end{align*}
	The second term in the last line can be further written as
	\begin{align*}
	&\partial^\bq\bp(\bx)'\widehat{\bQ}_m^{-1}\E_n[\bp(\bx_i)(\mu(\bx_i)-s^*(\bx_i))]\\
	&=\,-\partial^\bq\bp(\bx)'\widehat{\bQ}_m^{-1}
	\E_n\left[\bp(\bx_i)\mathscr{B}_{m,\mathbf{0}}(\bx_i)\right]\\
	&\qquad+\partial^\bq\bp(\bx)'\widehat{\bQ}_m^{-1}
	\E_n\Big[\bp(\bx_i)(\mu(\bx_i)-s^*(\bx_i)+\mathscr{B}_{m,\mathbf{0}}(\bx_i))\Big].
	\end{align*}
	By the same argument used to bound $\|\E_n[\bp(\bx_i)(\mu(\bx_i)-s^*(\bx_i))]\|_\infty$,
	\[\max_{1\leq k\leq K}\Big|\E_n\Big[p_k(\bx_i)(\mu(\bx_i)
	-s^*(\bx_i)+\mathscr{B}_{m, \mathbf{0}}(\bx_i))\Big]\Big|
	\lesssim_\P h^{m+\varrho+d}.
	\]
	Then by Assumption \ref{assumption: Local Basis} and Lemma \ref{lem: Q matrix},
	\[\|\partial^\bq\bp(\bx)'\widehat{\bQ}_m^{-1}
	\E_n[\bp(\bx_i)(\mu(\bx_i)-s^*(\bx_i)+\mathscr{B}_{m,\mathbf{0}}(\bx_i))]\|_{L_\infty(\mathcal{X})}
	\lesssim_\P h^{m+\varrho-[\bq]},
	\]
	which suffices to show the desired bias expansion.
	
	Suppose that Equation \eqref{eq: approx orth} holds. Then by Lemma \ref{lem: Q matrix} and Assumption \ref{assumption: Local Basis}, the second term in the bias expansion is $o_\P(h^{m-[\bq]})$. Thus the leading bias is further reduced to $\mathscr{B}_{m,\bq}(\bx)$. The proof of Lemma \ref{lem: bias} is complete.
\end{proof}

\subsection{Proof of Lemma \ref{lem: basis}}
\begin{proof}
	For $j=0,1$, the results immediately follow from Assumption \ref{assumption: Local Basis}, \ref{assumption: BC bases} and Lemma \ref{lem: Q matrix}. For $j=2$,
	\begin{align*}
	\sup_{\bx\in\mathcal{X}}\|\bgamma_{\bq, 2}(\bx)'\|_\infty
	&\leq
	\sup_{\bx\in\mathcal{X}}\|\bgamma_{\bq, 0}(\bx)'\|_\infty +
	\sup_{\bx\in\mathcal{X}}\|\bgamma_{\bq, 0}(\bx)'\bQ_{m, \tilde{m}}\bQ_{\tilde{m}}^{-1}\|_\infty +
	\sup_{\bx\in\mathcal{X}}\|\bgamma_{\bq, 1}(\bx)'\|_\infty\\
	&\lesssim h^{-d-[\bq]} + h^{-d-[\bq]}\|\bQ_{m, \tilde{m}}\|_\infty\|\bQ_{\tilde{m}}^{-1}\|_\infty.
	\end{align*}
	By Assumptions \ref{assumption: BC bases}, both $\bp(\cdot)$ and $\bpt(\cdot)$ are local bases. Then $\bQ_{m,\tilde{m}}$ has a finite number of nonzero elements in any row or column, and all elements in $\bQ_{m,\tilde{m}}$ are bounded by $Ch^d$ for some universal constant $C$. Thus
	$\|\bQ_{m,\tilde{m}}\|_\infty\lesssim h^d$,
	$\|\bQ_{m, \widetilde{m}}\|_1\lesssim h^d$ and 
	$\|\bQ_{m, \widetilde{m}}\|\lesssim h^d$. Then $\|\bgamma_{\bq, 2}(\bx)'\|_\infty\lesssim h^{-d-[\bq]}$.
	Similarly, $\|\bgamma_{\bq, 2}(\bx)\|\lesssim h^{-d-[\bq]}$. The lower bound on the $L_2$-norm follows by Assumption \ref{assumption: Local Basis}.
	Moreover, note that
    \begin{align*}
    &\|\bg{2}' - \bgamma_{\bq, 2}(\bx)'\|_\infty \\
    \leq&\, \|\bg{0}'-\bgamma_{\bq, 0}(\bx)'\|_\infty
    +\|\bg{0}'\widehat{\bQ}_{m, \tilde{m}}
    \widehat{\bQ}_{\tilde{m}}^{-1} - 
    \bgamma_{\bq, 0}\bQ_{m, \tilde{m}}\bQ_{\tilde{m}}^{-1}\|_\infty\\
    &+\|\bg{1}'-\bgamma_{\bq, 1}(\bx)'\|_\infty\\
    \leq&\, \|(\bg{0}'-\bgamma_{\bq, 0}(\bx)')\widehat{\bQ}_{m, \tilde{m}}
    \widehat{\bQ}_{\tilde{m}}^{-1}\|_\infty + 
    \|\bgamma_{\bq, 0}(\bx)'(\widehat{\bQ}_{m, \tilde{m}} -
    \bQ_{m, \tilde{m}})\widehat{\bQ}_{\tilde{m}}^{-1}\|_\infty \\ 
    &+ \|\bgamma_{\bq, 0}(\bx)')\bQ_{m, \tilde{m}}
    (\widehat{\bQ}_{\tilde{m}}^{-1}-\bQ_{\tilde{m}}^{-1})\|_\infty
    + h^{-d-[\bq]}\sqrt{\log n/(nh^d)},
    \end{align*}
    where the last line uses the results for $j=0,1$. Using the sparsity of $(\widehat{\bQ}_{m,\tilde{m}}-\bQ_{m, \tilde{m}})$ and the same argument for Lemma \ref{lem: Q matrix},
    $\|\widehat{\bQ}_{m,\tilde{m}}-\bQ_{m, \tilde{m}}\|_\infty
    \lesssim_\P h^d\sqrt{\log n/(nh^d)}$. Then the desired uniform bound follows from Lemma \ref{lem: Q matrix} and Assumption \ref{assumption: Local Basis}. The $L_2$-bound follows similarly.
    
    For $j=3$, notice that
    \begin{align*}
    \|\bgamma_{\bq, 3}(\bx)'\|_\infty \leq 
    &\,\|\bgamma_{\bq, 0}(\bx)'\|_\infty + \sum_{\bm{u}\in\Lambda_m}\Big\|\bgamma_{\bm{u}, 1}(\bx)'h_\bx^{m-[\bq]}B_{\bm{u}, \bq}(\bx)\Big\|_\infty \\
    &+\sum_{\bm{u}\in\Lambda_m}\Big\|\bgamma_{\bq, 0}(\bx)'\E[\bp(\bx_i)h_{\bx_i}^m B_{\bm{u},\mathbf{0}}(\bx_i)\partial^{\bm{u}}\bpt(\bx_i)']
    \bQ_{\tilde{m}}^{-1}\Big\|_\infty
    \end{align*}
    By Assumption \ref{assumption: Misspecification Error} and \ref{assumption: BC bases}, both $\bp(\cdot)$ and $\bpt(\cdot)$ are locally supported and all elements in $\bp(\bx_i)h_{\bx_i}^mB_{\bm{u}, \bm{0}}(\bx_i)\partial^{\bm{u}}
    \bpt(\bx)'$ are bounded by a universal constant. Using the argument given in Lemma \ref{lem: Q matrix},
    \[
    \begin{split}
    &\|\E[\bp(\bx_i)h_{\bx_i}^m B_{\bm{u},\mathbf{0}}(\bx_i)\partial^{\bm{u}}\bpt(\bx_i)']\|
    \lesssim h^d \qquad \text{and}\\
    &\|\E_n[\bp(\bx_i)h_{\bx_i}^m B_{\bm{u},\mathbf{0}}(\bx_i)\partial^{\bm{u}}\bpt(\bx_i)' -
    \E[\bp(\bx_i)h_{\bx_i}^m B_{\bm{u},\mathbf{0}}(\bx_i)\partial^{\bm{u}}\bpt(\bx_i)']\|
    \lesssim_\P h^d\sqrt{\frac{\log n}{nh^d}}.
    \end{split}
    \]
    Then the desired results follow from Assumption \ref{assumption: Local Basis} and Lemma \ref{lem: Q matrix}.
\end{proof}	

\subsection{Proof of Lemma \ref{lem: Avar}}
\begin{proof}
	For $j=0, 1$, the results directly follow from Assumption \ref{assumption: DGP}, \ref{assumption: Local Basis} and Lemma \ref{lem: Q matrix}.
	The proof for $j=2,3$ is divided into three steps.
	\smallskip
	
	\textit{Step 1:} We first establish the upper bounds on $\Omega_j(\bx)$.
	By Assumption \ref{assumption: DGP},
	\[\bSig{j}=\E[\e_i^2\bPi{j}(\bx_i)\bPi{j}(\bx_i)']\lesssim \E[\bPi{j}(\bx_i)\bPi{j}(\bx_i)'].
	\]
	To bound $\|\bSig{j}\|$, it suffices to give an upper bound on $\E[\bPi{j}(\bx_i)\bPi{j}(\bx_i)']$. By Assumption \ref{assumption: BC bases} and the same argument used in the proof of Lemma \ref{lem: Q matrix}, we have 
	$\|\E[\bPi{j}(\bx_i)\bPi{j}(\bx_i)]\|\lesssim h^d$. By Lemma \ref{lem: basis}, $\sup_{\bx\in\mathcal{X}}\|\bgamma_{\bq, j}(\bx)'\|\lesssim h^{-d-[\bq]}$. Thus $\sup_{\bx\in\mathcal{X}}\Omega_j(\bx)\lesssim h^{-d-2[\bq]}$.
	\smallskip
	
	\textit{Step 2:}
	Next, we show the lower bound on $\Omega_2$. Since $\sigma^2(\bx)\gtrsim 1$ uniformly over $\bx\in\mathcal{X}$, we have
	$\Omega_2(\bx)\gtrsim\bgamma_{\bq, 2}(\bx)'
	\E[\bPi{2}(\bx_i)\bPi{2}(\bx_i)']\bgamma_{\bq, 2}(\bx)$.
	Expanding the expression on the RHS of this inequality, we have a trivial lower bound:
	\begin{align*}
	&\bgamma_{\bq,2}(\bx)'\E[\bPi{2}(\bx_i)\bPi{2}(\bx_i)']\bgamma_{\bq,2}(\bx)\\
	=\,&\partial^\bq\bp(\bx)' \bQ_m^{-1}\partial^\bq\bp(\bx)
	+(\bgamma_{\bq,0}(\bx)'\bQ_{m,\tilde{m}}-\partial^\bq\bpt(\bx)')\bQ_{\tilde{m}}^{-1}
	(\bgamma_{\bq,0}(\bx)'\bQ_{m,\tilde{m}}-\partial^\bq\bpt(\bx)')'
	\\
	& -2\partial^\bq\bp(\bx)' \bQ_m^{-1}\bQ_{m,\tilde{m}}\bQ_{\widetilde{m}}^{-1}(\bgamma_{\bq,0}(\bx)'\bQ_{m,\tilde{m}}-\partial^\bq\bpt(\bx)')'\\
	=\,&\partial^\bq\bpt(\bx)'\bQ_{\tilde{m}}^{-1}\partial^\bq\bpt(\bx)+
	\Big[\partial^\bq\bp(\bx)' \bQ_m^{-1}(\bQ_m-\bQ_{m,\tilde{m}}\bQ_{\tilde{m}}^{-1}
	\bQ_{m,\tilde{m}}')\bQ_m^{-1}\partial^\bq\bp(\bx)\Big].
	\end{align*}
    By properties of projection operator, $(\bQ_m-\bQ_{m,\tilde{m}}\bQ_{\tilde{m}}^{-1}\bQ_{m,\tilde{m}}')$ is positive semidefinite. By Assumption \ref{assumption: BC bases} and Lemma \ref{lem: Q matrix}, $\partial^\bq\bpt(\bx)'\bQ_{\tilde{m}}^{-1}\partial^\bq\bpt(\bx)\gtrsim h^{-d-2[\bq]}$,
    and thus the desired lower bound is obtained.
	\smallskip
	
	\textit{Step 3:} 
	Now let us bound $\Omega_3(\bx)$ from below. Suppose condition (i) in Assumption \ref{assumption: BC bases} holds. Then there exists a linear map $\bUpsilon$ such that $\bPi{3}(\cdot)=\bUpsilon\bpt(\cdot)$.
	By Lemma \ref{lem: Q matrix} and Assumption \ref{assumption: DGP},
	$\Omega_3(\bx)\gtrsim\bgamma_{\bq, 3}(\bx)'
	\bUpsilon\E[\bpt(\bx_i)\bpt(\bx_i)']\bUpsilon'\bgamma_{\bq,3}(\bx)
	\gtrsim h^d\bgamma_{\bq,3}(\bx)'\bUpsilon\bUpsilon'\bgamma_{\bq,3}(\bx)$.
	Define $\mathbf{v}(\bx)':=\bgamma_{\bq,3}(\bx)'\bUpsilon$. Notice that for any function $s(\bx)\in\mathsf{span}(\bp(\cdot))$, there exists some $\mathbf{c}\in\mathbb{R}^{\tilde{K}}$ such that $s(\bx)=\bpt(\bx)'\mathbf{c}$. It follows that $\mathbf{v}(\bx)'\E[\bpt(\bx_i)s(\bx_i)]=\partial^\bq s(\bx)$.
	Then we have
	$$\|\mathbf{v}(\bx)\|\geq
	\frac{\|\bv(\bx)'\E[\bpt(\bx_i)s(\bx_i)]\|}
	{\|\E[\bpt(\bx_i)s(\bx_i)]\|}=
	\frac{\|\partial^\bq s(\bx)\|}
	{\|\E[\bpt(\bx_i)s(\bx_i)]\|}.$$
	Since the choice of $s(\bx)$ is arbitrary within the span of $\bp$, by Assumption \ref{assumption: Local Basis}(c) we can take a function in $\bp(\bx)$ to be $s(\bx)$ such that $\|\partial^\bq s(\bx)\|\geq Ch^{-[\bq]}|$ where $C$ is a constant independent of $\bx$ and $n$. Also, $\|\E[\bpt(\bx_i)s(\bx_i)]\|\lesssim h^d$. Hence $h^{-[\bq]-d}\lesssim\inf_{\bx\in\mathcal{X}}\|\bv(\bx)\|$. The desired bound follows.
	\smallskip
	
	\textit{Step 4:}
	Now suppose that condition (ii) in Assumption \ref{assumption: BC bases} holds. Again, since $\sigma^2(\bx)\gtrsim 1$ uniformly over $\bx\in\mathcal{X}$,
	$\Omega_3(\bx)\gtrsim\bgamma_{\bq, 3}(\bx)'\E[\bPi{3}(\bx_i)\bPi{3}(\bx_i)']\bgamma_{\bq,3}(\bx)$.
	
	Define
	$\mathscr{K}_h(\bx,\bx_1) := \bgamma_{\bq,3}(\bx)'\bPi{3}(\bx_1)$.
	It suffices to bound $\E_{\bx_1}[\mathscr{K}_h(\bx,\bx_1)^2]$,
	where $\E_{\bx_1}$ denotes the expectation with respect to the distribution of $\bx_1$.
	We write 
	$$\langle g_1,g_2\rangle_\mathcal{U}:=\int_\mathcal{U} g_1(\bx_1)g_2(\bx_1)f(\bx_1)d\bx_1$$ for the inner product of $g_1(\cdot)$ and $g_2(\cdot)$ with respect to the probability measure of $\bx_1$ on $\mathcal{U}\subset\mathcal{X}$.
	Clearly, $\langle g_1, g_2\rangle_\mathcal{X}=\E[g_1(\bx_1)g_2(\bx_1)]$. By Cauchy-Schwartz inequality,
	for $g\in L^2(\mathcal{U})$,
	\[\langle\mathscr{K}_h(\bx,\bx_1),g(\bx_1)\rangle_\mathcal{U}^2\leq\langle\mathscr{K}_h(\bx,\bx_1),\mathscr{K}_h(\bx,\bx_1)\rangle_\mathcal{U}\cdot\langle g(\bx_1),g(\bx_1)\rangle_\mathcal{U}.
	\]
	Given an evaluation point $\bx\in\mathcal{X}$, choose a polynomial
	$\varphi_h(\bx_1;\bx):=\frac{(\bx_1-\bx)^\bq}{h^{[\bq]}}$.
	Clearly, $\partial^\bq\varphi_h(\bx_1;\bx)=h^{-[\bq]}$. By Assumption \ref{assumption: BC bases}, the operator $\langle \mathscr{K}_h(\bx,\bx_1),\cdot\rangle_\mathcal{X}$ reproduces the $\bq$th derivative of $\varphi_h(\bx_1;\bx)$ at $\bx$, i.e.
	$\langle\mathscr{K}_h(\bx,\bx_1),\varphi_h(\bx_1;\bx)\rangle_\mathcal{X}=h^{-[\bq]}$. 
	
	In addition, we rectangularize $\Delta$ as described in the proof of Lemma \ref{lem: Q matrix}. Then for each $\bx\in\mathcal{X}$, we can choose a rectangular region that contains $\bx$ and consists of a fixed number of subrectangles in $\Delta^{\mathtt{rec}}$. Thus, the size of this region shrinks as $n\rightarrow\infty$. Specifically, let
	\[\mathcal{W}_\bx := \left\{\mathbf{\check{x}}: t_{\ell,\bx}-L\leq \check{x}_j\leq t_{\ell,\bx}+L,\, \ell=1,\ldots, d\right\},
	\]
	where $t_{\ell,\bx}$ is the closest point in $\Delta^{\mathtt{rec}}$ that is no greater than $x_\ell$, and $L$ is some fixed number to be determined which only depends on $d$, $m$ and $\tilde{m}$. If such a region spans outside of $\mathcal{X}$, take its intersection with $\mathcal{X}$.
	Then,
	$\langle \varphi_h(\bx_1;\bx),\varphi_h(\bx_1;\bx)\rangle_{\mathcal{W}_\bx}\lesssim h^d$,
	and
	\begin{equation*}
	\langle\mathscr{K}_h(\bx,\bx_1),\mathscr{K}_h(\bx,\bx_1)\rangle_\mathcal{X}
	\geq\langle\mathscr{K}_h(\bx,\bx_1),\mathscr{K}_h(\bx,\bx_1)\rangle_{\mathcal{W}_\bx}
	\gtrsim h^{-d}\,\langle\mathscr{K}_h(\bx,\bx_1),\varphi_h(\bx_1;\bx)\rangle_{\mathcal{W}_\bx}^2.
	\end{equation*}
	It suffices to show $|\langle\mathscr{K}_h(\bx,\bx_1),\varphi_h(\bx_1;\bx)\rangle_{\mathcal{X}\setminus\mathcal{W}_\bx}|$ can be made sufficiently small such that $$|\langle\mathscr{K}_h(\bx,\bx_1),\varphi_h(\bx_1;\bx)\rangle_{\mathcal{W}_\bx}|\gtrsim h^{-[\bq]}.$$
	
	By Lemma \ref{lem: Q matrix}, the elements of $h^d\bQ_m^{-1}$ and $h^d\bQ_{\tilde{m}}^{-1}$ exponentially decays when they get far away from the (block) diagonals. In view of Assumption \ref{assumption: BC bases}, with $\bx$ fixed, $\mathscr{K}_h(\bx,\bx_1)$ also exponentially decays as $\|\bx_1-
	\bx\|$ increases. Formally, write
	\begin{equation}\label{equation: script K decomp}
	\begin{split}
	\mathscr{K}(\bx, \bx_1)=&\,\partial^\bq\bp(\bx)'\bQ_m^{-1}\bp(\bx_1)
	+ \sum_{\bm{u}\in\Lambda_{m}}\partial^{\bm{u}}\bpt(\bx)' h_\bx^{m-[\bq]}B_{\bm{u},\bq}(\bx)\bQ_{\tilde{m}}^{-1}\bpt(\bx_1)\\
	&-\partial^\bq\bp(\bx)' \bQ_m^{-1}\E\Big[\bp(\bx_i)h_{\bx_i}^{\bm{u}}B_{\bm{u},\mathbf{0}}(\bx_i)\partial^{\bm{u}}\bpt(\bx_i)'\Big]\bQ_{\tilde{m}}^{-1}\bpt(\bx_1).
	\end{split}
	\end{equation}
	
	We temporarily change the meaning of subscripts of subrectangles in $\Delta^\mathtt{rec}$: let 
	$$\delta^\mathtt{rec}_{l_1\ldots l_d}:=\left\{\mathbf{\check{x}}: t_{\ell,l_\ell} \leq \check{x}_\ell \leq t_{\ell,l_\ell+1},\, \ell=1,\ldots, d \right\}$$
	with
	$\delta^{\mathtt{rec}}_\mathbf{0}$ denoting the subrectangle where $\bx$ is located, and index other subrectangles with $\delta^{\mathtt{rec}}_\mathbf{0}$ regarded as the ``origin''. 
	First notice that for any given point $\bx_0\in\mathcal{X}$, $\bp(\bx_0)$ and $\bpt(\bx_0)$ are two vectors containing a fixed number of nonzeros and all their elements are bounded by some universal constant. Their nonzero elements are obtained by evaluating those basis functions with local supports covering $\bx_0$. Moreover, $h^{-d}\E[\bp(\bx_i)h_{\bx_i}^{m}B_{\bm{u},\mathbf{0}}(\bx_i)\partial^{\bm{u}}\bpt(\bx_i)']$ also admits a block banded structure in the sense that only the products of basis functions with overlapping supports are nonzero, and all elements in this matrix is bounded by some universal constant. Hence for
	\[\mathbf{r}_1(\bx)'=h^{[\bq]}\partial^\bq\bp(\bx)' \bQ_m^{-1}\E\Big[\sum_{\bm{u}\in\Lambda_{m}}\bp(\bx_i)B_{\bm{u},\mathbf{0}}(\bx_i)\partial^{\bm{u}}\bpt(\bx_i)'\Big]\]
	and
	\[\mathbf{r}_2(\bx_1)'=h^d\bpt(\bx_1)'\bQ_{\tilde{m}}^{-1},\]
	we can find another two vectors $\bar{\mathbf{r}}_1(\bx)$ and $\bar{\mathbf{r}}_2(\bx_1)$ with strictly positive elements that are greater than the absolute values of the corresponding elements in $\mathbf{r}_1(\bx)$ and $\mathbf{r}_2(\bx_1)$, i.e., $\bar{\mathbf{r}}_1(\bx)$ and $\bar{\mathbf{r}}_2(\bx_1)$ are element-wise bounds on $\mathbf{r}_1(\bx)$ and $\mathbf{r}_2(\bx_1)$. Moreover, the elements of $\bar{\mathbf{r}}_1(\bx)$ and $\bar{\mathbf{r}}_2(\bx_1)$ are decaying exponentially with some rates $\vartheta_1, \vartheta_2\in(0,1)$ when they get far away from the positions of those basis functions in $\bp(\cdot)$ and $\bpt(\cdot)$ whose supports are around $\bx$ and $\bx_1$ respectively. Notice that for some constant $\vartheta\in(0, 1)$, $\{\sum_{l=i}^{\infty}\vartheta^l\}_{i=1}^{\infty}$ is also a geometric sequence. Therefore, the inner product between $\bar{\mathbf{r}}_1(\bx)$ and $\bar{\mathbf{r}}_2(\bx_1)$, i.e., the third term in \eqref{equation: script K decomp}, exponentially decays as $\|\bx_1-\bx\|$ increases. Similarly, the inner product between
	$\partial^\bq\bp(\bx)$ and  $\bQ_m^{-1}\bp(\bx_1)$ and that between $\sum_{\bm{u}\in\Lambda_m}\partial^{\bm{u}}\bpt(\bx)h_\bx^{m-[\bq]}B_{\bm{u},\bq}(\bx)$ and  $\bQ_{\tilde{m}}^{-1}\bpt(\bx_1)$ have the same property.
	Given these results, we have for some $\vartheta\in(0,1)$,
	\[\|h^{[\bq]+d}\mathscr{K}_h(\bx,\cdot)\|_{L_\infty(\delta_{l_1\ldots l_d})}\leq C\vartheta^{\sum_{\ell=1}^{d}|l_\ell|}.\]
	Meanwhile,
	$\|\varphi_h(\cdot;\bx)\|_{L_\infty(\delta_{l_1\ldots l_d}^{\mathtt{rec}})}\lesssim (|l_1|+1)^{q_1}\cdots(|l_d|+1)^{q_d}$.
	Let $\bm{l}=(l_1,\cdots, l_d)$ and $\mathbf{1}=(1, \cdots, 1)$. Denote $|\bm{l}|=(|l_1|,\ldots, |l_d|)$. The above results imply
	\[|\langle\mathscr{K}_h(\bx,\bx_1), \varphi_h(\bx_1;\bx)\rangle_{\delta^{\mathtt{rec}}_{l_1\ldots l_d}}|\lesssim h^{-[\bq]}(|\bm{l}|+\mathbf{1})^\bq\vartheta^{[|\bm{l}|]}.
	\]
	Then the desired result follows from the fact that
	$\sum_{[|\bm{l}|]=0}^{\infty}(|\bm{l}|+\mathbf{1})^\bq
	\vartheta^{[|\bm{l}|]}$
	exists. Therefore, we can choose $L$ large enough which is independent of $\bx$ and $n$ such that
	$|\langle\mathscr{K}_h(\bx,\bx_1),\varphi_h(\bx_1;\bx)\rangle_{\mathcal{W}_\bx}|\gtrsim h^{-[\bq]}$. Then the proof is complete.
\end{proof}

\subsection{Proof of Theorem 4.1}
In this section, we provide the proof of Theorem 4.1 in the main paper.

\begin{proof}
Regarding the integrated conditional variance,
\begin{align*}
&\quad \int_\mathcal{X}\V[\widehat{\partial^\bq\mu}_0(\bx)|\bX]w(\bx)\,d\bx
=\frac{1}{n}\tr\Big[\bSig{0}\int_\mathcal{X}\bgamma_{\bq,0}(\bx)\bgamma_{\bq,0}(\bx)'w(\bx)\,d\bx\Big]+o_\P\Big(\frac{1}{nh^{d+2[\bq]}}\Big)\\
&\leq\frac{1}{n}\lambda_{\max}\Big(\bQ_m^{-1}\E[\bp(\bx_i)\bp(\bx_i)'\sigma^2(\bx_i)]\bQ_m^{-1}\Big)
\tr\Big[\int_\mathcal{X}\partial^\bq\bp(\bx)\partial^\bq\bp(\bx)'w(\bx) d\bx\Big]+o_\P\Big(\frac{1}{nh^{d+2[\bq]}}\Big)\\
&\lesssim \frac{1}{nh^d}\tr\Big[\int_\mathcal{X}\partial^\bq\bp(\bx)\partial^\bq\bp(\bx)'w(\bx) d\bx\Big]
\lesssim \frac{1}{nh^{d+2[\bq]}},
\end{align*}	
where the first line holds by Lemma \ref{lem: Q matrix}, the second by Trace Inequality, and the last by the continuity of $w(\cdot)$ and Lemma \ref{lem: Q matrix}. Since $\sigma^2(\cdot)$ and $w(\cdot)$ are bounded away from zero, the other side of the bound follows similarly.

Regarding the integrated squared bias, we have
\begin{align}
&\int_{\mathcal{X}}\left(\E[\widehat{\partial^\bq\mu}_0(\bx)|\bX]-\partial^\bq\mu(\bx)\right)^2w(\bx)d\bx \nonumber \\
=&\int_{\mathcal{X}}\Big(\mathscr{B}_{m,\bq}(\bx)-
\bgamma_{\bq,0}(\bx)'\E[\bp(\bx_i)
\mathscr{B}_{m,\mathbf{0}}(\bx_i)]\Big)^2w(\bx)d\bx
+o_\P(h^{2m-2[\bq]}) \nonumber\\
=&\int_{\mathcal{X}}\mathscr{B}_{m,\bq}(\bx)^2w(\bx)d\bx+
\int_{\mathcal{X}}\Big(\bgamma_{\bq, 0}(\bx)'\E[\bp(\bx_i)\mathscr{B}_{m,\mathbf{0}}(\bx_i)]\Big)^2w(\bx)d\bx \nonumber\\
&-2\int_{\mathcal{X}}\mathscr{B}_{m,\bq}(\bx)
\bgamma_{\bq, 0}(\bx)'
\E[\bp(\bx_i)\mathscr{B}_{m,\mathbf{0}}(\bx_i)]w(\bx)d\bx+o_\P(h^{2m-2[\bq]}) \nonumber\\
=&: \mathrm{B}_1+\mathrm{B}_2-2\mathrm{B}_3+o_\P(h^{2m-2[\bq]}). \label{eq:3B}
\end{align}

Let $h_\delta$ be the diameter of $\delta$ and $\bt_\delta^*$ be an arbitrary point in $\delta$. Then
\begin{align*}
\mathrm{B}_1=&\sum_{\bm{u}_1,\bm{u}_2\in\Lambda_{m}}
\int_\mathcal{X}\Big[\partial^{\bm{u}_1}\mu(\bx)\partial^{\bm{u}_2}\mu(\bx)h_\bx^{2m-2[\bq]}B_{\bm{u}_1,\bq}(\bx)B_{\bm{u}_2,\bq}(\bx)\Big]w(\bx)d\bx\\
=&\sum_{\bm{u}_1,\bm{u}_2\in\Lambda_{m}}\,\sum_{\delta\in\Delta}\int_{\delta}\Big[h_{\delta}^{2m-2[\bq]}\partial^{\bm{u}_1}\mu(\bt_{\delta}^*)\partial^{\bm{u}_2}\mu(\bt_{\delta}^*)B_{\bm{u}_1,\bq}(\bx)
B_{\bm{u}_2,\bq}(\bx)\Big]w(\bt_{\delta}^*)d\bx+o(h^{2m-2[\bq]})\\
=&\sum_{\bm{u}_1,\bm{u}_2\in\Lambda_{m}}\,\sum_{\delta\in\Delta}
\Big\{h_{\delta}^{2m-2[\bq]}\partial^{\bm{u}_1}\mu(\bt_{\delta}^*)\partial^{\bm{u}_2}\mu(\bt_{\delta}^*)w(\bt_{\delta}^*)
\int_{\delta}B_{\bm{u}_1,\bq}(\bx)B_{\bm{u}_2,\bq}(\bx)d\bx\Big\}+o(h^{2m-2[\bq]})\\
\lesssim&\,h^{2m-2[\bq]}
\end{align*}
where the third line holds by the continuity of $\partial^{\bm{u}_1}\mu(\cdot)$, $\partial^{\bm{u}_2}\mu(\cdot)$, and $w(\cdot)$, and the last by Assumption \ref{assumption: quasi-uniform partition} and \ref{assumption: Misspecification Error}.
\begin{align*}
\mathrm{B}_2&=\tr\Big[ \bQ_m^{-1}\E[\bp(\bx_i)\mathscr{B}_{m,\mathbf{0}}(\bx_i)]
\E[\bp(\bx_i)'\mathscr{B}_{m,\mathbf{0}}(\bx_i)]\bQ_m^{-1}
\int_{\mathcal{X}}\partial^\bq\bp(\bx)\partial^\bq\bp(\bx)'w(\bx)d\bx\Big]\\
&\lesssim h^{-d-2[\bq]}\tr\Big[\E[\bp(\bx_i)\mathscr{B}_{m,\mathbf{0}}(\bx_i)]\E[\bp(\bx_i)'\mathscr{B}_{m,\mathbf{0}}(\bx_i)]\Big]
\lesssim h^{2m-2[\bq]}
\end{align*}
where the first inequality holds by Trace Inequality, Lemma \ref{lem: Q matrix}, and Assumption \ref{assumption: Local Basis}, and the second by Assumption \ref{assumption: Local Basis} and \ref{assumption: Misspecification Error}. Finally,
\[
|\mathrm{B}_3|
\leq\Big\|\int_{\mathcal{X}}\mathscr{B}_{m,\bq}(\bx)
\partial^\bq\bp(\bx)'w(\bx)d\bx\Big\|_\infty
\Big\|\bQ_m^{-1}\Big\|_\infty
\Big\|\E[\bp(\bx_i)\mathscr{B}_{m,\mathbf{0}}(\bx_i)]\Big\|_\infty
\lesssim h^{2m-2[\bq]}
\]
where the last inequality follows from Assumption \ref{assumption: Local Basis} and Lemma \ref{lem: Q matrix}.	
\end{proof}

\subsection{Proof of Theorem \ref{thm: IMSE-R}}
\begin{proof}
	We divide the proof into two steps.
	
	\textit{Step 1:} For the integrated variance,
	first define operators $\pmb{\mathscr{M}}(\cdot)$ and $\pmb{\mathscr{M}}_\bq(\cdot)$ as follows:
	\begin{equation*}
	\pmb{\mathscr{M}}(\phi):=\int_\mathcal{X}\bp(\bx)\bp(\bx)^\prime\phi(\bx)d\bx,\quad \text{and}\quad
	\pmb{\mathscr{M}}_\bq(\phi):=\int_\mathcal{X}\partial^\bq\bp(\bx)\partial^\bq\bp(\bx)^\prime\phi(\bx)d\bx.
	\end{equation*}
	Then,
	\begin{align*}
	&\int_\mathcal{X}\V[\widehat{\partial^\bq\mu}_0(\bx)|\bX]w(\bx)d\bx\\
	=&\frac{1}{n}\tr\Big[\bQ_m^{-1}\E[\bp(\bx_i)\bp(\bx_i)'\sigma(\bx_i)^2]\bQ_m^{-1}\int_\mathcal{X}\partial^\bq\bp(\bx)\partial^\bq\bp(\bx)'w(\bx)\,d\bx\Big]+o_\P\Big(\frac{1}{nh^{d+2[\bq]}}\Big)\\
	=&\frac{1}{n}\tr\Big[\pmb{\mathscr{M}}(f)^{-1}\pmb{\mathscr{M}}(\sigma^2f)\pmb{\mathscr{M}}(f)^{-1}\pmb{\mathscr{M}}_\bq(w)\Big]+o_\P\Big(\frac{1}{nh^{d+2[\bq]}}\Big).
	\end{align*}
	Moreover, define another operator $\pmb{\mathscr{D}}(\cdot)$:
	$\pmb{\mathscr{D}}(\phi):=\diag\{\phi(\tau_1),\phi(\tau_2),\cdots,\phi(\tau_K)\}$.
	Recall that $\btau_k$ is an arbitrary point in $\supp(p_k)$, for $k=1,\ldots, K$. Then we can write
	\begin{equation}\label{eq: imse matrix approx}
	\pmb{\mathscr{M}}(\phi)=\pmb{\mathscr{M}}(1)\pmb{\mathscr{D}}(\phi)-\pmb{\mathscr{E}}(\phi)
	\end{equation}
	where $\pmb{\mathscr{E}}(\phi)$ can be viewed as errors defined by Eq. \eqref{eq: imse matrix approx}. Similarly, write
	$$\pmb{\mathscr{M}}_\bq(\phi)=\pmb{\mathscr{M}}_\bq(1)\pmb{\mathscr{D}}(\phi)-\pmb{\mathscr{E}}_\bq(\phi).$$
	
	Then it directly follows that
	\begin{align*}
	&\pmb{\mathscr{M}}(f)^{-1}\pmb{\mathscr{M}}(\sigma^2f)=[\mathbf{I}-\pmb{\mathscr{U}}(f)]^{-1}[\pmb{\mathscr{D}}(f)^{-1}\pmb{\mathscr{D}}(\sigma^2f)-\pmb{\mathscr{L}}(f,\sigma^2f)],\quad\text{and}\\
	&\pmb{\mathscr{M}}(f)^{-1}\pmb{\mathscr{M}}_\bq(w)=[\mathbf{I}-\pmb{\mathscr{U}}(f)]^{-1}[\pmb{\mathscr{D}}(f)^{-1}\pmb{\mathscr{M}}(1)^{-1}\pmb{\mathscr{M}}_\bq(1)\pmb{\mathscr{D}}(w)-\pmb{\mathscr{L}}_\bq(f,w)]
	\end{align*}
	where
	\begin{align*}
	&\pmb{\mathscr{U}}(\phi):=\pmb{\mathscr{D}}(\phi)^{-1}\pmb{\mathscr{M}}(1)^{-1}\pmb{\mathscr{E}}(\phi),\\ &\pmb{\mathscr{L}}(\phi,\varphi):=\pmb{\mathscr{D}}(\phi)^{-1}\pmb{\mathscr{M}}(1)^{-1}\pmb{\mathscr{E}}(\varphi),\quad\text{and}\\
	&\pmb{\mathscr{L}}_\bq(\phi,\varphi):=\pmb{\mathscr{D}}(\phi)^{-1}\pmb{\mathscr{M}}(1)^{-1}\pmb{\mathscr{E}}_\bq(\varphi).
	\end{align*}
	The number of nonzeros in any row or any column of $\pmb{\mathscr{E}}(\phi)$ (and $\pmb{\mathscr{E}}_\bq(\phi)$) is bounded by some constant. In fact, as explained in the proof of Lemma \ref{lem: Q matrix}, it takes a (multi-layer) banded structure. If $\supp(p_k)\cap\supp(p_l)\neq\varnothing$, then by Assumption \ref{assumption: Local Basis} and the continuity of $f$, the $(k,l)$th element of $\pmb{\mathscr{M}}(f)$ can be approximated as follows:
	\begin{equation}\label{eq: imse Q element approx}
	\begin{split}
	  \int_\mathcal{X} p_k(\bx)p_l(\bx)f(\bx)d\bx
	=&f(\tau_k)\int_\mathcal{X} p_k(\bx)p_l(\bx)d\bx+o(h^d)\\
	=&f(\tau_l)\int_\mathcal{X} p_k(\bx)p_l(\bx)d\bx+o(h^d).
	\end{split}
	\end{equation}
	
	Moreover, since $\mathcal{X}$ is compact, $f$ is uniformly continuous, and then $\|\pmb{\mathscr{E}}(f)\|_1=o(h^d)$, $\|\pmb{\mathscr{E}}(f)\|_\infty=o(h^d)$, and $\|\pmb{\mathscr{E}}(f)\|=o(h^d)$.
	Since $\|\pmb{\mathscr{D}}(f)^{-1}\|\lesssim 1$ and $\|\pmb{\mathscr{M}}(1)^{-1}\|\lesssim h^{-d}$, we conclude $\|\pmb{\mathscr{U}}(f)\|=o(1)$. For $K$ large enough, we can make $\|\pmb{\mathscr{U}}(f)\|<1$, and thus
	$$[\mathbf{I}-\pmb{\mathscr{U}}(f)]^{-1}=\mathbf{I}+\pmb{\mathscr{U}}(f)+\pmb{\mathscr{U}}(f)^2+\cdots=\mathbf{I}+\pmb{\mathscr{W}}(f)$$
	where $\pmb{\mathscr{W}}(f):=\sum_{l=1}^{\infty}\pmb{\mathscr{U}}(f)^l$.
	Now we can write
	\begin{align*}
	&\tr\Big[\pmb{\mathscr{M}}(f)^{-1}\pmb{\mathscr{M}}(\sigma^2f)\pmb{\mathscr{M}}(f)^{-1}\pmb{\mathscr{M}}_\bq(w)\Big]\\
	=&\tr\Big[\Big(\mathbf{I}+\pmb{\mathscr{W}}(f)\Big)
	\Big(\pmb{\mathscr{D}}(f)^{-1}\pmb{\mathscr{D}}(\sigma^2f)-\pmb{\mathscr{L}}(f,\sigma^2f)\Big)
	\Big(\mathbf{I}+\pmb{\mathscr{W}}(f)\Big)\\
	&\qquad\;\times\Big(\pmb{\mathscr{D}}(f)^{-1}\pmb{\mathscr{M}}(1)^{-1}\pmb{\mathscr{M}}_\bq(1)\pmb{\mathscr{D}}(w)-\pmb{\mathscr{L}}_\bq(f,w)\Big)\Big]\\
	=&\tr\Big[\Big(\pmb{\mathscr{D}}(f)^{-1}\pmb{\mathscr{D}}(\sigma^2f)+\mathbf{E}_1\Big)
	\Big(\pmb{\mathscr{D}}(f)^{-1}\pmb{\mathscr{M}}(1)^{-1}\pmb{\mathscr{M}}_\bq(1)\pmb{\mathscr{D}}(w)+\mathbf{E}_2\Big)\Big]\\
	=&\tr\Big[\pmb{\mathscr{D}}(f)^{-1}\pmb{\mathscr{D}}(\sigma^2f)\pmb{\mathscr{D}}(f)^{-1}\pmb{\mathscr{M}}(1)^{-1}\pmb{\mathscr{M}}_\bq(1)\pmb{\mathscr{D}}(w)+\\
	&\hspace{3em}\mathbf{E}_1\pmb{\mathscr{D}}(f)^{-1}\pmb{\mathscr{M}}(1)^{-1}\pmb{\mathscr{M}}_\bq(1)\pmb{\mathscr{D}}(w)+\pmb{\mathscr{D}}(f)^{-1}\pmb{\mathscr{D}}(\sigma^2f)\mathbf{E}_2+\mathbf{E}_1\mathbf{E}_2\Big]
	\end{align*}
	where
	\begin{align*}
	&\mathbf{E}_1=-\pmb{\mathscr{L}}(f,\sigma^2f)+\pmb{\mathscr{W}}(f)\pmb{\mathscr{D}}(f)^{-1}\pmb{\mathscr{D}}(\sigma^2f)-\pmb{\mathscr{W}}(f)\pmb{\mathscr{L}}(f,\sigma^2f),\\
	&\mathbf{E}_2=-\pmb{\mathscr{L}}_\bq(f,w)+\pmb{\mathscr{W}}(f)\pmb{\mathscr{D}}(f)^{-1}\pmb{\mathscr{M}}(1)^{-1}\pmb{\mathscr{M}}_\bq(1)\pmb{\mathscr{D}}(w)-\pmb{\mathscr{W}}(f)\pmb{\mathscr{L}}_\bq(f,w).
	\end{align*}
	By assumptions in the theorem, $\vol(\delta_\bx)=\prod_{\ell=1}^d b_{\bx,\ell}=\prod_{\ell=1}^d\kappa_\ell^{-1}g_\ell(\bx)^{-1}+o(h^d)$. Hence
	\begin{align*}
	&\quad\tr\Big[\pmb{\mathscr{D}}(f)^{-1}\pmb{\mathscr{D}}(\sigma^2f)\pmb{\mathscr{D}}(f)^{-1}\pmb{\mathscr{M}}(1)^{-1}\pmb{\mathscr{M}}_\bq(1)\pmb{\mathscr{D}}(w)\Big]\\	&=\prod_{\ell=1}^d\kappa_\ell\tr\Big[\pmb{\mathscr{D}}(f)^{-1}
	\pmb{\mathscr{D}}(\sigma^2f)\pmb{\mathscr{D}}(f)^{-1}
	\pmb{\mathscr{M}}(1)^{-1}\pmb{\mathscr{M}}_\bq(1)
	\pmb{\mathscr{D}}(w)\pmb{\mathscr{D}}\Big(\prod_{\ell=1}^dg_\ell\Big)\pmb{\mathscr{D}}\Big(\prod_{\ell=1}^dg_\ell\Big)^{-1}
	\prod_{\ell=1}^d\kappa_\ell^{-1}\Big]\\
	&=\,\boldsymbol{\kappa}^{\mathbf{1}+2\bq}
	\left(\sum_{k=1}^K\Big[\frac{\sigma^2(\btau_k)w(\btau_k)}
	{f(\btau_k)}\prod_{\ell=1}^dg_\ell(\btau_k)\Big]
	\be_k'\pmb{\mathscr{M}}(1)^{-1}\boldsymbol{\kappa}^{-2\bq}
	\pmb{\mathscr{M}}_\bq(1)\be_k\vol(\delta_{\btau_k})\right)
	+o(\boldsymbol{\kappa}^{\mathbf{1}+2\bq})
	\end{align*}
	By Assumptions \ref{assumption: quasi-uniform partition}, \ref{assumption: Local Basis} and Lemma \ref{lem: Q matrix}, the summation in parenthesis is bounded from above and below.
	It remains to show all other terms are of smaller order. It directly follows from the same argument as that in the proof of \cite[Theorem 6.1]{Agarwal-Studden_1980_AoS} that the trace of the remaining terms is $o(\boldsymbol{\kappa}^{\mathbf{1}+2\bq})$.
	
	\textit{Step 2:} For the integrated squared bias, consider the  three leading terms $\mathrm{B}_1$, $\mathrm{B}_2$ and $\mathrm{B}_3$ defined in Equation \eqref{eq:3B}.
	For $\mathrm{B}_1$, by assumption in the theorem
	\begin{equation*}
	\mathscr{B}_{m,\bq}(\bx)=-\sum_{\bm{u}\in\Lambda_m}
	\partial^{\bm{u}}\mu(\bx)\Big(\prod_{\ell=1}^{d}
	\kappa_\ell^{-u_\ell+q_\ell}g_\ell(\bx)^{-u_\ell+q_\ell}\Big)\bb_{\bx}^{-\bm{u}+\bq}h_\bx^{m-[\bq]}B_{m,\bq}(\bx)+o(h^{m-[\bq]}).
	\end{equation*}
	Recall that $\boldsymbol{\kappa}=(\kappa_1,\ldots, \kappa_d)$. Define $\mathbf{g}(\bx):=(g_1(\bx), \ldots,g_d(\bx))$. Using the above fact and the same notation as in the proof of Theorem 4.1, we have
	\begin{align*}
	\mathrm{B}_1&=\sum_{\bm{u}_1,\bm{u}_2\in\Lambda_m}\int_{\mathcal{X}}\bigg[\partial^{\bm{u}_1}\mu(\bx)\partial^{\bm{u}_2}\mu(\bx)\frac{h_\bx^{2m-2[\bq]}B_{\bm{u}_1,\bq}(\bx)B_{\bm{u}_2,\bq}(\bx)}{\boldsymbol{\kappa}^{\bm{u}_1+\bm{u}_2-2\bq}\mathbf{g}(\bx)^{\bm{u}_1+\bm{u}_2-2\bq}\bb_\bx^{\bm{u}_1+\bm{u}_2-2\bq}}\bigg]w(\bx)d\bx+o(h^{2m-2[\bq]})\\
	&=\sum_{\bm{u}_1,\bm{u}_2\in\Lambda_m}\boldsymbol{\kappa}^{-(\bm{u}_1+\bm{u}_2-2\bq)}\bigg(\sum_{\delta\in\Delta}
	\bigg[\frac{\partial^{\bm{u}_1}g(t_\delta^*)\partial^{\bm{u}_2}g(t_\delta^*)w(t_\delta^*)}{\mathbf{g}(t_\delta^*)^{\bm{u}_1+\bm{u}_2-2\bq}}\bigg]\\
	&\hspace{12em}\times\bigg[\frac{h_\bx^{2m-2[\bq]}}{\bb_{\bx}^{\bm{u}_1+\bm{u}_2-2\bq}}\int_\delta B_{\bm{u}_1,\bq}(\bx)B_{\bm{u}_2,\bq}(\bx)d\bx\Big]\bigg)
	+o(h^{2m-2[\bq]})\\
	&=\sum_{\bm{u}_1,\bm{u}_2,\in\Lambda_m}\boldsymbol{\kappa}^{-(\bm{u}_1+\bm{u}_2-2\bq)}\eta_{\bm{u}_1,\bm{u}_2,\bq}\int_{\mathcal{X}}\frac{\partial^{\bm{u}_1}\mu(\bx)\partial^{\bm{u}_2}\mu(\bx)w(\bx)}{\mathbf{g}(\bx)^{\bm{u}_1+\bm{u}_2-2\bq}}d\bx+o(h^{2m-2[\bq]})
	\end{align*}
	where the last line holds by the integrability of $\partial^{\bm{u}_1}\mu(\bx)\partial^{\bm{u}_2}\mu(\bx)w(\bx)/\mathbf{g}(\bx)^{\bm{u}_1+\bm{u}_2-2\bq}$ over $\mathcal{X}$.
	
	For $\mathrm{B}_2$, first notice that
	\begin{equation*}
	\Big\|\E[\bp(\bx_i)\mathscr{B}_{m,\mathbf{0}}(\bx_i)]
	+\sum_{\bm{u}\in\Lambda_m}\boldsymbol{\kappa}^{-\bm{u}}\int_\mathcal{X}\Big(\bp(\bx)\partial^{\bm{u}}\mu(\bx)\mathbf{g}(\bx)^{-\bm{u}}\bb_\bx^{-\bm{u}}h_\bx^mB_{m,\mathbf{0}}(\bx)f(\bx)\Big)d\bx\Big\|_\infty=o(h^{m+d}),
	\end{equation*}
	implying that the errors given rise to by approximating $\bb_{\bx}$ are of smaller order. The integral in this approximation is a vector with typical elements given by
	\begin{align*}
	&\int_\mathcal{X}\Big(p_k(\bx)\partial^{\bm{u}}\mu(\bx)\mathbf{g}(\bx)^{-\bm{u}}\bb_\bx^{-\bm{u}}h_\bx^mB_{m,\mathbf{0}}(\bx)f(\bx)\Big)d\bx\\
	=&\,\partial^{\bm{u}}\mu(\btau_k)\mathbf{g}(\btau_k)^{-\bm{u}}f(\btau_k)\int_{\mathcal{X}}p_k(\bx)\bb_{\bx}^{-\bm{u}}h_\bx^mB_{m,\mathbf{0}}(\bx)d\bx+o(h^d).
	\end{align*}
	Recall that $\mathbf{v}_{\bm{u},\bq}$ defined in the theorem has the $k$th element equal to
	$$\frac{\partial^{\bm{u}}\mu(\btau_k)\sqrt{w(\btau_k)}}{\bkappa^\bq\mathbf{g}(\btau_k)^{\bm{u}-\bq}}\int_{\mathcal{X}}\partial^\bq p_k(\bx)\bb_{\bx}^{-(\bm{u}-\bq)}h_\bx^{m-[\bq]}B_{\bm{u},\bq}(\bx)d\bx.$$
	Given these results and Lemmas \ref{lem: Q matrix}, it follows that
	\begin{align*}
	\mathrm{B}_2&=\tr\Big[ \bQ_m^{-1}\E[\bp(\bx_i)\mathscr{B}_{m,\mathbf{0}}(\bx_i)]\E[\bp(\bx_i)^\prime\mathscr{B}_{m,\mathbf{0}}(\bx_i)]\bQ_m^{-1}\int_{\mathcal{X}}\partial^\bq\bp(\bx)\partial^\bq\bp(\bx)'w(\bx)d\bx\Big]\\
	&=\sum_{\bm{u}_1,\bm{u}_2,\in\Lambda_m}
	\boldsymbol{\kappa}^{-(\bm{u}_1+\bm{u}_2-2\bq)}
	\tr\Big(\bQ_m^{-1}\pmb{\mathscr{D}}(f/\sqrt{w})\mathbf{v}_{\bm{u}_2,\mathbf{0}}\mathbf{v}_{\bm{u}_1,\mathbf{0}}^\prime\pmb{\mathscr{D}}(f/\sqrt{w}) \bQ_m^{-1}\boldsymbol{\kappa}^{-2\bq}
	\pmb{\mathscr{M}}_\bq(w)\Big)\\
	&\hspace{30em}+o(h^{2m-2[\bq]})\\
	&=\sum_{\bm{u}_1,\bm{u}_2,\in\Lambda_m}
	\boldsymbol{\kappa}^{-(\bm{u}_1+\bm{u}_2-2\bq)}
	\Big(\mathbf{v}_{\bm{u}_1,\mathbf{0}}'\pmb{\mathscr{D}}(f/\sqrt{w})(\mathbf{I}+\pmb{\mathscr{W}})\pmb{\mathscr{D}}(f)^{-1}
	\pmb{\mathscr{M}}(1)^{-1}\boldsymbol{\kappa}^{-2\bq}\pmb{\mathscr{M}}_\bq(w)\\
	&\hspace{11em}\times\pmb{\mathscr{M}}(1)^{-1}\pmb{\mathscr{D}}(f)^{-1}(\mathbf{I}+\pmb{\mathscr{W}})'\pmb{\mathscr{D}}(f/\sqrt{w})\mathbf{v}_{\bm{u}_2,\mathbf{0}}\Big)+o(h^{2m-2[\bq]})\\
	&=\sum_{\bm{u}_1,\bm{u}_2,\in\Lambda_m}
	\boldsymbol{\kappa}^{-(\bm{u}_1+\bm{u}_2-2\bq)}
	\Big(\mathbf{v}_{\bm{u}_1,\mathbf{0}}'\pmb{\mathscr{D}}(1/\sqrt{w})\pmb{\mathscr{M}}(1)^{-1}\boldsymbol{\kappa}^{-2\bq}\pmb{\mathscr{M}}_\bq(1)
	\pmb{\mathscr{D}}(w)\\
	&\hspace{11em}\times\pmb{\mathscr{M}}(1)^{-1}\pmb{\mathscr{D}}(1/\sqrt{w})\mathbf{v}_{\bm{u}_2,\mathbf{0}}\Big)
	+o(h^{2m-2[\bq]})
	\end{align*}
	
	It should be noted that \eqref{eq: imse Q element approx} implies that the approximation given by \eqref{eq: imse matrix approx} is still valid if $\pmb{\mathscr{M}}(1)$ is pre-multiplied by $\pmb{\mathscr{D}}(\phi)$ instead of being post-multiplied if $\phi(\cdot)$ is continuous. Therefore, we can repeat the argument in Step 1 and further write the term in parenthesis in the last line as:
	\begin{align*}
	&\mathbf{v}_{\bm{u}_1,\mathbf{0}}'\pmb{\mathscr{D}}(1/\sqrt{w})\pmb{\mathscr{M}}(1)^{-1}\boldsymbol{\kappa}^{-2\bq}\pmb{\mathscr{M}}_\bq(1)\pmb{\mathscr{D}}(w)\pmb{\mathscr{M}}(1)^{-1}\pmb{\mathscr{D}}(1/\sqrt{w})\mathbf{v}_{\bm{u}_2,\mathbf{0}}\\
	=\,&\mathbf{v}_{\bm{u}_1,\mathbf{0}}'\pmb{\mathscr{D}}(1/\sqrt{w})\Big(\pmb{\mathscr{D}}(1/\sqrt{w})\pmb{\mathscr{M}}(1)\Big)^{-1}\boldsymbol{\kappa}^{-2\bq}\pmb{\mathscr{M}}_\bq(1)\Big(\pmb{\mathscr{M}}(1)\pmb{\mathscr{D}}(1/\sqrt{w})\Big)^{-1}\pmb{\mathscr{D}}(1/\sqrt{w})\mathbf{v}_{\bm{u}_2,\mathbf{0}}+o(1)\\
	=\,&\mathbf{v}_{\bm{u}_1,\mathbf{0}}'\bH_\mathbf{0}^{-1}\bH_\bq\bH_\mathbf{0}^{-1}\mathbf{v}_{\bm{u}_2,\mathbf{0}}+o(1)
	\end{align*}
	where $|\mathbf{v}_{\bm{u}_1,\mathbf{0}}'\bH_\mathbf{0}^{-1}\bH_\bq\bH_\mathbf{0}^{-1}\mathbf{v}_{\bm{u}_2,\mathbf{0}}|\lesssim 1$ by Assumption \ref{assumption: quasi-uniform partition}, \ref{assumption: Local Basis} and Lemma \ref{lem: Q matrix}.
	
	Finally, for $\mathrm{B}_3$, notice that
	\begin{align*}
	&\Big\|\int_\mathcal{X}\mathscr{B}_{m,\bq}(\bx)\partial^\bq\bp(\bx)'w(\bx)d\bx
	+\sum_{\bm{u}\in\Lambda_m}\int_{\mathcal{X}}\left(\frac{\partial^{\bm{u}}\mu(\bx)w(\bx)}{\boldsymbol{\kappa}^{\bm{u}-\bq}\mathbf{g}(\bx)^{\bm{u}-\bq}}\cdot\frac{h_\bx^{m-[\bq]}\partial^\bq\bp(\bx)^\prime B_{\bm{u},\bq}(\bx)}{\bb_{\bx}^{\bm{u}-\bq}}\right)d\bx\Big\|_\infty\\
	=&\,o(h^{m+d-2[\bq]}).
	\end{align*}
	Thus repeating the argument for $\mathrm{B}_2$, we have
	\begin{align*}
	\mathrm{B}_3&=\left(\int_{\mathcal{X}}\mathscr{B}_{m,\bq}(\bx)\partial^\bq\bp(\bx)'w(\bx)d\bx\right)\bQ_m^{-1}\E[\bp(\bx_i)\mathscr{B}_{m,\mathbf{0}}(\bx_i)]\\
	&=\sum_{\bm{u}_1,\bm{u}_2\in\Lambda_m}\bkappa^{-(\bm{u}_1+\bm{u}_2-2\bq)}\mathbf{v}_{\bm{u}_1,\bq}'\pmb{\mathscr{D}}(\sqrt{w})\bH_\mathbf{0}^{-1}\pmb{\mathscr{D}}(f)^{-1}\pmb{\mathscr{D}}(f/\sqrt{w})\mathbf{v}_{\bm{u}_2,\mathbf{0}}+o(h^{2m-2[\bq]})\\
	&=\sum_{\bm{u}_1,\bm{u}_2\in\Lambda_m}\boldsymbol{\kappa}^{-(\bm{u}_1+\bm{u}_2-2\bq)}\mathbf{v}_{\bm{u}_1,\bq}'\bH_\mathbf{0}^{-1}\mathbf{v}_{\bm{u}_2,\mathbf{0}}+o(h^{2m-2[\bq]})
	\end{align*}
	where $|\mathbf{v}_{\bm{u}_1,\bq}'\bH_\mathbf{0}^{-1}\mathbf{v}_{\bm{u}_2,\mathbf{0}}|\lesssim 1$ by Assumption \ref{assumption: quasi-uniform partition}, \ref{assumption: Local Basis} and Lemma \ref{lem: Q matrix}.
	Then the proof is complete.
\end{proof}

\subsection{Proof of Theorem 4.2}
\begin{proof}
	Continue the calculation in the proof of Theorem \ref{thm: IMSE-R}.
	For the integrated variance, when $\bq=\mathbf{0}$ and $\bp$ generates $J$ complete covers, $\pmb{\mathscr{M}}(1)^{-1}\pmb{\mathscr{M}}_\bq(1)=\mathbf{I}_K$ and hence
	$$\tr\Big[\pmb{\mathscr{D}}(f)^{-1}\pmb{\mathscr{D}}(\sigma^2f)
	\pmb{\mathscr{D}}(f)^{-1}\pmb{\mathscr{D}}(w)\Big]=\prod_{\ell=1}^d\kappa_\ell\times J\int_\mathcal{X}\frac{\sigma^2(\bx)w(\bx)}{f(\bx)}
	\prod_{\ell=1}^{d}g_\ell(\bx)\,d\bx+o(\prod_{\ell=1}^d\kappa_\ell).$$
	For the integrated squared bias, since the approximate orthogonality condition holds, both $\mathrm{B}_2$ and $\mathrm{B}_3$ are of smaller order, and the leading term in the integrated squared bias reduces to $\mathrm{B}_1$ only.
\end{proof}

\subsection{Proof of Lemma \ref{lem: pointwise linearization}}
\begin{proof}
	For $j=0,1$, the results directly follow from Assumption \ref{assumption: BC bases}, Lemma \ref{lem: bias}, and \cite[Lemma 4.1]{Belloni-Chernozhukov-Chetverikov-Kato_2015_JoE}.	
	For $j=2,3$,
	conditional on $\bX$, $R_{1n,\bq}(\bx)$ has mean zero, and its variance can be bounded as follows:
	\begin{align*}
	\V[R_{1n,\bq}(\bx)|\bX]&\lesssim \frac{1}{n}\left[\bg{j}'-\bgamma_{\bq,1}(\bx)'\right]
	\E[\bPi{j}(\bx_i)\bPi{j}(\bx_i)'] \left[\bg{j}-\bgamma_{\bq,j}(\bx)\right]\\
	&\lesssim \frac{1}{n}\|\bg{j}'-\bgamma_{\bq,j}(\bx)'\|^2
	\|\E[\bPi{j}(\bx_i)\bPi{j}(\bx_i)']\|\\
	&\lesssim_\P \frac{1}{n}\,h^{-2[\bq]-2d}\frac{\log n}{nh^d}\,h^d
	= \frac{\log n}{n^2h^{2d+2[\bq]}},
	\end{align*}
	where the third line holds by Lemma \ref{lem: basis} and the fact that $\E[\bPi{j}(\bx_i)\bPi{j}(\bx_i)']\lesssim h^d$ shown in the proof of Lemma \ref{lem: Avar}. Then by Chebyshev's inequality we conclude that
	$R_{1n,\bq}(\bx)\lesssim_\P \frac{\sqrt{\log n}}{nh^{d+[\bq]}}$.
	
	For the conditional bias $R_{2n, \bq}(\bx)$, we analyze the least squares bias correction and plug-in bias correction separately.
	For $j=2$, by construction,
	\begin{align*}
	\E[\widehat{\partial^\bq \mu}_2(\bx)|\bX]-\partial^\bq \mu(\bx)
	=\,&\left(\E[\widehat{\partial^\bq\mu}_1(\bx)|\bX]-\partial^\bq \mu(\bx)\right)-
	\partial^\bq\bp(\bx)'\widehat{\bQ}_m^{-1}
	\E_n[\bp(\bx_i)\mathfrak{B}_{\tilde{m},\mathbf{0}}(\bx_i)]\\
	=\,& O_\P(h^{m+\varrho-[\bq]})-\partial^\bq\bp(\bx)' \widehat{\bQ}_m^{-1}\E_n[\bp(\bx_i)\mathfrak{B}_{\tilde{m},\mathbf{0}}(\bx_i)]
	\end{align*}
	where $\mathfrak{B}_{\tilde{m},\mathbf{0}}(\bx_i)=
	\E[\widehat{\mu}_1(\bx_i)|\bX]-\mu(\bx_i)$ is the conditional bias of $\widehat{\mu}_1(\bx_i)$ and the last line follows from Lemma \ref{lem: bias}.
	
	Since
	$\|\E_n[\bp(\bx_i)\mathfrak{B}_{\tilde{m},\mathbf{0}}(\bx_i)]\|_\infty
	\leq\sup_{\bx\in\mathcal{X}}|\mathfrak{B}_{\tilde{m},\mathbf{0}}(\bx)|\,\|\E_n[|\bp(\bx_i)|]\|_\infty$,
	using the same proof strategy as that for Lemma \ref{lem: bias}, we have $\|\E_n[|\bp(\bx_i)|]\|_\infty\lesssim_\P h^d$. Also, by Lemma \ref{lem: bias}, $\sup_{\bx\in\mathcal{X}}|\mathfrak{B}_{\tilde{m},\mathbf{0}}(\bx)|\lesssim_\P h^{m+\varrho}$. Then by Lemma \ref{lem: basis}, the conditional bias of $\widehat{\partial^\bq \mu}_2(\bx)$ is $O_\P(h^{m+\varrho-[\bq]})$.
	
	Next, for $j=3$, using Lemma \ref{lem: bias}, we have
	\begin{align*}
	&\E[\widehat{\partial^\bq \mu}_3(\bx)|\bX]-\partial^\bq \mu(\bx)\\
	=\,&\E\Big[\widehat{\partial^\bq \mu}_0(\bx)+
	\sum_{\bm{u}\in\Lambda_m}\widehat{\partial^{\bm{u}}\mu}_1(\bx)h_\bx^{m-[\bq]}B_{\bm{u},\bq }(\bx)\\
	&\hspace{1.3in}+\partial^\bq\bp(\bx)'\widehat{\bQ}_m^{-1}\E_n[\bp(\bx_i)\widehat{\mathscr{B}}_{m,\mathbf{0}}(\bx_i)\Big|\bX\Big]
	-\partial^\bq \mu(\bx)\\
	=\,&\sum_{\bm{u}\in\Lambda_m}h_\bx^{m-[\bq]}B_{\bm{u},\bq}(\bx)
	\E[\widehat{\partial^{\bm{u}}\mu}_1(\bx)-\partial^{\bm{u}}\mu(\bx)|\bX]\\
	&+\partial^\bq\bp(\bx)'\widehat{\bQ}_m^{-1}\E_n\Big[\bp(\bx_i)
	\Big(\E[\widehat{\mathscr{B}}_{m,\mathbf{0}}(\bx_i)|\bX]-
	\mathscr{B}_{m,\mathbf{0}}(\bx_i)\Big)\Big] + O_\P(h^{m+\varrho-[\bq]})\\
	=\,&\partial^\bq\bp(\bx)'\widehat{\bQ}_m^{-1}\E_n\left[\bp(\bx_i)\Big(\E[\widehat{\mathscr{B}}_{m,\mathbf{0}}(\bx_i)|\bX]-
	\mathscr{B}_{m,\mathbf{0}}(\bx_i)\Big)\right]+O_\P(h^{m+\varrho-[\bq]}),
	\end{align*}
	where $\widehat{\mathscr{B}}_{m, \mathbf{0}}(\bx) = -\sum_{\bm{u}\in\Lambda_m}
	\widehat{\partial^{\bm{u}}\mu}_1(\bx)h_{\bx}^mB_{\bm{u}, \mathbf{0}}(\bx)$, and the last line follows from Assumption \ref{assumption: Misspecification Error}, \ref{assumption: BC bases} and Lemma \ref{lem: bias}. Also by Lemma \ref{lem: bias} and the fact that $B_{\bm{u},\mathbf{0}}(\cdot)$ is bounded, $\sup_{\bx\in\mathcal{X}}|\E[\widehat{\mathscr{B}}_{m,\mathbf{0}}(\bx)|
	\bX] - \mathscr{B}_{m,\mathbf{0}}(\bx)|\lesssim_\P h^{m+\varrho}$. The desired result immediately follows by the same argument as for $j=2$.	
\end{proof}

\subsection{Proof of Theorem 5.1} \label{proof: Theorem 5.1}
\begin{proof}
	For $j=0$, the result directly follows from \cite[Theorem 4.2]{Belloni-Chernozhukov-Chetverikov-Kato_2015_JoE} combined with our Lemma \ref{lem: pointwise linearization}.
	For $j=1,2,3$, by Lemma \ref{lem: Avar} and Lemma \ref{lem: pointwise linearization}, it suffices to show that 
	$\bgamma_{\bq, j}(\bx)'\G_n[\bPi{j}(\bx_i)\e_i]/\sqrt{\Omega_j(\bx)}$ weakly converges to the standard Normal distribution. First, by construction
	\[
	\V\Big[\frac{\bgamma_{\bq,j}(\bx)'}{\sqrt{\Omega_j(\bx)}}
	\G_n[\bPi{j}(\bx_i)\e_i]\Big]=1.
	\]
	Next, we write $a_{ni}:=\frac{\bgamma_{\bq, j}(\bx)'\bPi{j}(\bx_i)}{\sqrt{\Omega_j(\bx)}}$. For all $\vartheta>0$,
	\begin{align}
	\sum_{i=1}^{n}\E\left[a_{ni}^2\frac{\e_i^2}{n}\I\{|a_{ni}\e_i/\sqrt{n}|>\vartheta\}\right]&
	\leq\E\left[\E\left[a_{ni}^2\e_i^2\I\{|\e_i|>\vartheta\sqrt{n}/|a_{ni}|\}\Big|\bx_i\right]\right]\nonumber\\
	&\leq\E[a_{ni}^2]\cdot\sup_{\bx\in\mathcal{X}}
	\E\left[\e_i^2\I\{|\e_i|>\vartheta\sqrt{n}/|a_{ni}|\}\Big|\bx_i=\bx\right]\nonumber\\
	&\lesssim \frac{h^{-2[\bq]-2d}h^{d}}{h^{-d-2[\bq]}}
	\sup_{\bx\in \mathcal{X}}\E\left[\e_i^2\I\{|\e_i|>\vartheta\sqrt{n}/|a_{ni}|\}\Big|\bx_i=\bx\right]\nonumber\\
	&\lesssim \sup_{\bx\in \mathcal{X}}\E\left[\e_i^2\I\{|\e_i| > \vartheta\sqrt{n}/|a_{ni}|\}\Big|\bx_i=\bx\right], \label{eq: proof AN}
	\end{align}
	where the third line follows from Lemma \ref{lem: basis} and Lemma \ref{lem: Avar}. Since $|a_{ni}|\lesssim h^{-d/2}$ and $\log n/(nh^d)=o(1)$, it follows that $\sqrt{n}/|a_{ni}|\rightarrow\infty$ as $n\rightarrow\infty$. By the moment condition in the theorem, the upper bound in Eq. \eqref{eq: proof AN} goes to $0$ as $n\rightarrow\infty$, which completes the proof.	
\end{proof}

\subsection{Proof of Lemma \ref{lemma: uniform linearization}}
\begin{proof}	
	Consider the conditions in (i) hold. The proof is divided into two steps.
	
	\textit{Step 1:} We first bound $\sup_{\bx\in\mathcal{X}} |R_{1n,\bq}(\bx)|$ for $j=0,1,2,3$.		
	To simplify notation, we write $\bPi{j}(\bx_i)=(\pi_1(\bx_i),\ldots,\pi_{K_j}(\bx_i)'$ where $K_j=\dim(\bPi{j}(\cdot))$. We truncate the errors by an increasing sequence of constants $\{\vartheta_n: n\geq 1\}$ such that $\vartheta_n\asymp \sqrt{nh^d/\log n}$. Let $H_{ik}=\pi_k(\bx_i)(\e_i\I\{|\e_i|\leq\vartheta_n\}-
	\E[\e_i\I\{|\e_i|\leq\vartheta_n|\bx_i\}])$ and $T_{ik}=\pi_k(\bx_i)(\e_i\I\{|\e_i|>\vartheta_n\}-\E[\e_i\I\{|\e_i|>\vartheta_n|\bx_i\}])$. Regarding the truncated term, it follows from the truncation strategy, Assumption \ref{assumption: Local Basis}, and \ref{assumption: BC bases} that $|H_{ik}|\leq\vartheta_n$ and $\E[H_{ik}^2]\lesssim h^d$. By Bernstein's inequality, for every $t>0$,
	\begin{align}
	&\,\P\left(\max_{1\leq k\leq K_j}|\E_n[H_{ik}]|> h^d\sqrt{\log n/(nh^d)}\,t\right)\nonumber\\
	\leq&\,2\sum_{k=1}^{K_j}\exp\left\{-\frac{n^2h^{2d}h^{-d}\log n t^2/n}{C_1nh^d+C_2\vartheta_nnh^d\sqrt{\log n/(nh^d)}t}\right\}\nonumber\\
	\leq&\, C\exp\left\{\log n
	\Big(1-\frac{t^2}{C_1+C_2\vartheta_n\sqrt{\log n/(nh^d)}
		\,t}\Big)\right\},   \label{eq: uni-linear head}
	\end{align}
	which is arbitrarily small for $t$ large enough by the truncation strategy. By Lemma \ref{lem: basis},
	\begin{align*}
	&\quad\,\sup_{\bx\in \mathcal{X}}
	\Big|(\bg{j}'-\bgamma_{\bq,j}(\bx)')
	\E_n[\bPi{j}(\bx_i)(\e_i\I\{|\e_i|\leq\vartheta_n\}-
	\E[\e_i\I\{|\e_i|\leq\vartheta_n|\bx_i\}])]\Big|\\
	&\leq\sup_{\bx\in\mathcal{X}}\|\bg{j}'-
	\bgamma_{\bq,j}(\bx)'\|_\infty
	\|\E_n[\bPi{j}(\bx_i)(\e_i\I\{|\e_i|\leq\vartheta_n\}-
	\E[\e_i\I\{|\e_i|\leq\vartheta_n|\bx_i\}])]\|_\infty\\
	&\lesssim_\P h^{-[\bq]-d}\sqrt{\log n/(nh^d)}h^d\sqrt{\log n/(nh^d)}=h^{-[\bq]}\log n/(nh^d).
	\end{align*}
	
	Regarding the tails, let $\mathscr{K}_{ji}(\bx):=(\bg{j}'-\bgamma_{\bq,j}(\bx)')\bPi{j}(\bx_i)$. By Lemma \ref{lem: basis} and Assumption \ref{assumption: Local Basis}, we have
	$\sup_{\bx\in \mathcal{X}}|\mathscr{K}_{ji}(\bx)|\lesssim_\P h^{-d-[\bq]}\sqrt{\log n/(nh^d)}$.
	Let $\mathcal{A}_n(M)$ denote the event on which $\sup_{\bx\in\mathcal{X}}|\mathscr{K}_{ji}(\bx)|\leq Mh^{-d-[\bq]}\sqrt{\log n/(nh^d)}$ for some $M>0$, and $\I_{\mathcal{A}_n(M)}$ be an indicator function of $\mathcal{A}_n(M)$. Then by Markov's inequality, for any $t>0$,
	\begin{align}
	&\quad\,\P\Big(\sup_{\bx\in \mathcal{X}}
	\Big|\E_n[\I_{\mathcal{A}_n(M)}\mathscr{K}_{ji}(\e_i\I\{|\e_i|>\vartheta_n\}-\E[\e_i\I\{|\e_i|>\vartheta_n|\bx_i\}])]\Big|>
	\frac{t\log n}{nh^{d+[\bq]}}\Big)\nonumber \\
	&\lesssim \frac{Mh^{-d-[\bq]}\sqrt{\log n/(nh^d)}\E[|\e_i|\I\{|\e_i|>\vartheta_n\}]}{t h^{-[\bq]}
		\log n/(nh^d)}\nonumber\\
	&\leq \frac{M\sqrt{n}}{t\sqrt{h^d\log n}}
	\frac{\E[|\e_i|^{2+\nu}]}{\vartheta_n^{1+\nu}} \label{eq: uni-linear tail}
	\end{align}
	which is arbitrarily small for $t/M$ large enough by the additional moment condition specified in the lemma and the rate restriction. Since $\P(\mathcal{A}_n(M)^c)=o(1)$ as $M\rightarrow\infty$, simply let $t=M^2$ and $M\rightarrow\infty$, then the desired conclusion immediately follows.
	\smallskip
	
	\textit{Step 2:}
	Next, we bound $\sup_{\bx\in\mathcal{X}} |R_{2n,\bq}(\bx)|$.
	For $j=0,1$, the result directly follows from Lemma \ref{lem: bias}. For $j=2$, notice that
	the proof of Lemma \ref{lem: pointwise linearization} essentially establishes a bound on the uniform norm of $\partial^\bq\bp(\bx)' \bQ_m^{-1}\E_n[\bp(\bx_i)\mathfrak{B}_{\tilde{m},\mathbf{0}}(\bx_i)]$. The bound on $(\E[\widehat{\partial^\bq\mu}(\bx)|\bX]-\partial^\bq \mu(\bx))$ follows from Lemma \ref{lem: bias}. Then the desired bound on $R_{2n,\bq}(\bx)$ is obtained.
	
	For $j=3$, notice that
	by Assumption \ref{assumption: Misspecification Error} we can write $R_{2n,\bq}(\bx)$ explicitly as
	\begin{align*}
	&\bg{0}\Big\{\E_n[\bp(\bx_i)(\mu(\bx_i)-s^*(\bx_i))]-
	\E_n\Big[\bp(\bx_i)\sum_{\bm{u}\in\Lambda_m}h_{\bx_i}^m
	B_{\bm{u},\mathbf{0}}(\bx_i)\E\Big[\widehat{\partial^{\bm{u}}\mu}_1(\bx_i)\Big|\bX\Big]\Big]\Big\}\\
	&+\partial^\bq s^*(\bx)-
	\partial^\bq \mu(\bx)-\mathscr{B}_{m,\bq}(\bx)\\
	&+\mathscr{B}_{m,\bq}(\bx)+\sum_{\bm{u}\in\Lambda_m} h_\bx^{m-[\bq]}B_{\bm{u},\bq}(\bx)\partial^{\bm{u}}
	\bpt(\bx)'
	\widehat{\bQ}_{\tilde{m}}^{-1}\E_n[\bpt(\bx_i)\mu(\bx_i)].
	\end{align*}
    The second line is uniformly bounded by Assumption \ref{assumption: Misspecification Error}: $\sup_{\bx\in\mathcal{X}}|\partial^\bq s^*(\bx)-\partial^\bq \mu(\bx)-\mathscr{B}_{m,\bq}(\bx)|
	\lesssim h^{m+\varrho-[\bq]}$.
	The third line can be written as
	\begin{align*}
	&\mathscr{B}_{m,\bq}(\bx)+\sum_{\bm{u}\in\Lambda_m} h_\bx^{m-[\bq]}B_{\bm{u},\bq}(\bx)\partial^{\bm{u}}
	\bpt(\bx)'
	\widehat{\bQ}_{\tilde{m}}^{-1}\E_n[\bpt(\bx_i)\mu(\bx_i)]\\
	=\,&\sum_{\bm{u}\in\Lambda_m}h_\bx^{m-[\bq]}B_{\bm{u},\bq}(\bx)
	\Big(\partial^{\bm{u}}\bpt(\bx)'\widehat{\bQ}_{\tilde{m}}^{-1}\E_n[\bpt(\bx_i)\mu(\bx_i)]-\partial^{\bm{u}}\mu(\bx)\Big).
	\end{align*}
	By Assumption \ref{assumption: Misspecification Error},  $\sup_{\bx\in\mathcal{X}}|B_{\bm{u},\bq}(\bx)|
	\lesssim 1$. Moreover, as we have shown in the proof of Lemma \ref{lem: bias},
	$\sup_{\bx\in \mathcal{X}}
	|\partial^{\bm{u}}\bpt(\bx)'\widehat{\bQ}_{\tilde{m}}^{-1}\E_n[\bpt(\bx_i)\mu(\bx_i)]-\partial^{\bm{u}}\mu(\bx)|\lesssim_\P h^{\widetilde{m}-[\bm{u}]}$
	which suffices to show that the third line is $O_\P(h^{m+\varrho-[\bq]})$.
	
	In addition, using the proof strategy for Lemma \ref{lem: bias}, we have
	\[
	\begin{split}
	&\sup_{\bx\in\mathcal{X}}\Big|\bg{0}'\E_n[\bp(\bx_i)
	(\mu(\bx_i)-s^*(\bx_i)+\mathscr{B}_{m,\mathbf{0}}(\bx_i))]\Big|
	\lesssim_\P h^{m+\varrho-[\bq]}, \quad\text{and}\\
	&\sup_{\bx\in\mathcal{X}}\Big|\bg{0}'\E_n\Big[\bp(\bx_i)
	h_{\bx_i}^m B_{\bm{u},\mathbf{0}}(\bx_i)\Big(\E[
	\widehat{\partial^{\bm{u}}\mu}_1(\bx_i)|\bX]-\partial^{\bm{u}}\mu(\bx_i)\Big)\Big]\Big|\lesssim_\P h^{m+\varrho-[\bq]}.
	\end{split}
	\]
	
	Finally, if the conditions in (ii) hold, then it suffices to adjust the proof for $R_{1n,\bq}(\bx)$. We still use the same proof strategy, but let $\vartheta_n=\log n$. For the truncated term, it can be seen from  Bernstein's inequality (Equation \ref{eq: uni-linear head}) that the upper bound can be made arbitrarily small for $t$ large enough when $(\log n)^3/(nh^d)\lesssim 1$. On the other hand, when applying Markov's inequality to control the tail (Equation \eqref{eq: uni-linear tail}), we employ the exponential moment condition:
	\begin{align*}
	&\quad\,\P\left(\sup_{\bx\in \mathcal{X}}
	\Big|\E_n[\I_{\mathcal{A}_n(M)}\mathscr{K}_{ji}(\e_i\I\{|\e_i|
	>\vartheta_n\}-\E[\e_i\I\{|\e_i|>\vartheta_n|\bx_i\}])]\Big|>
	\frac{t\log n}{nh^{d+[\bq]}}\right)\\
	&\lesssim \frac{M\sqrt{n}}{t\sqrt{h^d\log n}}\E[|\e_i|\I\{|\e_i|>\vartheta_n\}]
	\leq \frac{M\sqrt{n}}{t\sqrt{h^d\log n}}
	\frac{\E[|\e_i|^3\exp(|\e_i|)]}{\vartheta_n^2\exp(\vartheta_n)}\\
	&\leq\frac{M}{t(\log n)^{5/2} \sqrt{nh^d}}\E[|\e_i|^3\exp(|\e_i|)]
	\end{align*}
	which is arbitrarily small for $t/M$ large enough. Thus the same bound on $R_{1n,\bq}$ is established.
	Then the proof is complete.
\end{proof}

\subsection{Proof of Theorem \ref{thm: uniform convergence}}
\begin{proof}
	Regarding the $L_2$ convergence, by Lemma \ref{lem: bias},
	\begin{align*}
	&\int_{\mathcal{X}}\Big(\widehat{\partial^{\bq}\mu}_0(\bx)
	-\partial^{\bq}\mu(\bx)\Big)^2\Big)w(\bx)d\bx\\
	=&\;\Big(\E_n[\bPi{0}(\bx_i)\e_i]'\Big)
	\Big(\int_{\mathcal{X}}\bg{0}\bg{0}'w(\bx)d\bx\Big) \Big(\E_n[\bPi{0}(\bx_i)\e_i]\Big)
	+O_\P(h^{2(m-[\bq])}).
	\end{align*}
	Notice that in the proof of Lemma \ref{lem: bias}, the uniform bound on the conditional bias does not require explicit expression of leading approximation error. Then by Lemma \ref{lem: Q matrix}, we have $\int_{\mathcal{X}}\bg{0}\bg{0}'w(\bx)d\bx\lesssim h^{-d-2[\bq]}$. Also, $\E[\|\E_n[\bPi{0}(\bx_i)\e_i]\|^2]\lesssim \E[\bPi{0}(\bx_i)'\bPi{0}(\bx_i)/n]\lesssim 1/n $. The desired $L_2$-convergence rate follows.
	
	Regarding the uniform convergence,	
	consider the case when the conditions of Lemma \ref{lemma: uniform linearization} hold. We use the same truncation strategy. Specifically, separate $\e_i$ into
	\[\e_i\I\{|\e_i|\leq\vartheta_n\}-\E[\e_i\I\{|\e_i|\leq\vartheta_n\}|\bx_i] \quad
	\text{and} \quad
	\e_i\I\{|\e_i|>\vartheta_n\}-\E[\e_i\I\{|\e_i|>\vartheta_n\}|
	\bx_i],
	\]
	where $\vartheta_n\asymp\sqrt{nh^d/\log n}$. By Lemmas \ref{lem: basis},
	$\sup_{\bx\in\mathcal{X}}|\bgamma_{\bq, j}(\bx)'\bPi{j}(\bx_i)|\lesssim_\P h^{-d-[\bq]}$.
	Then repeating the argument given in the proof of Lemma \ref{lemma: uniform linearization} for the truncated and tails respectively, we have
	\begin{equation*}
	\sup_{\bx\in\mathcal{X}}|\bgamma_{\bq,j}(\bx)'\E_n[\bPi{j}(\bx_i)
	\e_i]|
	\lesssim_\P h^{-[\bq]}\sqrt{\log n/(nh^d)}.
	\end{equation*}
	Moreover, $\bar{R}_{1n,\bq}=o(h^{-d/2-[\bq]}\sqrt{\log n/n})$ since $\log n/(nh^d)=o(1)$. In view of these bounds and Lemma \ref{lemma: uniform linearization}, the desired rate of uniform convergence follows.
	
	Finally, the same results can be proved under the conditions in (ii) of Lemma \ref{lemma: uniform linearization} if we let $\vartheta_n=\log n$ and assume 
	$(\log n)^3/(nh^d) \lesssim 1$.
\end{proof}

\subsection{Proof of Theorem \ref{thm: variance}}

\begin{proof}
	First consider the conditions in part (i) hold. Notice that for $j=0,1,2,3$,
	\begin{equation}\label{eq: pf thm var}
	\bSighat{j}-\bSig{j}=\E_n[(\widehat{\e}_{i,j}^2-\e_i^2)\bPi{j}(\bx_i)\bPi{j}(\bx_i)'] + \Big(\E_n[\e_i^2\bPi{j}(\bx_i)\bPi{j}(\bx_i)']-\bSig{j}\Big).
	\end{equation} 
	We then divide our proof into four steps.
	\smallskip
	
	\textit{Step 1:}
	For $\|\E_n[\bPi{j}(\bx_i)\bPi{j}(\bx_i)']\|$, by Assumption \ref{assumption: Local Basis}, \ref{assumption: BC bases}, and the same argument used in the proof of Lemma \ref{lem: Q matrix}, 
    $\|\E[\bPi{j}(\bx_i)'\bPi{j}(\bx_i)']\| \lesssim h^d$, and
    \[
    \|\E_n[\bPi{j}(\bx_i)\bPi{j}(\bx_i)']-
    \E[\bPi{j}(\bx_i)\bPi{j}(\bx_i)']\|
    \lesssim_\P h^d\sqrt{\log n/(nh^d)}. 
    \]
    By the triangle inequality and the fact that $\frac{\log n}{nh^d}=o(1)$, $\|\E_n[\bPi{j}(\bx_i)\bPi{j}(\bx_i)']
    \lesssim_\P h^d$.
    \smallskip	
	
	\textit{Step 2:}
	Next, we bound the second term in Equation \eqref{eq: pf thm var}. To simplify our notations, let $\mathbf{L}_j(\bx_i) :=\mathbf{W}_j^{-1/2}\bPi{j}(\bx_i)$ be the normalized basis, where $\mathbf{W}_j=\bQ_m$ for $j=0$, $\mathbf{W}_j=\bQ_{\tilde{m}}$ for $j=1$, and $\mathbf{W}_j=\diag\{\bQ_m, \bQ_{\tilde{m}}\}$ for $j=2,3$. Introduce a sequence of positive numbers:
	$M_n^2\asymp \frac{K^{1+1/\nu}n^{1/(2+\nu)}}{(\log n)^{1/(2+\nu)}}$.
	Then, we write
	\[\mathbf{H}_j(\bx_i)=\e_i^2\mathbf{L}_j(\bx_i)\mathbf{L}_j(\bx_i)'\I\{\|\e_i^2\mathbf{L}_j(\bx_i)\mathbf{L}_j(\bx_i)'\|\leq M_n^2\},\]
	and
	\[\mathbf{T}_j(\bx_i)=\e_i^2\mathbf{L}_j(\bx_i)\mathbf{L}_j(\bx_i)'\I\{\|\e_i^2\mathbf{L}_j(\bx_i)\mathbf{L}_j(\bx_i)'\|> M_n^2\}.\]
	Clearly, 
	\begin{align*}
	&\E_n[\mathbf{L}_j(\bx_i)\mathbf{L}_j(\bx_i)'\e_i^2]-
	\E[\mathbf{L}_j(\bx_i)\mathbf{L}_j(\bx_i)'\e_i^2]\\
	&= (\E_n[\mathbf{H}_j(\bx_i)]-
	\E[\mathbf{H}_j(\bx_i)])+
	(\E_n[\mathbf{T}_j(\bx_i)]-
	\E[\mathbf{T}_j(\bx_i)]).
	\end{align*}
	
	For the truncated terms, by definition, $\|\mathbf{H}_j(\bx_i)\|
	\leq M_n^2$. By the triangle inequality and Jensen's inequality, $\|\mathbf{H}_j(\bx_i)-
	\E[\mathbf{H}_j(\bx_i)]\|\leq 2M_n^2$. Also, by Assumption \ref{assumption: DGP},
	\begin{align*}
	\E[(\mathbf{H}_j(\bx_i)-\E[\mathbf{H}_j(\bx_i)])^2]&
	\leq\E[\e_i^4\|\mathbf{L}_j(\bx_i)\|^2\mathbf{L}_j(\bx_i)
	\mathbf{L}_j(\bx_i)'\I\{\|\e_i^2\mathbf{L}_j(\bx_i)
	\mathbf{L}_j(\bx_i)'\|\leq M_n^2 \}]\\
	&\leq M_n^2\E[\e_i^2\mathbf{L}_j(\bx_i)\mathbf{L}_j(\bx_i)'
	\I\{\|\e_i^2\mathbf{L}_j(\bx_i)\mathbf{L}_j(\bx_i)'\|\leq M_n^2\}]\\
	&\lesssim M_n^2\E[\mathbf{L}_j(\bx_i)\mathbf{L}_j(\bx_i)'],
	\end{align*}
	where the inequalities are understood in the sense of semi-definite matrices. Thus
	$\|\E[(\mathbf{H}_j(\bx_i)-\E[\mathbf{H}_j(\bx_i)])^2]\|
	\lesssim M_n^2$. Let
	$\vartheta_n=\sqrt{(\log n)^{\frac{\nu}{2+\nu}}/
		(n^{\frac{\nu}{2+\nu}}h^d)}$.
	By an inequality of \cite{Tropp2012_FCM} for independent matrices, for all $t>0$,
	\begin{align*}
	&\P\Big(\Big\|\E_n[\mathbf{H}_j(\bx_i)]
	-\E[\mathbf{H}_j(\bx_i)]\Big\|>\vartheta_n t\Big)\\
	&
	\leq\exp\Big(\log n-
	\frac{\vartheta_n^2 t^2/2}{M_n^2/n+M_n^2\vartheta_nt/(3n)}\Big)\\
	&\leq\exp\Big\{\log n \Big(1-\frac{t^2/2}{M_n^2\log n\,\vartheta_n^{-2}n^{-1}(1+\vartheta_nt/3)}\Big)\Big\}
	\end{align*}
	where $M_n^2\log n\,
	\vartheta_n^{-2}n^{-1} \asymp (\log n)^{\frac{1}{2+\nu}}/
	(n^{\frac{1}{2+\nu}}h^{d/\nu})=o(1)$ and $\vartheta_n=o(1)$. Hence
	\[\Big\|\E_n[\mathbf{H}_j(\bx_i)]
	-\E[\mathbf{H}_j(\bx_i)]\Big\|\lesssim_\P\vartheta_n=o(1).
	\]
	
	Regarding the tails, it directly follows from Lemma \ref{lem: Q matrix} that
	$\|\mathbf{T}_j(\bx_i)\|\lesssim h^{-d}\e_i^2\I\{\e_i^2
	\gtrsim M_n^2h^d\}$.
	Then by the triangle inequality, Jensen's inequality, and the assumption that $(2+\nu)$th moment of $\e_i$ is bounded,
	\begin{align*}
	\E\left[\Big\|\E_n[\mathbf{T}_j(\bx_i)]-
	\E[\mathbf{T}_j(\bx_i)]\Big\|\right]&\lesssim 2h^{-d}\E[\e_i^2\I\{|\e_i|\gtrsim M_n\sqrt{h^d}\}] \\
	&\lesssim \frac{2h^{-d(1+\nu/2)}
		\E[|\e_i|^{2+\nu}\I\{|\e_i|\gtrsim M_n\sqrt{h^d}\}]}
	{M_n^\nu} 
	\lesssim \vartheta_n.
	\end{align*}
	By Markov's inequality,
	$\|\E_n[\mathbf{T}_j(\bx_i)
	-\E[\mathbf{T}_j(\bx_i)]]\| \lesssim_\P \vartheta_n$. Since $\|\mathbf{W}_j^{1/2}\|\lesssim h^{d/2}$ and $\|\mathbf{W}_j^{-1/2}\|\lesssim h^{-d/2}$, we conclude that 
	$\|\E_n[\bPi{j}(\bx_i)\bPi{j}(\bx_i)'\e_i^2] -
	\bSig{j} \| \lesssim_\P h^d\vartheta_n = o_\P(h^d)$.
	\smallskip
	
	\textit{Step 3}: 
	The first term in Equation \eqref{eq: pf thm var} satisfies
	\begin{align*}
	&\,\|\E_n[(\widehat{\e}_{i,j}^2-
	\e_i^2)\bPi{j}(\bx_i)\bPi{j}(\bx_i)']\|\\
	\leq&\,\|\E_n[(\mu(\bx_i)-\widehat{\mu}_j(\bx_i))^2\bPi{j}(\bx_i)
	\bPi{j}(\bx_i)']\|+
	2\|\E_n[(\mu(\bx_i)-\widehat{\mu}_j(\bx_i))\e_i\bPi{j}(\bx_i)
	\bPi{j}(\bx_i)']\|\\
	\leq&\,\max_{1\leq i\leq n}|\mu(\bx_i)-\widehat{\mu}_j(\bx_i)|^2
	\|\E_n[\bPi{j}(\bx_i)\bPi{j}(\bx_i)']\|+\\
	&\,\max_{1\leq i\leq n}
	|\mu(\bx_i)-\widehat{\mu}_j(\bx_i)|
	\left(\|\E_n[\bPi{j}(\bx_i)\bPi{j}(\bx_i)']\|+
	\|\E_n[\bPi{j}(\bx_i)\bPi{j}(\bx_i)'\e_i^2]\|\right)
	\end{align*}
	where the last line follows from the fact that $2|a|\leq 1+a^2$. By Theorem \ref{thm: uniform convergence} and the results proved in Step 1 and 2, we have $\max_{1\leq i\leq n}
	|\mu(\bx_i)-\widehat{\mu}_j(\bx_i)|=R^{\mathtt{uc}}_{\mathbf{0},j}=
	o_\P(1)$, $\|\E_n[\bPi{j}(\bx_i)\bPi{j}(\bx_i)']\|\lesssim_\P h^d$ and $\|\E_n[\bPi{j}(\bx_i)\bPi{j}(\bx_i)'\e_i^2]\|\lesssim_\P h^d$. Hence	
	\begin{equation}\label{eq: pf sigma}
	\|\bSighat{j}-\bSig{j}\|\lesssim_\P h^d(R_{\mathbf{0}, j}^{\mathtt{uc}}+\vartheta_n)=o_\P(h^d)
	\end{equation}
	\smallskip
	
	\textit{Step 4:}
	Using all above results, we have
	\begin{align*}
	\quad\;\;\Big|\widehat{\Omega}_j(\bx)-\Omega_j(\bx)\Big|
	&=\,\|\bg{j}'\bSighat{j}\bg{j}-\bgamma_{\bq,j}(\bx)'\bSig{j}
	\bgamma_{\bq,j}(\bx)\|\\
	&\leq\, \|(\bg{j}'-\bgamma_{\bq,j}(\bx))'\bSighat{j}\bg{j}\|+
	\|\bgamma_{\bq,j}(\bx)'(\bSighat{j}-\bSig{j})\bg{j}\| \\
	&\quad+\|\bgamma_{\bq,j}(\bx)'\bSig{j}(\bg{j}-\bgamma_{\bq,j}(\bx))\|
	\end{align*}
	By Lemma \ref{lem: basis} and Equation \eqref{eq: pf sigma}, we have
	\[
	\sup_{\bx\in\mathcal{X}}\Big|\widehat{\Omega}_j(\bx)-\Omega_j(\bx)\Big|
	\lesssim_\P h^{-d-2[\bq]}(R_{\mathbf{0},j}^{\mathtt{uc}}+\vartheta_n)
	=o_\P(h^{-d-2[\bq]}).
	\]	
	
	Finally, when the conditions in part (ii) hold, we only need to adjust the proof in Step 2. Apply the same proof strategy with $M_n=\sqrt{Ch^{-d}}\log n$ and $\vartheta_n=\sqrt{\frac{(\log n)^3}{nh^d}}$. For the truncated term, since $M_n^2\log n/(\vartheta_n^2 n)\lesssim 1$, by Bernstein's inequality, $\|\E_n[\mathbf{H}_j(\bx_i)]-
	\E[\mathbf{H}_j(\bx_i)]\| \lesssim_\P \vartheta_n=o_\P(1)$. On the other hand, when applying Markov's inequality to bound the tail, we employ the stronger moment condition:
	\begin{align*}
	\E\Big[\|\E_n[\mathbf{T}_j(\bx_i)]-
	\E[\mathbf{T}_j(\bx_i)]\|\Big]&\lesssim 2h^{-d}\E[\e_i^2\I\{|\e_i|\geq M_n/\sqrt{Ch^{-d}}\}]\\
	&\lesssim \frac{h^{-3d/2}}{M_n\exp(M_n/\sqrt{Ch^{-d}})}\E[|\e_i|^3\exp(|\e_i|)]
	\lesssim \vartheta_n=o(1).
	\end{align*}
	
	Then the proof is complete.
\end{proof}

\subsection{Proof of Theorem \ref{thm: KMT coupling}}
\begin{proof}
By Lemma \ref{lem: hats off}, the remainders in the linearization of $\widehat{T}_j(\cdot)$ is $o_\P(r_n^{-1})$ in $\mathcal{L}^\infty(\mathcal{X})$, and hence we only need to show that $Z_j(\cdot)$ approximates $t_j(\cdot)$. 

Suppose that the conditions in (i) hold.
As a first step, we define $\mathscr{K}(\bx,\bx_i):=\frac{\bgamma_{\bq,j}(\bx)'\bPi{j}(\bx_i)}{\sqrt{\Omega_j(\bx)}}$ for each $j=0,1,2,3$, and will construct conditional coupling for $t_j(\cdot)$. One can verify conditions in Lemma 8.2 of the main paper (see Section 8.3) by using the properties of local bases considered in this paper, but here we directly follow the strategy given in the proof of Lemma 8.2, which will unify the proofs for $d=1$ and Haar basis with $d>1$. 

We need to construct a rearranged sequence $\{\bx_{i,n}\}_{i=1}^n$ of $\{\bx_i\}_{i=1}^n$. If $d=1$, we order them from the smallest to the largest: $x_{1,n}\leq\cdots\leq x_{n,n}$ (which are simply order statistics as in the proof of Lemma 8.2 of the main paper). For Haar basis with $d>1$ and $j=0$, we first classify $\{\bx_i\}_{i=1}^n$ into $\bar{\kappa}$ groups so that the points in the same group belong to the same cell (recall that $\bar{\kappa}$ is the number of cells in $\Delta$). Number the cells in an arbitrary way. Under the new ordering, the points in the first group is followed by the second, then the third, and so on. Ordering within each group is arbitrary. $\{\sigma_i\}_{i=1}^n$ and $\{\e_i\}_{i=1}^n$ are rearranged as $\{\sigma_{i,n}\}_{i=1}^n$ and $\{\e_{i,n}\}_{i=1}^n$ accordingly. 
Again, the ordering strategies described above only depend on the values of $\{\bx_i\}_{i=1}^n$, and thus conditional on $\bX$, the new sequence $\{\e_{i,n}\}_{i=1}^n$ is still independent. As explained in the proof of Lemma 8.2, we can find a sequence of i.i.d standard normal random variables $\{\zeta_{i,n}\}_{i=1}^n$ such that $\max_{1\leq l \leq n}|S_{l,n}|\lesssim_\P n^{\frac{1}{2+\bar{\nu}}}$ where $S_{l,n}=\sum_{i=1}^{l}(\e_{i,n}-\sigma_{i,n}\zeta_{i,n})$. We will write $\{\zeta_i\}_{i=1}^n$ for the sequence of $\zeta_{i,n}$'s rearranged based on the original ordering of $\{\bx_i\}_{i=1}^n$. Using summation by parts, 
\begin{align*}
&\,\sup_{\bx\in\mathcal{X}}\Big|\sum_{i=1}^{n}\mathscr{K}(\bx, \bx_{i,n})(\e_{i,n}-\sigma_{i,n}\zeta_{i,n})\Big|\\
=&\,\sup_{\bx\in\mathcal{X}}\Big|\mathscr{K}(\bx, \bx_{n,n})S_{n,n}
-\sum_{i=1}^{n-1}S_{i,n}\left(\mathscr{K}(\bx,\bx_{i+1,n})
-\mathscr{K}(\bx,\bx_{i,n})\right)\Big|\\
\leq&\,\sup_{\bx\in\mathcal{X}}\max_{1\leq i\leq n}|\mathscr{K}(\bx,\bx_i)||S_{n,n}|
+\sup_{\bx\in\mathcal{X}}\Big|
\frac{\bgamma_{\bq,j}(\bx)'}{\sqrt{\Omega_j(\bx)}}
\sum_{i=1}^{n-1}S_{i,n}(\bPi{j}(\bx_{i+1,n})-
\bPi{j}(\bx_{i,n}))\Big|\\
\leq&\,\sup_{\bx\in\mathcal{X}}\max_{1\leq i\leq n}|\mathscr{K}(\bx,\bx_i)||S_{n,n}|
+\sup_{\bx\in\mathcal{X}}
\bigg\|\frac{\bgamma_{\bq,j}(\bx)'}{\sqrt{\Omega_j(\bx)}}\bigg\|_\infty
\bigg\|\sum_{i=1}^{n-1}S_{i,n}(\bPi{j}(\bx_{i+1,n})-
\bPi{j}(\bx_{i,n}))\bigg\|_\infty.
\end{align*}
By Assumption \ref{assumption: Local Basis}, Lemma \ref{lem: basis}, and \ref{lem: Avar}, $\sup_{\bx\in\mathcal{X}}\max_{1\leq i\leq n}|\mathscr{K}(\bx,\bx_i)|\lesssim h^{-d/2}$ and
\begin{equation} \label{eq: uniform bound in thm 5.3}
\sup_{\bx\in\mathcal{X}}
\bigg\|\frac{\bgamma_{\bq,j}(\bx)'}{\sqrt{\Omega_j(\bx)}}\bigg\|_\infty
\lesssim h^{-d/2}.
\end{equation}
Also, if we write the $l$th element of $\bPi{j}(\cdot)$ as $\pi_{j,l}(\cdot)$, then
\begin{equation*}
\begin{split}
&\max_{1\leq l \leq K_j}
\Big|\sum_{i=1}^{n-1}\Big(\pi_{j,l}(\bx_{i+1,n})-\pi_{j,l}(\bx_{i,n})\Big)S_{l,n}\Big|\\
\leq& \max_{1\leq l \leq K_j}
\sum_{i=1}^{n-1}\Big|\pi_{j,l}(\bx_{i+1,n})-\pi_{j,l}(\bx_{i,n})\Big|
\max_{1\leq \ell \leq n}\Big|S_{\ell,n}\Big|
\lesssim_\P n^{\frac{1}{2+\bar{\nu}}},
\end{split}
\end{equation*}
where the last inequality follows from \cite[Corollary 5]{Sakhanenko_1991_SAM}, the ordering of $\{\bx_{i,n}\}_{i=1}^n$, Assumption \ref{assumption: Local Basis} and the fact that when $d>1$, each $\pi_{j,l}(\cdot)$ is a Haar basis function. This shows that there exists independent standard Normal random variables $\{\zeta_i\}_{i=1}^n$ such that
\[ \G_n[\mathscr{K}(\bx,\bx_i)\e_i]=_d z_j(\bx)+o_\P(r_n^{-1}) 
\quad \text{where} \quad z_j(\bx):=\G_n[\mathscr{K}(\bx,\bx_i)\sigma_i\zeta_i].
\]

Next, we convert the conditional coupling to an unconditional one. Notice that
\begin{equation*}
z_j(\bx)
=_{d|\bX}\;\frac{\bgamma_{\bq,j}(x)'}{\sqrt{\Omega_j(x)}}
\bSigbar{j}^{1/2}\bN_{K_j}
=\frac{\bgamma_{\bq,j}(\bx)'}{\sqrt{\Omega_j(\bx)}}
\bSig{j}^{1/2}\bN_{K_j} + \frac{\bgamma_{\bq,j}(\bx)'}{\sqrt{\Omega_j(\bx)}}
\left(\bSigbar{j}^{1/2}-\bSig{j}^{1/2}\right)\bN_{K_j}
\end{equation*}
where $\bSigbar{j}:= \E_n[\bPi{j}(\bx_i)\bPi{j}(\bx_i)'\sigma^2(\bx_i)]$, $\bN_{K_j}$ is a $K_j$-dimensional standard Normal vector (independent of $\bX$) and ``$=_{d|\bX}$'' denotes that two processes have the same conditional distribution given $\bX$.	
For the second term, we already have Equation \eqref{eq: uniform bound in thm 5.3}, and by the Gaussian maximal inequality (see \cite[Lemma 13]{Chernozhukov-Lee-Rosen_2013_ECMA}),
\[\E\left[ \left\|\left(\bSigbar{j}^{1/2}-\bSig{j}^{1/2}\right)
\bN_{K_j}\right\|_\infty\,\Big|\bX\right]
\lesssim \sqrt{\log n}
\,\left\|\bSigbar{j}^{1/2}-\bSig{j}^{1/2}\right\|.
\]
By the same argument given for Theorem \ref{thm: variance},
$\Big\|\bSigbar{j}-\bSig{j}\Big\|
\lesssim_\P h^d\sqrt{\log n/(nh^d)}$.
Then it follows from \cite[Theorem X.1.1]{Bhatia_2013_Book} that
\[\left\|\bSigbar{j}^{1/2}-\bSig{j}^{1/2}\right\|
\lesssim_\P h^{d/2}\left(\frac{\log n}{nh^d}\right)^{1/4}.
\]

For $j=0,1$, a sharper bound is available: by Theorem X.3.8 of \cite{Bhatia_2013_Book} and Lemma \ref{lem: Q matrix},
\[\|\bSigbar{j}^{1/2}-\bSig{j}^{1/2}\|\leq
\frac{1}{\lambda_{\min}(\bSig{j})^{1/2}}
\|\bSigbar{j}-\bSigma_{j}\|
\lesssim_\P h^{d/2}\sqrt{\log n/(nh^d)}.
\]
Thus using all these results, we have
\begin{equation*}
\E\Big[\sup_{\bx\in\mathcal{X}}\Big|\frac{\bgamma_{\bq,j}(\bx)'}{\sqrt{\Omega_j(\bx)}}\Big(\bSigbar{j}^{1/2}-
\bSig{j}^{1/2}\Big)
\bN_{K_j}\Big|\;\Big|\bX \Big]
\lesssim_\P h^{-d/2}\sqrt{\log n}
\left\|\bSigbar{j}^{1/2}-\bSig{j}^{1/2}\right\|=o_\P(r_n^{-1}),
\end{equation*}
where the last equality holds by the additional rate restriction given in the theorem (for $j=0,1$, no additional restriction is needed). By Markov inequality, this suffices to show that for any $\vartheta>0$,
\[
\P\left(\sup_{x\in\mathcal{X}}
\Big|z_j(\bx)-Z_j(\bx)\Big|>\vartheta\,\Big|\bX\right)= o_\P(r_n^{-1}).
\]
Since the conditional probability is bounded, by dominated convergence theorem, the desired result immediately follows.

When the conditions in (ii) hold, the proof remains the same except that we employ \cite[Theorem 1]{Sakhanenko_1985_Advances} to construct strong approximation to the partial sum process of $\{\e_{i,n}\}_{i=1}^n$.
\end{proof}

\subsection{Proof of Theorem \ref{thm: plug-in approx cond}}
\begin{proof}	
	Suppose that the conditions in (i) of Lemma \ref{lem: hats off} hold. The other case follows similarly.
	Let $z_j(\bx):=\G_n[\mathscr{K}(\bx,\bx_i)\sigma_i\zeta_i]$. First notice that for $j=0,1,2,3$,
	\[
	\begin{split}
	\sup_{\bx\in\mathcal{X}}\,|\widehat{z}_j(\bx)-z_j(\bx)|
	&=\sup_{\bx\in\mathcal{X}}\,\bigg|\frac{\bg{j}'}
	{\widehat{\Omega}_j(\bx)^{1/2}}\G_n[\bPi{j}(\bx_i)\zeta_i\hat{\sigma}(\bx_i)]
	-\frac{\bgamma_{\bq,j}(\bx)'}{\Omega_j(\bx)^{1/2}}
	\G_n[\bPi{j}(\bx_i)\zeta_i\sigma(\bx_i)]\bigg|\\
	&=:\sup_{\bx\in\mathcal{X}}\,|D_n(\bx)|.
	\end{split}
	\]
	Conditional on the data, $\{D_n(\bx),\bx\in\mathcal{X}\}$ is a Gaussian process with zero means. Let
	$\E^*$ denote the expectation with respect to the distribution of $\{\zeta_i\}_{i=1}^n$. For notational simplicity, define a norm $\|\cdot\|_{n,2}$ by $\|\mathbf{a}\|_{n,2}^2=n^{-1}\sum_{i=1}^{n}a_i^2$ for $\mathbf{a}\in\mathbb{R}^n$.
	Then by Cauchy-Schwarz inequality,
	\begin{align*}
	\E^*[D_n(\bx)^2]&=\frac{1}{n}\sum_{i=1}^n\bigg[
	\frac{\bg{j}'}{\widehat{\Omega}_j(\bx)^{1/2}}\bPi{j}(\bx_i)\hat{\sigma}(\bx_i)-\frac{\bgamma_{\bq,j}(\bx)'}{\Omega_j(\bx)^{1/2}}\bPi{j}(\bx_i)\sigma(\bx_i)
	\bigg]^2\\
	&=\frac{1}{n}\sum_{i=1}^{n}\bigg\{\frac{\bg{j}'
		-\bgamma_{\bq,j}(\bx)'}{\widehat{\Omega}_j(\bx)^{1/2}}
	\bPi{j}(\bx_i)\hat{\sigma}(\bx_i)\\
	&\quad+\bgamma_{\bq,j}(\bx)'(\widehat{\Omega}_j(\bx)^{-1/2}-
	\Omega_j(\bx)^{-1/2})\bPi{j}(\bx_i)\hat{\sigma}(\bx_i)\\
	&\quad+\frac{\bgamma_{\bq,j}(\bx)'}{\Omega_j(\bx)^{1/2}}\bPi{j}(\bx_i)
	(\hat{\sigma}(\bx_i)-\sigma(\bx_i))\bigg\}^2\\
	&\lesssim \bigg\|\frac{\bg{j}'-\bgamma_{\bq,j}(\bx)'}
	{\widehat{\Omega}_j(\bx)^{1/2}}\bPi{j}(\bx_i)\hat{\sigma}(\bx_i)\bigg\|_{n,2}^2\\
	&\quad+\Big\|\bgamma_{\bq,j}(\bx)'(\widehat{\Omega}_j(\bx)^{-1/2}-
	\Omega_j(\bx)^{-1/2})\bPi{j}(\bx_i)\hat{\sigma}(\bx_i)\Big\|_{n,2}^2\\
	&\quad+\Big\|\frac{\bgamma_{\bq,j}(\bx)'}{\Omega_j(\bx)^{1/2}}
	\bPi{j}(\bx_i)(\hat{\sigma}(\bx_i)-\sigma(\bx_i))\Big\|_{n,2}^2.
	\end{align*}
	
	By Lemmas \ref{lem: basis} and \ref{lem: Avar}, Theorem \ref{thm: variance}, and $\max_{1\leq i\leq n}
	|\hat{\sigma}^2(\bx_i)-\sigma^2(\bx_i)|=o_\P(1/(r_n\sqrt{\log n}))$, 
	\[
	\begin{split}
	\sup_{\bx\in\mathcal{X}}\E^*[D_n(\bx)^2]
	&\lesssim_\P \frac{\log n}{nh^d}+\Big(\frac{n^{\frac{2}{2+\nu}}(\log n)
	^{\frac{\nu}{2+\nu}}}{nh^d}+a_n\Big)+(\max_{1\leq i\leq n}|\hat{\sigma}(\bx_i)-\sigma(\bx_i)|)^2\\
	&=o_\P\Big(\frac{1}{r_n^2\log n}\Big)
	\end{split}
	\]
	where $a_n=h^{2m}$ for $j=0$ and $a_n=h^{2m+2\varrho}$ for $j=1,2,3$.
	Moreover, using the fact that for $\bx, \check{\bx}\in\mathcal{X}$,
	\begin{align*}
	&(\E^*[(D_n(\bx)-D_n(\check{\bx}))^2])^{1/2}\\
	&\lesssim\, \bigg\|\Big(\frac{\bg{j}_{\bq,j}(\bx)'}{\widehat{\Omega}_j(\bx)^{1/2}}
	-\frac{\widehat{\bgamma}_{\bq,j}(\check{\bx})'}{\widehat{\Omega}_j(\check{\bx})^{1/2}}\Big)\bPi{j}(\bx_i)\hat{\sigma}(\bx_i)\bigg\|_{n,2}\\
	&\qquad+\bigg\|\Big(\frac{\bgamma_{\bq,j}(\bx)'}{\Omega_j(\bx)^{1/2}}
	-\frac{\bgamma_{\bq,j}(\check{\bx})'}{\Omega_j(\check{\bx})^{1/2}}\Big)
	\bPi{j}(\bx_i)\sigma(\bx_i)\bigg\|_{n,2},
	\end{align*}
	we will show later in the proof that $\sup_{\delta\in\Delta}\sup_{\bx,\check{\bx}\in\clo(\delta)}
	(\E^*[(D_n(\bx)-D_n(\check{\bx}))^2])^{1/2}
	\lesssim_\P h^{-1}\|\bx-\check{\bx}\|$.
	Since $\sup_{\delta\in\Delta}\vol(\delta)\lesssim h^{d}$, and the number of cells in $\Delta$ is $\bar{\kappa}\asymp h^{-d}$, we apply Dudley's inequality (see the proof of \cite[Lemma 4.2]{Belloni-Chernozhukov-Chetverikov-Kato_2015_JoE} for more details of this inequality) to obtain
	$\E^*\Big[\sup_{\bx\in\mathcal{X}}|D_n(\bx)|\Big]=o_\P(r_n^{-1})$,
	which suffices to show for any $t>0$, $\P^*(\sup_{\bx\in\mathcal{X}}|D_n|>t r_n^{-1})=o_\P(1)$ by Markov inequality. Then, to show $z_j(\cdot)$ is approximated by $Z_j(\cdot)$, we only need to repeat the argument given in the proof of Theorem \ref{thm: KMT coupling}.
	
	In the end, we establish the desired Lipschitz bound on the basis and variance. For any $\bx,\check{\bx}\in\delta$, by Assumption \ref{assumption: Local Basis}, there are only a finite number of basis functions in $\bp$ and $\bpt$ which are active on $\delta$, and hence it follows from Lemma \ref{lem: Q matrix} that there exists some universal constant $C_1>0$ such that $\|\bgamma_{\bq,j}(\bx)-\bgamma_{\bq,j}(\check{\bx})\|\leq C_1h^{-[\bq]-1}h^{-d}\|\bx-\check{\bx}\|$. Similarly, by Lemma \ref{lem: basis} and Theorem \ref{thm: variance}, $\|\bg{j}-\widehat{\bgamma}_{\bq.j}(\check{\bx})\|\lesssim_\P h^{-[\bq]-1}h^{-d}\|\bx-\check{\bx}\|$, $|\Omega_j(\bx)-\Omega_j(\check{\bx})|\lesssim h^{-d-2[\bq]-1}\|\bx-\check{\bx}\|$, and $|\widehat{\Omega}_j(\bx)-\widehat{\Omega}_j(\check{\bx})|\lesssim_\P h^{-d-2[\bq]-1}\|\bx-\check{\bx}\|$.
	Then the proof is complete.
\end{proof}

\subsection{Proof of Theorem \ref{thm: wild bootstrap KMT}}
\begin{proof}	
	Let $\E^*$ be the expectation conditional on the data. Clearly, $\E^*[\widehat{\partial^\bq \mu}_j^*(\bx)]=\widehat{\partial^\bq \mu}_j(\bx)$. Therefore,
	\[\frac{\widehat{\partial^\bq \mu}_j^*(\bx)-\widehat{\partial^\bq \mu}_j(\bx)}
	{(\widehat{\Omega}^*_j(\bx)/n)^{1/2}}
	=\frac{\bg{j}'}{\widehat{\Omega}^*_j(\bx)^{1/2}}
	\G_n[\bPi{j}(\bx_i)\omega_i\widehat{\e}_{i,j}].
	\]
	Note that $\widehat{\e}_{i,j}=y_i-\widehat{\mu}_j(\bx_i)=\mu(\bx_i)-\widehat{\mu}_j(\bx_i)+\e_i$. Then
	\begin{align}
	&\frac{\widehat{\partial^\bq \mu}_j^*(\bx)-\widehat{\partial^\bq \mu}_j(\bx)}
	{(\widehat{\Omega}^*_j(\bx)/n)^{1/2}}
	=\frac{\bgamma_{\bq,j}(\bx)'}{\Omega_j(\bx)^{1/2}}
	\G_n[\bPi{j}(\bx_i)\widehat{\e}_{i,j}\omega_i]+o_\P(r_n^{-1}) \nonumber\\
	=\,&\frac{\bgamma_{\bq.j}(\bx)'}{\Omega_j(\bx)^{1/2}}
	\G_n[\bPi{j}(\bx_i)\e_i\omega_i] + \nonumber\\ &\frac{\bgamma_{\bq,j}(\bx)'}{\Omega_j(\bx)^{1/2}}
	\G_n[\bPi{j}(\bx_i)(\mu(\bx_i)-\widehat{\mu}_j(\bx_i))\omega_i]
	+o_\P(r_n^{-1}) \nonumber\\
	=\,&\frac{\bgamma_{\bq,j}(\bx)'}{\Omega_j(\bx)^{1/2}}
	\G_n[\bPi{j}(\bx_i)\e_i\omega_i]+o_\P(r_n^{-1})\quad\text{in}
	\quad\mathcal{L}^\infty(\mathcal{X}), \label{eq: wb leading term}
	\end{align}
	where the first and second equalities follow from a similar argument for Theorem \ref{thm: plug-in approx cond}, and the last line uses the rate of uniform convergence given in Theorem \ref{thm: uniform convergence}. Note that we employ different rate restrictions given in the theorem: for $j=1,2,3$, $h^{m+\varrho}\sqrt{\log n}\lesssim\sqrt{n}h^{d/2+m+\varrho}=o(r_n^{-1})$, whereas for $j=0$, $h^m\sqrt{\log n}\lesssim \sqrt{n}h^{d/2+m}=o(r_n^{-1})$. Moreover, either $\frac{\log n}{\sqrt{nh^d}}\lesssim\frac{n^{\frac{1}{2+\nu}}(\log n)^{\frac{1+\nu}{2+\nu}}}{\sqrt{nh^d}}=o(r_n^{-1})$, or 
	$\frac{\log n}{\sqrt{nh^d}}\lesssim\frac{(\log n)^2}{\sqrt{nh^d}}=o(r_n^{-1})$.
	Simply notice that $M_n\lesssim_\P 1$ means $\P(|M_n|>\vartheta_n)=o(1)$ for any $\vartheta_n\rightarrow\infty$. By Markov inequality, it immediately follows that $\P^*(|M_n|>\vartheta_n)=o_\P(1)$. Thus the above derivation still holds in $P$-probability if we replace $\P$ by $\P^*$.
	
	Then we only need to construct strong approximation to the leading term in Equation \eqref{eq: wb leading term}.
	Repeat the argument for Theorem \ref{thm: KMT coupling} conditional on the data. Regarding the conditional coupling step, we still rearrange terms according to the values of $\{\bx_i\}_{i=1}^n$ as described in the proof of Theorem \ref{thm: KMT coupling}. Now if $(2+\nu)$th moment of $\e_i$ is bounded, then conditional on the data, $\{\e_{i,n}\omega_{i,n}\}_{i=1}^n$ is independent and $\E_n[\E^*[|\e_i\omega_i|^{2+\nu}]]\lesssim_\P 1$. Thus, strong approximation to the partial sum process $\{\sum_{i=1}^{l}\e_{i,n}\omega_{i,n}: 1\leq l\leq n\}$ can be obtained by using \cite[Corollary 5]{Sakhanenko_1991_SAM}.
	When the conditions in (ii) of the theorem hold, \cite[Theorem 1]{Sakhanenko_1985_Advances} can be employed to construct the strong approximation to this partial sum process.
	Thus, conditional on the data, there exists a $K_j$-dimensional standard Normal vector $\mathbf{N}_{K_j}$ such that
	\[
	\widehat{z}^*(\cdot)=_{d^*}\frac{\bgamma_{\bq,j}(\cdot)'}{\Omega_j(\cdot)^{1/2}}
	\left(\E_n[\bPi{j}(\bx_i)\bPi{j}(\bx_i)'\e_i^2]\right)^{1/2}\bN_{K_j}+o_{\P^*}(r_n^{-1}).
	\]
	By Theorem \ref{thm: variance} and the argument given in the proof of Theorem \ref{thm: KMT coupling}, we can further show that the conditional Gaussian process on the RHS is further approximated by $Z_j(\cdot)$. Then the proof is complete.
	\end{proof}

\subsection{Proof of Theorem \ref{thm: wild bootstrap Yurinskii}}
\begin{proof}
	The proof is exactly the same as that for Theorem \ref{thm: wild bootstrap KMT} except that we apply the improved Yurinskii's inequality as in Theorem \ref{thm: Yurinskii coupling} conditional on the data.
	The conditions needed for coupling can be easily verified  
	since $\boldsymbol{\omega}$ is assumed to be independent of the data and bounded.
\end{proof}

\subsection{Proof of Theorem \ref{thm: confidence band}}
\begin{proof}
	The proof is similar to that given in Section 8.6 of the main paper, and thus omitted here.

\end{proof}

\subsection{Proof of Lemma \ref{lemma: B-Splines}}
\begin{proof}
	(a): Assumption \ref{assumption: Local Basis}(a) directly follows from the construction of tensor-product $B$-splines. Assumption \ref{assumption: Local Basis}(b) follows from \cite[Theorem 12.5]{Schumaker_2007_Book}.
	To show Assumption \ref{assumption: Local Basis}(c), note that for a univariate $B$-spline basis $\breve{\bp}_\ell(x_\ell)=(\breve{p}_1(x_\ell),\cdots,\\ \breve{p}_{K_\ell}(x_\ell))^\prime$, there exists a universal constant $C>0$ such that for any $\varsigma_\ell\leq m$, $1\leq l_\ell\leq K_{\ell}$, and $\delta\in\Delta$,
	$\|d^{\varsigma_\ell}\breve{p}_{l_\ell}(x_\ell)/dx_\ell^{\varsigma_\ell}\|_{L_\infty(\clo(\delta))}\lesssim h^{-\varsigma_\ell}$.
	Since there are only a fixed number of nonzero elements in $\bp$, we have for $[\bvsig]\leq m$,
	$\sup_{\delta\in\Delta}\sup_{\bx\in\clo(\delta)}\|\partial^{\bvsig}\bp(\bx)\|\lesssim h^{-[\bvsig]}$.
	To derive the other side of the bound, notice that the proof of \cite[Lemma 5.4]{Zhou-Wolfe_2000_SS} shows that for the univariate $B$-spline basis $\breve{\bp}_\ell(x_\ell)$, for any $\varsigma_\ell\leq m-1$, $x_\ell\in\mathcal{X}_\ell$,
	$\|d^{\varsigma_\ell}\breve{\bp}_\ell(x_\ell)/dx_\ell^{\varsigma_\ell}\|\gtrsim h^{-\varsigma_\ell}$.
	Since for any $x_\ell$, there are only $m$ nonzero elements in $\breve{\bp}_\ell(x_\ell)$, this suffices to show that for any $x_\ell\in\mathcal{X}_\ell$, there exists some $\breve{p}_{l_\ell}(x_\ell)$ such that
	$|d^{\varsigma_\ell}p_{l_\ell}(x_\ell)/dx_\ell^{\varsigma_\ell}|\gtrsim h^{-\varsigma_\ell}$.
	Then the lower bound directly follows from the construction of tensor-product $B$-splines.
	\smallskip
	
	(b):
	The proof of orthogonality between the constructed leading error and $B$-splines can be found in \cite{Barrow-Smith_1978_QAM}.
	Regarding the bias expansion, we first consider $\bvsig=\mathbf{0}$. Noticing that
	$$\mathscr{B}_{m,\mathbf{0}}(\bx)=
	-\sum_{\ell=1}^{d}\frac{\partial^{m}\mu(\bt_{\bx}^L)}{\partial x_\ell^{m}}\frac{b_{\ell,l_\ell}^{m}}{m!}B_{m}\Big(\frac{x_\ell-t_{\ell,l_\ell}}{b_{\ell,l_\ell}}\Big)+O(h^{m+\varrho})\; \text{ for }\bx\in \delta_{l_1\ldots l_d}$$
	where $B_m(\cdot)$ is the $m$th Bernoulli polynomial, we only need to focus on the first term on the RHS, denoted by $\bar{\mathscr{B}}_{m}(\bx)$. By construction, $\bar{\mathscr{B}}_{m}(\bx)$ is continuous on the interior of each subrectangle $\delta_{l_1\ldots l_d}$, and the discontinuity
	only takes place at boundaries of subrectangles. Let $J_{\mathbf{0}}$ denote the magnitude of a (generic) jump of $\bar{\mathscr{B}}_{m}(\bx)$. By Assumption \ref{assumption: DGP}, $J_{\mathbf{0}}$ is also a jump of $\bar{\mu}:= \mu+\bar{\mathscr{B}}_{m}$. We first bound it as  \cite{Barrow-Smith_1978_QAM} did in their proof. We introduce the following notation:
	\begin{enumerate}[label=\roman*.]
		\item $\boldsymbol{\tau}:=(\tau_1,\cdots,\tau_d)$ is a point on the boundary of a generic rectangle $\delta_{l_1\ldots l_d}$;
		\item $\boldsymbol{\tau}^-:=(\tau_1^-,\cdots,\tau_d^-)$ and $\boldsymbol{\tau}^+:=(\tau_1^+,\cdots,\tau_d^+)$ are two points close to $\boldsymbol{\tau}$
		but belong to two different subrectangles $\delta_{l_1\ldots l_d}^-:=\{\bx \colon t_{\ell,l_\ell}^-\leq x_\ell<t_{\ell,l_\ell+1}^-\}$ and $\delta_{l_1\ldots l_d}^+:=\{\bx \colon t_{\ell,l_\ell}^+\leq x_\ell<t_{\ell,l_\ell+1}^+\}$;
		\item $\bt^{L}_-$ and $\bt^{L}_+$ are the starting points of $\delta_{l_1\ldots l_d}^-$ and $\delta_{l_1\ldots l_d}^+$;
		\item $(b_{1,-},\cdots, b_{d,-})$ and $(b_{1,+},\cdots,b_{d,+})$ are the corresponding mesh widths of $\delta_{l_1\ldots l_d}^-$ and $\delta_{l_1\ldots l_d}^+$; 
		\item $\Xi:=\{\ell\colon\; \tau_\ell^--\tau_\ell$ and
		$\tau_\ell^+-\tau_\ell$ differ in signs$\}$.
	\end{enumerate}
	In words, the index set $\Xi$ indicates the directions in which we cross boundaries when we move from $\boldsymbol{\tau}_\ell^-$ to $\boldsymbol{\tau}_\ell^+$. To further simplify notation, we write $\bar{\mu}(\boldsymbol{\tau}^-):=\lim_{\bx\rightarrow\boldsymbol{\tau},\bx\in \delta_{l_1\ldots l_d}^-}\bar{\mu}(\bx)$ and $\bar{\mu}(\boldsymbol{\tau}^+):=\lim_{\bx\rightarrow\boldsymbol{\tau},\bx\in \delta_{l_1\ldots l_d}^+}\bar{\mu}(\bx)$. Then we have
	\begin{align*}
	J_{\mathbf{0}}&=|\bar{\mu}(\boldsymbol{\tau}^+)-\bar{\mu}(\boldsymbol{\tau}^-)|
	=\Big|\bar{\mathscr{B}}_{m}(\boldsymbol{\tau}^+)-\bar{\mathscr{B}}_{m}(\boldsymbol{\tau}^-))\Big|\\
	&=\sum_{\ell\in\Xi}(B_{m}(0)|/m!)
	\Big|\frac{\partial^m \mu(\bt^{L}_+)}{\partial x_\ell^{m}}b_{\ell,+}^{m}-\frac{\partial^{m}\mu(\bt^{L}_-)}{\partial x_\ell^{m}}b_{\ell,-}^{m}\Big|\\
	&=\sum_{\ell\in\Xi} (B_{m}(0)|/m!)
	\Big|\Big(\frac{\partial^{m}\mu(\bt^{L}_+)}{\partial x_\ell^{m}}-\frac{\partial^{m}\mu(\bt^{L}_-)}{\partial x_\ell^{m}}\Big)b_{\ell,+}^{m}+\frac{\partial^{m}\mu(\bt^{L}_-)}{\partial x_\ell^{m}}(b_{\ell,+}^{m}-b_{\ell,-}^{m})\Big|\\
	&\leq\sum_{\ell\in\Xi}(B_{m}(0)|/m!)
	\Big[O(h^{m+\varrho})+Ch^{m-1}|b_{\ell,+}-b_{\ell,-}|\Big]\\
	&\leq\sum_{\ell\in\Xi}(B_{m}(0)|/m!)
	\Big[O(h^{m+\varrho})+Ch^{m-1}O(h^{1+\varrho})\Big],
	\end{align*}
	where the fourth line follows from Assumption \ref{assumption: DGP}, and the last from the stronger quasi-uniformity condition specified in the Lemma. This suffices to show that $J_{\mathbf{0}}$ is $O(h^{m+\varrho})$.	
	
	Then we mimic the proof strategy used in \cite[Theorem 12.7]{Schumaker_2007_Book}. By \cite[Theorem 12.6]{Schumaker_2007_Book}, we can construct a bounded linear operator $\mathscr{L}[\cdot]$ mapping $\mathcal{L}_1(\mathcal{X})$ onto $\mathcal{S}_{\Delta,m}$ with $\mathscr{L}[s]=s$ for all $s\in\mathcal{S}_{\Delta,m}$. Specifically, $\mathscr{L}[\cdot]$ is defined as
	$$\mathscr{L}[\mu](\bx):=\sum_{l_1=1}^{K_1}\cdots\sum_{l_d=1}^{K_d}(\psi_{l_1\ldots l_d}\mu)p_{l_1\ldots l_d}(\bx)$$
	where $\{\psi_{l_1\ldots l_d}\}_{l_1=1,\ldots,l_d=1}^{K_1,\ldots,K_d}$ is the dual basis defined in \cite[Equation (12.24)]{Schumaker_2007_Book}. By multi-dimensional Taylor expansion, there exists a polynomial $\varphi_{l_1\ldots l_d}$ such that
	$\|\bar{\mu}-\varphi_{l_1\ldots l_d}\|_{L_{\infty}(\delta_{l_1\ldots l_d})}\lesssim h^{m+\varrho}$,
	and the degree of $\varphi_{l_1\ldots l_d}$ is no greater than $m-1$.
	Since $\mathscr{L}$ reproduces polynomials, we have
	\begin{align*}
	\|\bar{\mu}-\mathscr{L}[\bar{\mu}]\|_{L_{\infty}(\delta_{l_1\ldots l_d})}&\leq\|\bar{\mu}-\varphi_{l_1\ldots l_d}\|_{L_{\infty}(\delta_{l_1\ldots l_d})}+\|\mathscr{L}[\bar{\mu}-\varphi_{l_1\ldots l_d}]\|_{L_{\infty}(\delta_{l_1\ldots l_d})}\\
	&\leq C\|\bar{\mu}-\varphi_{l_1\ldots l_d}\|_{L_{\infty}(\delta_{l_1\ldots l_d})}\lesssim h^{m+\varrho}.
	\end{align*}
	Taking account of the jumps of $\bar{\mu}$ along boundaries, the approximation error of $\mathscr{L}[\bar{\mu}]$ is still $O(h^{m+\varrho})$. Evaluating the $L_\infty$ norm on all subrectangles, we conclude that there exists some $s^*\in\mathcal{S}_{\Delta,m}$ such that
	$\|\mu+\bar{\mathscr{B}}_{m}-s^*\|_{L_{\infty}(\mathcal{X})}\lesssim h^{m+\varrho}$.
	
	For other $\bvsig$, we only need to show that the desired result holds for $s^*=\mathscr{L}[\bar{\mu}]$. By construction of $\mathscr{L}$,
	\begin{equation}\label{eq: pf lemma spline}
	|\partial^{\bvsig}(\mathscr{L}[\bar{\mu}])|
	\leq\sum_{l_1=1}^{m+\kappa_1}\cdots\sum_{l_d=1}^{m+\kappa_d}
	|\psi_{l_1\ldots l_d}\bar{\mu}||\partial^{\bvsig} p_{l_1\ldots l_d}(\bx)|
	\leq Ch^{-[\bvsig]}\|\bar{\mu}\|_{L_\infty(\delta_{l_1\ldots l_d})} 
	\end{equation}
	where the last line follows from \cite[Theorem 12.5]{Schumaker_2007_Book}. Then we have
	\begin{align*}
	\|\partial^{\bvsig}\mu+\partial^{\bvsig}\bar{\mathscr{B}}_{m}-\partial^{\bvsig}(\mathscr{L}[\bar{\mu}])\|_{L_\infty(\delta_{l_1\ldots l_d})}&\leq\|\partial^{\bvsig}\mu+\partial^{\bvsig}\mathscr{B}_{m}^*-\partial^{\bvsig}\varphi_{l_1\ldots l_d}\|_{L_{\infty}(\delta_{l_1\ldots l_d})}\\
	&\quad+\|\partial^{\bvsig} (\mathscr{L}[\bar{\mu}-\varphi_{l_1\ldots l_d}])\|_{L_{\infty}(\delta_{l_1\ldots l_d})}\\
	&\leq O(h^{m+\varrho-[\bvsig]})+Ch^{-[\bvsig]}\|\bar{\mu}-\varphi_{l_1\ldots l_d}\|_{L_{\infty}(\delta_{l_1\ldots l_d})}\\
	&\lesssim h^{m+\varrho-[\bvsig]}
	\end{align*}
	where the second inequality follows from Taylor expansion and Equation \eqref{eq: pf lemma spline}. Moreover, by an argument similar to that for $J_{\mathbf{0}}$, the jump of $\partial^{\bvsig}\bar{\mathscr{B}}_{m}$ is $O(h^{m+\varrho-[\bvsig]})$. 
	\smallskip
	
	(c): By construction of $\bpt$, $\rho=1$. By the same argument as in part (a) and (b), $\bpt$ satisfies Assumption \ref{assumption: Local Basis} and \ref{assumption: Misspecification Error}. Finally, by definition of tensor-product splines, both $\bp$ and $\bpt$ reproduce polynomials of degree no greater than $m-1$. Then the proof is complete.  
\end{proof}

\subsection{Proof of Lemma \ref{lemma-w}}
\begin{proof}
	(a):
	Assumption \ref{assumption: Local Basis}(a) holds by the fact that the father wavelet is compactly supported and $\{\phi_{sl}\}$ is generated by translation and dilation. Assumption \ref{assumption: DGP}(b) follows from the fact that $\{\phi_{sl}\}$ is an orthonormal basis with respect to the Lebesgue measure. For Assumption \ref{assumption: Local Basis}(c), notice that
	$$\frac{d^{\varsigma_\ell}\phi(2^sx_\ell-l_\ell)}{d x_\ell^{\varsigma_\ell}}=2^{s\varsigma_\ell}\frac{d^{\varsigma_\ell}\phi(z)}{dz^{\varsigma_\ell}}\Big|_{z=2^sx_\ell-l_\ell}=b^{-\varsigma_\ell}\frac{d^{\varsigma_\ell}\phi(z)}{dz^{\varsigma_\ell}}\Big|_{z=2^sx_\ell-l_\ell}.$$
	Since the wavelet basis reproduces polynomials of degree no greater than $m-1$ and $\phi$ is assumed to have $q+1$ continuous derivatives, the desired bounds follow.
	\smallskip
	
	(b): We employ the strategy used in \cite{Sweldens-Piessens_1994_NM}, but extend their proof to the multidimensional case. First, we denote by $V^\ell_s$ the closure of the level-$s$ subspace spanned by $\{\phi_{sl}(x_\ell)\}$, and let $W^\ell_s$ be the orthogonal complement of $V^\ell_s$ in $V^\ell_{s+1}$. Then we write $\mathcal{V}_s:=\otimes_{\ell=1}^d V^\ell_s$ for the space spanned by the tensor-product level-$s$ father wavelets, and $\mathcal{W}_s$ denotes the orthogonal complement of $\mathcal{V}_s$ in $\mathcal{V}_{s+1}$. We use the fact that 
	$\mathcal{W}_s=\oplus_{i=1}^{2^d-1}\,\mathcal{W}_{s,i}$,
	where $\oplus$ denotes ``direct sum'', and each $\mathcal{W}_{s,i}$ takes the following form:
	$\mathcal{W}_{s,i}=\otimes_{\ell=1}^d Z_s^\ell$.
	Each $Z_s^\ell$ is either $V_s^\ell$ or $W_s^\ell$, but $\{Z^\ell_s\}_{\ell=1}^d$ cannot be identical to $\{V_s^\ell\}_{\ell=1}^d$. There are in total $(2^d-1)$ such subspaces.
	Accordingly, a typical element of a basis vector for $\mathcal{W}_s$ can be written as
	$$\bar{\psi}_{s\mathbf{l}\boldsymbol{\alpha}}(\bx)=\prod_{\ell=1}^{d}\left[\alpha_\ell\phi_{sl_\ell}(x_\ell)+(1-\alpha_\ell)\psi_{sl_\ell}(x_\ell)\right]$$  
	where $\mathbf{l}=(l_1,\ldots, l_d)$ and $\alpha_\ell=0\text{ or }1$, but $\boldsymbol{\alpha}=(\alpha_1, \ldots, \alpha_d)\neq(1,\ldots, 1)$. Then it directly follows from the properties of wavelet basis that for $\bar{\psi}_{s\mathbf{l}\boldsymbol{\alpha}}$, $s\geq m$,
	\begin{equation}\label{eq: vanishing moment}
	\langle\bx^{\bvsig},\bar{\psi}_{s\mathbf{l}\boldsymbol{\alpha}}(\bx)\rangle:=\int_\mathcal{X}\bx^{\bvsig}\bar{\psi}_{s\mathbf{l}\boldsymbol{\alpha}}(\bx)d\bx=0,\; \text{for }
	\bvsig \text{ such that } [\bvsig]\leq m, \text{ and  }\varsigma_\ell\neq m\;\forall \ell.
	\end{equation}
	
	Denote by $\mathscr{L}_s[\cdot]$ the orthogonal projection operator onto $\mathcal{W}_s$. Then the approximation error of the tensor-product wavelet space $\mathcal{V}_{s_n}$ can be written as
	\begin{align*}
	&\quad\;\sum_{s=s_n}^\infty\mathscr{L}_s[\mu](\bx)
	=\sum_{s=s_n}^{\infty}\sum_{\boldsymbol{\alpha}}\sum_{\mathbf{l}}
	\langle \mu(\check{\bx}),\bar{\psi}_{s\mathbf{l}\boldsymbol{\alpha}}(\check{\bx})
	\rangle \bar{\psi}_{s\mathbf{l}\boldsymbol{\alpha}}(\bx)\\
	&=\sum_{s=s_n}^{\infty}\sum_{\boldsymbol{\alpha}}\sum_{\mathbf{l}}\Big\langle \sum_{[\bvsig]\leq m}\partial^{\bvsig} \mu(\bx)\frac{(\check{\bx}-\bx)^{\bvsig}}{\bvsig!}+
	\vartheta_n(\check{\bx},\bx),\;\bar{\psi}_{s\mathbf{l}\boldsymbol{\alpha}}(\check{\bx})\Big\rangle\,\bar{\psi}_{s\mathbf{l}\boldsymbol{\alpha}}(\bx)
	\end{align*}
	where $\vartheta_n(\check{\bx},\bx)\lesssim \|\check{\bx}-\bx\|^{m+\varrho}$, and the inner product in the above equations are taken with respect to $\check{\bx}$ in terms of the Lebesgue measure. The index sets where $\boldsymbol{\alpha}$ and $\mathbf{l}$ live are described in Section \ref{section: examples} and the proof, and omitted in the above derivation for simplicity. By Assumption \ref{assumption: DGP} and Assumption \ref{assumption: Local Basis},
	\begin{equation*}
	\begin{split}
	&\sup_{\bx\in\mathcal{X}} \Big|\sum_{s=s_n}^{\infty}
	\sum_{\boldsymbol{\alpha}} \sum_{\mathbf{l}}
	\Big\langle\vartheta_n(\check{\bx},\bx),\;\bar{\psi}_{s\mathbf{l}\boldsymbol{\alpha}}(\check{\bx})\Big\rangle
	\bar{\psi}_{s\mathbf{l}\boldsymbol{\alpha}}(\bx)\Big|\\
	=&\sup_{\bx\in\mathcal{X}}\Big|\sum_{s=s_n}^{\infty}
	\Big(\frac{b}{2^{s-s_n}}\Big)^{m+\varrho}\sum_{\boldsymbol{\alpha}}\sum_{\mathbf{l}}
	\Big\langle \vartheta_n(\check{\bx},\bx)2^{s(m+\varrho)},\;
	\bar{\psi}_{s\mathbf{l}\boldsymbol{\alpha}}(\check{\bx})\Big\rangle
	\bar{\psi}_{s\mathbf{l}\boldsymbol{\alpha}}(\bx)\Big|
	\lesssim b^{m+\varrho}.
	\end{split}
	\end{equation*}
	Recall that $b=2^{-s_n}$.
	
	Regarding the leading terms \[\sum_{s=s_n}^{\infty}\sum_{\boldsymbol{\alpha}}
	\sum_{\mathbf{l}}\Big\langle \sum_{[\bvsig]\leq m}
	\partial^{\bvsig}\mu(\bx)\frac{(\check{\bx}-\bx)^{\bvsig}}{\bvsig!},\;\bar{\psi}_{s\mathbf{l}\boldsymbol{\alpha}}(\check{\bx})
	\Big\rangle \bar{\psi}_{s\mathbf{l}\boldsymbol{\alpha}}(\bx),
	\]
	it is clear that the coefficients of the wavelet basis can be viewed as a linear combination of the inner products of monomials and the mother wavelets themselves, and thus by Equation \eqref{eq: vanishing moment}, the leading error is of order $b^m$ and can be characterized as
	$$\mathfrak{B}_{m,\mathbf{0}}(\bx)=-\sum_{\bm{u}\in\Lambda_m}
	\frac{b^m}{\bm{u}!}\partial^{\bm{u}} \mu(\bx)B^{\mathtt{W}}_{\bm{u},\mathbf{0}}(\bx/b).$$
	$B^{\mathtt{W}}_{\bm{u},\mathbf{0}}$ is referred to as ``monowavelet'' in \cite{Sweldens-Piessens_1994_NM}. Here we extend it to the multidimensional case. Specifically, define a mapping
	\begin{align*}
	\varphi:\;\Lambda_m &\rightarrow \{1,\ldots, d\}\\
	\bm{u}&\mapsto \ell
	\end{align*}
	such that $\varphi(\bm{u})$th element of $\bm{u}$ is nonzero. We denote $\mathbf{l}_{-\ell}:=(l_1,\ldots, l_{\ell-1}, l_{\ell+1},\ldots, l_d)$ and 
	$\mathcal{L}_{s}^{-\ell}:=\Big\{\mathbf{l}_{-\ell}: l_{\ell'}\in\mathcal{L}_s, j'=\{1,\cdots,d\}\setminus\{\ell\}\Big\}$.
	Then define
	$$\varpi_{\bm{u},s}(\bx)=\sum_{l_{\varphi(\bm{u})}\in\mathcal{L}_s}\;\sum_{\mathbf{l}_{-\varphi(\bm{u})}\in
		\mathcal{L}_s^{-\varphi(\bm{u})}}
	c_m\psi(2^sx_{\varphi(\bm{u})}
	-l_{\varphi(\bm{u})})\prod_{\ell=1,\ldots d\atop \ell\neq\varphi(\bm{u})}\phi(2^sx_\ell-l_\ell)$$
	where $c_m:=\int_0^1 x^m\psi(x)\,dx$. Then $B^{\mathtt{W}}_{\bm{u},\mathbf{0}}(\cdot)$ can be expressed as
	\begin{equation}\label{eq: wavelet error}
	B^{\mathtt{W}}_{\bm{u},\mathbf{0}}(\bx)=\sum_{s=0}^{\infty}
	2^{-sm}\varpi_{\bm{u},s}(\bx)=:
	\sum_{s=0}^{\infty}\xi_{\bm{u},s}(\bx).
	\end{equation}
	Moreover, since the series in Equation \eqref{eq: wavelet error} converges uniformly and for $s\geq s_n$, $\varpi_{\bm{u},s}^*(\bx)$ is orthogonal to the tensor-product wavelet basis $\bp$ with respect to the Lebesgue measure, by Dominated Convergence Theorem, the approximate orthogonality condition holds.
	
	For other $\bvsig$, let
	\begin{equation*}
	\mathscr{B}_{m,\bvsig}(\bx)=-\sum_{\bm{u}\in\Lambda_m}
	\frac{b^{m-[\bvsig]}}{\bm{u}!}\partial^{\bm{u}}\mu(\bx)
	B^\mathtt{W}_{\bm{u},\bvsig}(\bx/b)
	\end{equation*}
	where $B^\mathtt{W}_{\bm{u},\bvsig}(\bx)=\partial^{\bvsig} B^\mathtt{W}_{\bm{u},\mathbf{0}}(\bx)$. Since for $\bvsig$ such that $[\bvsig]\leq \varsigma$, $\partial^{\bvsig}\phi$ and $\partial^{\bvsig}\psi$ are continuously differentiable,
	$\sum_{s=0}^{\infty}2^{-sm}\partial^{\bvsig}\varpi_{\bm{u},s}(\bx)$ converges uniformly, and hence we can interchange the differentiation and infinite summation. Therefore, $B^{\mathtt{W}}_{\bm{u},\bvsig}(\cdot)$ is well defined and continuously differentiable. Then the lipschitz condition on $B^{\mathtt{W}}_{\bm{u},\bvsig}(\cdot)$ in Assumption \ref{assumption: Misspecification Error} holds.
	
	Let $s^*$ be the orthogonal projection of $\mu$ onto $\mathcal{V}_{s_n}$. To complete the proof of part (b), it suffices to show
	$\|\partial^{\bvsig}\mu-\partial^{\bvsig} s^*
	+\mathscr{B}_{m,\bvsig}\|_{L_\infty(\mathcal{X})}
	\lesssim b^{m+\varrho-[\bvsig]}$.
	Given a resolution level $s_n$,
	\begin{align*}
	&\sum_{s=s_n}^{\infty}\sum_{\boldsymbol{\alpha}}
	\sum_{\mathbf{l}}\langle \mu(\check{\bx}),
	\bar{\psi}_{s\mathbf{l}\boldsymbol{\alpha}}(\check{\bx})\rangle
	\partial^{\bvsig}\bar{\psi}_{s\mathbf{l}\boldsymbol{\alpha}}(\bx)\\
	=&\sum_{s=s_n}^{\infty}\sum_{\boldsymbol{\alpha}}
	\sum_{\mathbf{l}} \Big\langle \sum_{[\bm{u}]\leq m}
	\partial^{\bm{u}}\mu(\bx)\frac{(\check{\bx}-\bx)^{\bm{u}}}{\bm{u}!}+\vartheta_n(\check{\bx},\bx),\;\bar{\psi}_{s\mathbf{l}\boldsymbol{\alpha}}(\check{\bx})\Big\rangle
	\partial^{\bvsig}\bar{\psi}_{s\mathbf{l}\boldsymbol{\alpha}}(\bx)\\
	=&\,b^{m-[\bvsig]}\sum_{s=s_n}^{\infty}
	\frac{2^{[\bvsig](s_n-s)}}{2^{m(s-s_n)}}
	\sum_{\boldsymbol{\alpha}}\sum_{\mathbf{l}}\frac{2^{sd}}{2^{-sm}}
	\Big\langle\sum_{[\bm{u}]\leq m}\partial^{\bm{u}} \mu(\bx)
	\frac{(\check{\bx}-\bx)^{\bm{u}}}{\bm{u}!}
	+\vartheta_n(\check{\bx},\bx),\\
	&\hspace{16em} 2^{-sd/2}\bar{\psi}_{s\mathbf{l}\boldsymbol{\alpha}}(\check{\bx})\Big\rangle\;\partial^{\bvsig}\Big(2^{-sd/2}\bar{\psi}_{s\mathbf{l}\boldsymbol{\alpha}}(\bx)\Big)
	\end{align*}
	By changing variables, the vanishing moments of the wavelet function, and the fact that geometric series converges, the last line uniformly converges to the $\bvsig$th derivative of the approximation error of $\mathcal{V}_{s_n}$, $\mathscr{B}_{m,\bvsig}(\cdot)$ is the leading error and the remainder is $O(b^{m+\varrho-[\bvsig]})$.
	\smallskip
	
	(c): By construction of $\bpt$, $\rho=1$. By the same argument as that for part (a)  and (b), $\bpt$ satisfies Assumption \ref{assumption: Local Basis} and \ref{assumption: Misspecification Error}. Finally, both $\bp$ and $\bpt$ reproduce polynomials of degree no greater than $m-1$. Thus Assumption \ref{assumption: BC bases} holds. The proof is complete.
\end{proof}

\subsection{Proof of Lemma \ref{lemma: Partitioning Estimators}}
\begin{proof}
	(a):
	By construction, each basis function $p_k(\bx)$ is supported by only one subrectangle, and there are only a fixed number of $p_k(\bx)$'s which are not identically zero on each subrectangle. Thus Assumption \ref{assumption: Local Basis}(a) is satisfied. In addition, given a generic subrectangle $\delta_{l_1\ldots l_d}$, store all basis functions supported on $\delta_{l_1\ldots l_d}$ in a vector $\bp_{l_1\ldots l_d}$. By \cite[Lemma A.3]{Cattaneo-Farrell_2013_JoE},
	$\bQ_{l_1\ldots l_d}:=\E[\bp_{l_1\ldots l_d}(\bx_i)\bp_{l_1\ldots l_d}(\bx_i)']\asymp \bI_{\dim(\bR(\cdot))}$,
	where $\bI_{\dim(\bR(\cdot))}$ is an identity matrix of size $\dim(\bR(\cdot))$. In fact, $\int_{\delta_{l_1\ldots l_d}}\bp_{l_1\ldots l_d}(\bx)\bp_{l_1\ldots l_d}(\bx)'\,d\bx$ is a finite-dimensional matrix with the minimum eigenvalue bounded from below by $Ch^d$ for some $C>0$. Hence for any $\mathbf{a}\in\mathbb{R}^{\dim(\mathbf{R}(\cdot))}$,
	$$\mathbf{a}^\prime\int_{\delta_{l_1\ldots l_d}}\bp_{l_1\ldots l_d}(\bx)\bp_{l_1\ldots l_d}(\bx)^\prime\,d\bx\, \mathbf{a}\geq Ch^d\mathbf{a}^\prime\mathbf{a}$$
	which suffices to show Assumption \ref{assumption: Local Basis}(b).
	
	To show Assumption \ref{assumption: Local Basis}(c), simply notice that given any $\bx\in\mathcal{X}$, there are only a fixed number of nonzero elements in $\partial^\bq\bp(\bx)$, and for any $k=1, \ldots, K$,
	$$\sup_{\delta\in\Delta}\sup_{\bx\in\clo(\delta)}|\partial^{\bvsig} p_k(\bx)|\lesssim h^{-[\bvsig]}\max_{[\boldsymbol{\alpha}]=m-1} \frac{\boldsymbol{\alpha}!}{(\boldsymbol{\alpha}-\bvsig)!}.$$
	Moreover, for any $\bx\in\mathcal{X}$, there exists some $p_k$ in $\bp$ such that for $[\bvsig]\leq m-1$,
	$|\partial^{\bvsig} p_k(\bx)|\gtrsim h^{-[\bvsig]}$.
	\smallskip
	
	(b):
	The result directly follows from the proofs of Lemma A.2 and Theorem 3 in \cite{Cattaneo-Farrell_2013_JoE}.
	The only difference is that we use shifted Legendre polynomials to re-express the approximating function $s^*(\bx)=\bp(\bx)^\prime\boldsymbol{\beta}^*$ and the leading error. Clearly, $\boldsymbol{\beta}^*$ is just a linear combination of coefficients of power series basis defined in their paper. The orthogonality between the approximation basis and the leading error holds by the property of Legendre polynomials and the fact that every basis function is locally supported on only one cell.
	\smallskip
	
	(c): By construction of $\bpt$, $\rho=1$. By the same argument as that for part (a) and (b), $\bpt$ satisfies Assumption \ref{assumption: Local Basis} and \ref{assumption: Misspecification Error}. Finally, when the degree of piecewise polynomials is increased, $\bpt$ spans a larger space containing the span of $\bp$, and both bases reproduce polynomials of degree no greater than $m-1$. Thus Assumption \ref{assumption: BC bases} holds.	
\end{proof}

\subsection{Proof of Lemma \ref{lemma-pp}}
\begin{proof}
	Assumption \ref{assumption: Local Basis}(a), \ref{assumption: Local Basis}(c) and \ref{assumption: Misspecification Error} directly follow from the construction of this basis and Taylor expansion restricted to a particular cell. For Assumption \ref{assumption: Local Basis}(b), given a generic cell $\delta$, by Assumption \ref{assumption: quasi-uniform partition}, we can find an inscribed ball of diameter $\ell_1$ that is proportional to $h_\bx$. Thus we can further find an inscribed rectangle with lengths equal to $\ell_2$ that is also proportional to $h_\bx$. Thus, by changing variables and the same argument as that for the basis defined on rectangular cells, we have Assumption \ref{assumption: Local Basis}(b) holds. The properties of $\bpt$ follow similarly as in Lemma \ref{lemma: Partitioning Estimators}.
\end{proof}

\subsection{Proof of Lemma \ref{lem: restriction matrix}}
\begin{proof}
	For the upper bound on the maximum eigenvalue, simply notice that all elements in $\bR$ is bounded by some constant $C$, and the number of nonzeros in any row or column of $\bR$ is bounded by some constant $L$. Then for any $\boldsymbol{\alpha}\in\mathbb{R}^{\vartheta_\iota}$ such that $\|\boldsymbol{\alpha}\|=1$, $\boldsymbol{\alpha}'\bR\bR'\boldsymbol{\alpha}\leq 
	L^2C^2\|\boldsymbol{\alpha}\|^2\lesssim 1$.
	
	For the other side of the bound, since $\mathbf{R}\mathbf{R}'$ is a symmetric block Toeplitz matrix, Szeg\"{o}'s theorem and its extensions state that the asymptotic behavior of Toeplitz or block Toeplitz matrices is characterized by the corresponding Fourier transformation of their entries. See \cite{Grenander-Szego_2001_book} for more details. Specifically,
	$\mathbf{R}\mathbf{R}'$ is transformed into the following matrix
		
	\scriptsize
	\begin{equation*}
	\pmb{\mathscr{T}}_{\bar{\iota}}(\omega)
	=\begin{bmatrix}
	\mathbf{r}_{11}'\mathbf{r}_{11}+2\mathbf{r}_{11}'\mathbf{r}_{21}\cos\omega
	&\cdots&\mathbf{r}_{11}'\mathbf{r}_{1\bar{\iota}}+(\mathbf{r}_{11}'\mathbf{r}_{2\bar{\iota}}+\mathbf{r}_{1\bar{\iota}}'\mathbf{r}_{21})\cos\omega\\
	\vdots&\ddots&\vdots\\
	\mathbf{r}_{1\bar{\iota}}'\mathbf{r}_{11}+(\mathbf{r}_{1\bar{\iota}}'\mathbf{r}_{21}+\mathbf{r}_{11}'\mathbf{r}_{2\bar{\iota}})\cos\omega&\cdots&\mathbf{r}_{1\bar{\iota}}'\mathbf{r}_{1\bar{\iota}}+2\mathbf{r}_{1\bar{\iota}}'\mathbf{r}_{2\bar{\iota}}\cos\omega
	\end{bmatrix},
	\end{equation*}\normalsize
	Using the representation of $\bR\bR'$ given in the discussion after the lemma, we can concisely write
	$\pmb{\mathscr{T}}_{\bar{\iota}}(\omega)=\mathbf{A}+(\mathbf{B}+\mathbf{B}')\cos\omega$.	
	By \cite[Equation (7)]{Gazzah-Regalia-Delmas_2001_IEEE}, we have
	\begin{equation*}
	\lambda_{\min}(\bR\bR^\prime)\rightarrow \underset{\omega\in[0,2\pi]}{\min}\;\lambda_{\min}(\pmb{\mathscr{T}}_{\bar{\iota}}(\omega))\quad \text{as}\quad \kappa\rightarrow\infty.
	\end{equation*}
	
	The minimum of the minimum eigenvalue function of $\pmb{\mathscr{T}}_{\bar{\iota}}(\omega)$ is attainable since each entry of $\pmb{\mathscr{T}}_{\bar{\iota}}(\omega)$ is a linear function of $\cos\omega$, and thus each coefficient of the corresponding characteristic polynomial is a continuous function of $\cos\omega$. By Theorem 3.9.1 of \cite{Tyrtyshnikov_2012_book}, there exist $\bar{\iota}$ continuous functions of $\cos\omega$ such that they are the roots of the characteristic polynomial, and thus the minimum eigenvalue is a continuous function of $\cos\omega$. In addition, since $\bR\bR^\prime$ is positive semidefinite, this function is nonnegative over $[-1,1]$.
	
	By construction, $\bR\bR^\prime$ is a real symmetric positive semidefinite matrix, and thus its eigenvalues are real and nonnegative. Moreover, given any fixed $\kappa$, $\bR\bR^\prime$ is positive definite since the restrictions specified in $\bR$ are non-redundant. Therefore, it suffices to show that the limit of the minimum eigenvalue sequence is bounded away from zero. The original problem is transformed into showing that the minimum eigenvalue of a finite-dimensional matrix $\pmb{\mathscr{T}}_{\bar{\iota}}(\omega)$ is strictly positive for any $\omega\in[0,2\pi]$.
	
	The next critical fact we employ is that the smallest eigenvalue as a function of a real symmetric matrix is concave \cite[see][Property 2.1]{Hiriart-Ye_1995_NM}. In our case, each entry is a linear function of $\cos\omega$, and thus $\lambda_{\min}(\pmb{\mathscr{T}}_{\bar{\iota}}(\omega))$ is concave with respect to $\cos\omega$. Therefore, the minimum of the smallest eigenvalue function could be attained only at two endpoints, i.e., when $\cos\omega=1$ or $\cos\omega=-1$.
	
	We start with the case in which $(m-1)$ continuity constraints are imposed at each knot, i.e., $\bar{\iota}=m-1$. Since each knot is treated the same way, the restriction matrix $\bR$ can be fully characterized by $\bar{\iota}$ row vectors. A typical restriction that the $\varsigma$th derivative ($0\leq \varsigma \leq \bar{\iota}-1$) is continuous at a knot can be represented by the following vector
	$$\Big(\underbrace{\bar{P}_0^{(\varsigma)}(1),\ldots,\bar{P}_{m-1}^{(\varsigma)}(1)}_{\text{left interval}},\,
	\underbrace{-\bar{P}_0^{(\varsigma)}(-1),\ldots,-\bar{P}_{m-1}^{(\varsigma)}(-1)}_{\text{right interval}}\Big)$$
	where we omit all zero entries and $\bar{P}_l^{(\varsigma)}(x)$, $0
	\leq l \leq m-1$, denotes the $\varsigma$th derivative of the \textit{normalized} Legendre polynomial of degree $l$. Generally, Legendre polynomial $P_l(x)$ can be written as
	$$P_l(x)=\frac{1}{2^l}\sum_{i=0}^{l}\binom{l}{i}^2(x-1)^{l-i}(x+1)^i,$$
	and they have the following properties: for any $l, l'\in\mathbb{Z}_+$
	$$P_l(1)=1,\quad P_l(-x)=(-1)^lP_l(x),\quad \int_{-1}^{1}P_l(x)P_{l'}(x)\,dx=\frac{2}{2l+1}\delta_{ll'}$$
	where $\delta_{ll'}$ is the Kronecker delta. Therefore, $\bar{P}_l(x)=\frac{\sqrt{2l+1}}{\sqrt{2}}P_l(x)$.
	
	Using these formulas,
	$$P_l^{(\varsigma)}(1)=2^{-\varsigma}\varsigma!\binom{l}{\varsigma}\binom{l+\varsigma}{\varsigma}=2^{-\varsigma}\frac{(l+\varsigma)(l+\varsigma-1)\cdots(l-\varsigma+1)}{\varsigma!}.$$
	Therefore,
	$\bar{P}_l^{(\varsigma)}(1)=\frac{\sqrt{2l+1}}{\sqrt{2}}2^{-\varsigma}\varsigma!\binom{l}{\varsigma}\binom{l+\varsigma}{\varsigma}$.
	In addition, since Legendre polynomials are symmetric or antisymmetric, we have
	$\bar{P}_l^{(\varsigma)}(-1)=(-1)^{l+\varsigma}\bar{P}_l^{(\varsigma)}(1)$.
	Thus we obtain an explicit expression for $\bR$.
	
	Next, $\bR\bR^\prime$ can be fully characterized by two matrices $\mathbf{A}$ and $\mathbf{B}$ as in \eqref{eq: RR'}. In what follows we use $A[\varsigma, \ell]$ to denote the $(\varsigma, \ell)$th element of $\mathbf{A}$. The same notation is used for $\mathbf{B}$ and $\pmb{\mathscr{T}}_{\bar{\iota}}(\omega)$. If we arrange restrictions by increasing $\varsigma$ (here we allow row and column indices to start from $0$), then we have
	$$A[\varsigma,\varsigma]=\sum_{u=\varsigma}^{\bar{\iota}}(2u+1)2^{-2\varsigma}\left[\varsigma!\binom{u}{\varsigma}\binom{u+\varsigma}{\varsigma}\right]^2,$$
	and for $\varsigma>\ell$
	$$A[\varsigma,\ell]=\begin{cases}
	0\quad &\varsigma+\ell\text{  is odd}\\
	\sum_{u=\varsigma}^{\bar{\iota}}(2u+1)2^{-\varsigma-\ell}\varsigma!\ell!\binom{u}{\varsigma}\binom{u+\varsigma}{\varsigma}\binom{u}{\ell}\binom{u+\ell}{\ell}\quad &\varsigma+\ell\text{  is even}
	\end{cases}.
	$$
	$\mathbf{B}$ can be expressed explicitly as well:
	\begin{align*}
	&B[\varsigma,\varsigma]=\sum_{u=\varsigma}^{\bar{\iota}}(-1)^{u+\varsigma+1}\frac{2u+1}{2}2^{-2\varsigma}\left[\varsigma!\binom{u}{\varsigma}\binom{u+\varsigma}{\varsigma}\right]^2,\quad \text{and}\\
	&B[\varsigma,\ell]=\begin{cases}
	\sum_{u=\varsigma}^{\bar{\iota}}(-1)^{u+\varsigma+1}\frac{2u+1}{2}2^{-\varsigma-\ell}\varsigma!\ell!\binom{u}{\varsigma}\binom{u+\varsigma}{\varsigma}\binom{u}{\ell}\binom{u+\ell}{\ell}\quad &\varsigma>\ell\\
	(-1)^{\ell+\varsigma}B[\ell,\varsigma]&\varsigma<\ell
	\end{cases}.
	\end{align*}
	
	Therefore, $(\varsigma, \ell)$th element of $\pmb{\mathscr{T}}_{\bar{\iota}}(\omega)$ is
	$$\mathscr{T}_{\bar{\iota}}(\omega)[\varsigma,\ell]=\begin{cases}
	0\quad&\ell+\varsigma\text{  is odd}\\
	A[\varsigma,\ell]+2B[\varsigma,\ell]\cos\omega&\ell+\varsigma\text{  is even}
	\end{cases}.
	$$
	
	When $\ell+\varsigma$ is even, the corresponding entry of $\pmb{\mathscr{T}}_{\bar{\iota}}(\omega)$ is nonzero. In addition, the summands in $A[\varsigma,\ell]$ and $2B[\varsigma,\ell]$ are the same in terms of absolute values and only differ in signs.
	Consider the case when $\cos\omega=1$ and $\varsigma\geq \ell$. There are several cases:
	
	\begin{enumerate}[label=(\roman*)]
		\item $\bar{\iota}$ is even, $\varsigma$ is even
		\begin{scriptsize}
		$$\mathscr{T}_{\bar{\iota}}(\omega)[\varsigma, \ell]=
		\sum_{u=\varsigma/2}^{(\bar{\iota}-2)/2}\Big(2(2u+1)+1\Big)2^{-\varsigma-\ell+1}\varsigma!\ell!\binom{2u+1}{\varsigma}\binom{2u+1+\varsigma}{\varsigma}\binom{2u+1}{\ell}\binom{2u+1+\ell}{\ell};
		$$
		\end{scriptsize}
		\item $\bar{\iota}$ is even, $\varsigma$ is odd
		\begin{scriptsize}
		$$\mathscr{T}_{\bar{\iota}}(\omega)[\varsigma,\ell]=
		\sum_{u=(\varsigma+1)/2}^{\bar{\iota}/2}\Big(2(2u)+1\Big)2^{-\varsigma-\ell+1}\varsigma!\ell!\binom{2u}{\varsigma}\binom{2u+\varsigma}{\varsigma}\binom{2u}{\ell}\binom{2u+\ell}{\ell};
		$$
	    \end{scriptsize}
		\item $\bar{\iota}$ is odd, $\varsigma$ is even
		\begin{scriptsize}
		$$\mathscr{T}_{\bar{\iota}}(\omega)[\varsigma, \ell]=
		\sum_{u=\varsigma/2}^{(\bar{\iota}-1)/2}\Big(2(2u+1)+1\Big)2^{-\varsigma-\ell+1}\varsigma!\ell!\binom{2u+1}{\varsigma}\binom{2u+1+\varsigma}{\varsigma}\binom{2u+1}{\ell}\binom{2u+1+\ell}{\ell};
		$$  
	    \end{scriptsize}
		\item $\bar{\iota}$ is odd, $\varsigma$ is odd
		\begin{scriptsize}
		$$\mathscr{T}_{\bar{\iota}}(\omega)[\varsigma,\ell]=
		\sum_{u=(\varsigma+1)/2}^{(\bar{\iota}-1)/2}\Big(2(2u)+1\Big)2^{-\varsigma-\ell+1}\varsigma!\ell!\binom{2u}{\varsigma}\binom{2u+\varsigma}{\varsigma}\binom{2u}{\ell}\binom{2u+\ell}{\ell}.
		$$
	    \end{scriptsize}
	\end{enumerate}
	
	When $\bar{\iota}$ is odd, $\pmb{\mathscr{T}}_{\bar{\iota}}(\omega)$ can be written as a Gram matrix $\pmb{\mathscr{T}}_{\bar{\iota}}(\omega)=\mathbf{GG}^\prime$ where
	$$\mathbf{G}=\begin{bmatrix}
	\bar{P}_1(1)&0&\bar{P}_3(1)&\cdots&0&\bar{P}_{m-1}(1)\\
	0&\bar{P}^{(1)}_2(1)&0&\cdots&\bar{P}^{(1)}_{m-2}(1)&0\\
	0&0&\bar{P}^{(2)}_3(1)&\cdots&0&\bar{P}^{(2)}_{m-1}(1)\\
	\vdots&\vdots&\vdots&&\vdots&\vdots\\
	0&0&0&\cdots&0&\bar{P}_{m-1}^{(m-2)}(1)
	\end{bmatrix}.
	$$
	Clearly, it is a row echelon matrix and has full row rank. Thus, the minimum eigenvalue of $\pmb{\mathscr{T}}_{\bar{\iota}}(\omega)$ is strictly positive.
	
	When $\bar{\iota}$ is even, $\pmb{\mathscr{T}}_{\bar{\iota}}(\omega)$ can be written as $\mathbf{GG}^\prime$ where
	$$\mathbf{G}=\begin{bmatrix}
	\bar{P}_1(1)&0&\bar{P}_3(1)&\cdots&\bar{P}_{m-2}(1)&0\\
	0&\bar{P}^{(1)}_2(1)&0&\cdots&0&\bar{P}^{(1)}_{m-1}(1)\\
	0&0&\bar{P}^{(2)}_3(1)&\cdots&\bar{P}^{(2)}_{m-2}(1)&0\\
	\vdots&\vdots&\vdots&&\vdots&\vdots\\
	0&0&0&\cdots&0&\bar{P}_{m-1}^{(m-2)}(1)
	\end{bmatrix}.
	$$
	The case of $\cos\omega=-1$ can be proved the same way. This suffices to show that the minimum eigenvalue of $\pmb{\mathscr{T}}_{\bar{\iota}}(\omega)$ as a function of $\cos\omega$ is strictly positive at the two endpoints, 1 and -1, thus completing the proof for $\bar{\iota}=m-1$.
	
	To complete the proof of the lemma, it remains to extend this result to the case in which fewer constraints are imposed. Compared with the case when $(m-1)$ constraints are imposed, some rows in the bigger restriction matrix are removed, and accordingly, $\bR\bR^\prime$ is a principle submatrix of the original one. By Cauchy Interlacing Theorem, the smallest eigenvalue of the principle submatrix must be no less than the smallest eigenvalue of the original matrix. Combining this fact with the results proved for $\bar{\iota}=m-1$, we have the minimum eigenvalue of $\bR\bR^\prime$ uniformly bounded away from $0$ when fewer restrictions are imposed, and then the proof is complete.
\end{proof}

\newpage
\newcounter{model}
\forloop{model}{1}{\value{model} < 8}{
	
	\begin{table}\renewcommand{\arraystretch}{1.1}
		\centering
		\caption{Pointwise Results, Model \themodel, Spline, Evenly-spaced, $n=1000$, 5000 Replications}
		\label{table:pointwise_uni\themodel}
		\resizebox{\textwidth}{!}{\input{simuls/output/table_pointwise_bs_es_\themodel.txt}}\\
	\end{table}\clearpage
	
}

\forloop{model}{1}{\value{model} < 8}{
	
	\begin{table}\renewcommand{\arraystretch}{1.1}
		\centering
		\caption{Pointwise Results, Model \themodel, Spline, Quantile-spaced, $n=1000$, 5000 Replications}
		\label{table:pointwise_qua\themodel}
		\resizebox{\textwidth}{!}{\input{simuls/output/table_pointwise_bs_qs_\themodel.txt}}\\
	\end{table}\clearpage
	
}


\begin{table}\renewcommand{\arraystretch}{1.1}
	\centering
	\caption{Tuning parameters, Spline, Evenly-spaced, $n=1000$, 5000 Replications}
	\label{table:kappa uni}
	\resizebox{!}{!}{
\begin{tabular}{lrrrrrrr}
\hline\hline
\multicolumn{1}{l}{}&\multicolumn{1}{c}{Min.}&\multicolumn{1}{c}{1st qu.}&\multicolumn{1}{c}{Median}&\multicolumn{1}{c}{Mean}&\multicolumn{1}{c}{3rd qu.}&\multicolumn{1}{c}{Max.}&\multicolumn{1}{c}{Std. dev.}\tabularnewline
\hline
{\bfseries Model 1, $\kappa_{\mathtt{IMSE}}=3.4$}&&&&&&&\tabularnewline
~~$\hat{\kappa}_{\mathtt{ROT}}$&$3$&$4$&$4$&$4.07$&$4$&$ 5$&$0.45$\tabularnewline
~~$\hat{\kappa}_{\mathtt{DPI}}$&$3$&$4$&$5$&$4.67$&$5$&$ 7$&$0.72$\tabularnewline
\hline
{\bfseries Model 2, $\kappa_{\mathtt{IMSE}}=7.2$}&&&&&&&\tabularnewline
~~$\hat{\kappa}_{\mathtt{ROT}}$&$5$&$6$&$6$&$6.07$&$6$&$ 7$&$0.26$\tabularnewline
~~$\hat{\kappa}_{\mathtt{DPI}}$&$6$&$7$&$8$&$7.79$&$8$&$10$&$0.57$\tabularnewline
\hline
{\bfseries Model 3, $\kappa_{\mathtt{IMSE}}=11.9$}&&&&&&&\tabularnewline
~~$\hat{\kappa}_{\mathtt{ROT}}$&$3$&$3$&$4$&$3.68$&$4$&$ 4$&$0.47$\tabularnewline
~~$\hat{\kappa}_{\mathtt{DPI}}$&$5$&$7$&$8$&$7.65$&$8$&$ 9$&$1.09$\tabularnewline
\hline
{\bfseries Model 4, $\kappa_{\mathtt{IMSE}}=4.4$}&&&&&&&\tabularnewline
~~$\hat{\kappa}_{\mathtt{ROT}}$&$2$&$3$&$3$&$2.94$&$3$&$ 3$&$0.24$\tabularnewline
~~$\hat{\kappa}_{\mathtt{DPI}}$&$4$&$4$&$4$&$4.00$&$4$&$ 5$&$0.02$\tabularnewline
\hline
{\bfseries Model 5, $\kappa_{\mathtt{IMSE}}=4.2$}&&&&&&&\tabularnewline
~~$\hat{\kappa}_{\mathtt{ROT}}$&$2$&$3$&$3$&$3.00$&$3$&$ 4$&$0.06$\tabularnewline
~~$\hat{\kappa}_{\mathtt{DPI}}$&$3$&$4$&$4$&$3.93$&$4$&$ 4$&$0.26$\tabularnewline
\hline
{\bfseries Model 6, $\kappa_{\mathtt{IMSE}}=3.3$}&&&&&&&\tabularnewline
~~$\hat{\kappa}_{\mathtt{ROT}}$&$2$&$2$&$2$&$2.10$&$2$&$ 3$&$0.31$\tabularnewline
~~$\hat{\kappa}_{\mathtt{DPI}}$&$3$&$3$&$3$&$3.00$&$3$&$ 4$&$0.04$\tabularnewline
\hline
{\bfseries Model 7, $\kappa_{\mathtt{IMSE}}=2.3$}&&&&&&&\tabularnewline
~~$\hat{\kappa}_{\mathtt{ROT}}$&$2$&$2$&$2$&$2.07$&$2$&$ 3$&$0.26$\tabularnewline
~~$\hat{\kappa}_{\mathtt{DPI}}$&$2$&$3$&$3$&$3.00$&$3$&$ 4$&$0.05$\tabularnewline
\hline
\end{tabular}
}
	\\
\end{table}\clearpage

\begin{table}\renewcommand{\arraystretch}{1.1}
	\centering
	\caption{Tuning parameters, Spline, Quantile-spaced, $n=1000$, 5000 Replications}
	\label{table:kappa qua}
	\resizebox{!}{!}{
\begin{tabular}{lrrrrrrr}
\hline\hline
\multicolumn{1}{l}{}&\multicolumn{1}{c}{Min.}&\multicolumn{1}{c}{1st qu.}&\multicolumn{1}{c}{Median}&\multicolumn{1}{c}{Mean}&\multicolumn{1}{c}{3rd qu.}&\multicolumn{1}{c}{Max.}&\multicolumn{1}{c}{Std. dev.}\tabularnewline
\hline
{\bfseries Model 1, $\kappa_{\mathtt{IMSE}}=3.4$}&&&&&&&\tabularnewline
~~$\hat{\kappa}_{\mathtt{ROT}}$&$3$&$4$&$4$&$4.07$&$4$&$ 5$&$0.45$\tabularnewline
~~$\hat{\kappa}_{\mathtt{DPI}}$&$3$&$4$&$5$&$4.69$&$5$&$ 7$&$0.72$\tabularnewline
\hline
{\bfseries Model 2, $\kappa_{\mathtt{IMSE}}=7.2$}&&&&&&&\tabularnewline
~~$\hat{\kappa}_{\mathtt{ROT}}$&$5$&$6$&$6$&$6.07$&$6$&$ 7$&$0.26$\tabularnewline
~~$\hat{\kappa}_{\mathtt{DPI}}$&$6$&$7$&$8$&$7.84$&$8$&$10$&$0.58$\tabularnewline
\hline
{\bfseries Model 3, $\kappa_{\mathtt{IMSE}}=11.9$}&&&&&&&\tabularnewline
~~$\hat{\kappa}_{\mathtt{ROT}}$&$3$&$3$&$4$&$3.68$&$4$&$ 4$&$0.47$\tabularnewline
~~$\hat{\kappa}_{\mathtt{DPI}}$&$5$&$7$&$8$&$7.66$&$8$&$ 9$&$1.05$\tabularnewline
\hline
{\bfseries Model 4, $\kappa_{\mathtt{IMSE}}=4.4$}&&&&&&&\tabularnewline
~~$\hat{\kappa}_{\mathtt{ROT}}$&$2$&$3$&$3$&$2.94$&$3$&$ 3$&$0.24$\tabularnewline
~~$\hat{\kappa}_{\mathtt{DPI}}$&$4$&$4$&$4$&$4.00$&$4$&$ 5$&$0.03$\tabularnewline
\hline
{\bfseries Model 5, $\kappa_{\mathtt{IMSE}}=4.2$}&&&&&&&\tabularnewline
~~$\hat{\kappa}_{\mathtt{ROT}}$&$2$&$3$&$3$&$3.00$&$3$&$ 4$&$0.06$\tabularnewline
~~$\hat{\kappa}_{\mathtt{DPI}}$&$3$&$4$&$4$&$3.93$&$4$&$ 4$&$0.26$\tabularnewline
\hline
{\bfseries Model 6, $\kappa_{\mathtt{IMSE}}=3.3$}&&&&&&&\tabularnewline
~~$\hat{\kappa}_{\mathtt{ROT}}$&$2$&$2$&$2$&$2.10$&$2$&$ 3$&$0.31$\tabularnewline
~~$\hat{\kappa}_{\mathtt{DPI}}$&$3$&$3$&$3$&$3.00$&$3$&$ 4$&$0.04$\tabularnewline
\hline
{\bfseries Model 7, $\kappa_{\mathtt{IMSE}}=2.3$}&&&&&&&\tabularnewline
~~$\hat{\kappa}_{\mathtt{ROT}}$&$2$&$2$&$2$&$2.07$&$2$&$ 3$&$0.26$\tabularnewline
~~$\hat{\kappa}_{\mathtt{DPI}}$&$2$&$3$&$3$&$3.00$&$3$&$ 4$&$0.04$\tabularnewline
\hline
\end{tabular}
}
	\\
\end{table}\clearpage

\forloop{model}{1}{\value{model} < 8}{
	
	\begin{table}\renewcommand{\arraystretch}{1.1}
		\centering
		\caption{Uniform Results, Model \themodel, Spline, $n=1000$, 5000 Replications.}
		\label{table:uniform\themodel}
		\resizebox{!}{!}{\input{simuls/output/table_uniform_bs_\themodel.txt}}\\
	\end{table}\clearpage
	
}

\forloop{model}{1}{\value{model} < 6}{
	
	\begin{table}\renewcommand{\arraystretch}{1.1}
		\centering
		\caption{Pointwise Results, Model \themodel, Wavelet, $n=1000$, 5000 Replications}
		\label{table:pointwise_wav\themodel}
		\resizebox{\textwidth}{!}{\input{simuls/output/table_pointwise_wav_\themodel.txt}}\\
	\end{table}\clearpage
	
}

\begin{table}\renewcommand{\arraystretch}{1.2}
	\centering
	\caption{Tuning parameters, Wavelets, $n=1000$, 5000 Replications}
	\label{table:kappa wav}
	\resizebox{.9\textwidth}{1.1\height}{
\begin{tabular}{lrrrrrrr}
\hline\hline
\multicolumn{1}{l}{}&\multicolumn{1}{c}{Min.}&\multicolumn{1}{c}{1st qu.}&\multicolumn{1}{c}{Median}&\multicolumn{1}{c}{Mean}&\multicolumn{1}{c}{3rd qu.}&\multicolumn{1}{c}{Max.}&\multicolumn{1}{c}{Std. dev.}\tabularnewline
\hline
{\bfseries Model 1, $s_{\mathtt{IMSE}}=2.4$}&&&&&&&\tabularnewline
~~$\hat{s}_{\mathtt{ROT}}$&$2$&$2$&$2$&$2.00$&$2$&$3$&$0.05$\tabularnewline
~~$\hat{s}_{\mathtt{DPI}}$&$2$&$3$&$3$&$2.81$&$3$&$3$&$0.39$\tabularnewline
\hline
{\bfseries Model 2, $s_{\mathtt{IMSE}}=3.5$}&&&&&&&\tabularnewline
~~$\hat{s}_{\mathtt{ROT}}$&$3$&$3$&$3$&$3.00$&$3$&$3$&$0.00$\tabularnewline
~~$\hat{s}_{\mathtt{DPI}}$&$3$&$3$&$3$&$3.03$&$3$&$4$&$0.16$\tabularnewline
\hline
{\bfseries Model 3, $s_{\mathtt{IMSE}}=4.2$}&&&&&&&\tabularnewline
~~$\hat{s}_{\mathtt{ROT}}$&$2$&$2$&$2$&$2.00$&$2$&$2$&$0.00$\tabularnewline
~~$\hat{s}_{\mathtt{DPI}}$&$3$&$3$&$3$&$3.15$&$3$&$4$&$0.35$\tabularnewline
\hline
{\bfseries Model 4, $s_{\mathtt{IMSE}}=2.7$}&&&&&&&\tabularnewline
~~$\hat{s}_{\mathtt{ROT}}$&$1$&$2$&$2$&$1.85$&$2$&$2$&$0.35$\tabularnewline
~~$\hat{s}_{\mathtt{DPI}}$&$2$&$2$&$2$&$2.00$&$2$&$2$&$0.00$\tabularnewline
\hline
{\bfseries Model 5, $s_{\mathtt{IMSE}}=2.6$}&&&&&&&\tabularnewline
~~$\hat{s}_{\mathtt{ROT}}$&$1$&$2$&$2$&$1.99$&$2$&$2$&$0.09$\tabularnewline
~~$\hat{s}_{\mathtt{DPI}}$&$2$&$2$&$2$&$2.00$&$2$&$2$&$0.00$\tabularnewline
\hline
\end{tabular}
}
	\\
\end{table}\clearpage

\forloop{model}{1}{\value{model} < 6}{
	
	\begin{table}\renewcommand{\arraystretch}{1.2}
		\centering
		\caption{Uniform Results, Model \themodel, Wavelet, $n=1000$, 5000 Replications.}
		\label{table:uniform_wav\themodel}
		\resizebox{0.5\textwidth}{.95\height}{\input{simuls/output/table_uniform_wav_\themodel.txt}}\\
	\end{table}\clearpage
	
}

\clearpage


\bibliography{Cattaneo-Farrell-Feng_2019_AoS--Bibliography}
\bibliographystyle{imsart-number}